\documentclass[11pt,reqno]{amsart}

\numberwithin{equation}{section} 
\usepackage{amsfonts,dsfont,scalerel,braket,booktabs}

\usepackage[usenames, dvipsnames]{color}
\newcommand{\mkr}{\color{black}}
\newcommand{\hl}{\color{black}}
\usepackage{amssymb}
\usepackage{subcaption}
\usepackage{bbm}
\usepackage{empheq}
\usepackage{fancyhdr}
\usepackage{tikz}
\usetikzlibrary{arrows,positioning} 
\usepackage{dsfont}
\usepackage{needspace}
\usepackage{booktabs}
\usepackage{framed}
\usepackage{sidecap}
\usepackage{bm} % for bold math symbols
    \usepackage{cases} 
    \usepackage{color}
\usepackage{graphicx, graphics}
\usepackage[ruled]{algorithm2e}
\usepackage[breaklinks=true,colorlinks=true,linkcolor=blue,urlcolor=blue,citecolor=blue]{hyperref}

\newtheorem{thm}{Theorem}
\newtheorem{lem}{Lemma}
\newtheorem{defi}{Definition}

\newtheorem{rem}{Remark}

\def\bt{\begin{thm}}
\def\et{\end{thm}}
\def\bl{\begin{lem}}
\def\el{\end{lem}}
\def\bd{\begin{defi}}
\def\ed{\end{defi}}
\def\bc{\begin{cor}}
\def\ec{\end{cor}}
\def\bp{\begin{proof}}
\def\ep{\end{proof}}
\def\br{\begin{rem}}
\def\er{\end{rem}}

\def\bpp{\begin{proof}}
\def\epp{\end{proof}}

\numberwithin{equation}{section} \oddsidemargin=-.0cm
\def\bi{\begin{itemize}}
\def\ei{\end{itemize}}
\def\ben{\begin{enumerate}}
\def\een{\end{enumerate}}
\def\be{\begin{equation}}
\def\ee{\end{equation}}
\def\bes{\begin{equation*}}
\def\ees{\end{equation*}}
\def\bea{\begin{equation} \begin{aligned}}
\def\eea{\end{aligned} \end{equation}}
\def\beas{\begin{equation*} \begin{aligned}}
\def\eeas{\end{aligned} \end{equation*}}
\def\d{\, \mathrm{d}}

\newcommand{\norm}[1]{\left\lVert#1\right\rVert}
\def\Forall{\text{ } \forall \:}
\def\bftau{{\boldsymbol{\tau}}}
\def\bfbeta{{\boldsymbol{\beta}}}
\def\tr{\mathrm{tr}}
\def\s{\mathfrak{s}}
\def\c{\mathfrak{c}}

\def\cH{\mathcal H}

\def\m{\mathfrak{m}}
\def\e{\epsilon}
\def\F{\mathcal{F}}
\def\R{\mathbf{R}}

\begin{document}

\title[Variational approach to closure of nonlinear dynamical systems]{Variational approach to closure of nonlinear dynamical systems: Autonomous case}


%%%%%%%%%%%%%%%%%%%%%%%%%%%%%%%%%%%%%%
%%%% JSP style for author & institute 
%\author{Micka\"el D. Chekroun \and Honghu Liu \and James C. McWilliams}
%\institute{M. D. Chekroun \at   Department of Earth and Planetary Sciences, Weizmann Institute, Rehovot 76100, Israel;
%Department of Atmospheric and Oceanic Sciences and Institute of Geophysics and Planetary Physics, University of California, Los Angeles, USA\\
%  \email{mchekroun@atmos.ucla.edu}
%    \and
% H. Liu \at Department of Mathematics, Virginia Tech, Blacksburg, Virginia 24061, USA \\
% \email{hhliu@vt.edu}
%  \and
%  J. C. McWilliams \at
%  Department of Atmospheric and Oceanic Sciences and Institute of Geophysics and Planetary Physics, University of California, Los Angeles, USA\\
%  \email{jcm@atmos.ucla.edu}
%}
%%%%%%%%%%%%%%%%%%%%%%%%%%%%%%%%%%%%%%%%

%%%%%%%%%%%%%%%%%%%%%%%%%%%%%%%%%%%%%%%%
%%% amsart style for author & institute
\author[Micka\"el D. Chekroun]{Micka\"el D. Chekroun}
\address[MDC]{Department of Earth and Planetary Sciences, Weizmann Institute, Rehovot 76100, Israel;
Department of Atmospheric and Oceanic Sciences and Institute of Geophysics and Planetary Physics, University of California, Los Angeles, USA}
\email{mchekroun@atmos.ucla.edu}
\author[Honghu Liu]{Honghu Liu}
\address[HL]{Department of Mathematics, Virginia Tech, Blacksburg, Virginia 24061, USA}
\email{hhliu@vt.edu}
\author[James C. McWilliams]{James C. McWilliams}
\address[JCM]{Department of Atmospheric and Oceanic Sciences and Institute of Geophysics and Planetary Physics, University of California, Los Angeles, USA}
\email{jcm@atmos.ucla.edu}

%%%%%%%%%%%%%%%%%%%%%%%%%%%%%%%%%%%%%%%%

%\date{April 30, 2019} % original submission date to JSP
\date{December 1, 2019} 

\begin{abstract}
A general approach for the derivation of nonlinear parameterizations of neglected scales is presented for nonlinear systems subject to 
an autonomous forcing. 
In that respect, dynamically-based formulas are derived subject to a free scalar parameter to be determined per mode to parameterize. 
For each high mode, this free parameter is obtained by minimizing a cost functional --- a parameterization defect --- depending on solutions from direct numerical simulation (DNS)  but over short training periods of length comparable to a characteristic recurrence or decorrelation time of the dynamics.

An important class of dynamically-based formulas, for our parameterizations to optimize, are obtained as parametric variations of manifolds approximating the invariant ones. To better appreciate the origins of the modified manifolds thus obtained, the standard approximation theory of invariant manifolds is revisited in Part I of this article. A special emphasis is put on backward-forward (BF) systems naturally associated with the original system,  whose asymptotic integration provides the leading-order approximation of invariant manifolds. 

Part II presents then (i) the modifications of these approximating manifolds based also on integration of the same BF systems but this time over a finite time $\tau$, and (ii) the variational approach aimed at making an efficient selection of $\tau$ per mode to  parameterize.  The parametric class of leading interaction approximation (LIA)  of the high modes obtained this way, is completed by another parametric class built from the quasi-stationary approximation (QSA); close to the first criticality, the QSA is an approximation to the LIA, but  it differs as one moves away from criticality. 

Rigorous results are derived that show that --- given a cutoff dimension --- the best manifolds that can be obtained through our variational approach, are manifolds which are  in general no longer invariant.
 The minimizers are objects, called the optimal parameterizing manifolds (PMs),  that are intimately tied to the conditional expectation of the original system, i.e.~the best vector field of the reduced state space resulting from averaging of the unresolved variables 
with respect to a probability measure conditioned on the resolved variables.

Applications to the closure of low-order models of Atmospheric Primitive Equations and Rayleigh-B\'enard convection are then discussed. 
The approach is finally illustrated --- in the context of the Kuramoto-Sivashinsky turbulence --- as providing efficient closures without slaving for a cutoff scale $k_\c$ placed within the inertial range and the reduced state space is just spanned by the unstable modes, without inclusion of  any stable modes whatsoever. The underlying optimal PMs obtained by our variational approach are far from slaving and  allow for  remedying the excessive backscatter transfer of energy to the low modes encountered by the LIA or the QSA parameterizations in their standard forms, when they are used at this cutoff wavelength.

\end{abstract}


\maketitle

\tableofcontents
%\linenumbers{}
\section{Introduction}

A number of theories have been proposed to explain the phenomenon of turbulence in fluid dynamics, but none has been universally accepted. Landau \cite{landau1959fluid} and Hopf \cite{hopf1948mathematical} suggested that turbulence is the result of an infinite sequence of bifurcations, each adding another independent period to a quasi-periodic motion of increasingly greater complexity. More recently, it has been shown numerically that the original quasiperiodic Landau's view of turbulence, with the amendment of the inclusion of stochasticity, may be well suited to describe certain turbulent behavior \cite{KCB18}, at least for the motion of large eddies. In the 1970's it has been theoretically argued and confirmed by many experiments  that dynamical systems may exhibit strange attractors which result in chaotic but deterministic behavior after a (very) few bifurcations have taken place. Ruelle and Takens \cite{ruelle1971nature} and others have suggested this as a mechanism underlying turbulence. 
In realistic physical problems one is seldomly able to carry out the mathematics beyond the first or second bifurcation, in particular regarding the derivation of reduced equations that capture effectively the amplitude and frequency content of the bifurcated solutions \cite{langford1979periodic,crawford1991introduction}. Noteworthy is normal form reduction that have been carried for degenerate singularities with simultaneous onset of co-existing and possibly many instabilities, but still close to first criticality 
\cite{coullet1983amplitude,arneodo1985asymptotic,elphick1987simple}.

It is typical of many bifurcation problems that, as the condition for instability is exceeded, increasingly many modes become unstable. This circumstance considerably complicates an effective reduction because it often corresponds to going through higher-order bifurcations to reach possibly chaos, for which a failure of the {\it slaving principle} of the unresolved variables onto the resolved ones --- mandatory for the success of standard reduction techniques --- is typically observed.

Center manifold techniques \cite{Vanderbauwhede89,crawford1991introduction,GH90} require such a slaving principle to provide an efficient reduction of the dynamics, and in that sense is reliable only in the vicinity of low-order bifurcations associated with the onset of instability. 
Center manifolds form a particular class of more general invariant manifolds associated with a fixed point, on which solutions obey {\it de facto} a slaving principle.  A comprehensive treatment of the computational aspects relative to the underlying parameterizations  can be found in \cite{haro2016parameterization}. The treatment in \cite{haro2016parameterization} is based on the so-called {\it parameterization method} \cite{cabre2003parameterization_I,cabre2003parameterization_II,cabre2005parameterization} itself built upon the invariance equation (see Eq.~\eqref{eq:invariance} below) and 
the associated cohomological equations that the sought (slaving) parameterization solves at different orders. 
The parameterization method allows for efficient computations for not only the case of invariant manifolds associated with fixed points, but also for the cases of invariant tori for autonomous or quasi-periodically forced systems, averaging and periodic diffeomorphisms \cite{Chekroun_al06}, invariant tori in Hamiltonian systems \cite{haro2016parameterization}, as well as  normally hyperbolic invariant tori. Other complementary approaches include e.g.~the Lyapunov-Schmidt reduction \cite{Golubitsky_al85,MW05} and the Lyapunov-Perron method \cite{Hen81,MW05}, as well as the usage of symmetries \cite{Golubitsky_al85,HI10}.

Despite the success for analyzing a broad class of bifurcations or detecting special solutions in dynamical systems such as quasi-periodic ones,  these methods relying on invariant manifold theory, have failed to prove their efficiency for reducing complicated behaviors resulting from the presence of chaos. In a certain sense, the ``story'' of the inertial manifold (IM) constitutes perhaps an epitome of this failure.  Despite appealing mathematical results showing existence of  IMs  for a broad class of dissipative systems \cite{Foias_al85CRAS,FST88,mallet1988inertial,CFNT89,Tem97}, and convergence error estimates when e.g.~slaving is not guaranteed to be satisfied (Approximate Inertial Manifold (AIM)) \cite{MT89,debussche1992construction,devulder1993rate,jones1994remark}, early promises \cite{Foias_al88KS,FMT88,JKT90,dubois1991subgrid,jolly1993bifurcation} have been challenged due to practical shortcomings pointed out for  efficient closure by IMs or AIMs for turbulent flows and route chaos \cite{daley1980development,foias1991approximate,jolly1991preserving,pascal1992nonlinear,graham1993computational,haugen1993spectral,garcia1995time}.

Essentially, the current  IM theory \cite{zelik2014inertial} predicts that the underlying slaving of the high modes to the low modes, holds when the cutoff wavenumber, $k_\c$,  is taken sufficiently far within the dissipative range, especially in ``strongly'' turbulent regimes that correspond e.g.~to the presence of many unstable modes. Still, as the AIM theory underlines, satisfactory closures may be expected to be derived for $k_\c$ corresponding to scales larger than what  the IM theory predicts. Nevertheless, as one seeks to further decrease $k_\c$   within the inertial range, standard AIMs fail typically in providing relevant closures  and one needs to rely on no longer a fixed  cutoff but instead a dynamic one so as to avoid energy accumulation on the cutoff level \cite{debussche1995nonlinear,dubois1998incremental,dubois1998dynamic}.

In general,  to aim at closing a given chaotic system at a fixed cutoff scale such that the neglected scales contain a non-negligible fraction of the energy\footnote{Such as ``cutting'' within the inertial range of turbulence.}, makes, a priori, the closure problem difficult to address. This difficulty is often manifested by either an  under- or over-parameterization of the small scales, i.e.~a deficient or excessive parameterization of the small-scale energy,  leading to an incorrect reproduction of the backscatter transfer of energy to the large scales \cite{kraichnan1976eddy,leith1990stochastic,piomelli1991subgrid,jansen2014parameterizing,berloff2015dynamically}. Thus, a deficiency in the (nonlinear) parameterization of the high modes leads to errors in the backscatter transfer of energy which is due to nonlinear interactions between the modes, especially those near the cutoff scale. 
We can speak of an inverse error cascade,  i.e.~errors in the modeling of the parameterized (small) scales that contaminate gradually the larger scales, and may spoil severely the closure skills for the resolved variables.

To remedy such a pervasive issue, it is thus reasonable, given a cutoff scale to seek for nonlinear parameterizations (manifolds) that minimize as much as possible a defect of parameterization in order to reduce spurious backscatter transfer of energy to the large scales. 
Obviously such manifolds should coincide with the invariant ones as one approaches towards the first bifurcation.

This latter point explains the two-part structure of our article.  We show here that an important class of dynamically-based formulas for our parameterizations are obtained as  parametric variations of manifolds approximating the invariant ones. To better appreciate the origins of the modified manifolds thus obtained, the standard approximation theory of invariant manifolds is revisited in Part I of this article. A special emphasis is put on backward-forward (BF) systems  naturally associated with the original system,  whose asymptotic integration provides the leading-order approximation of invariant manifolds. 

Part II presents then (i) the modifications of these approximating manifolds based also on integration of the same BF systems but this time over a finite time $\tau$, and (ii) the variational approach aimed at making an efficient selection of $\tau$ per mode to  parameterize, in order to minimize a parameterization defect.  The parametric class of {\it leading interaction approximation (LIA)}  of the high modes obtained this way, is completed by another parametric class built from the {\it quasi-stationary approximation (QSA)}; close to the first criticality, the QSA is an approximation to the LIA, but  differs as one moves away from criticality.

In this article our formulations are general, but our primary motivations are geophysical fluid dynamics, and our numerical illustrations are with simple systems of this type. With this in mind, we elaborate our approach for a broad class of ordinary differential equations (ODEs), that includes forced-dissipative systems of the form
\be\label{Eq_fundamental}
\frac{\d y} {\d t} =A y +B(y,y) +F, \;\; y\in \mathbb{C}^N.
\ee 
Here $A$ denotes a linear $N\times N$ matrix, $B$ a quadratic nonlinearity (as in the fluid advection operator) and $F$ a constant forcing, i.e.~autonomous. Such systems with complex entries arise e.g.~as equations for the perturbed variable around a mean state, when the latter are expressed in the eigenbasis $\{\boldsymbol{e}_j\}_{j=1}^N$ of the linearization at this mean state.

We decompose the phase space into the sum of the subspace, $E_\c$, of resolved variables (``coarse-scale''),  and the subspace, $E_\s$, of unresolved variables (``small-scale''). In practice $E_\c$ is spanned by the first few eigenmodes with dominant real parts (e.g.~unstable), and $E_\s$ by the rest.  Within this framework, 
and given a cutoff dimension, $m$ (i.e.~dim($E_\c$)=$m$), we consider for systems such as \eqref{Eq_fundamental} parametric families of nonlinear parameterizations of the form
 \bea\label{Eq_Htau}
 H_{\bftau}(\xi)&= \sum_{n\geq m+1} H_n(\tau_n,\xi) \boldsymbol{e}_n, \qquad \xi \in E_\c, \\
\boldsymbol{\tau}&=(\tau_{m+1},\cdots,\tau_N), \quad \tau_n\geq 0.
\eea
 
The purpose is to dispose of parameterizations that cover situations of slaving between the resolved and unresolved variables as well as 
situations for which slaving is not expected to occur (e.g.~far from criticality), as $\boldsymbol{\tau}$ is varied. In that respect, we aim at determining a family of parameterizations that include the leading-order approximation of invariant manifolds when the system is placed near the first bifurcation value. 
The theory of approximation of invariant manifolds revisited in Part I teaches us that such a family can be produced by 
finite time-integration of auxiliary BF systems derived from  Eq.~\eqref{Eq_fundamental}; see e.g.~\eqref{Eq_BF} and \eqref{Eq_BF_introb} below. 
This gives rise to the LIA class, for which taking the limit (under appropriate non-resonance conditions) of $H_n(\tau_n,\xi) $ as $\tau_n\rightarrow \infty$ provides the leading-order approximation of the invariant manifold; see Theorems \ref{thm:h1_general} and \ref{thm:h1} below.

We propose a variational approach to deal with situations far away from  criticality. 
It consists of determining the optimal $\tau_n$-value, $\tau_n^\ast$,  by minimizing (relevant) cost functionals that 
depend on solutions from  direct numerical simulation (DNS)  but over a training interval of length comparable to a characteristic recurrence or decorrelation time of the dynamics; see Secns.~ \ref{Sect_RBC} and \ref{Sec_KS_turbulence} below for applications.

Given a solution $y(t)$ of Eq.~\eqref{Eq_fundamental} available over an interval $I_T$ of length $T$, one such cost functional on which a substantial part of this article focuses on is given by the following parameterization defect
\be\label{Eq_minQnHn}
\mathcal{Q}_n(\tau_n,T)= \overline{\big|y_n(t) -H_{n}(\tau_n;y_\c (t)) \big|^2}.
\ee
Here $\overline{(\cdot )}$ denotes the time-mean over $I_T$ while $y_n(t)$ and $y_\c(t)$ denote the projections onto the high-mode $\boldsymbol{e}_n$ and the reduced state space $E_\c$ of $y(t)$, respectively. 
Our goal is then to optimize $\mathcal{Q}_n(\tau_n,T)$ by solving for each $ m+1 \leq  n \leq N$,
\be\label{Eq_minfunct}
\underset{\tau_n}\min \;  \mathcal{Q}_n(\tau_n,T).
\ee

This procedure corresponds to minimizing the variance of the residual error per high mode in case  $y_n$ and $H_n$ are zero-mean, 
and to minimizing the residual error as measured in a least-square sense, in the general case. 

Geometrically, as shown in Sec.~\ref{Sec_mode_adaptive} below, the graph of  $H_{\bftau}$ gives rise to a manifold $\mathfrak{M}_{\bftau}$ that satisfies 
\be\label{Eq_geom_interp0}
\overline{ \textrm{dist}(y(t),\mathfrak{M}_{\bftau})^2} \leq \sum_{n=m+1}^N \mathcal{Q}_n(\tau_n,T),
\ee 
{where $\textrm{dist}(y(t),\mathfrak{M}_{\bftau})$ denotes the distance of $y(t)$ (lying on the attractor) to the manifold $\mathfrak{M}_{\bftau}$.}

\begin{figure}[hbtp]
\centering
\includegraphics[height=0.33\textwidth, width=1\textwidth]{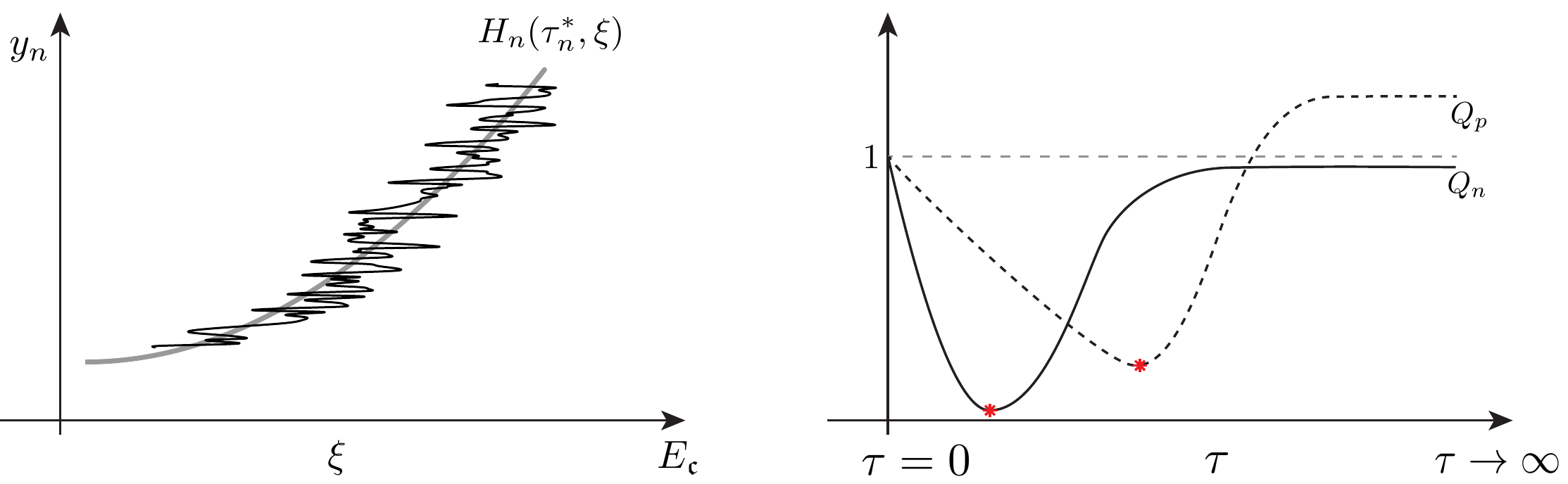}
\caption{{\footnotesize {\bf Left panel:} The optimal parameterization, $H_n(\tau_n^\ast,\xi)$, minimizing \eqref{Eq_minfunct} is shown (in gray). Here the dynamics (black curve) is transverse to it (i.e.~absence of slaving) while $H_n(\tau_n^\ast,\xi)$ provides the best (quadratic) parameterization in a least-square sense.  See Fig.~\ref{Visualization_BE} below for a concrete example in the case of a truncated Primitive Equation model due to Lorenz \cite{Lorenz80}. The parameter $\tau_n^\ast$ corresponds to the argmin of $Q_n$ (red asterisk) shown in the right panel. {\bf Right panel:}
Dependence on $\tau$ shown for two parameterization defects $Q_n$ and $Q_p$ given by \eqref{Eq_normalizedQ}, with $p, n \geq m+1$. The minimum is marked by a red asterisk. 
}}\label{Fig_intro}
\end{figure}

Thus minimizing each $\mathcal{Q}_n(\tau_n,T)$ (in the $\tau_n$-variable) is a natural idea to enforce closeness of $y(t)$ in a least-square sense to  the manifold $\mathfrak{M}_{\bftau}$.
The left panel in Fig.~\ref{Fig_intro} illustrates \eqref{Eq_geom_interp0} for the $y_n$-component: The  optimal parameterization, $H_n(\tau_n^\ast,\xi)$, minimizing \eqref{Eq_minfunct} is shown; it illustrates a situation where the dynamics is transverse to it (i.e.~absence of slaving) while $H_n(\tau_n^\ast,\xi)$ provides the best (quadratic) parameterization in a least-square sense.
 
In practice, the following {\it normalized parameterizing defect} (for the $n^{\textrm{th}}$ mode), $Q_n$, is a useful tool to compare the different parameterizations $H_{n}(\tau;\cdot)$ as $\tau$ is varied. It is defined as
\be\label{Eq_normalizedQ}
Q_n(\tau,T)=\frac{\overline{|y_n-H_{n}(\tau;y_\c )|^2}}{\overline{|y_n|^2}}.
\ee
%%%%%%%%%%%%
It provides a non-dimensional number to judge objectively of the quality of a parameterization.  If $Q_n(\tau,T)=0$ for each $n\geq m+1$, then $H_{\bftau}$ provides an exact slaving relation, and if $H_n=0$  i.e.~$H_{\bftau}\equiv 0$, corresponding to  a standard Galerkin approximation,  then $Q_n(\tau,T)=1$. Thus, the notion of (normalized) parameterizing defect allows us to bring another perspective on criticisms brought to the (approximate) inertial manifold theory \cite{heywood1993question,garcia1995time}: given a cutoff scale, if $Q_n(\tau,T)>1$ (over-parameterization) for several high modes, then a parameterization $H_{\bftau}$ may indeed lead to closure skills worse than those that would be obtained from a standard Galerkin scheme (cf.~$Q_p$ in Fig.~\ref{Fig_intro}; right). In other words, only a parameterization associated with a manifold that avoids such a situation  
is useful compared to a standard Galerkin scheme. This understanding alone is overlooked in the literature concerned with inertial manifolds and the like.  We call such a manifold a {\it parameterizing manifold (PM)}; see Definition \ref{def:PM}  for a precise characterization of a PM. 

 Minimizing the parameterization defects leads thus to an {\it optimal PM}, for the cost functionals $\mathcal{Q}_n$. 
We emphasize that each component $H_n$, of the parameterization $H_{\bftau}$ given in \eqref{Eq_Htau}, depends only  on $\tau_n$ (and not the other $\tau_p$'s for $p\neq n$), and thus
the cost functionals,  $\mathcal{Q}_n$, may be minimized independently from each other.  
%%%%%%%%%%%%%%%%%%%%%%%%%%%%%
 
The parametric dependence on $\bftau$ of $H_{\bftau}$ is of practical importance.
To understand this, let us consider for a moment a parameterization, $H_n$, given as  
a homogeneous quadratic polynomial of the $m$-dimensional $\xi$-variable with unknown coefficients (not depending on $\tau_n$).
To learn these coefficients via a standard regression would lead to $m(m-1)/2$ coefficients to estimate.  Instead, adopting the parametric formulation given in \eqref{Eq_minQnHn}, only the parameter $\tau$ needs to be learned (per high-mode) in case each coefficient of $H_n(\tau,\xi)$  is given by a function of $\tau$.
This way, we benefit from a significant reduction of the amount $N_T$ of snapshots $y(t_k)$ required from numerical integration of Eq.~\eqref{Eq_fundamental} to obtain robust parameterizations (in a statistical sense). Roughly speaking, if $N_T$ is smaller or comparable to $m(m-1)/2$, then learning the unknown (and arbitrary) coefficients of a homogeneous quadratic parameterization (not given under the parametric form \eqref{Eq_minQnHn}) is either undetermined or not robust statistically.

Explicit formulas for the coefficients of $H_n(\tau,\xi)$ are derived in Secns.~\ref{Sect_PM_with_forcing} and \ref{Sec_FMTtau} below. 
These formulas are dynamically-based in the sense that these coefficients 
 involve structural elements of the right-hand side (RHS) of Eq.~\eqref{Eq_fundamental} such as the eigenvalues $\beta_j$ of $A$, projections onto the $n^{\textrm{th}}$ high-mode of nonlinear interactions $B_{ij}^n$ between pairs of low eigenmodes $(\boldsymbol{e}_i, \boldsymbol{e}_j)$ of $A$ ($1 \leq i,j \leq m$), as well as possible nonlinear interactions between these modes and the forcing term.
   
For instance, for the LIA class, the coefficients  of the $H_n(\tau,\xi)$'s monomials are given by $D_{ij}^n (\tau)B_{ij}^n$ with 
 \bea \label{Eq_Dijn_intro}
&D_{ij}^n (\tau)=\frac{1-e^{-\tau \delta_{ij}^n}}{\delta_{ij}^n} , \qquad \tau>0,\\
& \textrm{with } \; \delta_{ij}^n=\beta_i+\beta_j -\beta_n.
 \eea
 We emphasize that at an heuristic level, the coefficient $D_{ij}^n (\tau)$ allows for balancing the denominator $\delta_{ij}^n$ by the numerator $1-e^{-\tau \delta_{ij}^n}$ when the former is small. Such compensating $\tau$-factors are in general absent from parameterizations built from invariant manifold or (approximate) inertial manifolds techniques.     
 
From the approximation theory of invariant manifolds revisited in Part I below,  one notes that $D_{ij}^n (\tau)$ is equal to $1/\delta_{ij}^n$ in the case of standard  approximation formulas of  invariant manifolds (Theorem \ref{thm:h1}), corresponding thus to the asymptotic case  $\tau\rightarrow \infty$ if $\delta_{ij}^n>0$. When adopting these approximation formulas outside their domain of applicability (i.e.~not for approximating an underlying invariant manifold), it corresponds typically to small   $\delta_{ij}^n$'s which without    
the compensating $\tau$-factors lead to an over-parameterization and an incorrect reproduction of the backscatter transfer of energy to the large scales.  This problem is typically encountered in invariant manifold approximation when small spectral gaps are present, regardless of whether the solution dynamics is simple or complicated; see {\hl the Supplementary Material} for a simple example. It turns out that, to seek for an optimal backward integration time $\tau$ actually helps alleviate this problem by introducing numerators balancing the small denominators present in standard LIA parameterizations such as provided by Theorem \ref{thm:h1} below.

At the same time, $\tau=0$ implies $D_{ij}^n (\tau)=0$, which corresponds to the  null parameterization, namely to a Galerkin approximation of dimension $m$.  Thus, minimizing the $Q_n$'s gives rise to an intermediate (and optimized) parameterization compared to a  Galerkin approximation ($H_n=0$) or an invariant manifold approximation ($Q_n=0$).

The right panel in Fig.~\ref{Fig_intro} shows a typical dependence on $\tau$ of the $Q_n$'s defined in \eqref{Eq_normalizedQ} for the LIA class. Similar dependences hold for the QSA class.  On a practical ground, the minimization problem \eqref{Eq_minfunct} is greatly facilitated by exploiting the explicit formulas of Secns.~\ref{Sect_PM_with_forcing} and \ref{Sec_FMTtau}.  An efficient minimization can be indeed operated by application of a simple gradient-descent algorithm in the real variable $\tau$, when the appropriate moments up to fourth order have been estimated;  see Appendix \ref{Sect_gradient_descent}.

We emphasize that the parameterization formulas of the LIA or QSA classes can be derived for dissipative nonlinear partial differential equations (PDEs) as well; see Sec.~\ref{Sec_KS_turbulence} below. The LIA class as rooted in the backward-forward method mentioned above was initially introduced for PDEs (possibly driven by a multiplicative linear noise)  in \cite[Chap.~4]{CLW15_vol2} and was applied to the closure 
of a stochastic Burgers equation in \cite[Chaps.~6 \& 7]{CLW15_vol2} and to optimal control in \cite{CL15}. 
The main novelty compared to these previous works is  the idea of optimizing per high mode the backward integration time, $\tau_n$, by minimization of the parameterization defect $Q_n$. Here, we also restrict ourselves to quadratic parameterizations that we prefer to optimize instead of computing higher-order terms that although being potentially useful make more cumbersome the numerical integration of the corresponding closure systems by adding too many extra terms in the RHS of the latter.

The justification of the variational approach proposed in this article relies on the ergodic 
theory of dissipative deterministic dynamical systems.  In that respect, given the flow $T_t$ associated with Eq.~\eqref{Eq_fundamental}, we assume in Part II of this article 
that $T_t$ possesses an invariant probability measure $\mu$, which is {\it physically relevant} \cite{eckmann_ruelle,collet2007concepts},  in the sense that time-average equals to ensemble average for trajectories emanating from Lebesgue almost every initial condition. More precisely, we say that the invariant measure, $\mu$, is physical if the following property holds for $y$ in a positive Lebesgue measure  set $B(\mu)$ (of $\mathbb{C}^N$) and for every continuous observable $\varphi:\mathbb{C}^N\rightarrow \mathbb{C}$
\be\label{Eq_phys_rev0}
\underset{T\rightarrow \infty}\lim \frac{1}{T} \int_0^{T}  \varphi(T_t( y)) \d t =\int \varphi( y) \d \mu ( y). 
\ee
This property assures that meaningful averages can be calculated and the statistics of the dynamical system can be investigated by the asymptotic distribution of orbits starting from Lebesgue almost every initial condition in e.g.~the basin of attraction $B(\mu)$ of the statistical equilibrium, $\mu$. 

It can be proven for e.g.~Anosov flows \cite{bowen1975ergodic}, partially hyperbolic systems \cite{alves2000srb}, Lorenz-like flows \cite{bonatti2000lorenz}, 
and observed experimentally for many others \cite{eckmann_ruelle,gallavotti1995dynamical,csg11,CGN18b} that a
common feature of (dissipative) chaotic systems is the transformation (under the action of the flow) of the initial
Lebesgue measure  into a probability measure with finer and finer scales, reaching asymptotically an invariant measure $\mu$ of Sinai-Ruelle-Bowen (SRB) type. This measure is singular with respect to the Lebesgue measure,  is supported by the local unstable manifolds contained in the global attractor or the non-wandering set \cite[Definition 6.14]{collet2007concepts}, and  if it has no zero Lyapunov exponents it satisfies \eqref{Eq_phys_rev0} \cite{young2002srb}.  This latter property is often referred to as the {\it chaotic hypothesis} that, roughly speaking, expresses an extension of the ergodic hypothesis to non-Hamiltonian systems \cite{gallavotti1995dynamical}.

At the core of our analysis, is the  disintegration $\mu_\xi$ of statistical equilibrium $\mu$ with respect to the resolved variable $\xi$ in $E_\c$; see~\cite[Sec.~3]{Chekroun2017a}. In our case,  the probability measure $\mu_\xi$ gives the conditional probability of the unresolved variables (in $E_\s$), contingent upon the value taken by the resolved variable $\xi$. Denoting by $y_\s(t)$ the high-mode projection of $y(t)$, Theorem \ref{Thm_variational-pb} below shows, under a natural boundedness assumption on the 2nd-order moments, that  
the optimal PM that minimizes the defect 
\be\label{QT_intro}
\mathcal{Q}_T(\Psi)=\overline{\norm{y_{\s}(t) -\Psi(y_{\c}(t))}^2}, 
\ee 
with $\Psi$ denoting a square-integrable mapping\footnote{With respect to the probability measure $\mathfrak{m}$ obtained as a projection of $\mu$ onto $E_\c$.} from $E_\c$ to $E_\s$,
is given, when $T \rightarrow \infty$, by   
\be\label{Def_h20}
 \Psi^\ast(\xi)=\int_{E_\s} \zeta \d \mu_{\xi}(\zeta), \qquad \xi \in E_\c.
\ee

This formula shows that the optimal PM corresponds actually to the manifold that maps to each resolved variable $\xi$ in $E_\c$,  the averaged value of the unresolved variable $\zeta$ in $E_\s$ as distributed according to the conditional probability measure $\mu_\xi$.  
In other words, the optimal PM provides the best manifold (in a least-square sense) that averages out the fluctuations of the unresolved variable. The closure system that consists of approximating the unresolved variables by this optimal parameterization provides then, when the high-mode to high-mode interactions are small,  the conditional expectation of the original system;  see Theorem \ref{Thm_variational-pb2} below. The latter provides the best vector field  of the reduced state space for which the effects of the unresolved variables are averaged out 
with respect to the probability measure $\mu_\xi$ on the space of unresolved variables, itself conditioned on the resolved variables. For slow-fast systems, in the limit of infinite time-scale separation, it is well-known that the slow dynamics is approximated (on bounded time scales) by the conditional expectation of the multiscale system \cite{kifer2001averaging,kifer2005another,pavliotis2008multiscale} and that slow trajectories may be obtained through a variational principle \cite{lebiedz2011variational}. Nevertheless, the conditional expectation may be useful to approximate other global features of the multiscale dynamics when time-scale separation is lacking. For instance, the low-frequency variability dynamics may be well approximated for chaotic systems that do not exhibit distinguished fast variables but rather episodic bursts of fast oscillations punctuated by slow oscillations for each variable; see \cite{CLM16_Lorenz9D} and Sec.~\ref{Sec_L9D_resurect} below.

The optimal PM, $\Psi^\ast$, comes with a normalized  parameterization defect defined by $Q_T(\Psi^\ast)=\mathcal{Q}_T(\Psi^\ast)/\overline{\|y_\s(t)\|^2}$, which satisfies necessarily (Theorem \ref{Thm_variational-pb})
\be
0\leq \underset{T\rightarrow \infty}\lim Q_T(\Psi^\ast) \leq 1.
\ee
This variational view on the  parameterization problem of the unresolved variables removes any sort of ambiguity that has surrounded the notion of (approximate) inertial manifold in the past. Indeed, within this paradigm shift,  given an ergodic invariant measure $\mu$ and a reduced dimension $m$, the optimal PM may have a parameterization defect very close to 1 and thus the best possible nonlinear parameterization one could ever imagine may not a priori do much better than a classical Galerkin approximation, and sometimes even worse.  To the opposite, the smaller $Q_T(\Psi^\ast)$ is (for $T$ large), the better the parameterization. All sort of nuances are actually admissible, even when the parameterization defect is just below unity; see \cite{CLM16_Lorenz9D}.

The parameterization defect analysis will be often completed by the evaluation of the {\it correlation parameterization}, $c(t)$ (see \eqref{Eq_corr_param}), that provides a measure of collinearity  between the parameterized variable $\Psi(y_{\c}(t))$ and the unresolved variable $y_{\s}(t)$, as time evolves. It allows thus for measuring how far from a slaving situation a given PM is on a more geometrical ground than with $Q_T$ (Sec.~\ref{Sec_variational_approach_Q}). As we will see in applications, the parameterization correlation  allows us, once an optimal PM has been determined, to select the dimension $m$ of the reduced state space according to the following criterium: $m$ should correspond to the lowest dimension of $E_\c$ for which the probability distribution function (PDF) of the corresponding {\it parameterization angle}, $\alpha(t)=\arccos (c(t)),
$ is the most skewed towards zero and the mode (i.e.~the value that appears most often) of this PDF is the closest to zero. The basic idea is that one should not only parameterize properly the statistical effects of the neglected scales but also avoid to lose their phase relationships with the retained scales \cite{mccomb2001conditional}. This is particularly important to derive closures that respect a certain phase coherence between the resolved and unresolved scales.

Although finite-time error estimates are easily accessible when PMs are used to derive   surrogate low-dimensional systems  in view of the optimal control of dissipative nonlinear PDEs (see e.g~\cite[Theorem 1 \& Corollary 2]{CL15}), error estimates that relate the parameterization defect to the ability of reproducing the original dynamics's long term statistics by a surrogate system  are difficult to produce for uncontrolled deterministic systems, in particular for chaotic regimes, due to the singular nature (with respect to the Lebesgue measure) of the invariant measure $\mu$ satisfying \eqref{Eq_phys_rev0}. In the stochastic realm, this invariant measure becomes smooth for a broad class of systems and the tools of stochastic analysis make the obtention of such estimates  more amenable albeit non trivial; see  \cite{chekroun2019grisanov}. Nevertheless, as discussed above, considerations from ergodic theory and conditional expectations are already insightful for the deterministic systems dealt with in this article.  They allow us to envision the addition of memory effects (non-Markovian terms) and/or stochastic parameterizations when a PM alone is not sufficient to provide an accurate enough closure. The addition of such ingredients are beyond the scope of this article, but are outlined in the Concluding Remarks (Sec.~\ref{Sec_concluding_rmk}) as a natural direction to extend the present work. The latter sets up a framework for determining, 
via dynamically-based formulas to optimize, approximations of the Markovian terms arising in the Mori-Zwanzig formalism \cite{Chorin_Hald-book,gottwald_crommelin_franzke_2017}; this formalism providing a conceptual framework to study the reduction of nonlinear autonomous systems.

The structure of this article is as follows. In Section \ref{Sect_UM_approx} we revisit the approximation formulas of invariant manifolds for equilibria. The leading-order approximation $h_k$ to these manifolds is obtained as the pullback limit of the high-mode part of the solution to an auxiliary backward-forward system (Theorem \ref{thm:h1_general}) and explicit formulas of $h_k$  are derived (Theorem \ref{thm:h1}). The resulting invariant manifold approximation formulas are applied to an El Ni\~no-Southern Oscillation ODE model in the  Supplementary Material, in the case of a  subcritical Hopf bifurcation. In Section \ref{Sect_PM_reduction}, we introduce the measure-theoretic framework in which our variational approach 
is formulated. Theorem \ref{Thm_variational-pb} characterizes the minimizers (optimal PMs) of the parameterization defect, and Theorem \ref{Thm_variational-pb2} shows that optimal PMs relate naturally to conditional expectations.  As a first application,  in Section \ref{Sec_L9D_resurect} the closure results of \cite{CLM16_Lorenz9D} concerning the low-order model atmospheric Primitive Equations of \cite{Lorenz80}, are enlightened by new insights introduced in this article.  Building upon the backward-forward systems of Section \ref{Sect_UM_approx}, we derive in  Section \ref{Sect_PM_formulas} parametric formulas of dynamically-based parameterizations aimed at being optimized.

Applications to the closure of a low-order model of Rayleigh-B\'enard convection are then discussed in Sec.~\ref{Sect_RBC}, for which a period-doubling regime and a chaotic regime are analyzed. 
In Section \ref{Sec_KS_turbulence} the approach is finally illustrated --- in the context of the Kuramoto-Sivashinsky turbulence --- as providing efficient closures without slaving and for cutoff scales placed well within the inertial range, keeping only the unstable modes in the reduced state space. It is shown that the variational approach introduced in this article allows for fixing the excessive backscatter transfer of energy to the low modes encountered by standard parameterizations. We conclude in Section \ref{Sec_concluding_rmk} by outlining future directions of research.

\newpage
\needspace{1\baselineskip}
\vspace{2cm}
%{\part*{\large \centerline{\bf Part I: Invariant manifold reduction revisited}}}
{\large \centerline{\bf Part I: Invariant manifold reduction revisited}}
%\vspace{1cm}

\section{Approximation formulas for invariant manifolds of nonlinear ODEs}  \label{Sect_UM_approx}
%%%%%%%%%%%%%%%%%%%%%%%%%%%%%%%%%%%%%%%%%%%%%%%%%%%%%
\subsection{Local invariant manifolds for equilibria: Validity and motivations for other parameterizations} \label{Sec_loc_invman}
Our framework takes place with autonomous systems of ordinary differential equations (ODEs) in $\mathbb{R}^N$ of the form: 
\be \label{Eq_ODEs_org}
\frac{\d Y}{\d t} = F(Y), 
\ee
for which  the vector field $F$ is assumed to be sufficiently smooth in the state variable $Y$.

Invariant manifold theory allows for the rigorous derivation of low-dimensional surrogate systems from which not only the system's qualitative behavior near e.g.~a steady state is preserved, but also quantitative features of the nonlinear dynamics are reasonably well approximated such as the solution's amplitude or possible dominant periods. This aspect of the theory is recalled below in {\hl the Supplementary Material}, for the unfamiliar reader.  

To set the ideas, assuming that $\overline{Y}$ is a steady state of the system \eqref{Eq_ODEs_org},
we rewrite the system \eqref{Eq_ODEs_org} in terms of the perturbed variable, $y = Y - \overline{Y}$, namely 
\bea \label{Eq_ODEs}
& \frac{\d y}{\d t} = A y + G(y), \text{ with } \\
& A = DF(\overline{Y}), \\
& G(y)= F(y + \overline{Y}) - A y,
\eea
where $DF(x)$ denotes the Jacobian matrix of $F$ at $x$.

From its definition, the nonlinear mapping, $G\colon \mathbb{R}^N \rightarrow \mathbb{R}^N$, satisfies
\be \label{Eq_G_properties}
G(0) = 0, \qquad \text{ and } \qquad  D G(0) = 0. 
\ee
As a consequence, $G(y)$ admits the following expansion for $y$ near the origin:
\begin{equation} \label{G Taylor}
G(y)= G_k(\underbrace{y, \cdots, y}_{k \text{ times}}) + O(\|y\|^{k+1}),
\end{equation}
where 
\bea \label{k-linear def}
G_k \colon \underbrace{ \mathbb{R}^N \times \cdots \times \mathbb{R}^N}_{k
\text{ times}} \rightarrow \mathbb{R}^N
\eea
denotes a homogenous polynomial of order $k\ge 2$. {\hl That is, $G_k$ is the homogeneous part of lowest degree.} Sometimes, $G_k(y)$ will be used as a compact notation for  $G_k(y,\,\cdots \, , y)$.

The spectrum of $A$ is denoted by $\sigma(A)$, i.e.
\be \label{eigen_A}
\sigma(A)=\{\beta_j  \in \mathbb{C} : j = 1,\cdots, N\},
\ee 
where the $\beta_j$s denote the eigenvalues of $A$ for which we have accounted for their algebraic multiplicity in the sense that if 
$\lambda$ is a root of multiplicity $p$ of the characteristic polynomial $\chi_A$, then e.g.
$\beta_1=\lambda,\cdots,\beta_p=\lambda$. The corresponding generalized eigenvectors are denoted by 
\be \label{eigenmode_A}
\{\boldsymbol{e}_j  \in \mathbb{C}^N : j = 1,\cdots, N\}.
\ee
%%%%%%%%%%%%%%%%%%%%%%%%%%%%%%%%%%%%%%%%%%%%%%%%%%%%%
The index in  \eqref{eigen_A} also accounts for an arrangement of the eigenvalues in   
lexicographical order, that is the eigenvalues are ordered so that their real parts decrease as the index increases, and for eigenvalues with the same real parts, they are arranged so that the imaginary parts decrease.

Taking into account this ordering, grouping the first $m$ eigenvalues of $A$,  and assuming
\be\label{cutoff_cond}
\mathrm{Re}(\beta_{m})\neq \mathrm{Re}(\beta_{m+1}), 
\ee
the spectrum of $A$ is decomposed as follows
\be
\sigma(A) =\sigma_{\c}(A)\cup \sigma_{\s}(A),
\ee
where
\be \label{Eq_unstablepart}
 \sigma_{\c}(A)=\{\beta_j, \; j = 1,\cdots, m\},
 \ee
and 
\be\label{Eq_stablepart}
 \sigma_{\s}(A)=\{\beta_j, \;   j = m+1,\cdots, N\}.
 \ee
Note that due to \eqref{cutoff_cond} and the aforementioned lexicographical order, we have
\be \label{Eq_gap}
\mathrm{Re}(\beta_m) >  \mathrm{Re}(\beta_{m+1}). 
\ee

This spectral decomposition implies a natural decomposition of $\mathbb{C}^N$:
\be \label{Eq_decomp_CN}
\mathbb{C}^N = E_{\c}  \oplus E_{\s},
\ee
in terms of the generalized eigenspaces  
\bea \label{eq:subspaces}
& E_{\c} = \mathrm{span}\{\boldsymbol{e}_j : j = 1, \cdots, m\}, \\
&  E_{\s} = \mathrm{span}\{\boldsymbol{e}_j: j = m+1, \cdots, N\}.
\eea
%%%
This spectral decomposition of $\mathbb{C}^N$ along with the corresponding canonical projectors $\Pi_{\c}$ and $\Pi_{\s}$ onto $E_{\c}$ and $E_{\s}$, respectively, are at the core of our dimension reduction of Eq.~\eqref{Eq_ODEs}.

The  theory of local invariant manifolds for equilibria says that the simple condition \eqref{Eq_gap} combined with the tangency condition \eqref{Eq_G_properties} about the nonlinear term $G$ ensure the existence of a local $m$-dimensional invariant manifold, namely  
a manifold obtained as the local graph over an open ball $ \mathfrak{B}$ in $E_{\c}$ centered at the origin, that is
 \be \label{critical manifolds}
\mathfrak{M}= \left \{ \xi+ h (\xi) : \xi \in \mathfrak{B} \subset E_{\c} \right \},
\ee
where $ h\colon  E_{\c} \rightarrow E_{\s}$ is a $C^1$-smooth manifold function such that $h(0) = 0$ and $D h(0) = 0$,
for which the following property holds:
\bi
\item[(i)] any solution $y(t)$ of Eq.~\eqref{Eq_ODEs} such that $y(t_0)$ belongs to $\mathfrak{M}$ for some $t_0$, stays on   $\mathfrak{M}$ over an interval of time $[t_0,t_0+\alpha)$, $\alpha>0$, i.e.
\be\label{Eq_inv_local}
y(t)=y_\c(t)+h(y_\c(t)), \; t\in [t_0,t_0+\alpha),
\ee
where $y_\c(t)$ denotes the projection of $y(t)$ onto the subspace $E_{\c}$. 
\ei
%%%%%
Additionally, if $\mathrm{Re}(\beta_{m+1})<0$ and $\mathrm{Re}(\beta_{m})\geq 0$, then the local invariant manifold is the so-called local center-unstable manifold and the following property holds
\bi
\item[(ii)] If there exists a trajectory $t\mapsto y(t)$ such that $y_\c(t)$ belongs to $ \mathfrak{B}$ for all $-\infty<t<\infty$, then the trajectory must lie on $\mathfrak{M}$.
\ei
Property (ii) implies that an invariant set $\Sigma$ of any type, e.g., equilibria, periodic orbits, invariant tori,
must lie in $\mathfrak{M}$ if its projection onto $E_\c$ is contained in  $\mathfrak{B}$, i.e.~if $\Pi_\c \Sigma \subset \mathfrak{B}$. Property \eqref{Eq_inv_local} holds then globally in time for the solutions that composed such invariant sets, and thus the knowledge of the $m$-dimensional variable, $y_\c(t)$, is sufficient to entirely determine any solution $y(t)$ that belongs to such an invariant set. Furthermore, $y_\c(t)$ is obtained as the solution of the following reduced $m$-dimensional problem
\be\label{Eq_reduced_absract}
 \frac{\d x}{\d t}= \Pi_\c A  x + \Pi_\c G(x + h(x)), \qquad x(0)=y_\c(0) \in \mathfrak{B},
\ee
which in turn characterizes the solution $y(t)$ in $\Sigma$, since the slaving relationship $y_\s(t)=h(y_\c(t))$ holds for any solution $y(t)$ that belongs to an invariant set $\Sigma$ for which $\Pi_\c \Sigma \subset \mathfrak{B}$.

More generally,  property $(i)$ allows for $y_\c(t)$ to leave the neighborhood $\mathfrak{B}$ for some time instance, $t$, and thus to violate the parameterization 
\eqref{Eq_inv_local} for $y(t)$, but does not exclude to have \eqref{Eq_inv_local} to hold again over another interval $[t_1,t_1+\alpha_1)$ as soon as $y(t_1)$  belongs to $\mathfrak{M}$.

Regarding the neighborhood $\mathfrak{B}$, the theory shows that it shrinks as the spectral gap, 
\bes
\gamma_m = \mathrm{Re}(\beta_{m})-\mathrm{Re}(\beta_{m+1}), 
\ees
gets small and the nonlinear term $G$ deviates quickly from the tangency condition as one moves away from the origin, leaving possible an (exact) parameterization only for solutions with sufficiently small amplitude.   Indeed, the existence of such a (local) exact parameterization or say in other words, of a local $m$-dimensional invariant manifold is subject to the following {\it spectral gap condition}:
\be\label{Eq_spectral_gapcond}
\gamma_m \geq C \mbox{Lip}(G\vert_{\mathcal{V}}),
\ee
where $\mbox{Lip}(G\vert_{\mathcal{V}})$ denotes the Lipschitz constant of the nonlinearity $G$, restricted to a neighborhood  $\mathcal{V}$ of the origin in $\mathbb{C}^N$ such that $\mathcal{V}\cap E_\c=\mathfrak{B}$, and $C >0$ is typically independent on $\mathcal{V}$. 
Due to the tangency condition \eqref{Eq_G_properties}, the condition \eqref{Eq_spectral_gapcond} always holds once $\mathcal{V}$ (and thus $\mathfrak{B}$)  is chosen sufficiently small.  The theory of local invariant manifolds makes thus sense if solutions with sufficiently small amplitudes lie in the neighborhood $\mathcal{V}$. This situation is encountered for many bifurcations, near criticality for which the system's linear part has modes that become unstable, although a condition on the asymptotic stability of the origin is often required to have a local attractor  that continuously unfolds from the origin as the bifurcation parameter is varied \cite[Theorem 6.1]{MW05}.        
In the context of e.g.~nonlinear oscillations that bifurcate from a steady state, local invariant manifolds provide exact parameterizations\footnote{As provided for instance by a center manifold or the unstable manifold of the origin.} of stable limit cycles near criticality in the case of a supercritical Hopf bifurcation, whereas it is the parameterization of the unstable limit cycle that emerges continuously from the steady state that is guaranteed to be exact, at least sufficiently close to criticality in the case of a subcritical Hopf bifurcation. {\hl In the Supplementary Material}, we show that the approximation formulas of Sec.~\ref{Sec_Leading_approx}, allow for approximating not only the unstable ``inner" unstable limit cycle but also the ``outer" stable limit cycle arising in an El Ni\~no-Southern Oscillation (ENSO) model via subcritical Hopf bifurcation.         

In any event, local invariant manifolds by their local nature, although useful in many applications do not allow for an  efficient dimension reduction of arbitrary or at least generic solutions. Attempts to extend the theory to a more 
global setting, have failed dramatically to systematically provide nonlinear parameterizations of type \eqref{Eq_inv_local} for a broader set of solutions, since, in general, the same type of spectral gap condition as \eqref{Eq_spectral_gapcond} is also encountered in such an endeavor. For instance, the theory of inertial manifolds is known to be conditioned on spectral gap conditions such as given by \eqref{Eq_spectral_gapcond} for which the Lipschitz constant is global or taken over a neighborhood $\mathcal{V}$ that contains the (projection onto $E_\c$ of the) global attractor.

Part II proposes a new framework to provide manifolds which are no-longer locally invariant --- and thus not subject to a spectral gap condition --- but still provide meaningful nonlinear parameterizations of nonlinear dynamics; these manifolds being called {\it parameterizing manifolds (PMs)}. Nevertheless, the calculation of PMs departs from the theory 
of approximation of local invariant manifolds which we revisit in the next section, before presenting the main, new, analytical ingredients in Sec.~\ref{Sect_PM_formulas}. 

The material presented in Sec.~\ref{Sec_Leading_approx} below will serve to derive (approximate) parameterizations for perturbed variable taken 
with respect to a mean state $\overline{Y}$, instead of a steady state; see Sec.~\ref{Sect_PM_with_forcing}. 
To set the ideas, we consider $F(Y)$ to be given by $L Y + B(Y,Y)$ with $L$ linear, and $B$ a quadratic homogeneous polynomial and symmetric, $B(X,Y)=B(Y,X)$. The equation for the perturbed variable $y$ then becomes
\be\label{Eq_pert_variable1}
\frac{\d y}{\d t} = (Ly + 2 B(y,\overline{Y})) +B(y,y) +B(\overline{Y},\overline{Y}),
\ee 
which adopting the notations of Eq.~\eqref{Eq_ODEs}, corresponds to 
$A=Ly + 2 B(y,\overline{Y})$ and $G(y)=B(y,y) + {\hl L \overline{Y}} + B(\overline{Y},\overline{Y})$.
Since $\overline{Y}$ is no longer a steady state, $G(0)\neq 0$,  and ${\hl L \overline{Y}} + B(\overline{Y},\overline{Y})$ is a time-independent forcing term. Thus the standard local invariant manifold theory for equilibria cannot be applied. 

Nevertheless, as shown in Sec.~\ref{Sect_PM_formulas} below, the theory underlying the derivation of approximation formulas for invariant manifolds is still relevant for their appropriate modification in view of providing approximate parameterizations in presence of forcing, once a  good representation of these formulas is adopted; see Theorem \ref{thm:h1_general} below for the representation of these approximation formulas (see \eqref{PB_rep}), and 
Sec.~\ref{Sect_PM_with_forcing} for the modified parameterizations in presence of forcing.

\subsection{Leading-order approximation of invariant manifolds} \label{Sec_Leading_approx}
This section is devoted to the derivation of analytic formulas for the approximation of the (local) invariant manifold function $h$ in \eqref{critical manifolds}. As shown below these formulas are easily obtained by relying only on the invariance property of  $\mathfrak{M}$, responsible for the {\it invariance equation} to be satisfied by $h$. We recall first the derivation of this fundamental equation; see also \cite[pp.~169-171]{Hen81} and \cite[VII.~A.~1]{crawford1991introduction}. {\hl For the existence of the invariant/center manifolds for ODEs, we refer to \cite{Vanderbauwhede89}.}

In that respect, note first that by applying respectively the projectors $\Pi_{\c}$ and $\Pi_{\s}$ on both sides of  Eq.~\eqref{Eq_ODEs} and by using that $A$ leaves  invariant the eigensubspaces $E_{\c}$ and $E_{\s}$, we obtain that Eq.~\eqref{Eq_ODEs} can be split as follows
\begin{subequations}
\begin{align}
\frac{\d y_{\c}}{\d t} = A_{\c} y_{\c} + \Pi_{\c}G( y_{\c} + y_{\s}), \label{eq:xc} \\
\frac{\d y_{\s}}{\d t} = A_{\s} y_{\s} + \Pi_{\s}G( y_{\c} + y_{\s}), \label{eq:xs}
\end{align}
\end{subequations}
with 
\be\label{Ac_As_def}
y_{\c} = \Pi_{\c} y \in E_{\c},\;  y_{\s} = \Pi_{\s} y \in E_{\s}, \; A_{\c} = \Pi_{\c} A\; \mbox{ and } A_{\s} = \Pi_{\s} A. 
\ee

Since $\mathfrak{M}$ is locally invariant, any solution $y(t)$ of Eq.~\eqref{Eq_ODEs} with initial datum on $\mathfrak{M}$ stays on $\mathfrak{M}$ as long as $y_{\c}(t)$  stays in $\mathcal{B}$ (where $\mathcal{B}$ is given in \eqref{critical manifolds}), i.e.
\be
y(t)=y_{\c}(t) + h(y_{\c}(t)),
\ee
provided  that $y_{\c}(t)$ lies in $\mathcal{B}$; see \eqref{Eq_inv_local}.

This implies, as long as $y_{\c}(t)$ belongs to $\mathcal{B}$, that $y_{\s}(t)=h(y_{\c}(t))$, which, when substituted into Eq.~\eqref{eq:xs}
gives
\be
\frac{\d h(y_{\c})}{\d t}  = A_{\s} h(y_{\c}) + \Pi_{\s}G( y_{\c} + h(y_{\c})). \label{eq:xs2}
\ee
On the other hand since $h$ is differentiable, we have by using Eq.~\eqref{eq:xc}, 
\be \label{eq:xs3}
\frac{\d h(y_{\c})}{\d t}= D h(y_{\c})\frac{\d y_{\c}}{\d t}  =D h(y_{\c}) [A_{\c} y_{\c} + \Pi_{\c}G( y_{\c} + h(y_{\c}))].
\ee

Then  \eqref{eq:xs2} and \eqref{eq:xs3} allow us to conclude that as long as  $y_{\c}(t)$  belongs to $\mathcal{B}$, $h$ evaluated along the corresponding ``segment'' of trajectory satisfies
\bea \label{eq:invariance0}
D  h(y_{\c}(t)) [A_{\c} y_{\c}(t) + \Pi_{\c} G(y_\c(t) + & h(y_{\c}(t))] -   A_{\s} h(y_{\c}(t))\\ 
&= \Pi_{\s} G(y_\c(t) + h(y_{\c}(t))),
\eea
which can be recast into the aforementioned  {\it invariance equation} to be satisfied by $h$, namely
\be \label{eq:invariance}
D  h(\xi) [A_{\c} \xi + \Pi_{\c} G(\xi + h(\xi))] -  A_{\s} h(\xi) = \Pi_{\s} G(\xi + h(\xi)),\; \xi\in  \mathcal{B}.
\ee
This functional equation is a nonlinear system of  first order PDEs that cannot be solved in closed form except in special cases.   However, one can solve Eq.~\eqref{eq:invariance} approximately by representing $h(\xi)$ as a formal power series. The solution is thus sought in terms of Taylor expansion in the $\xi$-variable and various numerical  techniques \textemdash\, based, e.g., on the resolution of the multilinear Sylvester equations associated with the invariance equation  \textemdash\,  have been proposed in the literature to find the corresponding coefficients  \cite{BK98,EvP04}. Once a power series approximation has been found, {\it a posteriori} error estimates can be checked by applying for instance \cite[Theorem 3, p. 5]{Car81}\footnote{According to this theorem, a candidate to a (truncated) Taylor expansion has to be first determined, and then it has to be checked to satisfy the invariance equation up to some order to ensure to be a genuine Taylor approximation; see also \cite[Thm.~6.2.3]{Hen81}.}.

For a broad class of systems, the leading-order approximation of $h$ can be efficiently and analytically calculated. It consists of dropping in Eq.~\eqref{eq:invariance} the terms involving nonlinear dependence on $h$. This operation leads to the following equation for the corresponding leading-order approximation $h_k$ (see, e.g., \cite{CLW15_vol1,Hen81}):
\be \label{h1_eqn}
D h_k (\xi) A_{\c} \xi  - A_{\s} h_k(\xi) =  \Pi_{\s} G_k(\xi),
\ee
where $G_k$ is the leading-order term in the Taylor expansion of $G$ about the origin; cf.~Eq.~\eqref{G Taylor}.

Easily checkable conditions on the eigenvalues of $A$, allows then for guaranteeing an analytic solution to Eq.~\eqref{h1_eqn}.  For instance, in the case $A$ is self-adjoint, it simply requires certain {\it cross non-resonance conditions} to be satisfied as stated in Theorem~\ref{thm:h1} below. Namely, for any given set of resolved modes for which their self-interactions (through the leading-order nonlinear term $G_k$) do not vanish when projected against an unresolved mode $\boldsymbol{e}_n$, it is required that some specific linear combinations of the corresponding eigenvalues dominate the eigenvalue associated with $\boldsymbol{e}_n$; see \eqref{NR} below. 

In the general case, when $A$ is not necessarily diagonal, the cross non-resonance condition is strengthened to the requirement that $\mathrm{Re}(\beta_{m+1}) < k \, \mathrm{Re}(\beta_{m})$ which ensures that 
the following Lyapunov-Perron integral $\mathfrak{I} \colon E_{\c} \rightarrow E_{\s}$,
\be \label{Eq_LP_integral}
\mathfrak{I}(\xi) = \int_{-\infty}^0 e^{-sA_{\s}} \Pi_{\s} G_k(e^{sA_{\c}}\xi) \d s,
\ee
is well defined and in fact provides a solution $h_k$ to Eq.~\eqref{h1_eqn}; see Theorem~\ref{thm:h1_general} below. This solutions provides actually the leading-order approximation of the (local) invariant manifold function $h$ if we assume furthermore that $\mathrm{Re}(\beta_{m+1}) <   \min\{ 2k \mathrm{Re}(\beta_{m}), 0\}$; see Theorem~\ref{thm:h1_general} again.

This Lyapunov-Perron integral itself possesses a flow interpretation: it is obtained as the pullback limit constructed from the solution of the following backward-forward auxiliary system 
\begin{subequations} \label{Eq_BF}
\begin{align}
& \frac{\mathrm{d} y^{(1)}_{\c}}{\d s} =  A_\c y^{(1)}_{\c}, && s \in[ -\tau, 0],    \label{BF1} \\
& \frac{\mathrm{d} y^{(1)}_{\s}}{\d s} = A_{\s} y_{\s}^{(1)}  +  \Pi_{\s} G_k\big(y^{(1)}_{\c}\big), &&  s \in [-\tau, 0], \label{BF2}\\
& \mbox{with } y^{(1)}_{\c}(s)\vert_{s=0} = \xi, \mbox{ and } y_{\s}^{(1)}(s)\vert_{s=-\tau}=0. 
\end{align}
\end{subequations}
Indeed, the solution to Eq.~\eqref{BF2} at $s=0$ is given by
\be\label{Eq_h1_tau}
h^{(1)}_\tau(\xi)=y^{(1)}_{\s}[\xi](0; -\tau) =\int_{-\tau}^0 e^{-sA_{\s}} \Pi_{\s} G_k(e^{sA_{\c}}\xi) \d s,
\ee
and taking the limit formally in \eqref{Eq_h1_tau} as $\tau\rightarrow \infty$, leads to $\mathfrak{I}$ given by \eqref{Eq_LP_integral}. 

The theorem below states more precisely the relationships between Eq.~\eqref{h1_eqn}, the Lyapunov-Perron integral \eqref{Eq_LP_integral}, and the solution to the backward-forward system \eqref{Eq_BF}.

\bt \label{thm:h1_general}

Consider Eq.~\eqref{Eq_ODEs}. Let the subspaces $E_\c$ and $E_{\s}$ be given by \eqref{eq:subspaces} and let $m$ be the dimension of $E_{\c}$. Assume \eqref{Eq_gap} and furthermore that 
\be \label{Eq_spectral_cond}
\mathrm{Re}(\beta_{m+1}) < k \, \mathrm{Re}(\beta_{m}),
\ee
where $k$ denotes the leading order of the nonlinearity $G$; cf.~\eqref{G Taylor}. 

Then, the Lyapunov-Perron integral 
\be\label{Eq_LyapPerron}
\mathfrak{I}(\xi) = \int_{-\infty}^0 e^{-sA_{\s}} \Pi_{\s} G_k(e^{sA_{\c}}\xi) \d s, \quad \xi \in E_{\c}, 
\ee
is well defined and is a solution to Eq.~\eqref{h1_eqn}. Moreover, $\mathfrak{I}$ is the pullback limit of the high-mode part of the solution to the backward-forward system \eqref{Eq_BF}:   
\be\label{PB_rep}
\mathfrak{I}(\xi)  = \lim_{\tau \rightarrow \infty} y^{(1)}_{\s}[\xi](0; -\tau),
\ee
where $y^{(1)}_{\s}[\xi](0; -\tau)$ denotes the solution to Eq.~\eqref{BF2} at $s=0$. 

Finally, if we assume furthermore that 
\be\label{Eq_gap_technical}
\mathrm{Re}(\beta_{m+1}) <   \min\{ 2k \mathrm{Re}(\beta_{m}), 0\},
\ee
then $\mathfrak{I}$ provides the leading-order approximation of the invariant manifold function $h$ in the sense that  
\be \label{Eq_leading_order_goal}
\|\mathfrak{I}(\xi) - h(\xi)\|_{E_{\s}} = o(\|\xi\|^k_{E_\c}), \quad \xi \in E_{\c}.
\ee
\et

\bp
First, we outline how condition \eqref{Eq_spectral_cond} combined with the fact that $G_k$ is a homogeneous polynomial of order $k$, ensure that  the Lyapunov-Perron integral $\mathfrak{I}$ is well defined. In that respect, we note first that natural estimates about $\|e^{t A_{\s}} \Pi_{\s}\|_{L(\mathbb{C}^N)}$ and $\|e^{t A_{\c}} \Pi_{\c}\|_{L(\mathbb{C}^N)}$ hold.  

This is essentially a consequence of \eqref{Eq_gap}. 
Indeed, any choice of real constants $\eta_1$ and $\eta_2$ such that 
\be \label{Eq_etas}
\mathrm{Re}(\beta_{m}) > \eta_1 > \eta_2 >  \mathrm{Re}(\beta_{m+1}),
\ee
ensures the existence of a constant $K > 0$ (depending on $\eta_1$ and $\eta_2$) such that the following estimates hold: 
\bea \label{Eq_dichotomy}
\|e^{tA_{\c}} \Pi_{\c}\|_{L(\mathbb{C}^N)} \le K e^{\eta_1 t}, \quad \forall t \le 0, \\
\|e^{tA_{\s}} \Pi_{\s}\|_{L(\mathbb{C}^N)} \le K e^{\eta_2 t}, \quad \forall t \ge 0.
\eea
The latter inequalities resulting essentially from the fact that $\|e^{tB} \|_{L(\mathbb{C}^N)}$ is bounded for $t\geq 0$ if $\mathrm{Re} \lambda <0$ for all $\lambda$ in $\sigma(B)$.

Since $G_k$ is a homogeneous polynomial of order $k$, there exists $C>0$ such that 
\be \label{Eq_Gk_est}
\|G_k(\xi)\| \le C \|\xi\|^k, \quad \forall \xi \in E_{\c}.
\ee
Now, by using \eqref{Eq_dichotomy} and \eqref{Eq_Gk_est},  we obtain for each $s \le 0$ that 
\beas
\|e^{-sA_{\s}} \Pi_{\s} G_k(e^{sA_{\c}}\xi)\| & \le K e^{-s \eta_2} \|G_k(e^{sA_{\c}}\xi)\|  \\
 & \le  C K e^{-s \eta_2} \|e^{sA_{\c}}\xi\|^k \\
 & \le  C K^2 e^{-s (\eta_2-k\eta_1)} \|\xi\|^k.
\eeas
Assumption \eqref{Eq_spectral_cond} allows us to choose $\eta_1$ and $\eta_2$ in \eqref{Eq_etas} such that $\eta_2-k\eta_1 < 0$ which in turns leads to 
\bea  \label{Eq_I_est}
\left\|\int_{-\infty}^0 e^{-sA_{\s}} \Pi_{\s} G_k(e^{sA_{\c}}\xi) \d s \right \|  & \le \int_{-\infty}^0 \|e^{-sA_{\s}} \Pi_{\s} G_k(e^{sA_{\c}}\xi)\| \d s \\
& \le C K^2\|\xi\|^k \int_{-\infty}^0 e^{-s (\eta_2-k\eta_1)} \d s \\
& =  \frac{C K^2\|\xi\|^k}{k\eta_1 - \eta_2}, \qquad \forall \xi \in E_{\c}. 
\eea
We have thus shown that $\mathfrak{I}$ is well defined. 

We show next that $\mathfrak{I}$ satisfies Eq.~\eqref{h1_eqn}. To do so, for any $\xi$ in $E_{\c}$ we introduce the following function 
\bea
\psi \colon (-\infty, 0] & \rightarrow E_{\s} \\
& t \mapsto \mathfrak{I}(e^{tA_{\c}} \xi) = \int_{-\infty}^t e^{(t-s)A_{\s}} \Pi_{\s} G_k(e^{sA_{\c}}\xi) \d s.
\eea

On one hand, by differentiating $\psi(t) = \int_{-\infty}^t e^{(t-s)A_{\s}} \Pi_{\s} G_k(e^{sA_{\c}}\xi) \d s$, we obtain
\be
\frac{\d \psi}{\d t} = \Pi_{\s} G_k(e^{t A_{\c}}\xi) + A_{\s} \int_{-\infty}^t e^{(t-s)A_{\s}} \Pi_{\s} G_k(e^{sA_{\c}}\xi) \d s.
\ee
On the other, using that $\psi(t) =  \mathfrak{I}(e^{tA_{\c}} \xi)$, we have  
\be
\frac{\d \psi}{\d t} = D\mathfrak{I}(e^{tA_{\c}} \xi) A_{\c} e^{tA_{\c}} \xi. 
\ee
It follows then that 
\be
D\mathfrak{I}(e^{tA_{\c}} \xi) A_{\c} e^{tA_{\c}} \xi = \Pi_{\s} G_k(e^{t A_{\c}}\xi) + A_{\s} \int_{-\infty}^t e^{(t-s)A_{\s}} \Pi_{\s} G_k(e^{sA_{\c}}\xi) \d s, \qquad \forall t \le 0. 
\ee
Set $t=0$ in the above equality, we then obtain 
\bes
D\mathfrak{I}(\xi) A_{\c} \xi = \Pi_{\s} G_k(\xi) + A_{\s} \int_{-\infty}^0 e^{-s A_{\s}} \Pi_{\s} G_k(e^{sA_{\c}}\xi) \d s, \;\; \forall \xi \in E_\c,
\ees
which is equivalent to 
\bes
D\mathfrak{I}(\xi) A_{\c} \xi - A_{\s} \mathfrak{I}(\xi) = \Pi_{\s} G_k(\xi), \qquad \forall \xi \in E_{\c}. 
\ees
We have thus verified that $\mathfrak{I}$ is a solution to Eq.~\eqref{h1_eqn}. 
 
Recall from Eq.~\eqref{Eq_h1_tau} that the high-mode part of the solution to the backward-forward system \eqref{Eq_BF} is given (at $s=0$) by:
\be 
y^{(1)}_{\s}[\xi](0; -\tau) =\int_{-\tau}^0 e^{-sA_{\s}} \Pi_{\s} G_k(e^{sA_{\c}}\xi) \d s,
\ee
By using the same type of estimates as in \eqref{Eq_I_est}, it is easy to show that the limit, $\lim_{\tau \rightarrow \infty} y^{(1)}_{\s}[\xi](0; -\tau)$, exists and it is equal to $\mathfrak{I}(\xi)$. 

The leading-order approximation property stated in \eqref{Eq_leading_order_goal} under the assumption \eqref{Eq_gap_technical} is a direct consequence of the general result \cite[Corollary 7.1]{CLW15_vol1} proved for stochastic evolution equations in infinite dimension, driven by a multplicative white noise which thus applies to our finite dimensional and deterministic setting. Indeed, to apply \cite[Corollary 7.1]{CLW15_vol1}, we are only left with the checking of constants $\eta_1$ and $\eta_2$ for which  \cite[condition (7.1)]{CLW15_vol1} is verified, namely
\be
\eta_{\s} < \eta_2 < \eta_1 < \eta_{\c}, \quad \eta_2 < 2k\eta_1 < 0,
\ee
with $\eta_{\s} = \mathrm{Re}(\beta_{m+1})$ and $\eta_{\c} = \mathrm{Re}(\beta_{m})$ here. One can readily check that this condition is guaranteed under the assumptions \eqref{Eq_gap} and  \eqref{Eq_gap_technical}. Indeed, if $\mathrm{Re}(\beta_{m+1}) <   2k \mathrm{Re}(\beta_{m}) < 0$, we just need to choose 
\bes
\eta_1 = \mathrm{Re}(\beta_{m}) - \epsilon \mbox{ and  } \eta_2 = \mathrm{Re}(\beta_{m+1}) +  \epsilon,
\ees 
with sufficiently small positive $\epsilon$; and if $\mathrm{Re}(\beta_{m+1}) <  0 <  2k \mathrm{Re}(\beta_{m})$, we just need to choose $\eta_1= - \epsilon$ and $\eta_2 = \mathrm{Re}(\beta_{m+1})+  \epsilon$ with again $\epsilon$ sufficiently small. 
\ep

The next Theorem shows, under a slightly relaxed spectral condition (see \eqref{NR} below), that if the matrix $A$ is assumed to be diagonal, then even when the Lyapunov-Perron integral  
\eqref{Eq_LyapPerron} is no longer defined, a solution $h_k$ to Eq.~\eqref{h1_eqn} can still be derived and that this solution possesses even an explicit expression.

This expression consists of an expansion in terms of the eigenvectors  $\boldsymbol{e}_n$ lying in the eigenspace $E_\s$, and whose coefficients are homogeneous polynomials of order $k$ in the $\xi$-variable lying in eigenspace $E_\c$; the coefficients of these polynomials being themselves expressed in terms of ratios between the linear combinations of eigenvalues of $A$ and the  corresponding eigenmodes interactions through the leading-order nonlinear term $G_k$; see \eqref{h1_part2}.    More precisely, we have

\bt \label{thm:h1}
Consider Eq.~\eqref{Eq_ODEs}. Let the subspaces $E_\c$ and $E_{\s}$ be given by \eqref{eq:subspaces} and let $m$ be the dimension of $E_{\c}$. Assume \eqref{Eq_gap} and that the matrix $A$ is diagonal under its eigenbasis $\{\boldsymbol{e}_j \in \mathbb{C}^N : j = 1,\cdots, N\}$. We denote by $\{\boldsymbol{e}_j^\ast, j=1,\cdots,N\}$ the eigenvectors of the conjugate transpose $A^*$.

Recalling that $G_k$ denotes the leading-order homogeneous polynomial in the expansion of $G$ (see \eqref{G Taylor}), 
let us assume furthermore that the eigenvalues $\beta_j$ of $A$ satisfies the following cross non-resonance condition:
\begin{equation}  \label{NR} \tag{NR}
\begin{aligned} 
& \Forall \, (i_1, \cdots, i_k ) \in \mathcal{I}^k,  \ n \in \{ m+1, \cdots, N\},  \text{ it holds that} \\
 & \Bigl (\langle G_k(\boldsymbol{e}_{i_1}, \cdots, \boldsymbol{e}_{i_k}), \boldsymbol{e}_n^\ast \rangle \neq 0 \Bigr) \Longrightarrow  \biggl ( \sum_{j=1}^{k} \beta_{i_j} - \beta_n \neq 0 \biggr),  
\end{aligned}
\end{equation}
where $\mathcal{I}= \{1, \cdots, m\}$, and $\langle \cdot, \cdot \rangle$ denotes the inner product on $\mathbb{C}^N$ defined by 
\be \label{Eq_inner_product}
\langle a, b \rangle = \sum_{i= 1}^N a_i \overline{b_i}, \qquad a, b \in \mathbb{C}^N.
\ee

Then, a solution to Eq.~\eqref{h1_eqn} exists, and is given by
\be \label{h1_part1}
h_k(\xi)  = \sum_{n=m+1}^N h_{k,n}(\xi) \boldsymbol{e}_n,  \quad  \, \xi=(\xi_1, \cdots, \xi_m) \in E_{\c},
\ee
where $h_{k,n}(\xi)$ is a homogeneous polynomial of degree $k$ in the variables $\xi_1, \cdots$, $\xi_m$ given by
\be \label{h1_part2}
h_{k,n}(\xi) = \sum_{(i_1, \cdots, i_k )\in \mathcal{I}^k} \frac{\langle G_k(\boldsymbol{e}_{i_1}, \cdots, \boldsymbol{e}_{i_k}), \boldsymbol{e}_n^\ast \rangle}{\sum_{j = 1}^k \beta_{i_j} - \beta_n} \xi_{i_1} \cdots \xi_{i_k}.
\ee
\et

\needspace{1\baselineskip}
\br \label{rmk:h1}

\hspace*{2em}  \vspace*{-0.4em}
\bi

\item[(i)] The formulas \eqref{h1_part1}--\eqref{h1_part2} for the case of real and symmetric matrices, are known; see e.g.~\cite[Appendix A]{MW14}. The result presented in Theorem~\ref{thm:h1} extends nevertheless  these formulas to cases for which $A$ is diagonalizable in $\mathbb{C}$, allowing in particular for an arbitrary number of complex conjugate eigenpairs.  The case when the neutral/unstable modes correspond to a single complex conjugate pair has been dealt with in \cite[Appendix A]{MW14}. Even in this special case, our formulas are in contradistinction simpler than those given in \cite[Eq. (A.1.15)]{MW14}. This is due to the use of generalized eigenvectors adopted here and the method of proof of Theorem \ref{thm:h1} which relies on the calculation of spectral elements of the homological operator $\mathcal{L}_A$ naturally associated with Eq.~\eqref{h1_eqn}; see \eqref{eq:ad_A} below.

\item[(ii)] The case of eigenvalues of higher-order multiplicity is more involved.   
The presence of Jordan blocks makes indeed the derivation of general analytic formulas challenging but still possible by the method used in the derivation of the formulas \eqref{h1_part1}--\eqref{h1_part2}.  Communication about these formulas will be pursued elsewhere.

\item[(iii)] By only assuming the \eqref{NR} condition, the solution to Eq.~\eqref{h1_eqn} given by the formulas \eqref{h1_part1}--\eqref{h1_part2} is not necessarily unique.  This situation happens for instance when we have a $k$-uple $(i_1, \cdots, i_k)$ and an index $n$ for which $\langle G_k(\boldsymbol{e}_{i_1}, \cdots, \boldsymbol{e}_{i_k}), \boldsymbol{e}_n^\ast \rangle = 0$ while $\sum_{j=1}^{k} \beta_{i_j} - \beta_n = 0$. In this case, we can add to any solution $h_k$ to Eq.~\eqref{h1_eqn} a monomial $c x_{i_1} \cdots x_{i_k}$ with any scalar coefficient $c$ and get another solution; see \eqref{h1_eqn_expansion}--\eqref{Gk_expansion} below.

\item[(iv)]
Note that if the \eqref{NR} condition is strengthened to 
\begin{equation}  \label{NR2}  %\tag{NR2}
\begin{aligned} 
& \Forall \, (i_1, \cdots, i_k ) \in \mathcal{I}^k,  \ n \in \{ m+1, \cdots, N\},  \text{ it holds that} \\
 & \Bigl (\langle G_k(\boldsymbol{e}_{i_1}, \cdots, \boldsymbol{e}_{i_k}), \boldsymbol{e}_n^\ast \rangle \neq 0 \Bigr) \Longrightarrow \biggl ( \sum_{j=1}^{k} \mathrm{Re}(\beta_{i_j}) - \mathrm{Re}(\beta_n) > 0 \biggr), 
\end{aligned}
\end{equation}
then the expression of $h_k$ given by \eqref{h1_part1}--\eqref{h1_part2} results directly from the expression of
Lyapunov-Perron integral $\mathfrak{I}$. Indeed, 
\bea
\mathfrak{I}(\xi) &= \int_{-\infty}^0 e^{-sA_{\s}} \Pi_{\s} G_k\Big(\sum_{i=1}^m e^{\beta_i s}\xi_i \boldsymbol{e}_i \Big) \d s\\
&= \int_{-\infty}^0 \sum_{j=m+1}^N e^{-s\beta_j} \Big \langle G_k\Big(\sum_{i=1}^m e^{\beta_i s}\xi_i \boldsymbol{e}_i \Big),\boldsymbol{e}_{n}  \Big \rangle \boldsymbol{e}_{n} \d s
\eea
i.e.
\be\label{Exp_J}
\mathfrak{I}(\xi)=  \sum_{j=m+1}^N  \sum_{(i_1, \cdots, i_k )\in \mathcal{I}^k}  \, \Big \langle G_k\Big(\boldsymbol{e}_{i_1}, \cdots, \boldsymbol{e}_{i_k} \Big),\boldsymbol{e}_{n}^\ast  \Big \rangle  \xi_{i_1} \cdots \xi_{i_k} \boldsymbol{e}_{n} \int_{-\infty}^0 e^{(\beta_{i_1} + \cdots + \beta_{i_k}- \beta_j)s} \d s,
\ee
recalling that $G_k(u)$ denotes $G_k(u,\cdots,u)$, a homogeneous polynomial or order $k$. 
The condition \eqref{NR2} ensures that the integrals in \eqref{Exp_J} are well-defined, leading to  
\eqref{h1_part1}--\eqref{h1_part2} after integration.

Of course, by assuming only \eqref{NR} instead of \eqref{NR2}, the Lyapunov-Perron integral may not be well defined anymore. But as shown below, the solution to Eq.~\eqref{h1_eqn} still exists, and is given again by \eqref{h1_part1}--\eqref{h1_part2}. 

\item[(v)] Finally, it is worth mentioning that cross non-resonance conditions of the form  
\bes
\sum_{j=1}^{k}\beta_{i_j} - \beta_n \neq 0, \Forall \, (i_1, \cdots, i_k ) \in \mathcal{I}^k,  \ n \in \{ m+1, \cdots, N\},
\ees
is also encountered for the study of normal forms on an invariant manifolds; see, e.g. \cite[Sect.~3.2.1]{Haro}, \cite[Thm.~2.4]{Faria06} and also \cite[Thm.~3.1]{Bibikov79}.

\ei

\er

\bpp[{\bf Proof of Theorem~\ref{thm:h1}}]

The proof is inspired by Lie algebra techniques used in the derivation
of normal forms for ODEs (see, e.g., \cite[Chap.~5]{Arnold88} and \cite[Chap.~1]{Bibikov79}). We proceed in three steps.

\medskip
{\bf Step 1}. We seek a solution to Eq.~\eqref{h1_eqn} as a mapping 
$h_{k} : E_\c \rightarrow E_\s$ that admits the following expansion:
\be  \label{h1_expansion_goal}
h_{k}(\xi) = \sum_{n = m+1}^N \left( \sum_{(i_1, \cdots, i_k) \in \mathcal{I}^k} \Psi^n_{i_1, \cdots, i_k}(\xi) \right) \boldsymbol{e}_n, \quad  \, \xi=(\xi_1, \cdots, \xi_m) \in E_{\c}.
\ee
Here, for each $(i_1, \cdots, i_k) \in \mathcal{I}^k$, the function $\Psi^n_{i_1, \cdots, i_k}(\xi)$ is a complex-valued homogeneous polynomial of degree $k$ given by
\be \label{Psi_def0}
\Psi^n_{i_1, \cdots, i_k}(\xi) = \Gamma^n_{i_1, \cdots, i_k} \xi_{i_1} \cdots \xi_{i_k}.
\ee
The task is then to determine the coefficients $\Gamma^n_{i_1, \cdots, i_k}$ (in $\mathbb{C}$) by using Eq.~\eqref{h1_eqn}.

\medskip

{\bf Step 2}.  In that respect, we introduce the following homological operator $\mathcal{L}_{A}$:
\be \label{eq:ad_A}
\mathcal{L}_{A}[\phi](\xi) = D \phi(\xi) A_{\c} \xi - A_\s \phi(\xi), \qquad \xi \in E_{\c},
\ee
where $\phi \colon E_\c \rightarrow E_\s$ is a smooth function.

A key observation consists of noting that the $E_\s$-valued function, $\xi\mapsto  \Psi^n_{i_1, \cdots, i_k}(\xi)\boldsymbol{e}_n$, provides an eigenfunction of $\mathcal{L}_{A}$  corresponding to the eigenvalue  $\sum_{j = 1}^k \beta_{i_j} - \beta_n$, in other words that the following identity holds
\be \label{eigen_Phi}
\mathcal{L}_{A}[\Psi^n_{i_1, \cdots, i_k}(\xi) \boldsymbol{e}_n](\xi)  = \left [ \sum_{j = 1}^k \beta_{i_j} - \beta_n \right] \Psi^n_{i_1, \cdots, i_k}(\xi) \boldsymbol{e}_n.
\ee

In order to check \eqref{eigen_Phi}, we first calculate $D \phi(\xi) A_{\c} \xi$ when $\phi(\xi)=\Psi^n_{i_1, \cdots, i_k}(\xi) \boldsymbol{e}_n$.  In that respect, denoting by $e^n_j$ the $j^{\mathrm{th}}$ component of $\boldsymbol{e}_n$, the Jacobian matrix $D [\Psi^n_{i_1, \cdots, i_k}(\xi) \boldsymbol{e}_n]$, given by the following $N\times m$ matrix,
\be
D [\Psi^n_{i_1, \cdots, i_k}(\xi) \boldsymbol{e}_n]   =  \begin{pmatrix}
\frac{\partial \Psi^n_{i_1, \cdots, i_k}(\xi)}{\partial \xi_1} e^n_1  & \cdots  & \cdots \frac{\partial \Psi^n_{i_1, \cdots, i_k}(\xi)}{\partial \xi_m} e^n_1 \\
\vdots & \vdots & \vdots \\
\frac{\partial \Psi^n_{i_1, \cdots, i_k}(\xi)}{\partial \xi_1} e^n_N  & \cdots  & \cdots \frac{\partial \Psi^n_{i_1, \cdots, i_k}(\xi)}{\partial \xi_m} e^n_N \\
\end{pmatrix},
\ee
possesses the following representation 
\bea
D [\Psi^n_{i_1, \cdots, i_k}(\xi) \boldsymbol{e}_n]& = \boldsymbol{e}_{n} \Big( \frac{\partial \Psi^n_{i_1, \cdots, i_k}(\xi)}{\partial \xi_1}, \cdots, \frac{\partial \Psi^n_{i_1, \cdots, i_k}(\xi)}{\partial \xi_m} \Big) \\
& = \Gamma^n_{i_1, \cdots, i_k} \boldsymbol{e}_n \boldsymbol{B}(\xi).
\eea
where $\boldsymbol{B}(\xi)=(B_1(\xi), \cdots, B_m(\xi))$ is an $m$-dimensional row vector whose components are given for any  $j$  in $\{1, \cdots, m\}$ by 
\be\label{Exp_Bj}
B_j(\xi) = \frac{\partial }{\partial \xi_j} \big(\xi_{i_1} \cdots \xi_{i_k}\big) =
 \begin{cases}
{\displaystyle p\xi_j^{p-1}\prod_{\substack{i_\ell \neq j}} \xi_{i_\ell}}, & \text{ if  $j \in \{i_1,\cdots, i_k\}$}, \\
0, & \text{otherwise},
\end{cases} 
\ee
where $p$ denotes the number of indices in the set $\{i_1,\cdots, i_k\}$ that equal $j$. 

Thus,
\be 
 D [\Psi^n_{i_1, \cdots, i_k}(\xi) \boldsymbol{e}_n] A_{\c} \xi  = \Gamma^n_{i_1, \cdots, i_k} \boldsymbol{e}_n \boldsymbol{B}(\xi) A_{\c} \xi.
\ee
which leads to 
\be\label{eq:DP_identity2}
D [\Psi^n_{i_1, \cdots, i_k}(\xi) \boldsymbol{e}_n] A_{\c} \xi   = \Gamma^n_{i_1, \cdots, i_k} \boldsymbol{e}_n \boldsymbol{B}(\xi) \left( \beta_1 \xi_1, \cdots, \beta_m \xi_m \right)^\tr,
\ee
since $A$ is assumed to be diagonal.

By noting that the product $ \boldsymbol{B}(\xi) \left( \beta_1 \xi_1, \cdots, \beta_m \xi_m \right)^\tr $ is nothing else that 
$\sum_{j = 1}^k \beta_j \xi_{i_1} \cdots \xi_{i_k},$ and recalling the expression of $\Psi^n_{i_1, \cdots, i_k}(\xi)$ in \eqref{Psi_def0}, we infer from \eqref{eq:DP_identity2} that 
\bea
D [\Psi^n_{i_1, \cdots, i_k}(\xi) \boldsymbol{e}_n] A_{\c} \xi =\sum_{j = 1}^k \beta_{i_j} \Psi^n_{i_1, \cdots, i_k}(\xi) \boldsymbol{e}_n.
\eea

On the other hand, 
\be
A_\s  \Psi^n_{i_1, \cdots, i_k}(\xi) \boldsymbol{e}_n = \beta_n \Psi^n_{i_1, \cdots, i_k}(\xi) \boldsymbol{e}_n,
\ee
and recalling the definition of $\mathcal{L}_A$ in \eqref{eq:ad_A}, the identity \eqref{eigen_Phi} follows.

\medskip 

{\bf Step 3}.  By using the expansion of  $h_{k}(\xi)$ given by \eqref{h1_expansion_goal} in Eq.~\eqref{h1_eqn}, and by using the fact that $\Psi^n_{i_1, \cdots, i_k}(\xi)\boldsymbol{e}_n$ are eigenvectors of the homological operator $\mathcal{L}_{A}$ with eigenvalue $\sum_{j = 1}^k \beta_{i_j} - \beta_n$ (cf.~\eqref{eigen_Phi}), we get  
\beas
\sum_{n = m+1}^N \Bigl[ \sum_{(i_1, \cdots, i_k) \in \mathcal{I}^k} \Bigl ( \sum_{j = 1}^k \beta_{i_j} - \beta_n\Bigr) \Psi^n_{i_1, \cdots, i_k}(\xi) \Bigr] \boldsymbol{e}_n = \Pi_{\s}G_k(\xi).
\eeas
Recalling from \eqref{Psi_def0} that $\Psi^n_{i_1, \cdots, i_k} = \Gamma^n_{i_1, \cdots, i_k} \xi_{i_1} \cdots \xi_{i_k}$, we obtain
\bea  \label{h1_eqn_expansion}
\sum_{n = m+1}^N \Bigl[ \sum_{i_1, \cdots, i_k \in \mathcal{I}^k} \Bigl ( \sum_{j = 1}^k \beta_{i_j} - \beta_n\Bigr) \Gamma^n_{i_1, \cdots, i_k} \xi_{i_1} \cdots \xi_{i_k} \Bigr] \boldsymbol{e}_n = \Pi_{\s}G_k(\xi).
\eea
 
At the same time, since $G_k$ is a homogeneous polynomial of order $k$ and $\xi = \sum_{i = 1}^m \xi_i \boldsymbol{e}_i$, we obtain
\bea \label{Gk_expansion}
\Pi_{\s}G_k(\xi) & = \sum_{n =m+1}^N \langle G_k(\xi), \boldsymbol{e}_n^\ast \rangle \boldsymbol{e}_n \\
& = \sum_{n =m+1}^N \sum_{(i_1, \cdots, i_k) \in \mathcal{I}^k} \xi_{i_1} \cdots \xi_{i_k} \langle G_k(\boldsymbol{e}_{i_1}, \cdots, \boldsymbol{e}_{i_k}), \boldsymbol{e}_n^\ast \rangle \boldsymbol{e}_n.
\eea

By using the above identity in \eqref{h1_eqn_expansion}, we obtain the following formulas for the coefficients $\Gamma^n_{i_1, \cdots, i_k}$ in \eqref{Psi_def0}: 
\be \label{Psi_coef}
\Gamma^n_{i_1, \cdots, i_k} = \frac{\langle G_k(\boldsymbol{e}_{i_1}, \cdots, \boldsymbol{e}_{i_k}), \boldsymbol{e}_n^\ast \rangle}{\sum_{j = 1}^k \beta_{i_j} - \beta_n}.
\ee
The formula of $h_k$ given in \eqref{h1_part1}--\eqref{h1_part2} is thus derived by combining \eqref{h1_expansion_goal}, \eqref{Psi_def0} and \eqref{Psi_coef}. The proof is complete. 
\epp

\subsection{Analytic formulas for higher-order approximations} \label{Sect_G2_case}
We discuss briefly here simple considerations to derive higher-order approximations of an invariant manifold. The approach relies on the use of a power series expansion of the manifold function $h$ in the invariance equation \eqref{eq:invariance}. However, instead of keeping all the monomials at a given degree arising from this expansion, we filter out terms that carries significantly less energy compared with those that are kept. This elimination procedure relies on the assumption that the projected ODE dynamics onto the resolved subspace $E_{\c}$ contains most of the energy; an assumption which is often met in practical applications concerned with invariant manifold reduction. To present the idea in a simple setting, we consider below the case for which $G(y) = G_2(y,y) + G_3(y,y,y)$ and a cubic approximation is sought.

When $G = G_2 + G_3$, the leading-order approximation of $h$ is $h_2$ given by \eqref{h1_part1}--\eqref{h1_part2} with $k=2$. Recall also $h_2$ satisfies \eqref{h1_eqn}.
To determine the approximation of order $3$, we replace $h$ in the invariance equation \eqref{eq:invariance} by $h^{\mathrm{app}}= h_2 + \psi$, where $\psi$ represents the homogeneous cubic terms in the power expansion of $h$, to be determined. By identifying all the terms of order two, we recover \eqref{h1_eqn} with $k=2$ to be satisfied for $h_2$, and by identifying all the terms of order three, we obtain the following equation for $\psi$: 
\be \label{Eq_invariance_psi}
D \psi(\xi) A_{\c} \xi  - A_{\s} \psi(\xi) =  -D h_2(\xi) \Pi_{\c} G_2(\xi) + \Pi_{\s} G_2(\xi, h_2(\xi)) + \Pi_{\s} G_2(h_2(\xi), \xi) + \Pi_{\s} G_3(\xi). 
\ee

Notice that the LHS of \eqref{Eq_invariance_psi} is  $\mathcal{L}_A\psi$, and that the RHS  is a homogeneous cubic polynomial in the $\xi$-variable. If most of the energy of the ODE dynamics is contained in the low modes, one gets that the energy carried by $y_\s$ is much smaller than $\|y_\c\|^2$. It is then reasonable to expect that the energy carried by $h_2(\xi)$ is much smaller than $\|\xi\|^2$ for $\xi = y_\c(t)$ as $t$ varies. This energy consideration implies that on the RHS of \eqref{Eq_invariance_psi}, the term $\Pi_{\s} G_3(\xi)$ dominates the other three terms provided that $\|G_2(y, y)\|/\|y\|^2$ is on the same order of magnitude as $\|G_3(y, y,y)\|/\|y\|^3$. Thus, it is reasonable to seek for a good approximation of  $\psi$ by simply solving the equation: 
\be
D  h_3(\xi) A_{\c} \xi  - A_{\s} h_3(\xi) =  \Pi_{\s} G_3(\xi). 
\ee
Note that this is exactly \eqref{h1_eqn} with $k=3$. In virtue of Theorem \ref{thm:h1}, the existence of $h_3$ is guaranteed under the non-resonance condition \eqref{NR}, and $h_3$ is given by \eqref{h1_part1}--\eqref{h1_part2}. We denote this cubic parameterization by 
\bea \label{Eq_Phi}
\Phi(\xi)&= h_2(\xi) + h_3(\xi) \\
 & = \sum_{n=m+1}^N \bigg( \sum_{(i_1, i_2)\in \mathcal{I}^2} \frac{\langle G_2(\boldsymbol{e}_{i_1}, \boldsymbol{e}_{i_2}), \boldsymbol{e}_n^\ast \rangle}{\beta_{i_1} + \beta_{i_2} - \beta_n} \xi_{i_1} \xi_{i_2} 
+ \hspace{-.5cm}\sum_{(i_1, i_2, i_3)\in \mathcal{I}^3} \frac{\langle G_3(\boldsymbol{e}_{i_1}, \boldsymbol{e}_{i_2}, \boldsymbol{e}_{i_3}), \boldsymbol{e}_n^\ast \rangle}{\beta_{i_1} + \beta_{i_2} + \beta_{i_3} - \beta_n} \xi_{i_1} \xi_{i_2}\xi_{i_3}
\bigg)\boldsymbol{e}_n,
\eea
with $\mathcal{I}= (1, \cdots, m)$.  See the Supplementary Material for an application to the derivation of effective reduced models able to capture a subcritical Hopf bifurcation arising in an ENSO model.

In what precedes, we considered the case $G$ of order 3, and determined approximations  of order 3. 
We could nevertheless, seek for higher-order approximations of invariant manifolds, independently of the nonlinearity to be of high-order or not. For instance if $G(y)=B(y,y)$, i.e. quadratic, we outline hereafter how recursive solutions to a hierarchy of homological equations arise naturally once we look for   higher-order approximations. 

In that respect, we introduce some notations. 
We denote by $\mbox{Poly}_k(E_\c;E_\s)$ (resp.~$\mbox{Poly}_k(E_\c;E_\c)$) the space of vectors in $E_\s$ (resp.~$E_\c$) whose components are homogeneous polynomials of order $k$ in the $E_\c$-variable.   Given a polynomial  $\mathcal{P}$ in $\mbox{Poly}_k(E_\c;E_\s)$ or in $\mbox{Poly}_k(E_\c;E_\c)$, the symbol $\big[ \mathcal{P}(\xi) \big]_k$ represents the collection of terms of order $k$ in $\mathcal{P}$.
  
By seeking a solution, $\Psi$, to the invariance equation Eq.~\eqref{eq:invariance} under the form,
 \be
 \Psi(\xi)=\sum_{k\geq 2} \Psi_k(\xi), \; \Psi_k \in  \mbox{Poly}_k(E_\c;E_\s).
 \ee 
we infer that the $\Psi_k$'s satisfy the following recursive 
 homological equations given by
 \be\label{Eq_general_homol}
 \mathcal{L}[\Psi_k](\xi)= \Big[\Pi_{\s}B(\Phi_{<k}(\xi),\Phi_{<k}(\xi))\Big]_k-\sum_{\ell=2}^{k-1} D\Psi_{k-\ell+1}(\xi)  \Big[\Pi_\c B(\Phi_{<\ell}(\xi),\Phi_{<\ell}(\xi))\Big]_{\ell}  
 \ee
 where $\Phi_{<\ell}(\xi)$ denotes 
 \be
 \Phi_{<\ell}(\xi)=\xi +\sum_{j=2}^{\ell-1} \Psi_j(\xi).
 \ee
Note that with the convention $\sum_{2}^1 \equiv0$, we recover the first homological equation, namely
 \be \label{h1_eqnb}
\mathcal{L}[\Psi_2](\xi) =  \Pi_{\s} B(\xi,\xi).
\ee
In other words $\Psi_2=h_2$.  We refer to \cite{haro2016parameterization} for a detailed account regarding the rigorous and computational aspects for the determination of solutions to Eq.~\eqref{Eq_general_homol}.   \cite[Chap.~11]{kuehn2015multiple} contains also a detailed survey of algorithms to compute numerically invariant manifolds for fast-slow systems.

\newpage
\needspace{1\baselineskip}
\vspace{2em}
{\large \centerline{{\bf Part II: Variational approach to closure}}}

\section{Optimal parameterizing manifolds} \label{Sect_PM_reduction}
\subsection{Variational formulation}\label{Sec_variational_approach_Q}

\subsubsection{Parameterizing manifolds (PM) and parameterization defect}\label{PM_sec}
 A cornerstone of our approach presented below is the 
notion of {\it parameterizing manifold} (PM) that we recall below from \cite{CL15,CLW15_vol2,CLM16_Lorenz9D}.
Our framework takes place in finite dimension as in Part I, however here we consider more general systems of the form 
\be\label{Eq_ODE_gen}
\frac{\d y}{\d t} = A y + G(y) +F, \qquad y\in \mathbb{C}^N,
\ee
where $F$ denotes a time-independent forcing in $\mathbb{C}^N$, $A$ is a $N\times N$ matrix with complex entries,  while $G$ is assumed to be a smooth nonlinearity for which we do not assume $G(0)=0$ anymore.   
In practice  Eq.~\eqref{Eq_ODE_gen} can be thought as derived in the perturbed variable from an original system, for which  $A$ is either the Jacobian matrix at a mean state ($F\neq 0$) or at a steady state ($F=0$), although the concepts presented below do not restrict to such situations.    Hereafter we assume that $A,F$ and $G$ are such that classical solutions (at least $C^1$) exist and that the corresponding initial value problem possesses a unique solution, at least for initial data taken in an open domain $\mathcal{D}$ of $\mathbb{C}^N$.  Dynamically-based formulas to design PMs for Eq.~\eqref{Eq_ODE_gen} are given in Secns.~\ref{Sect_PM_with_forcing} and \ref{Sec_FMTtau} below. For the moment we recall the definition of a PM, and introduce the notion of parameterization defect that will be used for the optimization of PMs\footnote{Note however that other cost functionals may be considered at this stage; see Sec.~\ref{Sec_FMTtau} below.}.

\begin{defi}\label{def:PM}
 Let $T > 0$ and $0\leq t_1 <t_2 \leq \infty$. Let $y$  be a solution to Eq.~\eqref{Eq_ODE_gen}, and $\Psi \colon E_{\c} \rightarrow E_{\s}$ be a continuous mapping satisfying the  following energy inequality for all $t$ in $[t_1,t_2)$
  \begin{align}\label{inequ:2}
   \int_t^{t+T} \norm{y_{\s}(s) - \Psi(y_{\c}(s))}{}^2 \d s  < \int_t^{t+T} \norm{y_{\s}(s)}{}^2 \d s,
  \end{align}
  where $y_{\c}(s)=\Pi_{\c} y(s)$ and $y_{\s}(s)=\Pi_{\s}y(s)$, with $\Pi_\c$ and $\Pi_\s$ that denote the canonical projectors onto $E_\c$ and $E_\s$, respectively ($E_\c$ and $E_\s$ being defined in \eqref{eq:subspaces}).

  Then, the manifold, $  \mathfrak{M}_\Psi$, defined as the graph of $\Psi$, i.e.
   \begin{align}
  \mathfrak{M}_{\Psi}=\{ \xi+ \Psi(\xi)\; | \; \xi \in E_{\c}\},
 \end{align}
is a finite-horizon parameterizing manifold associated with the system of ODEs \eqref{Eq_ODE_gen}, over the time interval $[t_1,t_2)$.
The time-parameter $T$ measuring the length of the ``finite-horizon'' is independent on $t_1$ and $t_2$.  If \eqref{inequ:2} holds for $t_2=\infty$, then $ \mathfrak{M}_{\Psi}$ is simply called a finite-horizon parameterizing manifold, and if it holds furthermore for all $T$, it is called a parameterizing manifold (PM).
\end{defi}

Given a parameterization $\Psi$ of the unresolved variables (in $E_\s$) in terms of the resolved ones (in $E_\c$), a natural non-dimensional number, the {\it parameterization defect}, is defined as
\begin{align}\label{Eq_PD}
Q_T(t,\Psi)=\frac{   \int_t^{t+T} \norm{y_{\s}(s) -\Psi(y_{\c}(s))}{}^2 \d s }{   \int_t^{t+T} \norm{y_{\s}(s)}{}^2 \d s }, \qquad t \in [t_1,t_2).
\end{align}
Sometimes, the dependence on $t$ will be secondary, and by making $t=t_1$ in \eqref{Eq_PD} with $t_1$ sufficiently large so that  for instance transient dynamics has been removed,  we will denote  $Q_T(t,\Psi)$ simply by $Q_T(\Psi)$.
In any event,  either $Q_T(t,\Psi)$ or $Q_T(\Psi)$ allows us to compare objectively two manifolds in their ability to parameterize the variables that lie in the subspace $E_\s$ by those that lie in the subspace $E_{\c}$. Clearly a situation corresponding to an exact slaving of the variables in $E_\s$ by those in $E_\c$ as encountered in the invariant manifold theory revisited in Part I, corresponds to $Q_T(\Psi)\equiv 0$ for any solution $y$ that lies on the invariant manifold,   $\mathfrak{M}_{\Psi}$, associated with the parameterization $\Psi$. If furthermore  $\mathfrak{M}_{\Psi}$ attracts e.g.~exponentially any trajectory like in the case of an inertial manifold, then $Q_T(\Psi)\rightarrow 0$, as $T \rightarrow \infty$ whatever the solution $y$. 
 
A standard $m$-dimensional Galerkin approximation based on the modes in $E_\c$ (with dim$(E_\c)=m$), corresponds  to $\Psi=0$ and thus to   $Q_T(\Psi)\equiv 1$. 
Thus, 
\bes
 \mathfrak{M}_\Psi \mbox{ is a PM if and only if } Q_T(\Psi)<1 \; \mbox{ for all }T>0.  
\ees
Clearly, given a parameterization $\Psi$, it may happen that the corresponding parameterization defect $Q_T(\Psi)$ fluctuates from solutions to solutions, and depends also substantially on the time interval $[t_1,t_2)$ over which the initial time $t$ is taken to compute the integrals in \eqref{Eq_PD}, as well as the horizon $T$.    

Nevertheless, given a set of solutions of interest, a horizon $T$, an interval $[t_1,t_2)$, and a set dimension of the reduced state space (i.e.~dim($E_\c$)$=m$), one is naturally inclined for seeking for parameterizations, $\Psi$, that come with the smallest parameterization defect. In other words, we aim at solving the following minimization problem
\be\label{Eq_Jfunc}
\underset{\Psi \in \mathcal{E}}\min \int_{t}^{t+T} \bigl \|y_\s(s) - \Psi(y_\c(s)) \bigr \|^2 \, \d s, 
\ee
for which  $\mathcal{E}$ denotes a space of parameterizations that makes not only tractable the determination of a minimizer, but also that is not too greedy in terms of data. This latter requirement comes from important practical considerations.  For instance, for high-dimensional systems (e.g.~$N$ of about few hundred thousands), one has typically $y(t)$ available over a relatively small interval of time, and thus if e.g.~$m\sim N/100$ and 
the choice of $\mathcal{E}$ is too naive, such as homogeneous polynomials in the $E_\c$-variable, with arbitrary coefficients, one might easily face an overfitting problem in which too many coefficients have to be determined while not enough snapshots of $y(s)$ are available over $[t, t+T]$. Section \ref{Sect_PM_formulas} below shows that the backward-forward system \eqref{Eq_BF} provides a space  $\mathcal{E}$ of dynamically-based parameterizations that allow to bypass this difficulty as the coefficients to be determined are dependent only on a scalar parameter, the backward integration time $\tau$ in  \eqref{Eq_BF}.

These practical considerations are central in our approach but before providing their details,  we consider in the next section other important theoretical questions. These questions deal with the existence (and uniqueness) of minimizers to \eqref{Eq_Jfunc} on one hand, and with the characterization of the closure system that is reached once  \eqref{Eq_Jfunc} is solved, on the other. 
Thus, we show in Sec.~\ref{Sec_cond_expec}  below that, under assumptions of ergodicity, reasonable  for a broad class of forced-dissipative nonlinear systems such as arising in fluid dynamics, the minimization problem \eqref{Eq_Jfunc} possesses a unique solution, as $T\rightarrow \infty$; see Theorem \ref{Thm_variational-pb} and also \cite[Theorem A.1 and Remark 4.1]{CLM16_Lorenz9D}. 
We call the corresponding minimizer, the {\it optimal parameterizing manifold}.   We conclude finally by showing that an optimal PM, once used as a substitute of the unresolved variables, leads to a reduced system in $E_\c$ that gives the conditional expectation of the original system, i.e.~the best vector field of the reduced state space resulting from averaging of the unresolved variables 
with respect to a probability measure conditioned on the resolved variables; see Theorem \ref{Thm_variational-pb2} below.

We emphasize that PMs have already demonstrated their utility in other applications. For instance, PMs have shown their usefulness for the effective determination of surrogate low-dimensional systems  in view of the optimal control of dissipative nonlinear PDEs. 
 In this case, rigorous error estimates show that parameterization defects arise naturally in the efficient model reduction of optimal control problems (see \cite[Thm.~1 and Cor.2]{CL15}) as furthermore supported by detailed numerical results (see \cite[Sec.~5.5]{CL15} and \cite{Chekroun_al16post}). Speaking roughly, these estimates show that the smaller is the parameterization defect, the better a low-dimensional controller designed from the surrogate system, behaves.  
  Error estimates that relate the parameterization defect to the ability of reproducing the original dynamics' long term statistics by a surrogate system  are difficult to produce for uncontrolled deterministic systems, in particular for chaotic regimes such as considered hereafter in Secns.~\ref{Sect_RBC} and \ref{Sec_KS_turbulence}, due to the singular nature (with respect to the Lebesgue measure) of the underlying invariant measure. In the stochastic realm, this invariant measure becomes smooth for a broad class of systems and the tools of stochastic analysis make the obtention of such estimates  more amenable albeit non trivial; see  \cite{chekroun2019grisanov}. Nevertheless, considerations from ergodic theory and conditional expectations are already insightful for the deterministic systems dealt with in this article as explained in Sec.~\ref{Sec_cond_expec} below.

\subsubsection{Parameterization correlation and angle}\label{Sec_corr_param}
Given a parameterization $\Psi$ that is not trivial (i.e.~$\Psi\neq 0$), we define the {\it parameterization correlation} as, 
\be\label{Eq_corr_param}
c(t)= \frac{\mathrm{Re}\langle \Psi(y_{\c}(t)), y_{\s}(t) \rangle}{\|\Psi(y_{\c}(t))\| \; \|y_{\s}(t)\|}.
\ee
It provides a measure of collinearity  between the parameterized variable $\Psi(y_{\c}(t))$ and the unresolved variable $y_{\s}(t)$, as time evolves. In case of exact slaving, $y_{\s}(t)=\Psi(y_{\c}(t))$ and thus $c(t)\equiv 1$.

The parameterization correlation, $c(t)$, is another key quantity in our approach. 
Speaking roughly, we aim for not only at finding a PM with the smallest parameterization defect but also with  
a parameterization correlation, $c(t)$, to be as much close to one as possible. The basic idea is to find parameterizations that approximate as much as possible an ideal slaving situation, for regimes in which slaving does not hold necessarily. 

\begin{figure}[hbtp]
\centering
   \includegraphics[height=0.35\textwidth, width=.75\textwidth]{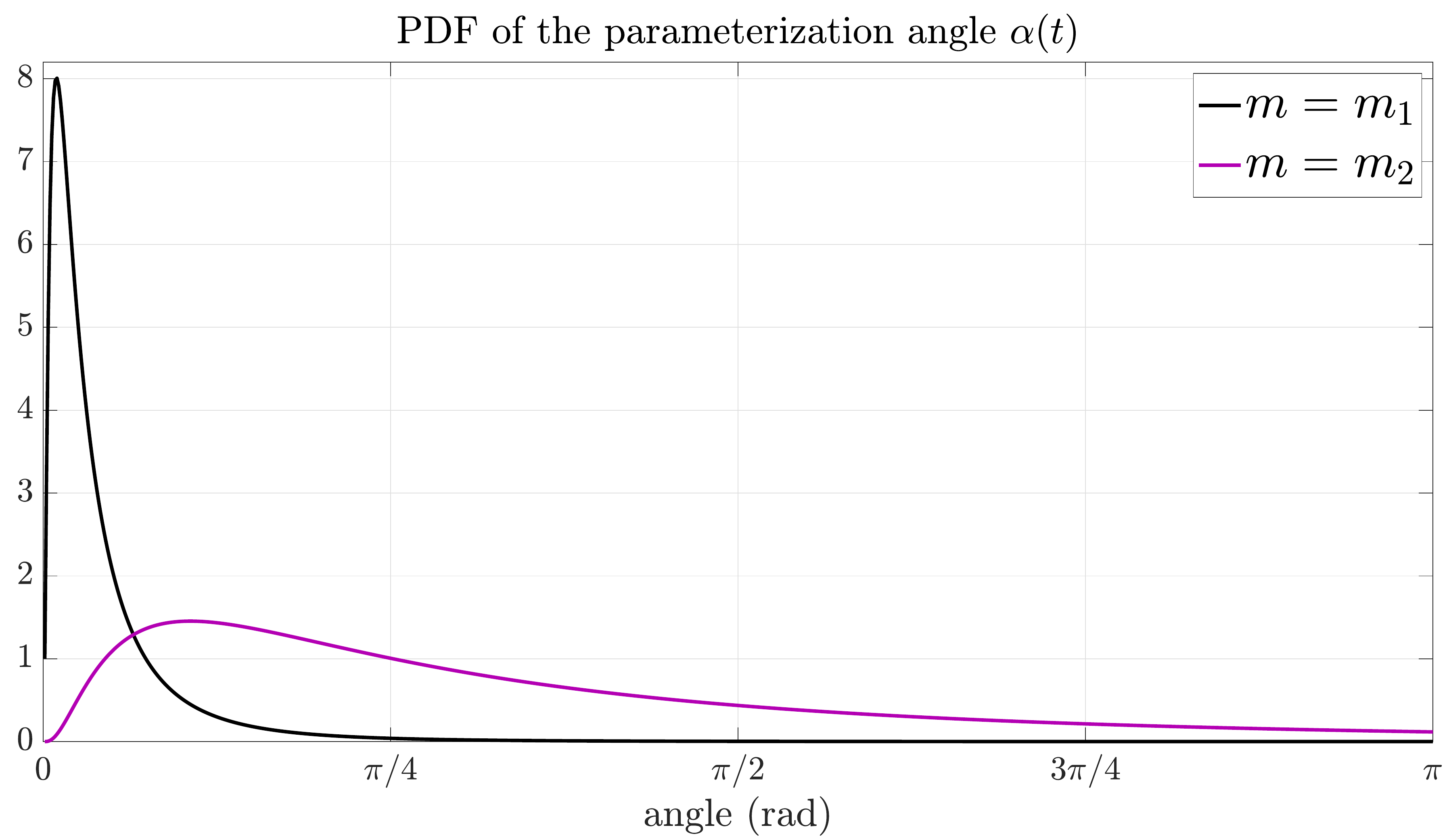}
 \caption{{\footnotesize {\bf Effect of the reduced dimension $m$: Schematic.} This effect is schematically shown here on the PDF of the parameterization angle $\alpha(t)$.  Here a case corresponding to $m_1>m_2$, is depicted: $m_1$ is large enough to be a successful PM while $m_2$ is not.}}
\label{Fig_intro_angle}
\end{figure}

In particular, the parameterization correlation  allows us, once an optimal PM has been determined, to select the dimension $m$ of the reduced phase space according to the following criterium: $m$ should correspond to the lowest dimension of $E_\c$ for which the probability distribution function (PDF) of the corresponding {\it parameterization angle}, 
\be\label{Eq_alpha}
\alpha(t)=\arccos (c(t)),
\ee 
is the most skewed towards zero and the mode of this PDF (i.e.~the value that appears most often) is the closest to zero; see Fig.~\ref{Fig_intro_angle}.    

As a rule of thumb, we aim at finding PMs, $\Psi$, such that: 
\bi
\item[1.] The parameterization defect, $Q_T(\Psi)$, is as small as possible, and
\item[2.] The PDF of the parameterization angle $\alpha(t)$ is skewed towards zero as much as possible, and its mode (i.e.~the value that  appears most often) is close to zero. 
\ei
We illustrate in Sections \ref{Sec_L9D_resurect} and \ref{Sect_RBC} below that, when breakdown of slaving principle occurs, these rules manifest as a natural framework to diagnose and select a parameterization.   Nevertheless as the dimension of the original problem gets large, one may have to inspect a modewise version of $Q_T$ (as discussed in Sec.~\ref{Sec_mode_adaptive}) as well as of $\alpha(t)$; see Sec. \ref{Sec_More_Results} for the latter. In any case, the idea is that one should not only parameterize properly the statistical effects of the neglected scales but also avoid to lose their phase relationships with the retained scales \cite{mccomb2001conditional}. This is particularly important to derive closures that respect a certain phase coherence between the resolved and unresolved scales.

\subsection{Optimal parameterizing manifold and conditional expectation}\label{Sec_cond_expec} 
We present in this section the main results that serve as a foundational basis for the applications discussed hereafter. 
We denote by $X$ the vector field associated with Eq.~\eqref{Eq_ODE_gen}
i.e.~
\be\label{Eq_vector_field}
X(y)=A y + G(y) +F, \quad \mbox{for all} \; y \in \mathbb{C}^N. 
\ee

To simplify the presentation, we assume this vector field  to be 
sufficiently smooth and dissipative on $\mathbb{C}^{N}$, such that the corresponding flow, $T_t$, is well-defined.
We assume, furthermore,  that $T_t$ possesses an invariant probability measure $\mu$, which is {\it physically relevant} \cite{eckmann_ruelle,collet2007concepts},  in the sense that the following property holds for $y$ in a positive Lebesgue measure  set $B(\mu)$ (of $\mathbb{C}^N$) and for every continuous observable $\varphi:\mathbb{C}^N\rightarrow \mathbb{C}$
\be\label{Eq_phys_rev}
\underset{T\rightarrow \infty}\lim \frac{1}{T} \int_0^{T}  \varphi(T_t( y)) \d t =\int \varphi( y) \d \mu ( y). 
\ee
This property assures that meaningful averages can be calculated and the statistics of the dynamical system can be investigated by the asymptotic distribution of orbits starting from Lebesgue almost every initial condition in e.g.~the basin of attraction, $B(\mu)$, of the statistical equilibrium $\mu$.

Recall that, like all probability measures invariant under $T_t,$ an invariant measure that satisfies \eqref{Eq_phys_rev} is supported by the global attractor $\mathcal{A}$ when the latter exists; e.g.~\cite[Lemma 5.1]{chekroun_glatt-holtz}.  In the case a global attractor is not known to exist, an invariant measure has its support in the {\it non-wandering set}, $\Lambda$; see \cite[Remark 1.4, p. 197]{FMRT01}.

It can be proven for e.g.~Anosov flows \cite{bowen1975ergodic}, partially hyperbolic systems \cite{alves2000srb}, Lorenz-like flows \cite{bonatti2000lorenz}, 
and observed experimentally for many others \cite{eckmann_ruelle,gallavotti1995dynamical,csg11,CGN18b} that a
common feature of (dissipative) chaotic systems is the transformation (under the action of the flow) of the initial
Lebesgue measure  into a probability measure with finer and finer scales, reaching asymptotically an invariant measure $\mu$ of Sinai-Ruelle-Bowen (SRB) type. This measure is singular with respect to the Lebesgue measure,  is supported by the local unstable manifolds contained in $\mathcal{A}$ or in $\Lambda$ \cite[Def.~6.14]{collet2007concepts}, and if it has no zero Lyapunov exponents it satisfies \eqref{Eq_phys_rev} \cite{young2002srb}.  This latter property is often referred to as the {\it chaotic hypothesis} that, roughly speaking, expresses an extension of the ergodic hypothesis to non-Hamiltonian systems \cite{gallavotti1995dynamical}. We work thus hereafter within this hypothesis and we  assume  furthermore that \eqref{Eq_phys_rev} holds for $\varphi$ that lies in the space of integrable function, $L_\mu^1 (\mathbb{C}^N)$, with respect to the invariant measure $\mu$.

%%%%%%%%%%%%%%%%%%%%%%%%%%%%%%%%%%%%
Having clarified the ergodic framework within which we will frame our variational approach, we 
consider now a high-mode parameterization of the form
\be\label{Param_psi}
\Psi(\xi) = \sum_{n = m+1}^N \Psi_n (\xi) \boldsymbol{e}_n, \; \; \xi \in E_\c,
\ee
with the $\boldsymbol{e}_n$'s denoting the eigenmodes  of the linear part, $A$,  that span the subspace $E_\s$. The regularity assumption made on $\Psi$ is clarified hereafter; see Theorem \ref{Thm_existence_lim}. In practice, $\Psi$ does not need to cover the whole range $[m+1,N]$ and some $\Psi_n$ may be zero.

We denote by $\m$ the push-forward of the measure $\mu$ by the projector $\Pi_\c$ onto $E_\c$, namely  
\be\label{Definition_m}
\m(B)=\mu(\Pi_\c^{-1}(B)), \quad B\in \mathcal{B}(E_\c),
\ee
where $\mathcal{B}(E_\c)$ denotes the family of Borel sets of $E_\c$; i.e.~the family of sets that can be formed from open sets (for the topology on $E_\c$ induced by the norm $\|\cdot \|_{E_\c}$) through the operations of countable union, countable intersection, and relative complement.

In what follows (see Sec.~\ref{Sect_PM_formulas}), given a solution $y(t)$ that emanates from $y_0$ in $B(\mu)$, we also consider the parameterization defect, $\mathcal{Q}_n$, associated with the parameterization $\Psi_n$ of the $n^{\rm th}$-eigenmode, namely 
\be\label{Eq_QnT}
\mathcal{Q}_n(T)=\frac{1}{T}\int_{0}^T \Big|\langle y_\s(t), \boldsymbol{e}_n^\ast\rangle -\Psi_n (y_\c(t))\Big|^2 \d \, t,
\ee
where we recall that $\{\boldsymbol{e}_j^\ast\}_{j=1}^N$ denotes the eigenvectors of the conjugate transpose $A^\ast$.

In the case $\{\boldsymbol{e}_n\}$ forms an orthonormal basis of $\mathbb{C}^N$, {\hl namely when $A$ is a Hermitian matrix}, we have due to the Parseval's identity, 
\be\label{Eq_Q_T}
\mathcal{Q}_T(\Psi)=\frac{1}{T}\int_{0}^T \norm{y_\s(t) -\Psi (y_\c(t))}^2 \d \, t=\sum_{n=m+1}^N \mathcal{Q}_n(T).
\ee
However this equality does not always hold, in general.  Indeed, by writing $y_\s(t) = \sum_{n=m+1}^N y_n(t) \boldsymbol{e}_{n}$ with $y_n(t) = \langle y_\s(t), \boldsymbol{e}_n^\ast\rangle$, we {\mkr remark that} 
\bes
\norm{y_\s(t) -\Psi (y_\c(t))}^2 = \sum_{n_1, n_2=m+1}^N \bigg \langle \Big(y_{n_1}(t) - \Psi_{n_1} (y_\c(t))\Big)  \boldsymbol{e}_{n_1},  \Big(y_{n_2}(t) - \Psi_{n_2} (y_\c(t)) \Big) \boldsymbol{e}_{n_2} \bigg \rangle,
\ees
and the {\mkr latter identity is reduced to} $\sum_{n=m+1}^N |y_{n}(t) - \Psi_{n_1} (y_\c(t))|^2$ {\mkr when  $\langle \boldsymbol{e}_{j}, \boldsymbol{e}_{k} \rangle = \delta_{j, k}$} for all $j,k = m+1, \cdots, N$.

Thus, solving \eqref{Eq_Jfunc} is not {\hl always} equivalent to solving the following family of variational problems
\be\label{Eq_Jfunc_n}
\underset{\Psi_n \in \mathcal{E}}\min \int_{0}^T \Big|\langle y_\s(t), \boldsymbol{e}_n^\ast\rangle -\Psi_n (y_\c(t))\Big|^2 \d \, t, \qquad m+1 \leq n \leq N.
\ee
As we will see, for practical reasons we will often prefer to solve \eqref{Eq_Jfunc_n} rather than  \eqref{Eq_Jfunc}; see Sec.~\ref{Sec_mode_adaptive} below. Nevertheless, the existence and uniqueness of minimizers for either \eqref{Eq_Jfunc_n} or  \eqref{Eq_Jfunc}, are dealt with in the same way. Hereafter, we present the latter only in the case of  \eqref{Eq_Jfunc} (allowing for the simplification of certain statements) and leave to the reader the corresponding statements and proofs in the case of the minimization problems  \eqref{Eq_Jfunc_n}.

 In that respect, we select the space of parameterizations, $\mathcal{E}$, to be the Hilbert space constituted by $E_\s$-valued functions of the resolved variables $\xi$ in $E_\c$, that are square-integrable with respect to $\mathfrak{m}$, namely 
\be\label{E_class}
\mathcal{E}=L^2_{\m}(  E_\c; E_\s)=\Big\{\Psi:   E_\c\rightarrow E_\s\,\, \mathrm{measurable} \; \mbox{and} \; \mbox{such} \; \mbox{that} \; \int_{ E_\c} \|\Psi(\xi) \|^2\d \m(\xi) <\infty\Big\}.
\ee

Our approach to minimize, $\mathcal{Q}_T(\Psi)$ (in $\mathcal{E}$), and to identify parameterizations for which the normalized parameterization defect 
\be\label{Eq_normalized_defect}
Q_T(\Psi)=\mathcal{Q}_T(\Psi) \langle \norm{y_\s}^2\rangle_T^{-1},
\ee
satisfies
\be\label{Eq_limQ}
0<\underset{T\rightarrow \infty}\lim Q_T(\Psi)<1, 
\ee
relies substantially on the general disintegration theorem of probability measures; see e.g.~\cite[p.~78]{dellacherie1978probabilities}.
In \eqref{Eq_normalized_defect}, we have denoted by $\langle \norm{y_\s}^2\rangle_T$ the time-mean of $y_\s$ over $[0,T]$.
The disintegration theorem states that given a probability measure $\mu$ on $\mathbb{C}^N$, a vector subspace $V$ of $\mathbb{C}^N$, and a Borel-measurable mapping $\mathfrak{p}: \mathbb{C}^N \rightarrow  V$, then there exists a uniquely determined family of probability measures $\{\mu_{x}\}_{x\in V}$ such that, for $\m$-almost  all \footnote{i.e.~up to an exceptional set of null measure with respect to $\m$.} $x$ in $V$, $\mu_{x}$ is  concentrated on the pre-image $\mathfrak{p}^{-1} (\{x\})$ of $x$, i.e.~$\mu_{x} \left(\mathbb{C}^N\setminus \mathfrak{p}^{-1} (\{x\}) \right) = 0$,
and such that for every Borel-measurable function $\phi: \mathbb{C}^N \rightarrow \mathbb{C}$,
\be\label{Eq_desint}
\int \phi (y) \d \mu (y) = \int_{V} \Big(\int_{y\in \mathfrak{p}^{-1} (\{x\})} \phi(y) \d\mu_{x} (y) \Big)\d \m (x).
\ee
Here $\m$ denotes the {\it push-forward} in $V$ of the measure $\mu$  by the mapping $\mathfrak{p}$, i.e.~$\mathfrak{m}$ is given by \eqref{Definition_m} where $\Pi_\c$ is replaced by $\mathfrak{p}$. Note that when $\mathfrak{p}$ is the projection onto $V$, the probability measure $\mu_x$ is the conditional probability of the unresolved variables, contingent upon the value of the resolved variable to be $x$; see also \cite[Supporting Information]{Chek_al14_RP}. 

Hereafter, we apply this theorem with the reduced phase space, $V$, to be the subspace of the resolved variables, $E_\c$, and the mapping $\mathfrak{p}$ to be the projector $\Pi_\c$ onto $E_\c$. In this case, a decomposition analogous  to \eqref{Eq_desint} holds for the measure $\mu$ itself, namely  
\be\label{Eq_desint2}
\mu(B\times F)=\int_{F} \mu_{\xi} (F) \d \m (\xi), \qquad B\times F \in \mathcal{B}(E_\c)\otimes\mathcal{B}(E_\s).
\ee
First, we state a result identifying natural conditions under which, $\underset{T\rightarrow \infty}\lim\mathcal{Q}_T(\Psi)$ exists. 

%%%%%%%
\bt\label{Thm_existence_lim}
Assume that Eq.~\eqref{Eq_ODE_gen} admits an invariant probability measure $\mu$ satisfying \eqref{Eq_phys_rev} and that
the unresolved variable $\zeta$ in $E_\s$ has a finite energy in the sense that % for all $\xi$ in $E_\c$,
\be\label{Eq_finite-energy}
 \int \norm{\zeta}^2 \d \mu <\infty.
\ee

If $\Psi$ lies in $L^2_{\mathfrak{m}}(E_\c,E_\s)$, then for a.e.~solution $y(t)$ of Eq.~\eqref{Eq_ODE_gen} that emanates from an initial datum $y_0$ in the basin of attraction $B(\mu)$, the limit $\underset{T\rightarrow \infty}\lim\mathcal{Q}_T(\Psi)$ exists, and is given by 
\be\label{Eq_limit_QnA}
\underset{T\rightarrow \infty}\lim\mathcal{Q}_T(\Psi)=\int_{(\xi,\zeta)\in E_\c \times E_\s} \|\zeta -\Psi(\xi)\|^2 \d \mu. 
\ee
\et
\bp
This theorem is a direct  consequence of the ergodic property \eqref{Eq_phys_rev} applied to the observable 
\be
\varphi(\xi,\zeta)= \|\zeta -\Psi(\xi)\|^2.
\ee  
Indeed, first, let us note that $\varphi(\xi,\zeta)=\|\zeta\|^2-2\langle \zeta,\Psi(\xi)\rangle + \|\Psi(\xi)\|^2$ satisfies 
\be\label{Eq_toto}
\int \varphi(\xi,\zeta) \d\mu \leq  \int \|\zeta\|^2 \d\mu_\xi (\zeta) + \int \|\Psi(\xi)\|^2 \d \mathfrak{m} +\int (\| \zeta\|^2 +\| \Psi(\xi)\|^2) \d \mu,
\ee 
by application of \eqref{Eq_desint2} and the Fubini's theorem for the two first integrals in the RHS of \eqref{Eq_toto}, and of the Cauchy-Schwarz and Young inequalities for the third integral.  Another application of \eqref{Eq_desint2} and the Fubini's theorem for this latter integral shows that  $\varphi$ lies in $L^1_{\mu}(\mathbb{C}^N)$, since $\Psi$  belongs to $L^2_{\mathfrak{m}}(E_\c,E_\s)$ and \eqref{Eq_finite-energy} holds.
\ep

We are now in position to show the existence of a unique minimizer to the minimization problem
\be\label{Eq_variational-pb0}
\underset{\Psi \in \mathcal{E}}\min\bigg(\underset{T\rightarrow \infty}\lim\mathcal{Q}_T(\Psi)\bigg),
\ee
i.e.~to ensure the existence of an optimal manifold minimizing the parameterization defect. The minimizer is also characterized; see \eqref{Def_h2} below.  An earlier version of such results may be found in \cite[Theorem A.1]{CLM16_Lorenz9D} for the special case of a truncated Primitive Equation model due to  Lorenz \cite{Lorenz80}. The general case is dealt with below. 
%%%%%%%%%%%%%%
\bt\label{Thm_variational-pb}
Assume that the assumptions of Theorem \ref{Thm_existence_lim} hold.
Then the  minimization problem
\be\label{Eq_variational-pb}
\underset{\Psi \in \mathcal{E}}\min \int_{(\xi,\zeta)\in E_\c \times E_\s} \norm{\zeta -\Psi(\xi)}^2 \d \mu,
\ee
possesses a unique solution in $\mathcal{E}=L^2_{\m}(E_\c,E_\s)$ whose argmin  is given 
by
\be\label{Def_h2}
 \Psi^\ast(\xi)=\int_{E_\s} \zeta \d \mu_{\xi}(\zeta), \qquad \xi \in E_\c.
\ee
Furthermore
\be\label{Eq_PM_check}
\underset{T\rightarrow \infty}\lim\mathcal{Q}_T(\Psi^\ast) \leq \underset{T\rightarrow \infty}\lim\mathcal{Q}_T(\Psi), \;\; \forall \; \Psi \in L^2_{\m}(E_\c,E_\s).
\ee
\et

\bp
The proof is a direct consequence of the disintegration theorem applied to the ergodic measure $\mu$. Let us introduce the following Hilbert space of $E_\s$-valued functions 
\be
L^2_{\mu}( E_\c\times E_\s; E_\s)=\Big\{f: E_\c\times E_\s\rightarrow E_\s, \,\, \mathrm{measurable \, and \, s.t.} \; \int_{E_\c\times E_\s} \|f(\xi,\zeta)\|^2  \d \mu (\xi,\zeta) <\infty\Big\}.
\ee

Let us define the expectation $\mathbb{E}_\mu(g)$ with respect to the invariant measure $\mu$ by
\be\label{Def_expec}
\mathbb{E}_\mu(g)=\int_{E_\c\times E_\s} g(\xi,\zeta) \d \mu (\xi,\zeta), \qquad g\in L^2_{\mu}( E_\c\times E_\s; E_\s).
\ee

By applying to the ambient Hilbert space  $L^2_{\mu}( E_\c\times E_\s; E_\s)$, the standard projection theorem onto closed convex sets \cite[Theorem 5.2]{brezis_book}, one defines (given $\Pi_\c$) the conditional expectation $\mathbb{E}_\mu[g\vert \Pi_\c]$ of $g$ as the unique function in $\mathcal{E}$ that satisfies the inequality 
\be\label{Eq_best_app}
\mathbb{E}_\mu[\|g-\mathbb{E}_\mu[g| \Pi_\c]\|^2] \leq \mathbb{E}_\mu[\| g-\Psi \|^2], \; \mbox{for all } \Psi \in\mathcal{E}.
\ee

The general disintegration theorem of probability measures, applied to $\mu$ (see \eqref{Eq_desint}),
provides the following  explicit representation of the conditional expectation 
\be
\mathbb{E}_\mu[g| \Pi_\c]=\int_{E_\s} g(\xi,\zeta) \d\mu_\xi(\zeta),
\ee
with $\mu_\xi$ denoting the disintegrated measure of $\mu$ in \eqref{Eq_desint2}.

Now let us take $g(\xi,\zeta)=\zeta$,  then 
\be
\mathbb{E}_\mu[\zeta| \Pi_\c]=\Psi^\ast, 
\ee
with $\Psi^\ast$ defined by \eqref{Def_h2}. 
 We have then 
\be
\norm{\Psi^\ast (\xi)}^2 \leq \int \norm{\zeta}^2 \d \mu_\xi(\zeta),
\ee
and by using \eqref{Eq_desint} we have
\be
\int \norm{\Psi^\ast (\xi)}^2 \d \m (\xi) \leq  \int \norm{\zeta}^2 \d \mu.
\ee
This inequality shows that  $\Psi^\ast$ lies in $L^2_{\m}(E_\c,E_\s)$ due to assumption \eqref{Eq_finite-energy}.

We have then from \eqref{Eq_best_app}, 
\be
\mathbb{E}_\mu[\|\zeta-\Psi^\ast\|^2] \leq \mathbb{E}_\mu [\| \zeta-\Psi \|^2], \; \mbox{ for all }\Psi \in \mathcal{E}.
\ee

By recalling that 
\be
\mathbb{E}_\mu[\|\zeta-\Psi^\ast\|^2] =\int_{E_\s\times E_\s } \|\zeta-\Psi^\ast(\xi)\|^2 \d \mu (\xi,\zeta)=\int \|\zeta-\Psi^\ast(\xi)\|^2 \d \mu(\xi,\zeta),
\ee
one obtains then,  by applying respectively \eqref{Eq_phys_rev} to $\varphi=\|\zeta-\Psi^\ast\|^2$ and $\varphi=\|\zeta-\Psi\|^2$, that for all $\Psi $ in $\mathcal{E}$,
\be
\underset{T\rightarrow \infty}\lim\; \frac{1}{T} \int_0^T \|y_\s(t)-\Psi^\ast(y_\c(t))\|^2 \d t \leq \underset{T\rightarrow \infty}\lim \; \frac{1}{T} \int_0^T \|y_\s(t)-\Psi(y_\c(t))\|^2 \d t.
\ee
The proof is complete. 
\ep
%%%%%%%%%%%%%%%%%%%%%%%%%%%%%%%%%%%%%%
The manifold obtained as the graph of $\Psi^\ast$ given by \eqref{Def_h2} will be called the {\it optimal PM}. 
Formula \eqref{Def_h2} shows that the optimal PM corresponds actually to the manifold that maps to each resolved variable $\xi$ in $E_\c$,  the averaged value of the unresolved variable $\zeta$ in $E_\s$ as distributed according to the conditional probability measure $\mu_\xi$.  
In other words, the optimal PM provides the best manifold (in a least-square sense) that averages out the fluctuations of the unresolved variable. 

By making $\Psi\equiv0$ in \eqref{Eq_PM_check}, this optimal PM  comes with a (normalized) parameterization defect \eqref{Eq_normalized_defect} that satisfies necessarily 
\be
0\leq \underset{T\rightarrow \infty}\lim Q_T(\Psi^\ast) \leq 1.
\ee
This variational view on the  parameterization problem of the unresolved variables removes any sort of ambiguity that has surrounded the notion of (approximate) inertial manifold in the past. Indeed, within this paradigm shift,  given an ergodic invariant measure $\mu$ and a reduced dimension $m$ (defining thus a projector $\Pi_\c$), the optimal PM may have a parameterization defect very close to 1 and thus the best possible nonlinear parameterization one could ever imagine cannot {\it a priori} do much better than a classical Galerkin approximation, and sometimes even worse.  To the opposite, the smaller $Q_T(\Psi^\ast)$ is (for $T$ large), the best the parameterization. All sort of nuances are actually admissible, even when the parameterization defect is just below unity; see \cite{CLM16_Lorenz9D} and Sec.~\ref{Sec_L9D_resurect} below.  

 We emphasize that although the theory presented in this section has been shaped for asymptotic values of $T$, in practice 
we will be instead interested to seek for optimal PMs learned over a training length as short as possible (to rely on as few as possible DNS snapshots). In that respect, it is where the parametric families of dynamically-based parameterizations derived in Sec.~\ref{Sect_PM_formulas} below (and relying on Part I) become useful. We will indeed show that by applying these formulas in practice, we are able to derive optimal PMs trained over short training intervals of length comparable to a characteristic recurrence or decorrelation time of the dynamics; see Secns.~ \ref{Sect_RBC} and \ref{Sec_KS_turbulence} below.

\needspace{1\baselineskip}
\br\label{Rmk_problems_to_overcome}
\hspace*{2em}  %\vspace*{0.1ex}
\bi
\item[(i)] The ergodic property \eqref{Eq_phys_rev} can be relaxed into weaker forms such as considered in e.g. \cite{FMRT01,chekroun_glatt-holtz}.   These relaxed versions hold for a broad class of dissipative systems including systems of ODEs and even PDEs, as long as a global attractor exists \cite[Theorem 2.2]{chekroun_glatt-holtz}. However these weaker forms 
do not guarantee the existence of the limit in  \eqref{Eq_limit_QnA} and the latter would be replaced instead by a notion of generalized limit  involving e.g.~averaging over accumulations points.  The  statistical equilibrium  $\mu$ is then not guaranteed to be unique. 

Nevertheless, bearing these changes in mind, the proof presented above can be easily adapted and the conclusion of Theorem \ref{Thm_variational-pb} remains valid with however a form of optimality that is now subject to the choice of the statistical equilibrium. Within this ergodic framework, several optimal parameterizing manifolds may co-exist but for each statistical equilibrium there is only one optimal parameterizing manifold. The same is true if a global attractor $\mathcal{A}$ is not guaranteed to exist: $\mathcal{A}$ must be replaced by the non-wandering set $\Lambda$, and the optimal PM is unique for trajectories sampled according to the statistical equilibrium $\mu$.

\item[(ii)]  With the nuances brought up in (i) above,  Theorem \ref{Thm_variational-pb} applies thus to any relevant Galerkin truncations of systems of PDEs arising in fluid dynamics; see \cite{CLM16_Lorenz9D} and Sec.~\ref{Sec_L9D_resurect} below for an application to a 9D Galerkin truncation of the Primitive Equations of the atmosphere due to Lorenz \cite{Lorenz80}.
\item[(iii)]  Theorem \ref{Thm_variational-pb} is fundamental for understanding and interpretation but is of little interest for computing the optimal PM in practice, except in specific problems for which $\mu$ is known explicitly (see e.g.~\cite[Sec.~4]{Chekroun2017a}) or can be approximated semi-analytically \cite{majda2001mathematical,majda2003systematic}; see also \cite{vanden2003fast} for an  
alternative approach to estimate numerically $\mu_\xi$ in the context of slow-fast systems. In Section \ref{Sect_PM_formulas} below we introduce instead explicit dynamically-based parameterizations that, once optimized according to a mode-adaptive approach, provide an efficient way to determine PMs that although suboptimal (for \eqref{Eq_variational-pb}) will be shown to be skillful for closure in practice; see Secns.~\ref{Sect_RBC} and  \ref{Sec_KS_turbulence} below.      
\ei
\er

We  have then the following result relating the conditional expectation to the optimal PM. 
We state this theorem in the case of quadratic interactions, motivated by applications in fluid dynamics; see also \cite[Sec.~4.3]{CLM16_Lorenz9D} and Sec.~\ref{Sec_L9D_resurect} below, for an illustration.  
\bt\label{Thm_variational-pb2}
Under the conditions of Theorem \ref{Thm_variational-pb} if $G$ is a quadratic nonlinearity $B$ in Eq.~\eqref{Eq_ODE_gen}, 
the conditional expectation, $\mathbb{E}_\mu[X| \Pi_\c]$, satisfies 
\be\label{Eq_cond_exp_optiPM}
\mathbb{E}_\mu[X| \Pi_\c](\xi)=A_\c \xi +\Pi_\c B(\xi,\xi) +\Pi_\c \big(B(\xi,\Psi^\ast(\xi)) +B(\Psi^\ast(\xi),\xi)\big)+F_\c+\eta(\xi), \; \xi \in E_\c,
\ee
where $X$ is the vector field given by \eqref{Eq_vector_field}, $\Psi^\ast$ is the optimal PM guaranteed by Theorem \ref{Thm_variational-pb}, and $\eta$ is given by
\be
\eta(\xi) =\int_{\zeta \in E_\s} \Pi_\c B(\zeta,\zeta)\d\mu_\xi(\zeta).
\ee 
Thus in the case $\eta=0$, the optimal PM, $\Psi^\ast$, provides the conditional expectation $\mathbb{E}_\mu[X| \Pi_\c]$, i.e. 
\be
\mathbb{E}_\mu[X| \Pi_\c](\xi)=A_\c \xi +\Pi_\c B(\xi,\xi) +\Pi_\c \big(B(\xi,\Psi^\ast(\xi)) +B(\Psi^\ast(\xi),\xi)\big)+F_\c.
\ee
\et

\bp
Expanding $X(\xi +\zeta)$ (with $(\xi,\zeta)$ in $E_
\c \times E_\s$) and integrating with respect to the disintegrated probability measure, $\mu_\xi$, we get (by using that $\int \d \mu_\xi =1$)
\bea
\mathbb{E}_\mu[X| \Pi_\c](\xi)&=A_\c \xi +\Pi_\c B(\xi,\xi)+F_\c+\eta(\xi)+\int \bigg( \Pi_\c \big(B(\xi,\zeta) +B(\zeta,\xi)\big)  \bigg)\d \mu_\xi(\zeta),\\
&=A_\c \xi +\Pi_\c B(\xi,\xi)+F_\c +\eta(\xi) +\Pi_\c B\bigg(\xi,\int \zeta \d \mu_\xi(\zeta)\bigg) +\Pi_\c B\bigg(\int \zeta \d \mu_\xi(\zeta),\xi\bigg),
\eea
which given the expression of $\Psi^\ast$ in \eqref{Def_h2}, gives \eqref{Eq_cond_exp_optiPM}. 
\ep

\subsection{Inertial manifolds and optimal PMs}\label{Sec_IMs}
To avoid any confusion, we clarify the distinction between the concept of an inertial manifold (IM) and that of an optimal {\hl parameterizing manifold (PM)}. First of all, an IM is a particular case of an asymptotic PM since when an inertial manifold $\Psi$ exists, $Q_T(\Psi)=0$ for all $T$ sufficiently large. 
We list below some  important points to better appreciate the differences between the two concepts.  
\bi
\item[(i)] When an IM, $\Psi$, exists,  then $\Psi=\Psi^\ast$ in \eqref{Def_h2} with 
$\mu_\xi$ being the Dirac mass (in $E_\s$) concentrated on $\Psi(\xi)$, i.e.~$\mu_\xi=\delta_{\Psi(\xi)}$. Furthermore in this case, the probability distribution $p_\alpha$ of the parameterization angle, $\alpha(t)$ given by \eqref{Eq_alpha}, is given by the Dirac mass $\delta_0$ (on the real line) concentrated at $0$.
\item[(ii)]  Working with the eigenbasis of  the linear part of Eq.~\eqref{Eq_ODE_gen} and assuming that an IM exists, let $m_\ast$ denote the minimal dimension of the reduced state space required for an IM to exist. 
If $m=\mbox{dim}(E_\c)<m_\ast$ then there  is no inertial manifold but a PM still exists in general as supported by Theorem \ref{Thm_existence_lim}.
One may wonder however whether more can be said when $m<m_\ast$. 

This is where the parameterization defect, $Q_T$, and the parameterization angle, $\alpha(t)$, provide useful mutual informations.
Typically when $m<m_\ast$, seeking for a manifold that minimizes $Q_T$ allows for parameterizing optimally (in a least square sense) the statistical effects of the neglected scales in terms of those retained. However one should keep in mind to avoid losing the phase relationships between the resolved and unresolved scales, and in that sense the distribution  $p_\alpha$ should not be too spread.  For systems  with a high-dimensional global attractor one may need to inspect a modewise version of $Q_T$ (as discussed in Sec.~\ref{Sec_mode_adaptive} below) as well as of $\alpha(t)$ for the design of the nonlinear parameterization; see Sec. \ref{Sec_More_Results} for the latter in the context of 1D Kuramoto-Sivashinsky turbulence.
\ei

Thus, even for systems that admit an IM, an optimal PM often provides an efficient closure based on much fewer modes compared to an inertial form. Such an observation about efficient reduced dimension is known by the  practitioner familiar with the notion of approximate inertial manifold (AIM). An AIM provides a manifold such that the attractor lies within a neighborhood of it that shrinks as the reduced dimension $m$ is increased \cite{MT89,debussche1992construction,devulder1993rate}. Nevertheless, as the reduced dimension is set too low, a given AIM may suffer from e.g.~an over-parameterization of the small scales resulting into dramatic errors backscattering to the large scales; see Sec.~\ref{Sec_KS_turbulence}. This is because the AIM approach does not address the question of finding an optimal manifold that minimizes the parameterization defect while keeping the reduced dimension as low as possible. This is the focus of the PM approach proposed in this article which is thus, in essence, variational rather than concerned with the rate of convergence with $m$ as in standard AIM theory.

\subsection{A reduced-order Primitive Equation example: PM and breakdown of slaving principles}\label{Sec_L9D_resurect}
The conditional expectation is related to the optimal PM according to Theorem \ref{Thm_variational-pb2}, making thus the optimal PM an essential ingredient for the closure problem.  
 Depending on the problem at hand, the conditional expectation provides e.g.~the reduced equations that filter out the fast gravity waves from truncated Primitive Equations (PE) of the atmosphere; see \cite{CLM16_Lorenz9D}. 
Truncations corresponding to $\eta=0$ in \eqref{Eq_cond_exp_optiPM}, i.e.~when the high-high interactions do not contribute to the low-mode dynamics, is particularly favorable for the conditional expectation to provide such a filtering property. As shown numerically in \cite{CLM16_Lorenz9D}, the conditional expectation provides indeed such a ``low-pass filter'' closure 
for the truncated PE proposed by Lorenz in 1980 \cite{Lorenz80}, when a critical Rossby number, $\epsilon^\ast$, is crossed. 
We reproduce hereafter some of these numerical results and provide new, complementary understanding based on the theory of PMs such as discussed in this article. 

The model of \cite{Lorenz80}, when rescaled following \cite{CLM16_Lorenz9D},  becomes
\bea \label{Eq_L9D_rescale}
\e^2 a_i \frac{\d X_i}{\d t} &= \e^3 a_i b_iX_jX_k - \e^2 c(a_i - a_k) X_j Y_k 
+ \e^2 c(a_i - a_j) Y_j X_k \\
&\hspace{6cm}- 2\e c^2Y_jY_k - \e^2 N_0 a_i^2 X_i + a_i(Y_i - Z_i), \\
a_i \frac{\d Y_i}{\d t} &=  -  \e a_kb_k X_jY_k - \e a_jb_j Y_jX_k 
+ c(a_k-a_j)Y_jY_k -a_iX_i-N_0a_i^2Y_i, \\
\frac{\d Z_i}{\d t} &=  - \e b_kX_j(Z_k-H_k) -\e  b_j(Z_j-H_j)X_k 
+ cY_j(Z_k-H_k) \\
&\hspace{6cm}- c(Z_j-H_j)Y_k + g_0 a_iX_i-K_0a_iZ_i + \F_i. \\
\eea
The above equations are written for each cyclic permutation of the set of indices $(1, 2, 3)$, namely, for  
\be\label{Eq_ijk}
(i, j, k) \in \{ (1,2,3), (2,3,1), (3,1,2)\}.
\ee
We refer to \cite{CLM16_Lorenz9D} for a detailed description of this model and its parameters. For our purpose, it is sufficient to know that the time, $t$, is an $\mathcal{O}(1)$-slow
time, and that $X_i$'s,$Y_i$'s, and $Z_i$'s are $\mathcal{O}(1)$-amplitudes for the divergent
velocity potential, streamfunction, and dynamic height, respectively.
In this setting $N_0$ and $K_0$ are rescaled damping coefficients in
the slow time.  The $\F_i$'s are $\mathcal{O}(1)$ control parameters
that, in combination with variations of $\e$, can be used to affect
regime transitions/bifurcations. In a general way, $\e$, can be
identified with the Rossby number.

Solutions of higher-order accuracy in $\e > 0$ that are entirely slow
in their evolution are, by definition, balanced solutions, and
\cite{Gent_McWilliams82} showed by construction several examples of
explicitly specified, approximate balanced models. One of these, the
Balance Equations (BE), was conspicuously more accurate than the
others when judged in comparison with apparently slow solutions of \eqref{Eq_L9D_rescale}. 
The BE approximation consists of a parameterization of the $X_i$'s and $Z_i$'s variables, in terms of 
the $Y_i$'s variables.  The $\boldsymbol{Z}$-component of this parameterization has an explicit expression. 
The $\boldsymbol{X}$-component of this parameterization, denoted by $\Phi$, is however obtained implicitly, by solving a system of differential-algebraic equations derived 
from Eq.~\eqref{Eq_L9D_rescale} under a balance assumption that consists of replacing the dynamical equation for the $X_i$'s by algebraic relations.   Eventually, we arrive at a 3D reduced system of ODEs, simply called the BE, and that takes the form
\be\label{BE_eq}
a_i \frac{\d Y_i}{\d t} =  -  \e a_kb_k \Phi_j(\boldsymbol{Y})Y_k - \e a_jb_j Y_j\Phi_k(\boldsymbol{Y})
+ c(a_k-a_j)Y_jY_k -a_i\Phi_i(\boldsymbol{Y})-N_0a_i^2Y_i,
 \ee
with $(i,j,k)$ as in  \eqref{Eq_ijk}. We refer to  \cite[Sec.~3.1]{CLM16_Lorenz9D} for a derivation.

%%%%%%%%%%%%%%%%%%%%%%%%%%%%%%%%%%%%%%%%%%%%%%%%%%%%%%%%%%
\begin{figure}[hbtp]
   \centering
   \includegraphics[height=0.4\textwidth, width=.9\textwidth]{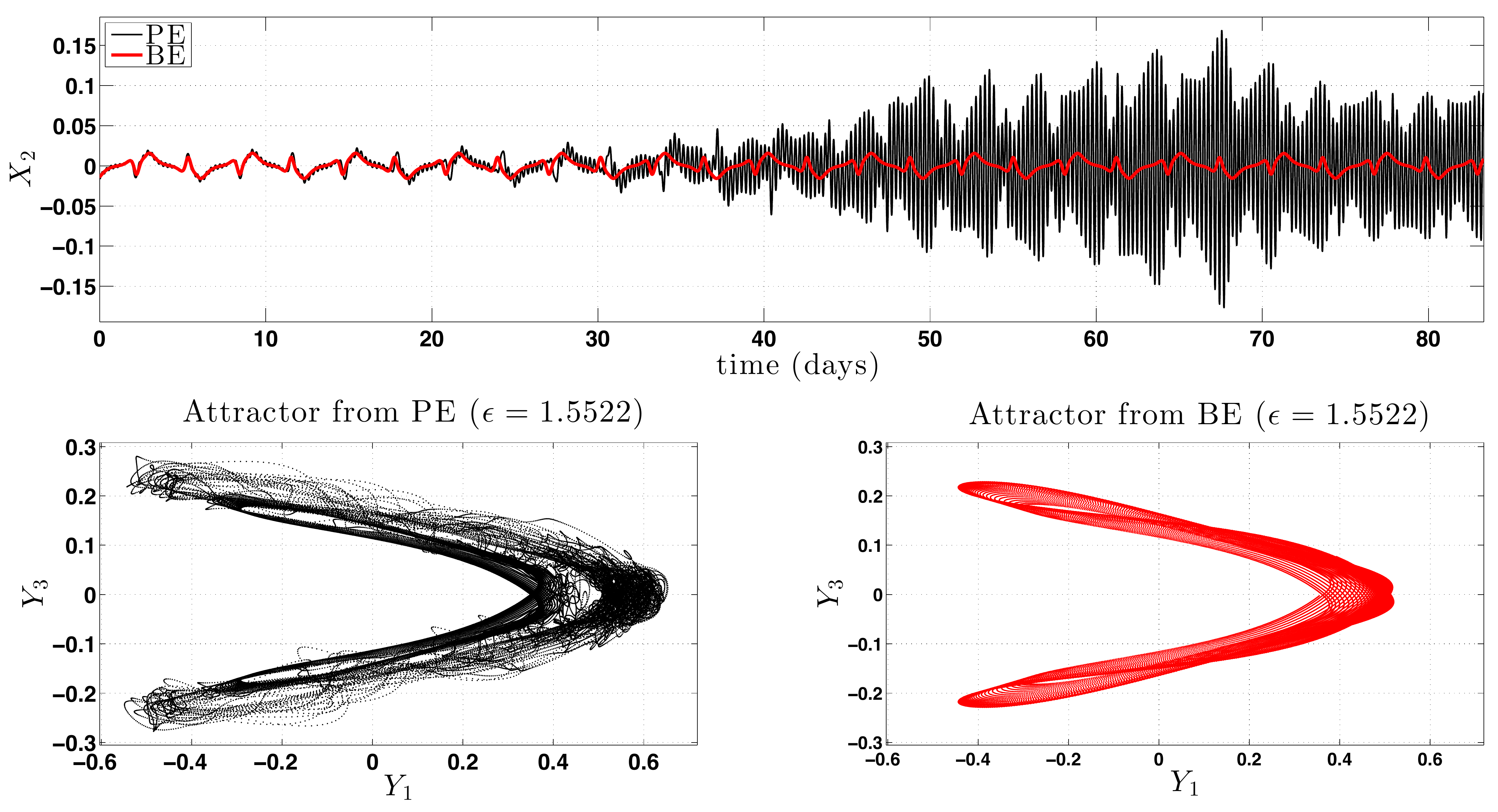}
  \caption{{\footnotesize {\bf Attractor comparison between PE and BE ({\hl reproduced from \cite{CLM16_Lorenz9D}, with permission from Elsevier}).} A slow-variable projection of the global attractor associated with Eq.~\eqref{Eq_L9D_rescale} (lower-left panel) and its approximation obtained from the BE reduced model (lower-right panel). Even in presence of energetic bursts of fast oscillations in the fast variables (here such an episode is shown in the upper panel for the $X_2$-variable (black curve)), the BE model \eqref{BE_eq} is able to capture the coarse-grained topological features of the projected attractor onto the slow variables. This is because the BE manifold provides a good approximation of the optimal PM given in \eqref{Def_h2} that averages here out (optimally) the fast oscillations.}}   \label{red-regime_216}
\vspace{-.2cm}
\end{figure}
%%%%%%%%%%%%%%%%%%%%%%%%%%%%%%%%

For certain Rossby numbers for which energetic bursts of fast oscillations occur in the course of time (occurring for $\epsilon >\epsilon^\ast$),  Chekroun {\it et al.} \cite{CLM16_Lorenz9D} have shown that the underlying BE manifold (associated with the BE parameterization of the $\boldsymbol{X}$- and $\boldsymbol{Z}$-variables), provides a very good approximation of the optimal PM for this problem, and thus of the conditional expectation in virtue of Theorem \ref{Thm_variational-pb2}, i.e.~the best approximation in the $\boldsymbol{Y}$-variable for which the ``fast'' $\boldsymbol{X}$- and $\boldsymbol{Z}$-variables  are averaged out.  In other words, the BE \eqref{BE_eq} provides a nearly optimal reduced vector field that averages out the fast oscillations contained in the $\boldsymbol{Y}$-variable.  Figure \ref{red-regime_216}, reproduced from \cite{CLM16_Lorenz9D}, illustrates this feature for the model \eqref{Eq_L9D_rescale}. The lower-right panel shows that the BE reduced model  is able to capture the coarse-grained topological features of the projected attractor onto the ``slow'' variables, $Y_1$ and $Y_3$, when compared with the projection onto the same variables of the attractor associated with the full Eq.~\eqref{Eq_L9D_rescale}. For the rest of this section we will use the BE as if it were the optimal PM. All the results presented hereafter correspond to $\epsilon =1.5522 >\epsilon^\ast$; see  \cite{CLM16_Lorenz9D}.

\begin{figure}[htbp]
		\includegraphics[width=.9\textwidth,height=0.6\textwidth]{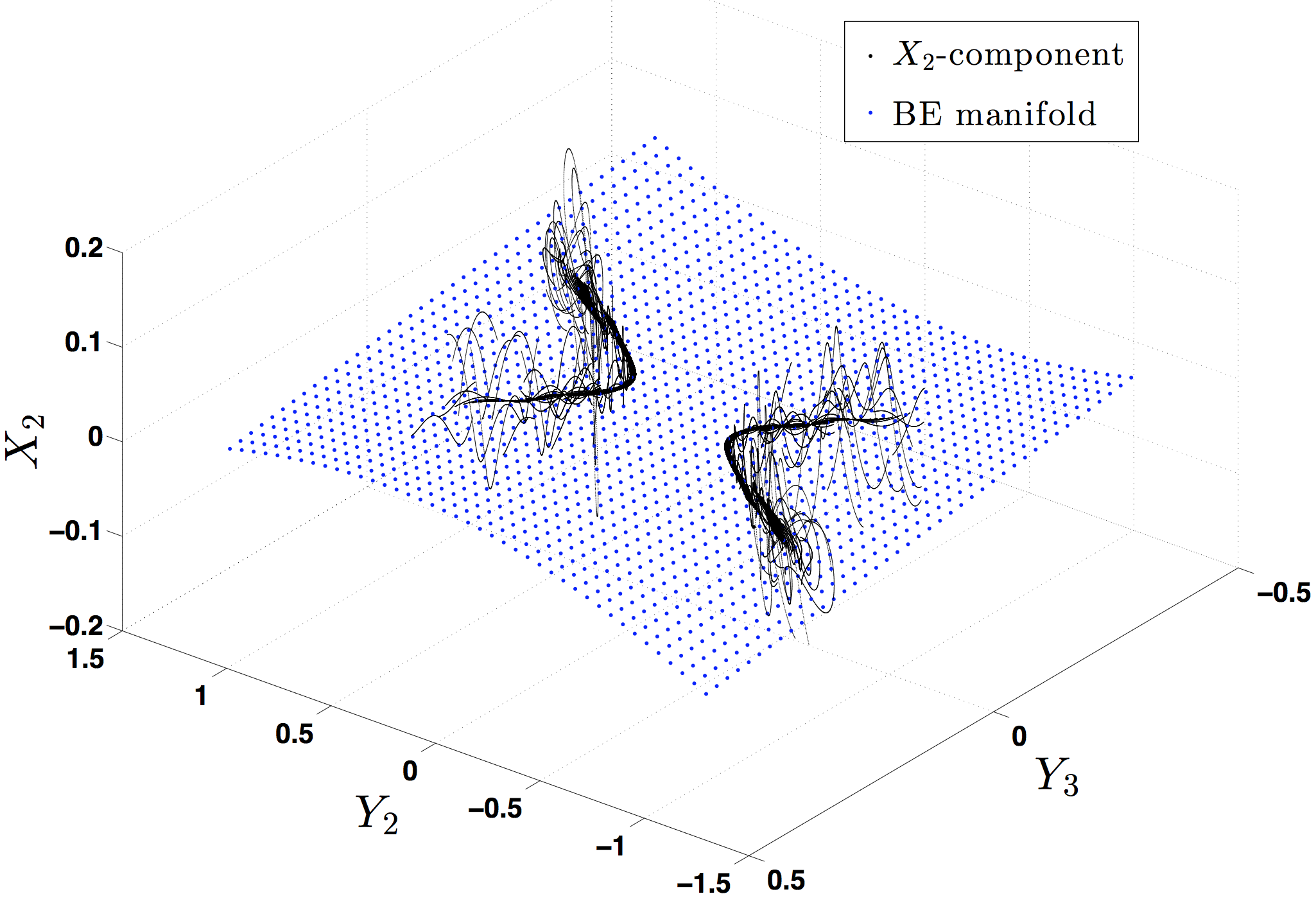}
	\caption{\footnotesize {\bf The BE manifold for the $X_2$-variable.}  Note that in order to obtain this representation, the $Y_1$-variable, involved also in the BE parameterization $\Phi$ along with the $Y_2$- and $Y_3$-variables,  has been set to its most probable value. The  black curve shows the resulting $X_2$-variable obtained after solving Eq.~\eqref{Eq_L9D_rescale} while the blue dots correspond to the BE parameterization $\Phi$ involved in \eqref{BE_eq}.}
	\label{Visualization_BE}
\end{figure}

The underlying BE manifold is a 6D manifold obtained as graph of a 6D-valued mapping of a 3D-variable ($\boldsymbol{Y}$), and as such 
only slices can be represented in 3D. Such a slice is shown in Fig.~\ref{Visualization_BE}. More exactly, it shows the  
$X_2$-variable as parameterized by the slow $Y_2$- and $Y_3$-variables. Note that in order to obtain this representation, the $Y_1$-variable, involved also in the BE parameterization $\Phi$ along with the $Y_2$- and $Y_3$-variables,  has been set to its most probable value conferring to  Fig.~\ref{Visualization_BE} a certain ``typicalness.'' This being kept in mind, the slice thus obtained of the BE manifold (and shown in Fig.~\ref{Visualization_BE}) will be simply called the BE manifold, for simplifying the discourse.

%%%%%%%%%%%%%%%%%%%%%%%%%%%%%%%%%%%%%%%%%%%%%%
%%%%%%%%%%%%%%%%%%%%%%%

\begin{figure}[htbp]
\centering
		\includegraphics[width=.85\textwidth,height=0.35\textwidth]{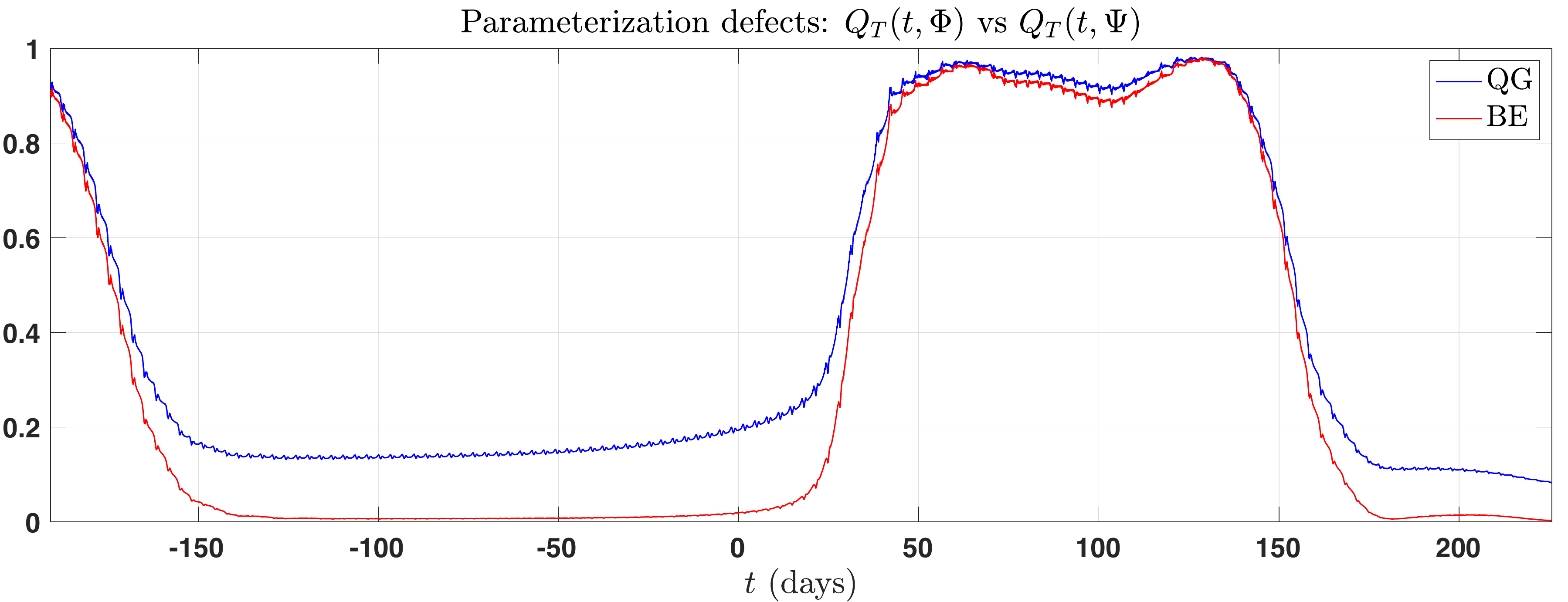}
	\caption{\footnotesize {\bf Parameterization defects of the BE manifold $\Phi$ and the QG manifold $\Psi$.}  Here the parameterization defects as given by \eqref{Eq_PD}, $Q_T(t,\Phi)$ (red curve) and $Q_T(t,\Psi)$ (blue curve), are computed for the BE manifold, $\Phi$, and  for the QG manifold $\Psi$ \cite[Eq.~(4.22)]{CLM16_Lorenz9D}; each with $T=80$ (for the rescaled system \eqref{Eq_L9D_rescale}) which corresponds to 10 days in the time-variable of the original Lorenz model \cite{Lorenz80}.}
	\label{Fig_BE_defect}
\end{figure}

As evidenced in Fig.~\ref{Visualization_BE}, a PE solution on the attractor --- as observed through the $X_2$-variable --- possesses an intricate transversal component to the BE manifold that seems to exclude its  
parameterization by a smooth manifold, whereas, at the same time, a substantial portion of the trajectory lies very close to the BE manifold. It is this latter portion of the dynamics that is well captured by the BE manifold and that allows for approximating the aforementioned conditional expectation. Here Fig.~\ref{Visualization_BE} reveals thus simple geometric features (not identified in \cite{CLM16_Lorenz9D}), which are responsible for the BE to provide in the space of slow variables, a vector field that approximates the PE dynamics. It does so by filtering out the (fast) oscillations contained in the PE solutions; {\mkr the fast dynamics} corresponding, in this representation, to the transversal part of the dynamics. Indeed, a closer inspection reveals that this transversal part of the dynamics corresponds exactly to the aforementioned burst of fast oscillations.    
This is confirmed by computing the parameterization defect. In that respect, Figure \ref{Fig_BE_defect} shows the parameterization defect $t\mapsto Q_T(t,\Phi)$ (given by \eqref{Eq_PD}) of the BE manifold $\Phi$ for a time horizon set to $T=80$ (for the rescaled system \eqref{Eq_L9D_rescale}) which corresponds to 10 days in the time-variable of the  original Lorenz model \cite{CLM16_Lorenz9D}. Figure \ref{Fig_BE_defect} shows that 
$Q_T(t,\Phi)$ oscillates, as $t$ evolves, between values right above zero and right below one (red curve). The rising of values taken by $Q_T(t,\Phi)$ occurs over time windows for which the parameterized $\boldsymbol{X}$-variable contains a significant fraction of the total energy, such as manifested by bursts of fast oscillations in the $X_2$-variable shown in the upper panel of Fig.~\ref{red-regime_216} between 40 and 80 days.  To the contrary, when the PE solutions get very close to the BE manifold, the dynamics is almost slaved to this manifold and   
$Q_T(t,\Phi)\approx 0$.

%%%%%%%%%%%%%%%%%%%%%%%%%%%%%%%%%%%%%%%%%%%%%%
%%%%%%%%%%%%%%%%%%%%%%%%%%

Complementarily, the parameterization  defect $Q_T(t,\Psi)$ has been computed for the standard Quasigeostrophic (QG) manifold  \cite[Eq.~(4.22)]{CLM16_Lorenz9D} that can be derived for $\epsilon=0$ and is associated with the famous quadratic Lorenz
system \cite{lorenz1963deterministic}; see \cite[Sec.~4.2]{CLM16_Lorenz9D}. Here again a similar behavior is observed for 
$Q_T(t,\Psi)$ (blue curve in Fig.~\ref{Fig_BE_defect}) with the noticeable difference that $Q_T(t,\Psi)$ stays further away from zero than $Q_T(t,\Phi)$ does,  as $t$ evolves.

%%%%%%%%%%%%%%%%%%%%%%%%%%%%%%%%%%%%%%%%%%%%%%
\begin{figure}[htbp]
		\includegraphics[width=.98\textwidth,height=0.48\textwidth]{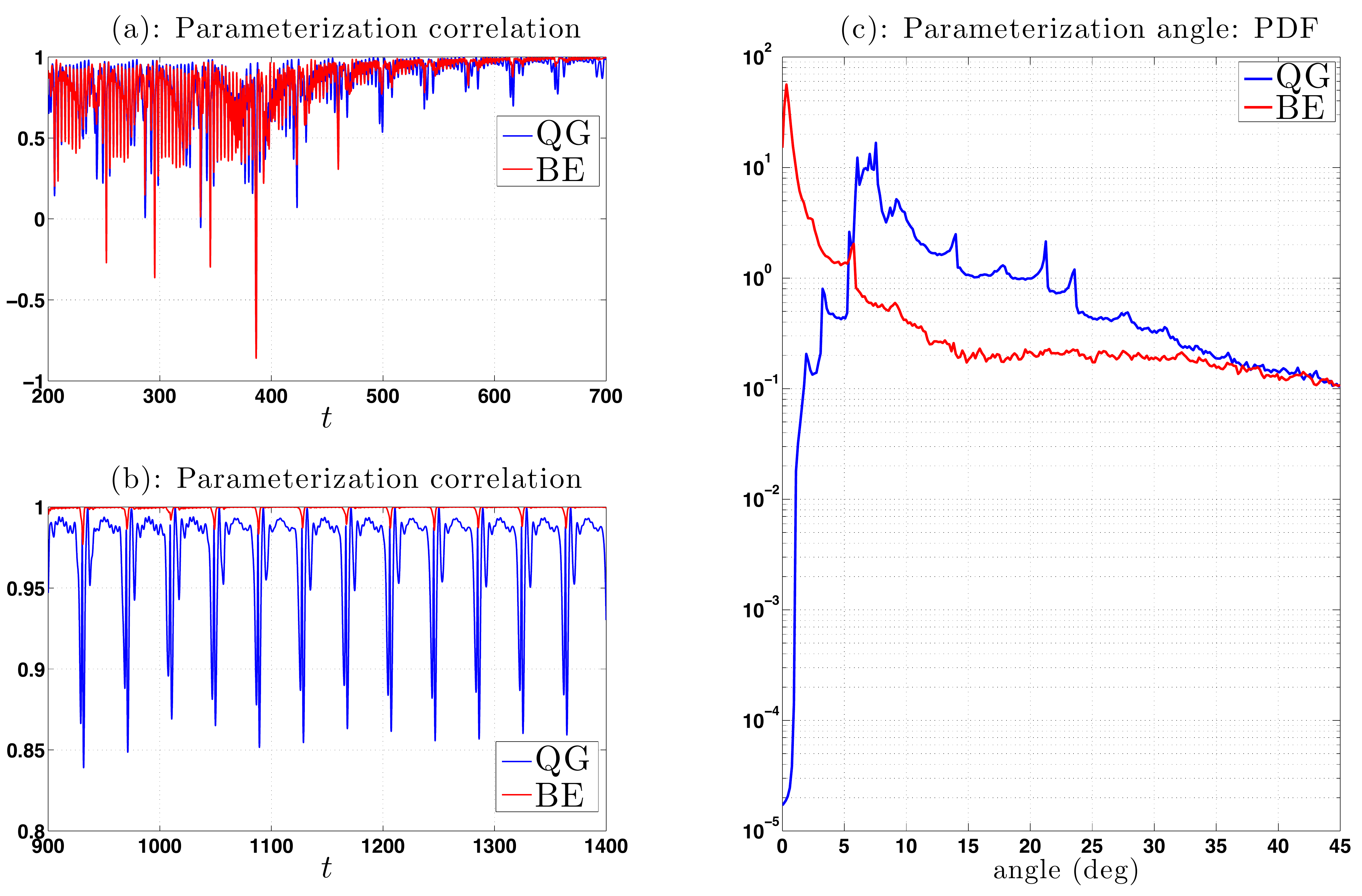}
	\caption{\footnotesize {\bf Parameterization correlation and angle.} The parameterization correlation, $c(t)$ given by \eqref{Eq_corr_param}, is shown for the BE manifold ($\Psi=\Phi$, red curve) and the QG manifold (``$\Psi=$QG manifold,'' blue curve), over two consecutive time windows for panels (a) and (b); {\mkr the range of fluctuations over the 2nd window (panel (b)) being smaller to the range shown in the 1st window (panel (a))}. The time-episode shown in panel (a) corresponds to the presence of energetic bursts of fast oscillations in the solutions ($Q_T \approx 1$ for the BE), whereas panel (b) corresponds to a time-episode devoid of such oscillations ($Q_T \approx 0$ for the BE). The PDFs of the corresponding parameterization angle $\alpha(t)$ given by \eqref{Eq_alpha}, estimated after long integration of  Eq.~\eqref{Eq_L9D_rescale}, are shown in panel (c).}\label{Fig_Param_angle_BE}
\end{figure}
%%%%%%%%%%%%%%%%%%%%%%%%%%

The parameterization correlation, $c(t)$ given by \eqref{Eq_corr_param}, has been also computed for the BE and the QG manifolds. 
The results are shown in Panels (a) and (b) of Fig.~\ref{Fig_Param_angle_BE}, over  different time intervals. Although when an episode of fast (gravity waves) oscillations occurs in the PE solutions, the parameterization correlation can deviate substantially from 1 for the BE and QG manifolds (panel (a)),  the parameterization correlation gets, comparatively, much closer to 1 for the BE than for the QG manifold over time intervals for which  the slow, Rossby waves dominate the dynamics (panel (b)). This phenomenon is confirmed statistically at the level of the probability distribution for the corresponding parameterization angle, $\alpha(t)=\arccos(c(t))$. The PDF of the latter is much more skewed towards zero for the BE manifold than for the QG manifold supporting thus, at a quantitative level, the visual rendering of Fig.~\ref{Visualization_BE} which suggests  that a substantial portion of the PE trajectory lies very close to the BE manifold. More precisely, Fig.~\ref{Fig_Param_angle_BE}-(c) shows that the mode of the PDF of $\alpha(t)$ (i.e.~the value that appears most often) for the BE manifold is located very close to zero, whereas $\alpha(t)$ almost never reaches such a level of proximity to zero for the QG manifold.     
In that sense, the BE manifold is a manifold that is close to be locally invariant in the sense of (i) of Sec.~\ref{Sec_loc_invman}, that is a slaving relationship like \eqref{Eq_inv_local} almost takes place over time, while being brutally violated from time to time (transversal part of the PE dynamics to the manifold; see Fig.~\ref{Visualization_BE}).

Thus the BE manifold provides an example of a manifold that  is close to be locally invariant and that provides a (nearly optimal) PM.
However, nothing excludes the existence of dynamics that although getting very close to a given manifold over certain time windows  (almost slaving situation), experiences excursions far away from it so often that in average the parameterization defect gets greater than one, making this manifold to be a non-parameterizing one. Situations for which the dynamics lies in the vicinity of a given manifold (without large excursions) is also a favorable context for this manifold to be a PM; see Sec.~\ref{Sect_chaotic-RB} below for such an example.

Noteworthy are also the tails of the PDFs of the parameterization angle $\alpha(t)$ for both, the BE and  QG manifolds, which do not drop off suddenly as $\alpha$ increases: this is symptomatic of the fact that the PE solutions get frequently far away from these manifolds as time evolves. As a comparison, we refer to Sec.~\ref{Sect_chaotic-RB} below for an example of parameterization angle $\alpha$ whose PDF drops suddenly as  $\alpha$ increases.      

Although enlightening, this example of (excellent) approximation of the optimal PM (and thus of the conditional expectation) that the BE manifold provides, exploits specific aspects of the problem at hand, encapsulated in the very derivation of the BE manifold. The question of efficient dynamically-based formulas for the approximation of an optimal PM in a general context, thus remains. The next section addresses this issue.

\section{Parameterizing manifolds and mode-adaptive minimization: Dynamically-based formulas} \label{Sect_PM_formulas}

 In this section we derive dynamically-based formulas for designing parameterizing manifolds in practice. 
The formulas derived in Sec.~\ref{Sect_PM_with_forcing} below take their origin in the pullback representation \eqref{PB_rep} (in Theorem \ref{thm:h1_general}) and the associated backward-forward system \eqref{Eq_BF} that arise in the approximation theory of invariant manifolds revisited in Part I. The parametric class of leading interaction approximation (LIA)  of the high modes obtained this way is completed by another parametric class built from the quasi-stationary approximation (QSA) in Sec.~\ref{Sec_FMTtau}; close to the first criticality, the QSA is an approximation to the LIA, but  differs as one moves away from criticality. 
We also make precise hereafter the corresponding minimization problems to solve in order to optimize our parameterizations in practice, within a mode-adaptive optimization procedure (Sec.~\ref{Sec_mode_adaptive}).

\subsection{Backward-forward method:  General considerations} \label{Sect_BF_method_basics}
%%%%%%%%%%%%%%%%%%%%%%%%%%%%%%%%%%%%%%%%%%%%%
We first  show that the parameterization $h^{(1)}_\tau$ given in \eqref{Eq_h1_tau}, as obtained by finite-time integration of the backward-forward system \eqref{Eq_BF}, satisfies an equation analogous to Eq.~\eqref{h1_eqn} satisfied by $h_k$. 
\begin{lem}\label{important_lemma}
The manifold function  $h^{(1)}_\tau$ defined by \eqref{Eq_h1_tau} satisfies the following system of first order quasilinear PDEs:
\be \label{eq:invariance_h1_tau}
\mathcal{L}_A [h] (\xi)= \Pi_{\s} G_k(\xi) - e^{\tau A_{\s}} \Pi_{\s} G_k(e^{-\tau A_{\c}} \xi).
\ee
with $\mathcal{L}_A [h] (\xi)=Dh(\xi) A_\c \xi- A_\s h(\xi)$ and $A_\c$, $A_\s$ defined in \eqref{Ac_As_def}. 
\end{lem}

\bp
In \eqref{Eq_h1_tau}, by replacing $\xi$ with $e^{t A_{\c}} \xi$, we get
\bea
\Phi(t)=h^{(1)}_\tau (e^{t A_{\c}} \xi)  &= \int_{-\tau}^0 e^{-s A_{\s}} \Pi_{\s} G_k(e^{s A_{\c}} e^{t A_{\c}} \xi) \d s \\
& = \int_{-\tau}^0 e^{-s A_{\s}} \Pi_{\s} G_k(e^{(s +t) A_{\c}} \xi) \d s \\
& = \int_{t-\tau}^t e^{-(s'-t) A_{\s}} \Pi_{\s} G_k(e^{s' A_{\c}} \xi) \d s'.
\eea
We obtain then
\bea \label{dPhi_eq1}
\frac{\d \Phi(t)}{\d t} & = \Pi_{\s} G_k( e^{t A_{\c}} \xi) - e^{\tau A_{\s}}  \Pi_{\s} G_k( e^{(t-\tau) A_{\c}} \xi) \\
& \quad + A_{\s} \int_{t-\tau}^t e^{-(s'-t) A_{\s}} \Pi_{\s} G_k(e^{s' A_{\c}} \xi) \d s' \\
& = \Pi_{\s} G_k( e^{t A_{\c}} \xi) - e^{\tau A_{\s}}  \Pi_{\s} G_k( e^{(t-\tau) A_{\c}} \xi) + A_{\s} \Phi(t).
\eea

On the other hand, we also have
\be \label{dPhi_eq2}
\frac{\d \Phi(t)}{\d t} = [Dh^{(1)}_\tau (e^{t A_{\c}} \xi)] A_{\c} e^{t A_{\c}} \xi.
\ee

Equation \eqref{eq:invariance_h1_tau} follows by equating the RHSs of \eqref{dPhi_eq1} and \eqref{dPhi_eq2} and by taking the limit $t \rightarrow 0$. 
\ep
%%%%%%%%%%%%%%%%%%%%%%%%%%%%%%%%%%%%%%%%%%%%%%%%%%%
This lemma provides the equation satisfied by the parameterization $h^{(1)}_\tau$ given by \eqref{Eq_h1_tau}. However this parameterization is built from the backward-forward system \eqref{Eq_BF} associated with Eq.~\eqref{Eq_ODEs} that does not include forcing terms, unlike for more  general systems of ODEs such as Eq.~\eqref{Eq_ODE_gen} dealt with in Sec.~\ref{Sect_PM_reduction}.

To extend the parameterization $h^{(1)}_\tau$ to systems that include forcing terms, we thus naturally seek for solution 
of the backward-forward system associated with Eq.~\eqref{Eq_ODE_gen}, namely 
\begin{subequations} \label{Eq_BF_quad0}
\begin{align}
& \frac{\mathrm{d} y^{(1)}_{\c}}{\d s} =  A_\c y^{(1)}_{\c} + \Pi_{\c} F, && s \in[ -\tau, 0],    \label{BF1_quad0} \\
& \frac{\mathrm{d} y^{(1)}_{\s}}{\d s} = A_\s y_{\s}^{(1)}  +  \Pi_{\s} G_k\big(y^{(1)}_{\c}\big) + \Pi_{\s} F, &&  s \in [-\tau, 0], \label{BF2_quad0}\\
& \mbox{with } y^{(1)}_{\c}(s)\vert_{s=0} = \xi, \mbox{ and } y_{\s}^{(1)}(s)\vert_{s=-\tau}=0.
\end{align}
\end{subequations}
Here $\Pi_\s=\mbox{Id}_{\mathbb{C}^N}-\Pi_\c$ with $\Pi_{\c}$ denoting the  canonical projector onto the eigensubspace, $E_\c$, spanned by the dominant eigenmodes of $A$.

By going through similar calculations than for the proof of Lemma \ref{important_lemma}, 
the high-mode solution of \eqref{Eq_BF_quad0}, $y^{(1)}_{\s}[\xi](0; -\tau)$, denoted here by $\Psi^{(1)}_\tau (\xi)$, satisfies then 
\bea\label{Eq_h1_taugen}
 \mathcal{L}_A [\Psi^{(1)}_\tau] (\xi)+D\Psi^{(1)}_\tau(\xi) \Pi_\c F= \Pi_{\s} G_k(\xi) - e^{\tau A_{\s}} \Pi_{\s}& G_k(S_F(-\tau)\xi)\\
 &  +(\mbox{Id}- e^{\tau A_{\s}})  \Pi_{\s} F,
\eea
with 
\be\label{Def_SF}
S_F(t) \xi =e^{t A_\c} \xi -A_\c^{-1}(\mbox{Id}-e^{t A_\c})\Pi_{\c} F.
\ee
Obviously $\Psi^{(1)}_\tau=h^1_\tau$ when $F\equiv 0$. 

In practice, in order to find an explicit expression of the parameterization $\Psi^{(1)}_\tau$, one prefers to solve \eqref{Eq_BF_quad0} rather than solving Eq.~\eqref{Eq_h1_taugen} directly. Note that we could have adopted the same strategy for deriving the formulas of Theorem \ref{thm:h1}, i.e.~by solving the backward-forward system \eqref{Eq_BF} in this case. 

The manifold $\mathfrak{M}_{\tau}$ associated with $\Psi^{(1)}_\tau$ possesses a natural geometric interpretation.  Given a solution $y(t)$ of Eq.~\eqref{Eq_ODE_gen} and denoting by  $U_\tau y_\c(t)$ the lift of $y_\c(t)$ onto the manifold $\mathfrak{M}_{\tau}$, i.e.~ $U_\tau y_\c(t)=y_\c(t)+\Psi^{(1)}_{\tau}(y_\c(t))$, we obtain
\bea
\overline{ \textrm{dist}(y(t),\mathfrak{M}_{\tau})^2} \leq \overline{\|y(t)-U_\tau y_\c (t)\|^2}=\overline{\|y_\s(t) -\Psi^{(1)}_{\tau}(y_\c(t))\|^2},
\eea
where the overbar denotes the time average over $[0,T].$
In other words,
\bea\label{Eq_interp}
\overline{ \textrm{dist}(y(t),\mathfrak{M}_{\tau})^2} \leq \mathcal{Q}_T(\Psi^{(1)}_{\tau}),
\eea
with $\mathcal{Q}_T$ that denotes the parameterization defect 
\be\label{Eq_Q_Tbis}
\mathcal{Q}_T(\Psi^{(1)}_{\tau})=\frac{1}{T}\int_{0}^T \norm{y_\s(t) -\Psi^{(1)}_{\tau}(y_\c(t))}^2 \d \, t.
\ee

Thus, we understand a practical advantage in restricting ourself to the $\Psi^{(1)}_\tau$-class of parameterizations instead of the more general $\mathcal{E}$-class considered in \eqref{E_class}. Indeed, once an explicit expression for $\Psi^{(1)}_\tau$ is derived, it allows us to greatly simplify the minimization problem involved in Theorem \ref{Thm_variational-pb}, by replacing it with the minimization in the scalar variable $\tau$ of the cost functional  $\mathcal{Q}_T$ given by \eqref{Eq_Q_Tbis}. Although the corresponding minimizer is a priori suboptimal compared to the more general minimization problem \eqref{Eq_variational-pb}, we will see in applications that it provides in various instances an efficient parameterization.

Furthermore, based on \eqref{Eq_interp}, minimizing $\mathcal{Q}_T(\Psi^{(1)}_{\tau})$ in the $\tau$-variable has the following useful interpretation: it forces, within the $\Psi^{(1)}_{\tau}$-parametrization class, the manifold $\mathfrak{M}_{\tau}$ to get the closest  to the trajectory $y(t)$, in a least-square sense. As mentioned earlier, an alternative approach, the AIM approach, has been proposed in the literature, but the latter is {\it asymptotic} in essence rather than the PM approach presented here which is {\it variational}.  The AIM approach consists indeed of seeking for a family of manifolds, $\mathcal{M}_m$,  for which $\overline{\textrm{dist}(u(t),\mathcal{M}_m)}$ vanishes to zero as $m=\textrm{dim}(\mathcal{M}_m)\rightarrow \infty$; see e.g.~\cite{Temam88,Temam_IMfamily89,titi1990approximate,debussche1992construction}.  In contradistinction, the PM approach consists for a given reduced dimension, $m$, of seeking for a manifold $\mathfrak{M}$ within a certain parametric class of dynamically-based parameterizations, for which  $\overline{ \textrm{dist}(u(t),\mathfrak{M})}$ is minimized. 

Thus, given a reduced dimension, $m$, seeking for the best approximation within a parameterization class is at the core of the PM approach and, as shown in Sec.~\ref{Sect_PM_reduction}, is quintessential to address closure problems, in the sense that it relates naturally to the conditional expectation i.e.~to the best closure that can be derived out of nonlinear parameterizations alone; see Theorem \ref{Thm_variational-pb2}.

\br\label{Estimates_Q_literature}
Given the limitations on our ability to estimate the norms, it is in general hard to derive sharp estimates of 
$Q_T(\Psi_{\bftau}^{(1)})$. Nevertheless, some related estimates have been produced  
about $\overline{ \textrm{dist}(y(t),\mathcal{M})^2}/\|y(t)\|^2$, for the 2D Navier-Stokes equations   \cite{foias1991approximate,chae1992ensemble} when $\mathcal{M}$ denotes the manifold associated with the quasi-stationary approximation; see  \eqref{Eq_QSA} below.  
\er

\subsection{Mode-adaptive optimization}\label{Sec_mode_adaptive} 

Although the minimization in the scalar variable $\tau$ of the cost functional $\mathcal{Q}_T$  in \eqref{Eq_Q_Tbis} is more appealing  than solving the general minimization problem  \eqref{Eq_variational-pb}, we may suffer from the fact that the parameter $\tau$ to be optimized, is chosen globally, irrespectively e.g.~to the content of energy of a particular high mode to parameterize.   To better account for the distribution of energy across the modes,  we propose instead to optimize parameterizations of the form
\be\label{Eq_Phi_tau_opt0}
\Phi^{(1)}_{\bftau}(\xi)= \sum_{n = m+1}^N \Phi_n (\tau_n, \boldsymbol{\beta}, \xi) \boldsymbol{e}_n, \;\; \bftau=(\tau_{m+1},\cdots,\tau_N),
\ee
in the multivalued $\bftau$-variable. We emphasize that each  parameterization $\Phi_n$  depends only  on $\tau_n$ (and not the other $\tau_p$'s for $p\neq n$), and thus each  $\Phi_n$ may be optimized independently from each other.

This way, we are left for each of the  $n^{\textrm{th}}$  mode, with a parameterization  to optimize, $\Phi_n (\tau_n, \boldsymbol{\beta}, \xi)$,  that is a scalar function of the scalar variable $\tau_n$.  
Following Sec.~\ref{Sect_BF_method_basics} and assuming $A$ diagonalizable (in $\mathbb{C}^N$), we obtain $\Phi_n (\tau_n, \boldsymbol{\beta}, \xi)$, for each $m+1\leq n \leq N$, as the high-mode part $y^{(1)}_{n}$ of the solution (at $s=0$) to  the backward-forward system
\begin{subequations} \label{Eq_BF_introb}
\begin{align}
& \frac{\mathrm{d} y^{(1)}_{\c}}{\d s} =  A_\c y^{(1)}_{\c} + \Pi_{\c} F, && s \in[ -\tau_n, 0],    \label{BF1_introa} \\
& \frac{\mathrm{d} y^{(1)}_{n}}{\d s} = \beta_n y_{n}^{(1)}  +  \Pi_{n} G_k\big(y^{(1)}_{\c}\big) + \Pi_{n} F, &&  s \in [-\tau_n, 0], \label{BF2_introb}\\
& \mbox{with } y^{(1)}_{\c}(s)\vert_{s=0} = \xi, \mbox{ and } y_{n}^{(1)}(s)\vert_{s=-\tau_n}=0,
\end{align}
\end{subequations}
in which the RHS in Eq.~\eqref{BF2_quad0} has been replaced by $\beta_n y_{n}^{(1)}  +  \Pi_n G_k\big(y^{(1)}_{\c}\big) + \Pi_n F$. 
Here $\Pi_{n} X=\langle X, \boldsymbol{e}_n^\ast\rangle$, for any $X$ in $\mathbb{C}^N$.

Explicit formulas of the $\Phi_n (\tau_n, \boldsymbol{\beta}, \xi)$'s are given in Sec.~\ref{Sect_PM_with_forcing} below when $G_k$ is a  quadratic nonlinearity. We show hereafter that minimizing for each $n$ the parameterization defect naturally associated with $\Phi_n$ leads to an optimal parameterization, $\Phi^{(1)}_{\bftau}$, with a clear geometrical interpretation. To do so --- given a fully resolved solution $y(t)$ of the underlying $N$-dimensional ODE system \eqref{Eq_quadODEs} available over a training interval $[0,T]$ --- we consider for each $n\geq m+1$, the parameterization defect
\be\label{Q_n_defined}
\mathcal{Q}_n(\tau_n,T)=\frac{1}{T}\int_0^T \big| \Pi_n y(t)- \Phi_n(\tau, \boldsymbol{\beta}, y_\c(t))\big|^2 \d t,
\ee
with $y_\c(t)=\Pi_{\c} y(t)$.

Denoting by $\mathfrak{M}_{\bftau}$ the manifold associated with the parameterization  $\Phi^{(1)}_{\bftau}$ given by \eqref{Eq_Phi_tau_opt0}, 
we have 
\bea
\overline{ \textrm{dist}(y(t),\mathfrak{M}_{\bftau})^2}& \leq \overline{\bigg\|y(t) -\bigg(y_\c(t)+\sum_{n\geq m+1} \Phi_n(\tau_n,\bfbeta, y_\c(t)) \boldsymbol{e}_n)\bigg)\bigg\|^2}\\
& \leq \overline{\bigg\|\sum_{n\geq m+1} (\Pi_n y(t)-\Phi_n(\tau_n,\bfbeta, y_\c(t)))\boldsymbol{e}_n\bigg\|^2}.
\eea
Taking the eigenvectors of $A$ to be normalized, we are thus left, thanks to the triangular inequality, with the following estimate
\be\label{Eq_geom_interp}
\overline{ \textrm{dist}(y(t),\mathfrak{M}_{\bftau})^2} \leq \underset{n\geq m+1}\sum  \overline{\bigg|\Pi_n y(t)-\Phi_n(\tau_n,\bfbeta, y_\c(t))\bigg|^2 }=\underset{n\geq m+1}\sum \mathcal{Q}_n(\tau_n,T).
\ee
Thus minimizing each $\mathcal{Q}_n(\tau_n,T)$ (in the $\tau_n$-variable) is a natural idea to enforce closeness of $y(t)$ in a least-square sense to  
the corresponding manifold $\mathfrak{M}_{\bftau}$.  Note that we could have chosen  to minimize $\mathcal{Q}_T$  as given in \eqref{Eq_Q_Tbis} but with   $\Phi^{(1)}_{\bftau}$ replacing $\Psi^{(1)}_{\tau}$. The resulting minimization would become however 
more challenging in high-dimension as it would require to  minimize $\mathcal{Q}_T(\Phi^{(1)}_{\bftau})$ in the multdimensional variable $\bftau$. Except when the basis $\{\boldsymbol{e}_j\}_{j=1}^N$ is orthonormal (see \eqref{Eq_Q_T}), the two approaches 
are not equivalent, i.e.~minimizing $\mathcal{Q}_T(\Phi^{(1)}_{\bftau})$ in the vector $\bftau$, vs minimizing $\mathcal{Q}_n(\tau_n,T)$ in the scalar $\tau_n$ for each $n\geq m+1$. We opted for the latter as a simple algorithm can be proposed to minimize  $\mathcal{Q}_n$ efficiently; see Appendix \ref{Sect_gradient_descent}.  Nevertheless, even in this scalar case, a certain care must be paid, as the mapping $\tau \mapsto \mathcal{Q}_n(\tau,T)$ is not guaranteed to be convex; see Sec.~\ref{Sect_RBC}. Furthermore, depending on the dynamics (and the training interval $[0,T]$) local minima may appear that require also a special care in order to properly design an efficient parameterization for the problem at hand; see Remark \ref{Remark_training_interval} below.

\subsection{Parametric Leading-Interaction Approximation}\label{Sect_PM_with_forcing}
In this section, we focus on the case of quadratic nonlinear interactions under constant forcing, for which we derive parameterization formulas by solving the backward-forward systems \eqref{Eq_BF_introb} (for $G_k$ quadratic) presented in Sec.~\ref{Sec_mode_adaptive} above. 
Our approach allows for deriving parameterizations that take into account interactions between the forcing components and  the nonlinear terms, at the leading order.  As already pointed out in Sec.~\ref{Sec_mode_adaptive}, these parameterizations are conditioned on the 
choice of a finite collection $\bftau$ of scalar parameters. For these reasons we will refer to  $\Phi_{\bftau}^{(1)}$ given by \eqref{Eq_Phi_tau_opt} as the 
{\it parametric Leading-Interaction Approximation (LIA)}. As $\bftau$ varies, the corresponding class of parameterizations will be referred to  as the $\Phi^{(1)}$-class or simply the LIA class.

The ODE system considered here is of the form: 
\be \label{Eq_quadODEs}
\frac{\d y}{\d t} = A y +  B(y,y) + F,   \;\;  y \in \mathbb{C}^N,
\ee
where $A$ is an $N\times N$ matrix with complex entries, $B$ denotes quadratic nonlinear interactions with complex coefficients, and $F$ is a constant forcing term in $\mathbb{C}^N$.

Given the spectral elements $(\beta_j,\boldsymbol{e}_j)$ of the matrix $A$ that we assume diagonalizable (in $\mathbb{C}^N$), we decompose the state space into resolved and unresolved subspaces as follows
\be % \label{Eq_decomp_CN}
\mathbb{C}^N = E_{\c}  \oplus E_{\s},
\ee
where
\bea %\label{eq:subspaces}
E_{\c} &= \mathrm{span}\{\boldsymbol{e}_i : i = 1, \cdots, m\}, \\
E_{\s} &= \mathrm{span}\{ \boldsymbol{e}_i: i = m+1, \cdots, N\},
\eea
see also \eqref{eigen_A}--\eqref{eq:subspaces}.

We define the projection of a vector $X$ in $\mathbb{C}^N$ onto $\boldsymbol{e}_j$ as follows
\be\label{Proj_formula}
\Pi_{j} X =\langle X, \boldsymbol{e}_j^\ast\rangle,
\ee
with $\{\boldsymbol{e}_j^\ast\}$ denoting the eigenvectors of the conjugate transpose, $A^\ast$.
The projectors $\Pi_\c$ is then explicitly given by
\be\label{Def_Pic}
\Pi_\c X=\sum_{j=1}^m (\Pi_j X) \boldsymbol{e}_j \mbox{ and } A_\c=\textrm{diag}(\beta_1,\cdots,\beta_m). 
\ee

Recall that according to the convention \eqref{cutoff_cond} (of Sec.~\ref{Sec_loc_invman}) made throughout this article, 
the reduced state space $E_\c$ is spanned by modes that come either as conjugate pairs or as a real eigenvector. 
As a result, $\Pi_\c X$ is real if $X$ is real.

For each given unresolved mode $\boldsymbol{e}_n$ ($n\geq m+1$), a parameterization $y^{(1)}_{n}$ of the corresponding unresolved variable 
\be
Y_n = \Pi_n y,
\ee
is obtained from the following backward forward system: 
\begin{subequations} \label{Eq_BF_quad}
\begin{align}
& \frac{\mathrm{d} y^{(1)}_{\c}}{\d s} =  A_\c y^{(1)}_{\c} + \Pi_{\c} F, && s \in[ -\tau, 0],    \label{BF1_quad} \\
& \frac{\mathrm{d} y^{(1)}_{n}}{\d s} = \beta_n y_{n}^{(1)}  +  \Pi_{n} B\big(y^{(1)}_{\c}, y^{(1)}_{\c} \big) + \Pi_{n} F, &&  s \in [-\tau, 0], \label{BF2_quad}\\
& \mbox{with } y^{(1)}_{\c}(s)\vert_{s=0} = \xi \in E_\c, \mbox{ and } y_{n}^{(1)}(s)\vert_{s=-\tau}=0.
\end{align}
\end{subequations}

Note that the solution to  \eqref{BF1_quad} is given by: 
 \be
 y^{(1)}_{\c}(t) = e^{A_{\c} t} \xi - \int_{t}^0 e^{A_{\c} (t-s)}\Pi_{\c} F \d s, \qquad t \in [-\tau, 0],
 \ee
 which admits the following explicit expression:
 \be \label{Eq_uc_soln}
 y^{(1)}_{\c}(t) = \sum_{j=1}^m \Big( e^{\beta_j t} \xi_j + \gamma_j(t) \Pi_jF \Big)  \boldsymbol{e}_j,
 \ee
 where
 \be
 \gamma_j(t) = \begin{cases}
 \frac{\exp(\beta_j t) - 1}{\beta_j}, & \text{if $\beta_j \neq 0$}, \\
 t, & \text{otherwise}. 
 \end{cases} 
 \ee

 The solution to \eqref{BF2_quad} is given by: 
 \be\label{Eq_preformula}
y^{(1)}_{n}[\xi](t) = \int_{-\tau}^t e^{\beta_n(t-s)} \Pi_n B( y^{(1)}_{\c}(s),  y^{(1)}_{\c}(s)) 
 \d s +  \int_{-\tau}^t e^{\beta_n(t-s)}  \Pi_{n} F \d s, \qquad t \in [-\tau, 0],
 \ee
which leads to the following parameterization for the high mode $e_n$:
\be \label{eq_Phi_n_abstract}
\Phi_n(\tau, \xi) = \int_{-\tau}^0 e^{-\beta_n s} \Pi_n B( y^{(1)}_{\c}(s),  y^{(1)}_{\c}(s)) 
 \d s +  \int_{-\tau}^0 e^{-\beta_n s }  \Pi_{n} F \d s. 
\ee
By using \eqref{Eq_uc_soln} in the nonlinear term $\Pi_n B( y^{(1)}_{\c}(s),  y^{(1)}_{\c}(s))$ and expanding this term, the first integral $I$ in the RHS of \eqref{eq_Phi_n_abstract} becomes after simplification 
\bea
I= \sum_{i, j = 1}^m  U_{i,j}^n(\tau, \boldsymbol{\beta}) B_{i,j}^n F_{i}F_{j} +\sum_{i, j = 1}^m & V_{i,j}^n(\tau, \boldsymbol{\beta}) F_{j} (B^n_{i,j} + B^n_{j,i})  \xi_{i} \\
&+ \sum_{i, j = 1}^m D_{i,j}^n(\tau, \boldsymbol{\beta}) B_{i,j}^n \xi_{i} \xi_{j},
\eea
where 
\be\label{B_ijn}
B_{i,j}^n=\langle B(\boldsymbol{e}_{i}, \boldsymbol{e}_{j}), \boldsymbol{e}_n^\ast \rangle,
\ee
the coefficients $D_{i,j}^n(\tau, \boldsymbol{\beta})$  of the quadratic terms (in the $\xi$-variable) are given by 
\bea \label{Eq_D_term0}
D_{i, j}^n(\tau,\bfbeta)=  \begin{cases}
\frac{1 - \exp\big(-(\beta_{i} + \beta_{j} - \beta_{n})\tau\big)}{\beta_{i} + \beta_{j} - \beta_{n}}, & \text{if $\beta_{i} + \beta_{j} - \beta_{n}\neq 0$,} \\
\tau, & \text{otherwise},
\end{cases}
\eea
while the coefficients in the constant and linear terms are given respectively by
\be
U_{i,j}^n(\tau, \boldsymbol{\beta}) = \begin{cases}
 \frac{1}{\beta_{i} \beta_{j}}\Big(D_{i,j}^n(\tau, \boldsymbol{\beta})  
 - \frac{1 - \exp(-\tau(\beta_{i} - \beta_{n}))}{\beta_{i} - \beta_{n}} \\
\qquad\qquad - \frac{1 - \exp(-\tau(\beta_{j} - \beta_{n}))}{\beta_{j} - \beta_{n}} 
 - \frac{1 - \exp(\tau \beta_{n})}{\beta_{n}} \Big), & \text{if $\beta_{i} \neq 0$ and $\beta_{j} \neq 0$}, \vspace{2ex}\\
 \frac{1}{\beta_{i}}\Big(\frac{\tau \exp(-\tau(\beta_{i} - \beta_{n}))}{\beta_{i} - \beta_{n}} - \frac{1 - \exp(-\tau(\beta_{i} - \beta_{n}))}{(\beta_{i} - \beta_{n})^2} \\
\qquad\qquad + \frac{\tau \exp(\tau \beta_{n})}{\beta_n} + \frac{1-\exp(\tau \beta_{n})}{(\beta_{n})^2} \Big), & \text{if $\beta_{i} \neq 0$ and $\beta_{j} = 0$}, \vspace{2ex}\\
 \frac{1}{\beta_{j}}\Big(\frac{\tau \exp(-\tau(\beta_{j} - \beta_{n}))}{\beta_{j} - \beta_{n}} - \frac{1 - \exp(-\tau(\beta_{j} - \beta_{n}))}{(\beta_{j} - \beta_{n})^2} \\
\qquad\qquad + \frac{\tau \exp(\tau \beta_{n})}{\beta_n} + \frac{1-\exp(\tau \beta_{n})}{(\beta_{n})^2} \Big), & \text{if $\beta_{i} = 0$ and $\beta_{j} \neq 0$}, \vspace{2ex}\\
- \frac{(\tau)^2 \exp(\tau \beta_{n})}{\beta_{n}} - \frac{2}{\beta_n} \Big(
\frac{\tau \exp(\tau \beta_{n})}{ \beta_{n}} + \frac{1-\exp(\tau \beta_{n})}{(\beta_{n})^2} \Big), & \text{if $\beta_{i} = 0$ and $\beta_{j} = 0$},
 \end{cases}
\ee
and 
\be
V_{i,j}^n(\tau, \boldsymbol{\beta}) = \begin{cases}
 \frac{1 - \exp(-\tau(\beta_{i} + \beta_{j} - \beta_{n}))}{\beta_{j}(\beta_{i} + \beta_{j} - \beta_{n})} - \frac{1 - \exp(-\tau(\beta_{i} - \beta_{n}))}{\beta_{j}(\beta_{i} - \beta_{n})}, & \text{if $\beta_{j} \neq 0$}, \vspace{2ex}\\
 \frac{\tau \exp(-\tau(\beta_{i} - \beta_{n}))}{\beta_{i} - \beta_{n}} - \frac{1 - \exp(-\tau(\beta_{i} - \beta_{n}))}{(\beta_{i} - \beta_{n})^2}, & \text{otherwise}.
 \end{cases}
\ee
%%%%%%%%%%%%%%%%%%%%%%%%%%%%%%%%%%%%%%%%%%%%%%%%%%%%%%
By adding $\int_{-\tau}^0 e^{\beta_n(t-s)}  \Pi_{n} F \d s$ to the constant and linear terms in $I$, we can form
\bea \label{Formula_GammaF}
\Gamma_n(F,\bfbeta,\tau,\xi) =  \sum_{i, j = 1}^m  U_{i,j}^n(\tau, \boldsymbol{\beta}) B_{i,j}^n F_{i}F_{j}+ \sum_{i, j = 1}^m   V_{i,j}^n(\tau, \boldsymbol{\beta})  F_{j} (B^n_{i,j} &+ B^n_{j,i})   \xi_{i} \\
& - \frac{1 - e^{\tau \beta_n}}{\beta_n} \Pi_{n} F,
\eea
leading thus to
\be \label{Eq_Phi_tau}
\boxed{\Phi_n(\tau, \boldsymbol{\beta}, \xi) =\Gamma_n(F,\bfbeta,\tau,\xi) + \sum_{i, j = 1}^m D_{i,j}^n(\tau, \boldsymbol{\beta}) B_{i,j}^n \xi_{i} \xi_{j}. }
\ee
The optimal $\tau$ value for each of the unresolved mode is obtained by minimizing the corresponding parameterization defect $\mathcal{Q}_n$ defined in \eqref{Q_n_defined}. In other words, given a fully resolved solution $y(t)$ of the underlying $N$-dimensional ODE system \eqref{Eq_quadODEs} available over a training interval $[0,T]$ (after possible removal of transient dynamics), we solve for each $m+1 \leq n \leq N$ the following minimization problem 
%%%%%%%%%%%%%%%%%%%%%%%%%%%%%%%%%%%%%%%%%%%%%%%
\begin{numcases}{\label{Min_formulation_h1_quadcase}}
     \; \underset{\tau}\min \int_0^T \big| \Pi_n y(t)- \Phi_n(\tau, \boldsymbol{\beta}, \Pi_{\c} y(t))\big|^2 \d t, \label{PD_h1_quadcase}\\
 \;\textrm{where } \Phi_n(\tau,\boldsymbol{\beta}, \xi)  \textrm{ is given by \eqref{Eq_Phi_tau}.}\nonumber 
    \end{numcases}
The resulting minimizers $\tau_n^\ast$ whose collection is denoted by $\bftau^\ast$, allows us then to define the following 
optimal parameterization within the LIA class
\be\label{Eq_Phi_tau_opt}
\Phi^{(1)}_{\bftau^\ast}(\xi)= \sum_{n = m+1}^N \Phi_n (\tau_n^\ast, \boldsymbol{\beta}, \xi) \boldsymbol{e}_n.
\ee
In what follows we will sometimes denote by LIA($\bftau$), the parameterization  $\Phi^{(1)}_{\bftau}$ (see \ref{Eq_Phi_tau_opt}) with $\Phi_n$ given by  \eqref{Eq_Phi_tau}. 

Although providing in general only a suboptimal solution to the more general family of minimization problems \eqref{Eq_Jfunc_n} discussed in Sec.~\ref{Sec_variational_approach_Q}, we will refer to the optimal LIA, $\Phi^{(1)}_{\bftau^\ast}$, as the optimal PM
when the context is clear;  see Sec.~ \ref{Sect_RBC} below.
As mentioned above, Appendix~\ref{Sect_gradient_descent} presents a simple gradient-descent method  to determine efficiently, the $\tau_n^\ast$'s (and thus $\bftau^\ast$) in practice; as pointed out above, see however Remark \ref{Remark_training_interval} below in the presence of local minima.

\br\label{Rmk_LIA_h2}
Note that for $F=0$, and when $\beta_i+\beta_j>\beta_n$,  the LIA class includes the leading-order approximation, $h_2$, given by 
\eqref{h1_part1}-\eqref{h1_part2} (with $k=2$) of the invariant manifold dealt with in Sec.~\ref{Sec_Leading_approx}, in the sense that then for  all  $\xi$ in $E_\c$,
\be
\lim_{\bftau \rightarrow \infty} \Phi^{(1)}_{\bftau}(\xi)=h_2(\xi). 
\ee 
Furthermore $\Phi^{(1)}_{\bftau} \equiv 0$ when $\bftau=0$, i.e.~the LIA class contains Galerkin approximations of dimension $m=\textrm{dim}(E_\c)$. 
\er

\br
Note that in  the expression of $\Phi_n$ given by \eqref{Eq_Phi_tau},  the term $\Gamma_n(F,\bfbeta,\tau,\xi)$ takes into account interactions between the 
low-mode components of the forcing, $F$, as well as cross-interactions between the low-mode components of $F$ and the low-mode variable $\xi$ in $E_\c$.   It also includes the $n^{th}$ high-mode  component of the forcing.  

We emphasize that these formulas can be derived for PDEs as well, as rooted in the backward-forward method recalled above and initially introduced for PDEs (possibly driven by a multiplicative linear noise)  in \cite[Chap.~4]{CLW15_vol2}; see also \cite[Sec.~3.2]{CL15}. The main novelty compared to  \cite[Chap.~4]{CLW15_vol2} is  the idea of optimizing, high-mode by high-mode, the backward integration time, $\tau_n$, of Eq.~\eqref{Eq_BF_quad}, by minimization of the parameterization defect $Q_n$.
\er

\br\label{Rem_Conj_pair}
Note that when $\beta_{n+1}=\overline{\beta_n}$, we have $\boldsymbol{e}_{n+1}^\ast=\overline{\boldsymbol{e}_n^\ast}$ and therefore 
$\Pi_{n+1} X=\overline{\Pi_n X}$  when $X$ is real according to \eqref{Proj_formula}. 
Furthermore when $B( u^{(1)}_{\c}(s),  u^{(1)}_{\c}(s))$ and $F$ are real, we have according to   
\eqref{eq_Phi_n_abstract}, that $\Phi_{n+1}=\overline{\Phi_{n}}$ when evaluated on a real vector $\xi$ of $E_\c$. 
\er 

%%%%%%%%%%%%%%%%%%%%%%%%%%%%%%%%%%%%%%%%%%%%%%%%%%%%%
\subsection{ Parametric Quasi-Stationary Approximation and another cost functional}\label{Sec_FMTtau}
Other cost functionals than $\mathcal{Q}_n(\tau_n,T)$ could have been considered to seek for optimal LIA. 
For instance, 
\be\label{Eq_Defect_bis}
\mathcal{J}_n(\tau,T; \Phi_n)=\bigg| \overline{\Big [\Pi_n y(t)\Big]^2}-\overline{\Big[\Phi_n(\tau,\bfbeta, y_\c(t)))\Big]^2}\bigg|.
\ee
Here $\overline{(\cdot)}$ denotes a time-averaging over an interval of length $T$. 
The minimization of the $\mathcal{J}_n$'s leads in general to different optimal LIA compared to the one obtained by solving the minimization problems \eqref{Min_formulation_h1_quadcase}.

If the mean value of $y_n(t)$ is zero, minimizing $\mathcal{Q}_n$ consists of minimizing the variance of the residual error, i.e.~$\overline{|y_n-f(\tau,y_\c)|^2}$, for a given parameterization $f(\tau,\cdot)$. By construction, minimizing $\mathcal{J}_n$ consists instead of minimizing the residual error of the variance approximation, i.e.~$|\overline{|y_n|^2}-\overline{|f(\tau,y_\c)|^2}|$. The latter cost functional better accounts for the distribution of energy across the modes; see Sec.~\ref{Sec_More_Results} for an illustration.

Although a geometric interpretation like \eqref{Eq_geom_interp} is not available for such a cost functional, 
minimizing 
\eqref{Eq_Defect_bis} leads in general to a better reproduction of the energy budget across the high modes. For this reason, the cost functional \eqref{Eq_Defect_bis} will be adopted for certain applications; see Sec.~\ref{Sec_KS_turbulence} below.

While the LIA class may be preferred when forcing terms are present (especially when e.g.~only the low modes are forced), another class of parameterization 
is particularly suited to systems that do not include forcing terms.  Still, in presence of such terms this other class may be relevant 
in certain applications (when e.g.~only the high modes are forced) and thus we present hereafter the derivation of the corresponding formulas that take into account (constant) forcing as for LIA.  

This class is rooted in the following {\it Quasi-Stationary approximation (QSA)} for Eq.~\eqref{Eq_quadODEs}
\be\label{Eq_balance0}
\Pi_\s A z+ \Pi_\s B(\xi,\xi)+\Pi_\s F=0, \quad\; \xi  \in E_\c, \; z \in  E_\s. 
\ee
The QSA arises in homogeneous turbulence theory \cite{FMT88}; see Remark \ref{Rmk_QSA_backscatter} below. It consists of neglecting the terms $\Pi_\s[B(y_\s,y_\c) + B(y_\s,y_\s)]$ in virtue of the energy content of the small structures being small, and following a suggestion of Kraichnan balancing $\d y_\s /\d t$ with $\Pi_\s B(y_\c,y_\s)$, i.e., with the advection of small eddies by large eddies; see \cite{foias1991approximate}.

After solving \eqref{Eq_balance0}, the QSA parameterization is then obtained as $z=K(\xi)$ with $K$ given by 
\be\label{Eq_QSA}
K(\xi)=(-A_\s)^{-1} (\Pi_{\s}B(\xi,\xi)+\Pi_\s F). 
\ee 
 In contrast, the standard LIA  is obtained by solving the backward-system \eqref{Eq_BF_quad0} asymptotically, and the parameterization LIA($\bftau$) is obtained after solving the backward-systems \eqref{Eq_BF_quad}.

Similar to what precedes, we use a dynamic version of  Eq.~\eqref{Eq_balance0} to get access to a parametric family of dynamically-based parameterizations such that $K$ belongs to this family, as in Remark \ref{Rmk_LIA_h2} regarding the LIA class that includes $h_2$.  By assuming $A$ diagonal (in $\mathbb{C}$), 
we consider thus for $\tau>0$  
\bea\label{FMT_tau}
&\frac{\mathrm{d} z_{n}}{\d s} = \beta_n z_n  +  \Pi_{n} B\big(\xi, \xi \big) + \Pi_{n} F,\\
& z_n(-\tau)=0.
\eea
Solving Eq.~\eqref{FMT_tau} for each $n$, leads then to the following high-mode parameterization
\be\label{Eq_Psin}
\Psi_n(\tau,\boldsymbol{\beta},\xi)=\delta_n(\tau)\bigg( \sum_{i, j = 1}^m B_{i j }^n \xi_{i} \xi_{j}+\Pi_{n} F\bigg),
\ee
with $B_{i j }^n$ given by \eqref{B_ijn} and where
\bea \label{Eq_D_term00}
\delta_n(\tau)=  \begin{cases}
\beta_n^{-1}(e^{\beta_{n}\tau}-1), & \text{if $\beta_{n}\neq 0$,} \\
\tau, & \text{otherwise}.
\end{cases}
\eea
We arrive then at the following {\it parametric} QSA or simply denoted QSA$(\bftau)$:
\be\label{QSA_tau}
\Psi_\bftau(\xi)=\sum_{n=m+1}^N \Psi_n(\xi,\boldsymbol{\beta},\xi) \boldsymbol{e}_n. 
\ee
In particular, if $\beta_n<0$ for all $n\geq m+1$, since $\delta_n(\tau)\underset{\tau\rightarrow \infty}\longrightarrow-\beta_n^{-1}$, then for  all  $\xi$ in $E_\c$,
\be\label{K_limit_QSA}
\lim_{\bftau \rightarrow \infty}\Psi_\bftau(\xi)=K(\xi), 
\ee
with $K$ given by \eqref{Eq_QSA}. Furthermore $\Psi_{\bftau} \equiv 0$ when $\bftau=0$, i.e.~the QSA class contains also Galerkin approximations of dimension $m=\textrm{dim}(E_\c)$. 

In Sec.~\ref{Sec_KS_turbulence} below, we show applications of this parameterization class (called the QSA class),  from which 
the {\it optimal QSA} is determined by solving for each $m+1 \leq n \leq N$ the following minimization problem 
%%%%%%
\begin{numcases}{\label{Min_formulation_h1_quadcase_b}}
     \; \underset{\tau}\min \; \bigg| \overline{\Big [\Pi_n y(t)\Big]^2}-\overline{\Big[\Psi_n(\tau,\bfbeta, y_\c(t)))\Big]^2}\bigg|.\label{PD_FMT_quadcase}\\
 \;\textrm{where } \Psi_n(\tau,\boldsymbol{\beta}, \xi)  \textrm{ is given by \eqref{Eq_Psin}.}\nonumber 
    \end{numcases}
 The algorithm presented in Appendix~\ref{Sect_gradient_descent} to solve \eqref{Min_formulation_h1_quadcase}, can be easily adapted to solve \eqref{Min_formulation_h1_quadcase_b} (after smoothing) and thus to determine the minimizers $\tau_n^\ast$; the details are left to the reader.

As recalled above, Remark \ref{Rmk_LIA_h2} emphasizes that  the leading-order approximation  $h_2(\xi)$ (given by \eqref{Eq_LyapPerron} with $G_k=B$) of the invariant manifold dealt with in Sec.~\ref{Sec_Leading_approx} may be obtained as a limit LIA($\bftau$): here \eqref{K_limit_QSA} shows that the standard QSA, $K(\xi)$, may also be obtained as a limit of QSA($\bftau$). It is noteworthy that the theory of approximation of invariant manifolds shows that these two limiting objects, $h_2(\xi)$ and $K(\xi)$, are actually related. More precisely, \cite[Lemma 4.1]{CLW15_vol2} shows that near the first criticality and when $F=0$,  the QSA and the leading-order approximation $h_2(\xi)$, are linked according to the following approximation relation 
\be
h_2(\xi)=(-A_\s)^{-1} \Pi_{\s}B(\xi,\xi) + O(  \|\xi\|^2), \quad \Forall \xi \in E_\c.
\ee
Thus when $F=0$, one should not expect much difference between the parameterizations LIA($\bftau$) and QSA($\bftau$) for large values of $\bftau$ (and under the appropriate conditions on the $\beta_k$'s).  

However, if $\bftau$ has components with small values, differences are expected to occur between the corresponding LIA($\bftau$) and QSA($\bftau$) parameterizations. 
To better appreciate these differences, let us introduce the function $f(\tau) = p^{-1}(1-e^{-p \tau})$ and note that $f(\tau)=\delta_n(\tau)$ when $p= -\beta_n$ and that $f(\tau) =D_{ij}^n(\tau)$ (given by \eqref{Eq_D_term0}) when $p=\beta_i + \beta_j - \beta_n$.
Thus when $F=0$ the LIA and QSA classes differ only by these coefficients. 

To simplify, let us assume that the eigenvalues of $A$ are real and that  $E_\c $ contains all and only the unstable modes. In this case,  $p= \beta_i + \beta_j - \beta_n$ is always bigger than $p= -\beta_n$. Now if we assume furthermore that $p>0$ (in either case) we have  
\be
0 \leq f(\tau) < p^{-1},
\ee
and therefore due to  \eqref{Eq_Psin} and \eqref{Eq_Phi_tau} (with $F=0$), the range of the coefficient in front of each monomial is  larger for $\Psi_n(\tau,\xi)$ than for $\Phi_n(\tau,\xi)$, in this case. 
This allows in practice for  $\Psi_n(\tau,\xi)$ to span a larger range of values which in turn may lead to smaller values of $\mathcal{Q}_n$ or $\mathcal{J}_n$. The situation described here is exactly what happens for the closure problem considered below in Sec.~\ref{Sec_KS_turbulence} within the context of Kuramoto-Sivashinsky turbulence, when one sets the cutoff wavenumber to be the highest wavenumber among the unstable modes.  
As we will show in Sec.~\ref{Sec_KS_turbulence} for different turbulent regimes, the QSA($\bftau$) when optimized (either for  $\mathcal{Q}_n$ or $\mathcal{J}_n$) provides a drastic improvement compared to the standard QSA, $K(\xi)$, for such cutoff scales. 

\br\label{Rmk_QSA_backscatter}
As mentioned right after \eqref{Eq_balance0}, the QSA is a well-known parameterization in homogeneous turbulence and has been rigorously proved to provide an AIM in \cite{FMT88} for the 2D Navier-Stokes equations.  The QSA also arises in atmospheric turbulence in the so-called nonlinear normal-mode initialization\cite{machenhauer1977dynamics,baer1977complete,tribbia1979nonlinear,leith1980nonlinear,daley1980development,ghil1991data,daley1993atmospheric}; see \cite{debussche1991inertial} for rigorous results. Nevertheless, when the cutoff wavelength is too low within the inertial range it is known that the standard QSA  suffers from over-parameterization leading then to errors in the backscatter transfer of energy, i.e.~errors in the modeling of the parameterized (small) scales that contaminate gradually the larger scales. We show in Sec.~\ref{Sec_KS_turbulence}, in the context of KS turbulence that 
by solving the minimization problems \eqref{Min_formulation_h1_quadcase_b}, the optimal QSA fixes this problem remarkably.  
\er

%%%%%%%%%%%%%%%%%%%%%%%%%%%%%%%%%%%%%%%%%%%%%%%%%%%
%
%                                 RBC applications
%
%%%%%%%%%%%%%%%%%%%%%%%%%%%%%%%%%%%%%%%%%%%%%%%

\section{Applications to a reduced-order Rayleigh-B\'enard system} \label{Sect_RBC}
In this section, we apply the PM approach --- as presented in its practical aspects in Sec.~\ref{Sect_PM_formulas} --- to   
a Galerkin system of nine nonlinear ODEs examined in \cite{Reiterer_al98}  and obtained from a triple Fourier expansion to the Boussinesq equations governing thermal convection in a 3D spatial domain.

The PM approach is applied to two parameter regimes for this 9D Rayleigh-B\'enard (RB) convection system: (i) a regime located right after the first period-doubling bifurcation occurring for this system (Sec.~\ref{Sect_period-doubling-RB}), and (ii) a regime corresponding to chaotic dynamics  that takes place right after the period-doubling cascade (Sec.~\ref{Sect_chaotic-RB}).

We show hereafter for both cases, that, given a reduced state space,  $E_\c$, the dynamically-based parameterization, LIA($\bftau$),  of Sec.~\ref{Sect_PM_with_forcing} when optimized in the $\bftau$-variable,  by minimizing\footnote{while maximizing,  in certain circumstances, the parameterization correlation, $c(t)$, given by \eqref{Eq_corr_param}; see Sec.~\ref{Sect_period-doubling-RB}.} the parameterization defects \eqref{Min_formulation_h1_quadcase},  provides efficient low-dimensional closures of the original  RB system.

To prepare the numerical results of Secns.~\ref{Sect_period-doubling-RB} and \ref{Sect_chaotic-RB},  we first recall the 9D RB system and give the details of its  LIA($\bftau$)-closure  in Sec.~\ref{Sect_9D_system_and_closure}. We emphasize that the closures are determined in each case with respect to a mean state $\overline{\boldsymbol{C}}$, leading in particular to equations for the perturbed variable, $\boldsymbol{C}-\overline{\boldsymbol{C}}$, of the form \eqref{Eq_pert_variable1}.

\subsection{Optimal PM closure} \label{Sect_9D_system_and_closure}
Like \cite{Reiterer_al98}, our study below deals with three-dimensional cells with square planform in dissipative Rayleigh-B\'enard convection. In that respect, the 9D RB system derived in \cite[Section 2]{Reiterer_al98} takes the form:
%%%%%%%%%%%%%%%%%%%%%%%%%%%%%%%%%%%%
\bea \label{Eq_9DRBC}
\dot{C_1}&=-\sigma b_1 C_1 -C_2 C_4 +b_4 C_4^2+b_3 C_3 C_5 -\sigma b_2 C_7, \\
\dot{C_2}&=-\sigma C_2 +C_1 C_4 -C_2 C_5+ C_4 C_5 -\frac{\sigma}{2} C_9, \\
\dot{C_3}&=-\sigma b_1 C_3 +C_2 C_4 -b_4 C_2^2-b_3 C_1 C_5 +\sigma b_2 C_8, \\
\dot{C_4}&=-\sigma C_4 -C_2 C_3 - C_2C_5+C_4 C_5 +\frac{\sigma}{2}C_9, \\
\dot{C_5}&=-\sigma b_5 C_5 + \frac{1}{2} C_2^2 - \frac{1}{2} C_4^2, \\
\dot{C_6}&=- b_6 C_6 + C_2 C_9 - C_4C_9, \\
\dot{C_7}&=- b_1 C_7 -rC_1 + 2 C_5 C_8 - C_4 C_9, \\
\dot{C_8}&=- b_1 C_8  + rC_3 -2 C_5 C_7 +   C_2 C_9, \\
\dot{C_9}&=- C_9 -r C_2 +r C_4 -2 C_2 C_6 +2 C_4 C_6+  C_4 C_7 -C_2 C_8.
\eea
%%%%%%%%%%%%%%%%%%%%%%%%%%%%%%%%%%%%
Here $\sigma$ denotes the Prandtl number, and $r$ denotes the reduced Rayleigh number defined to be the ratio between the Rayleigh number $R$ and its critical value $R_c$ at which the convection sets in. The coefficients $b_i$'s are given by
\bea
& b_1= \frac{4(1+a^2)}{1+2a^2},  && b_2=\frac{1+2a^2}{2(1+a^2)}, && b_3=\frac{2(1-a^2)}{1+a^2}, \\
& b_4=\frac{a^2}{1+a^2}, && b_5=\frac{8a^2}{1+2a^2}, &&  b_6=\frac{4}{1+2a^2},
\eea
with $a=\frac{1}{2}$ being the critical horizontal wavenumber of the square convection cell. 

With the purpose to derive a closure for Eq.~\eqref{Eq_9DRBC}, we first put Eq.~\eqref{Eq_9DRBC} into the following compact form: 
\be
\dot{\boldsymbol{C}} = A \boldsymbol{C} + B(\boldsymbol{C},\boldsymbol{C}),
\ee
where $\boldsymbol{C}=(C_1, \cdots C_9)^{\mathrm{tr}}$, $A$ is the $9\times9$ matrix given by
\be
A= \begin{pmatrix}
-\sigma b_1 & 0 & 0 & 0 & 0 & 0 & -\sigma b_2 & 0 & 0 \\  % row 1
0 & -\sigma & 0 & 0 & 0 & 0 & 0 & 0 & -\frac{\sigma}{2}  \\ % row 2
0 & 0 &-\sigma b_1& 0 & 0 & 0 & 0 & \sigma b_2 &  0  \\ % row 3
0 & 0 & 0 & -\sigma & 0 & 0  & 0 & 0 & \frac{\sigma}{2} \\ % row 4
0 & 0 & 0 & 0 & -\sigma b_5 & 0 &  0  & 0 & 0 \\ % row 5
0 & 0 & 0 & 0 & 0 & - b_6 & 0  & 0  & 0 \\ % row 6
-r & 0 & 0 & 0 & 0 & 0 & - b_1  & 0  & 0 \\ % row 7
0 & 0 & r & 0 & 0 & 0 & 0 & - b_1  & 0 \\ % row 8
0 & -r & 0 & r & 0 & 0 & 0 & 0  & -1 % row 9
\end{pmatrix}, 
\ee
and the quadratic nonlinearity $B$ is defined by 
\bea\label{B_term_RBC}
B(\boldsymbol{\phi},\boldsymbol{\psi}) =  \begin{pmatrix}
-\phi_2 \psi_4 +b_4 \phi_4 \psi_4 +b_3 \phi_3 \psi_5 \\
\phi_1 \psi_4 -\phi_2 \psi_5+ \phi_4 \psi_5 \\
\phi_2 \psi_4 -b_4 \phi_2\psi_2-b_3 \phi_1 \psi_5 \\
-\phi_2 \psi_3 - \phi_2\psi_5+\phi_4 \psi_5 \\
\frac{1}{2} \phi_2 \psi_2 - \frac{1}{2} \phi_4 \psi_4 \\
\phi_2 \psi_9 - \phi_4\psi_9 \\
2 \phi_5 \psi_8 - \phi_4 \psi_9 \\
-2 \phi_5 \psi_7 +   \phi_2 \psi_9 \\
-2 \phi_2 \psi_6 +2 \phi_4 \psi_6+  \phi_4 \psi_7 
\end{pmatrix} 
\eea
for any $\boldsymbol{\phi}=(\phi_1, \cdots, \phi_9)^\mathrm{tr}$ and $\boldsymbol{\psi}=(\psi_1, \cdots, \psi_9)^\mathrm{tr}$ in $\mathbb{C}^9$. 

We consider next fluctuations defined with respect to a mean state. 
In that respect, we subtract from $\boldsymbol{C}(t)= (C_1(t), \cdots, C_9(t))$ its mean value $\overline{\boldsymbol{C}}$, which is estimated, in practice, from simulation of Eq.~\eqref{Eq_9DRBC} on the same training interval $T$ than used to optimize our parameterizations hereafter. 
The corresponding ODE system for the fluctuation variable, $\boldsymbol{D} = \boldsymbol{C} - \overline{\boldsymbol{C}}$, is then given by:
\be \label{Eq_RBC_fluct}
\frac{\mathrm{d} \boldsymbol{D}}{\mathrm{d} t} = L \boldsymbol{D} + B(\boldsymbol{D},\boldsymbol{D}) + {\hl A \overline{\boldsymbol{C}}} +  B(\overline{\boldsymbol{C}}, \overline{\boldsymbol{C}}),
\ee
with 
\be\label{RB_Linear_part_perturbed}
L \boldsymbol{D}=A \boldsymbol{D} + B(\overline{\boldsymbol{C}},\boldsymbol{D}) + B(\boldsymbol{D},\overline{\boldsymbol{C}}). 
\ee

Denote the spectral elements of the matrix $L$ by $\{(\beta_j, \boldsymbol{e}_j) \; : \; 1 \le j \le 9\}$ and those of $L^{\ast}$ by $\{(\beta^*_j, \boldsymbol{e}^*_j) \; : \; 1 \le j \le 9\}$.  By taking the expansion of $\boldsymbol{D}$ under the eigenbasis of $L$,
\be
\boldsymbol{D}= \sum_{j=1}^9 y_j \boldsymbol{e}_j \quad \text{ with } \quad y_j = \langle\boldsymbol{D}, \boldsymbol{e}^*_j \rangle,
\ee
and assuming that $L$ is diagonal under its eigenbasis, we rewrite Eq.~\eqref{Eq_RBC_fluct} in the variable $\boldsymbol{y}=(y_1, \cdots, y_9)^\mathrm{tr}$ as follows: 
\bea \label{Eq_RBC_eigenbasis}
\dot{y}_j = \beta_j y_j + \sum_{k,\ell = 1}^9 \langle B(\boldsymbol{e}_k, \boldsymbol{e}_\ell), \boldsymbol{e}^*_j \rangle y_k y_{\ell} + \langle {\hl A \overline{\boldsymbol{C}}} + B(\overline{\boldsymbol{C}}, \overline{\boldsymbol{C}}), \boldsymbol{e}^*_j \rangle, \quad j = 1, \cdots, 9.
\eea

Now we take the reduced state space $E_\c$ to be spanned by the first $m$ eigenvectors of $A$ for some $m < 9$, where the eigenvalues are ranked according to the ordering \eqref{Eq_gap} adopted here from Sec.~\ref{Sec_loc_invman}, i.e.~the modes are ordered according to their linear rate of growth/decay.
For each  $m+1\leq n \leq 9 $, we approximate the (unresolved) variable $y_n$ by the parameterization $\Phi_n(\tau_n^\ast,\boldsymbol{\beta},\cdot)$ obtained from \eqref{Eq_Phi_tau} after minimization of \eqref{Min_formulation_h1_quadcase}, given a training interval of length $T$ that will be specified hereafter depending on the context.  

The resulting $m$-dimensional optimal PM closure (in the LIA class) reads then 
\bea \label{Eq_RBC_reduced}
\dot{x}_j & = \beta_j x_j + \sum_{k,\ell = 1}^m \langle B(\boldsymbol{e}_k, \boldsymbol{e}_\ell), \boldsymbol{e}^*_j \rangle x_k x_\ell 
 \\
& \quad +  \sum_{k = 1}^m \sum_{\ell=m+1}^9 \Big( \langle B(\boldsymbol{e}_\ell, \boldsymbol{e}_k), \boldsymbol{e}^*_j \rangle + \langle B(\boldsymbol{e}_k, \boldsymbol{e}_\ell), \boldsymbol{e}^*_j \rangle\Big) x_k \Phi_\ell(\tau_\ell^\ast,\boldsymbol{\beta}, x_1,\cdots, x_m) \\
& \quad +  \sum_{k,\ell=m+1}^9 \langle B(\boldsymbol{e}_\ell, \boldsymbol{e}_k), \boldsymbol{e}^*_j \rangle  \Phi_k(\tau_k^\ast, \boldsymbol{\beta}, x_1,\cdots, x_m)  \Phi_\ell(\tau_\ell^\ast, \boldsymbol{\beta},x_1,\cdots, x_m)  \\
& \quad + \langle {\hl A \overline{\boldsymbol{C}}} + B(\overline{\boldsymbol{C}}, \overline{\boldsymbol{C}}), \boldsymbol{e}^*_j \rangle, \quad j = 1, \cdots, m.
\eea

Once the optimal PM closure \eqref{Eq_RBC_reduced} is solved, an approximation, $\boldsymbol{C}^{\text{PM}}(t)$, of the solution $\boldsymbol{C}(t)$ to the original system \eqref{Eq_9DRBC} is obtained as follows,
\be \label{Eq_soln_reconst_9DRBC}
\boldsymbol{C}^{\text{PM}}(t) = \sum_{j=1}^m x_j(t) \boldsymbol{e}_j  + \sum_{n=m+1}^9 \Phi_n(\tau_n^\ast, \boldsymbol{\beta}, x_1(t),\cdots, x_m(t)) \boldsymbol{e}_n  + \overline{\boldsymbol{C}}. 
\ee

\subsection{Closure in a period-doubling regime} \label{Sect_period-doubling-RB}
As the reduced Rayleigh number  $r$ increases, the first period-doubling bifurcation for Eq.~\eqref{Eq_9DRBC} occurs at approximately $r = 13.97$, and the dynamics becomes chaotic at approximately $r=14.22$ after successive periodic-doubling bifurcations. We have set $r= 14.1$ to examine how the PM approach operates in a period-doubling regime. As a benchmark, for the same reduced dimension, $m$, as used for the optimal PM closure \eqref{Eq_RBC_reduced},  we determine the reduced system of the form \eqref{Eq_reduced_absract} in which $h$ is replaced by the approximation $h_2$ given by  \eqref{h1_part1}-\eqref{h1_part2} (with $k=2$) in Theorem~\ref{thm:h1}, i.e.~the parameterization that provides the leading-order  approximation of the local invariant manifold for an equilibrium. Applying the ideas of Sec.~\ref{Sec_loc_invman} to Eq.~\eqref{Eq_9DRBC}, the calculations of  $h_2$ are made about a steady state of Eq.~\eqref{Eq_9DRBC}, taken here to be the closest steady state $\overline{Y}$ to the mean state, $\overline{\boldsymbol{C}}$.  If one denotes by $F$ the RHS of Eq.~\eqref{Eq_9DRBC}, the linear part $A$ in \eqref{Eq_ODEs} is then taken to be given by $DF(\overline{Y})$. 

Thus, denoting by $(\lambda_j,\boldsymbol{f}_j)$ the spectral elements of $DF(\overline{Y})$ and those of $\big(DF(\overline{Y})\big)^{\ast}$ by $(\lambda^*_j, \boldsymbol{f}^*_j)$, the following reduced system based on the  invariant manifold approximation $h_2$,
\bea \label{Eq_RBC_reduced_h1}
\dot{z}_j & = \lambda_j z_j + \sum_{k,\ell = 1}^m \langle B(\boldsymbol{f}_k, \boldsymbol{f}_\ell), \boldsymbol{f}^*_j \rangle z_k z_\ell 
 \\
& \quad +  \sum_{k = 1}^m \sum_{\ell=m+1}^9 \Big( \langle B(\boldsymbol{f}_\ell, \boldsymbol{f}_k), \boldsymbol{f}^*_j \rangle + \langle B(\boldsymbol{f}_k, \boldsymbol{f}_\ell), \boldsymbol{f}^*_j \rangle\Big) z_k h_{2,\ell}( z_1,\cdots, z_m) \\
& \quad +  \sum_{k,\ell=m+1}^9 \langle B(\boldsymbol{f}_\ell, \boldsymbol{f}_k), \boldsymbol{f}^*_j \rangle  h_{2,k}(z_1,\cdots, z_m)  h_{2,\ell}(z_1,\cdots, z_m),  \quad j = 1, \cdots, m,
\eea
serves us as a benchmark. Here $h_{2,n}$ ($6\leq n \leq 9$) is given by \eqref{h1_part2} in which $G_k$ is replaced by $B$ given by \eqref{B_term_RBC} and the $(\beta_j,\boldsymbol{e}_j)$'s replaced by the $(\lambda_j,\boldsymbol{f}_j)$'s. 

From the solution $z(t)=(z_1(t),\cdots,z_m(t))^{\textrm{tr}}$ of the reduced system \eqref{Eq_RBC_reduced_h1}, the following approximation of $\boldsymbol{C}(t)$ is then obtained,
\be\label{Eq_CIM}
\boldsymbol{C}^{\text{IM}}(t) = \sum_{j=1}^m z_j(t) \boldsymbol{f}_j  + \sum_{n=m+1}^9 h_{2,n}( z_1(t),\cdots, z_m(t)) \boldsymbol{f}_n +\overline{Y}. 
\ee
%%%%%%%%%%%%%%%%%%%%%%%%%%%%%%%%%%%%%%%%%
For the numerical results presented hereafter, the reduced state space $E_\c$ is taken to be spanned by the first five eigenmodes, i.e.~by setting $m=5$ in this section. To determine our optimal PM closure, we used the quadratic parameterization, $\Phi_n(\tau,\cdot)$ given by \eqref{Eq_Phi_tau}, in order to parameterize each of the modes $\boldsymbol{e}_n$ with $6 \leq n\leq 9$.
For each $6 \leq n\leq 9$, each of this parameterization is optimized in the $\tau$-variable by minimizing the parameterization defect 
\be\label{toto_Qn}
Q_n(\tau,T; t_0)= \int_{t_0}^{t_0+T} \big| \Pi_n y(t)- \Phi_n(\tau, \boldsymbol{\beta}, \Pi_{\c} y(t))\big|^2 \d t,
\ee
for some $t_0$ chosen so that transient dynamics has been removed. 
%%%%%%%%%%%%
Since the dynamics to emulate by a closure is here periodic, we selected $T=3 T_p/4$, where $T_p$ ($\approx 17.25$) corresponds to the period of the solution to the 9D RB system \eqref{Eq_9DRBC} in order to do not use all the available information about the periodic orbit. 
Other choices could have been made for the training interval such as $T=T_p/2$. Note that we observed that the choice of $t_0$ plays a key role here.  As discussed in Remark \ref{Remark_training_interval} below, depending on $t_0$ the global minimizer $\tau_n^\ast$ of $Q_n$ here, does not provide necessarily the best parameterization within the $\Phi_n$-class, and one may have to rely on the parameterization correlation $c(t)$  (see \eqref{Eq_corr_param}) to discriminate between other local minimizers of $Q_n$. The results presented below corresponds to a time origin, $t_0$, for which the global minimizer of the $Q_n$'s lead to the best parameterization within the $\Phi_n$-class.

Despite the aforementioned $t_0$-dependence, for the sake of keeping the notations as concise as possible, the dependence on $t_0$ will not be made apparent for the numerical results presented below.  
This being said, whatever the length $T$ of the training interval, we have used the same training interval $[t_0,t_0+T]$ to estimate the mean state, $\overline{\boldsymbol{C}}$, than used for evaluating the cost functionals  $Q_n$ in \eqref{toto_Qn}.
 
The mean state, $\overline{\boldsymbol{C}}$, plays a key role in the determination of the closure as it determines the linear part $L$ defined in \eqref{RB_Linear_part_perturbed}, and thus the spectral elements $(\beta_j,\boldsymbol{e}_j)$ arising in the formulation of the parameterizations, $\Phi_n(\tau,\cdot)$ (see \eqref{Eq_Phi_tau}), and of the corresponding closure \eqref{Eq_RBC_reduced}. Numerically, a fourth-order Runge-Kutta method is used to solve Eq.~\eqref{Eq_RBC_eigenbasis} with a time-step size taken to be $\delta t = 5 \times 10^{-3}$ to determine a numerical approximation of $y(t)$. The minimization algorithm for the parameterization defect described in Appendix \ref{Sect_gradient_descent} is used to find the minimizer $\tau_n^\ast$ of $Q_n(\tau,T)$. In that respect, the trapezoid rule is used to approximate the integrals involved in \eqref{Eq_Q_elements_continuous}. 

%%%%%%%%%%%%%%%%%%%%Period Doubling%%%%%%%%%%%%%%%%%%%
\begin{figure}[hbtp]
   \centering
\includegraphics[width=.95\textwidth, height=0.45\textwidth]{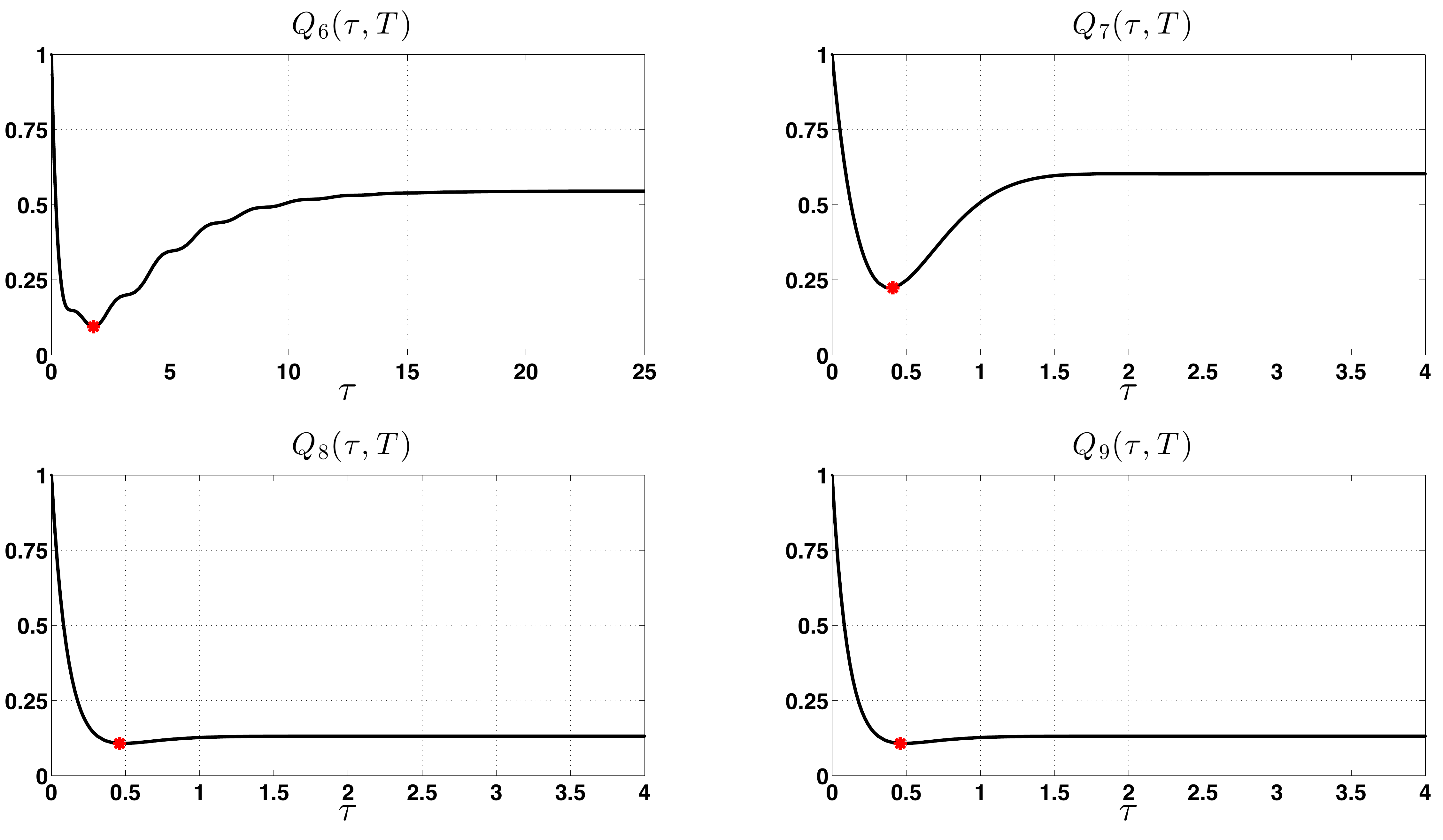}
  \caption{{\footnotesize {\bf $Q_n(\tau,T)$ vs $\tau$  for Eq.~\eqref{Eq_9DRBC} for $r=14.1$ (period-doubling regime) and $m=5$.} For each parameterized mode shown here, the minimum is marked by a red dot.}}   
  \label{fig:RB9D_Qplot_perioddoubling}
\end{figure}
%%%%%%%%%%%

The mapping $\tau \mapsto Q_n(\tau,T)$ is shown 
 in Fig.~\ref{fig:RB9D_Qplot_perioddoubling}  from $n=6$ to $n=9$ and exhibits a non-convex behavior for each $n$, although this behavior is more pronounced for $n=6$ and $n=7$.  The minimizer $\tau_n^\ast$ found by the algorithm of  Appendix \ref{Sect_gradient_descent} corresponds to the abscissa of the red dot shown in each of the panels.  
 Among the parameterized modes, the minima of $Q_n$ that are the most clearly distinguishable occur for the ``adjacent'' modes --- $\boldsymbol{e}_6$ and $\boldsymbol{e}_7$ ---  located   
next to the cutoff dimension, i.e.~for the modes whose real part of the corresponding eigenvalues is the closest (from below) to the real part of $\beta_5$. Nevertheless we emphasize that the ``wavy'' shape of the graph of $Q_6(\tau,T)$ may experience noticeable changes when $t_0$ varies. These changes may be manifested by the emergence of local minima that can modify substantially the global minimizer and thus affect the determination of the optimal PM; a sensitivity issue that can be fixed by the calculation of $c(t)$ given by \eqref{Eq_corr_param}; see Remark \ref{Remark_training_interval}.

Thus, the minimization of the $Q_n$'s possibly completed by the analysis of the parameterization correlation, $c(t)$, allows us to determine the optimal PM, $\Phi^{(1)}_{\bftau^\ast}$, for Eq.~\eqref{Eq_RBC_eigenbasis}  and $E_\c=\mbox{span}\{\boldsymbol{e}_1,\cdots,\boldsymbol{e}_5\}$.  For our choice of $t_0$, the global minima of the $Q_n$'s provide the optimal PM. The values of the parameterization defects for this optimal PM are then given by, 
 $Q_6(\tau_6^\ast,T)=9.5\times 10^{-2}$, $Q_7(\tau_7^\ast,T)=2.2 \times 10^{-1}$ and $Q_8(\tau_8^\ast,T)=Q_9(\tau_9^\ast,T)=1.1\times 10^{-1}$. 
By comparison, for the invariant manifold approximation the parameterization defects (with $h_{2,n}$ replacing $\Phi_n$ in \eqref{toto_Qn}) are given by $Q_6(h_2)=1.8\times 10^{-1}$,    $Q_7(h_2)=2.2$   and  $Q_8(h_2)=Q_9(h_2)=8.2\times 10^{-1}$. 
Note that in both cases, $Q_8=Q_9$, since here $\beta_9=\overline{\beta_8}$ (and $\lambda_9=\overline{\lambda_8}$) and the corresponding parameterizations are just conjugate to each other; see Remark \ref{Rem_Conj_pair}. 

These values of the parameterization defects should be put in perspective with the energy budget for a better appreciation of the exercise of parameterization conducted here.  Table \ref{Table_Frac_energy1} summarizes how the energy is distributed (in average) among the modes, over the training interval $[0,T]$.  The distribution of energy is explained in part (but not only) by  the spectral decomposition and ordering \eqref{Eq_gap} adopted here from Sec.~\ref{Sec_loc_invman}, i.e.~the modes are ordered according to their linear rate of growth/decay.
 In our case, it turns out that Eq.~\eqref{Eq_RBC_eigenbasis} is a genuine forced-dissipative system in which the $\beta_j$'s have all their real parts negative. Thus the ordering is here from the least to the most stable ones; the least stable modes ($\boldsymbol{e}_1$ and $\boldsymbol{e}_2$) containing most of the energy.

It is noteworthy that it is exactly (and only) for mode $\boldsymbol{e}_7$ --- the mode that contains the smallest fraction of energy --- that the parameterization defect $Q_7(h_2)$ for $h_2$ is above 1, leading to an over parameterization for this mode.  Despite the small fraction of energy contained in a given mode, it is known that an over parameterization of such a mode can lead to an overall misperformance of the associated closure. 

\begin{figure}[hbtp]
   \centering
\includegraphics[width=0.95\textwidth, height=0.65\textwidth]{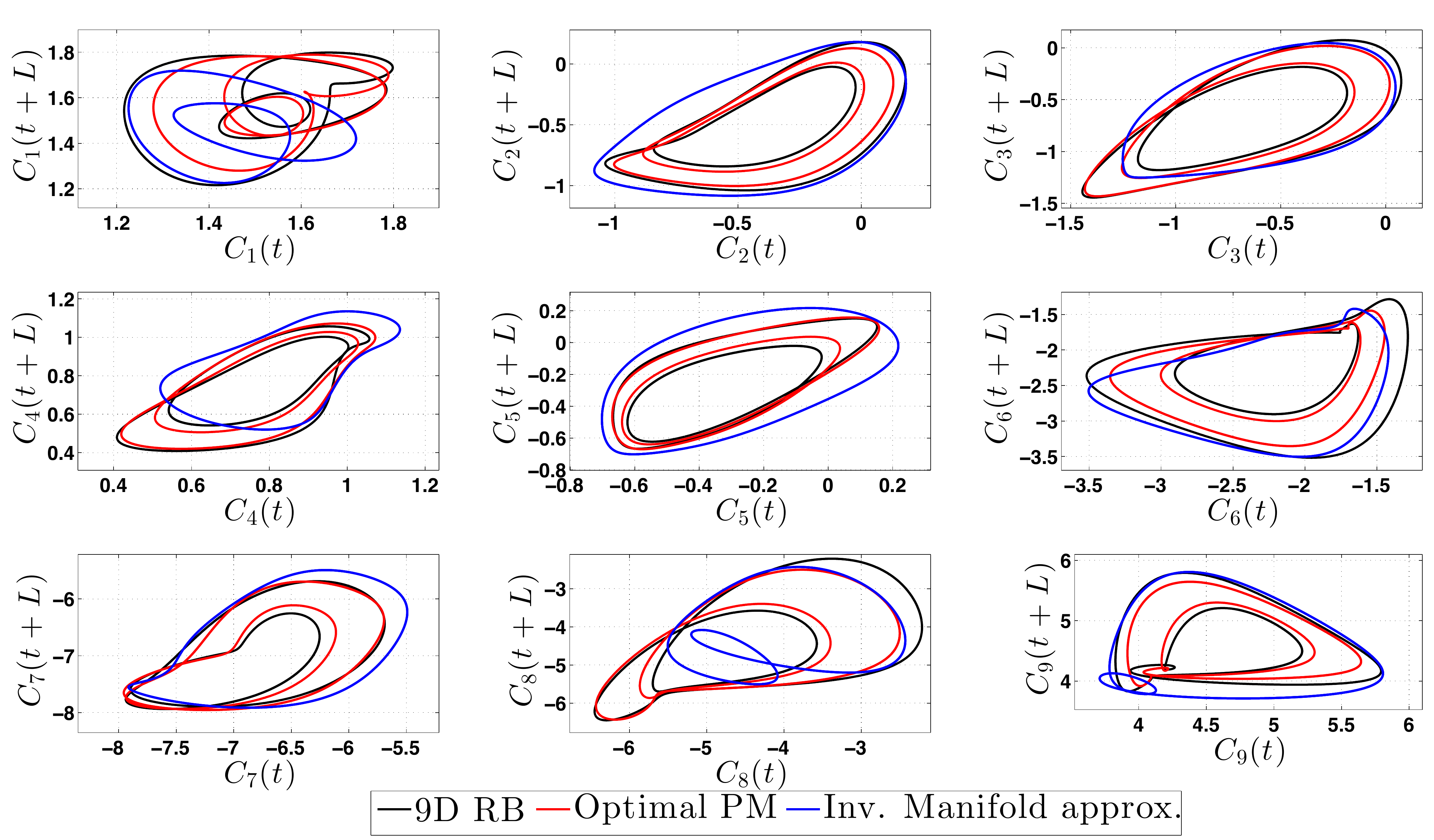}
\caption{{\footnotesize {\bf Attractor approximation for $r =14.1$ and $m=5$.} Here the attractor is projected onto the delay coordinates, $(C_j(t),C_j(t+L))$ ($1\leq j \leq 9$), for the original 9D RB system (black curve).  Here $L=1$. The approximation $\boldsymbol{C}^{\text{PM}}$ given by \eqref{Eq_soln_reconst_9DRBC} and obtained from the 5D optimal PM closure \eqref{Eq_RBC_reduced} is shown by the red curve. 
The approximation $\boldsymbol{C}^{\text{IM}}$ given by \eqref{Eq_CIM} and obtained from the 5D reduced system  \eqref{Eq_RBC_reduced_h1}  based on the invariant manifold approximation $h_2$,  is shown by the blue curve.}}   
\label{Fig_attractor_RB9D_perioddoubling}
\end{figure}

In contradistinction, $Q_7(\tau_7^\ast,T)$ is of same order of magnitude than the $Q_n$'s for modes $\boldsymbol{e}_6$, $\boldsymbol{e}_8$ and $\boldsymbol{e}_9$. As a result,  the optimal PM, $\Phi^{(1)}_{\bftau^\ast}$, provides comparatively, a much more efficient closure than 
when the parameterization $h_2$ is used. Figure \ref{Fig_attractor_RB9D_perioddoubling} shows for instance that in terms of attractor reconstruction, the approximation $\boldsymbol{C}^{\text{IM}}(t)$ given by \eqref{Eq_CIM} and obtained from the 5D reduced system  \eqref{Eq_RBC_reduced_h1}  based on  $h_2$ (blue curve), fails --- compared to its counterpart $\boldsymbol{C}^{\text{PM}}(t)$ obtained from the 5D optimal PM closure \eqref{Eq_RBC_reduced} (red curve) --- in capturing, within the embedded phase space,  the intricate behavior of the original model's periodic orbit  (black curve).  

%%%%%%%%%%%%%%%%%%%%%%%%%%%%%%%%%%%
%%%%%%%%%%%%%%%%%%%%%%%%%%%%%%%%%%%%%%%%%

%%%%%%%%%%%%%%%%%%%%%
\begin{table}[h] 
\caption{Averaged fraction of energy over $[t_0,t_0+T]$: Period-doubling regime}
\label{Table_Frac_energy1}
\centering
\begin{tabular}{ccccccccc}
\toprule\noalign{\smallskip}
   $\boldsymbol{e}_1$ &  $\boldsymbol{e}_2$ & $\boldsymbol{e}_3$   &  $\boldsymbol{e}_4$ &  $\boldsymbol{e}_5$ &  $\boldsymbol{e}_6$ & $\boldsymbol{e}_7$&  $\boldsymbol{e}_8$&  $\boldsymbol{e}_9$    \\ 
\noalign{\smallskip}\hline\noalign{\smallskip}
 $42.14\%$ & $42.14\%$ & $1.81\%$ & $3.87\%$ & $3.87\%$ & $4.27\%$ & $0.20\%$ & $0.86\%$ & $0.86\%$\\
\noalign{\smallskip} \bottomrule 
\end{tabular}
\end{table}
%%%%%%%%%%%%%%%%%%%%%%%%%%

\begin{figure}[hbtp]
   \centering
\includegraphics[width=0.9\textwidth, height=0.4\textwidth]{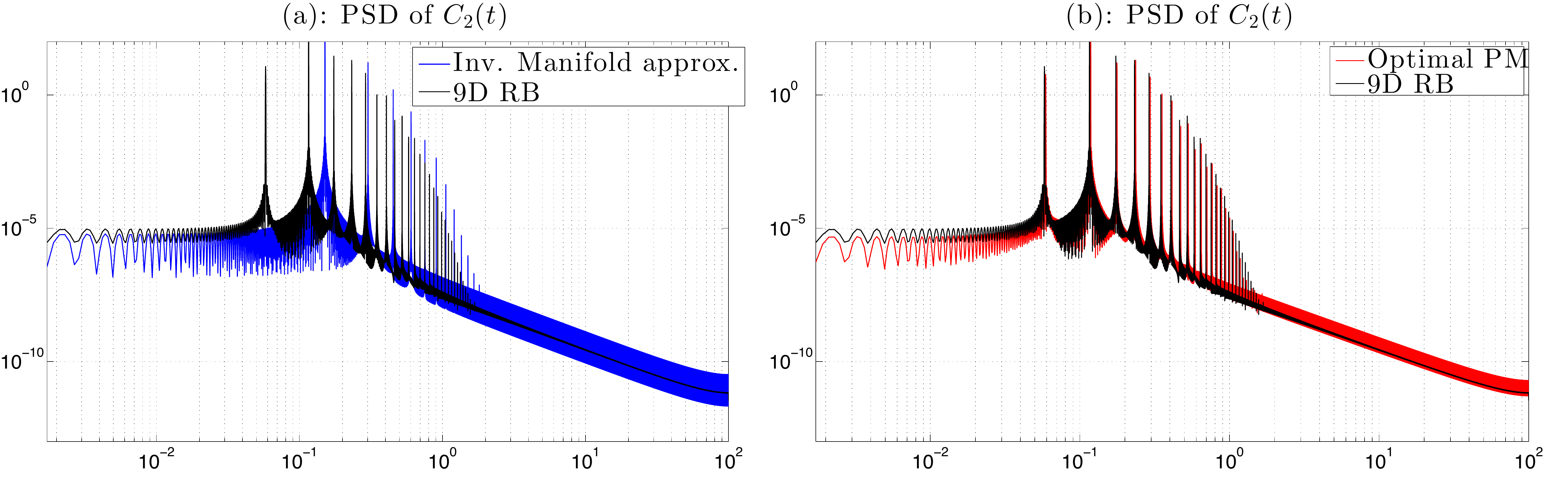}
\caption{{\footnotesize {\bf PSD approximation for $r =14.1$ and $m=5$.} Here the PSDs are estimated for $C_2(t)$ obtained from the original 9D RB system (black curve --- panels (a) and (b)), for $C_2^{\text{PM}}(t)$ obtained from the 5D optimal PM closure \eqref{Eq_RBC_reduced} (red curve -- panel (b)), and for  $C_2^{\text{IM}}(t)$ obtained from the 5D reduced system  \eqref{Eq_RBC_reduced_h1} based on invariant manifold approximation   (blue curve -- panel(a)). A semi-log scale is used for panels (a) and (b).}}   
\label{Fig_PSD_RB9D_period_doubling}
\end{figure}

A closer examination of the power spectral density (PSD) reveals that $\boldsymbol{C}^{\text{IM}}(t)$ fails in reproducing the dominant frequency and its subharmonics, whereas $\boldsymbol{C}^{\text{PM}}(t)$ captures them almost perfectly; compare panel (a) and (b) of Fig.~\ref{Fig_PSD_RB9D_period_doubling}.  The length of simulation $T_f$ for the original dynamics and the   5D optimal PM closure \eqref{Eq_RBC_reduced} used for the estimation of these PSDs is $T_f=1000$. Recall that for the latter, such results 
are obtained by optimizing the parameterization defects on a training interval of length $T$ equals only to three fourth of the period $T_p$  of the original dynamics, demonstrating thus good skills at least in the frequency domain.   
Similar skills than those shown in Fig.~\ref{Fig_PSD_RB9D_period_doubling}  for $C_2(t)$, hold for the other system's components.  

%%%%%%%%%%%%%%%%%%%%%%%%%%%%%%%%%%%
%%%%%%%%%%%%%%%%%%%%%%%%%%%%%%%%%%%%%%%%%

As progressing through the period-doubling cascade, the inability of the invariant manifold approximation, $h_2$, in reproducing the main  features of the RB system's solutions, is getting even worse, in particular right after the onset of chaos. The next section shows that the reduced systems \eqref{Eq_RBC_reduced}, to the contrary, provide still low-dimensional efficient closures (when driven by the appropriate optimal PM) for such chaotic regimes. 

\br\label{Remark_training_interval}
Depending on $t_0$ (after removal of transient), the global minimizer $\tau_n^\ast$ of $Q_n$, does not provide necessarily the best parameterization within the $\Phi_n$-class, and one may have to rely on the parameterization correlation $c(t)$  (see \eqref{Eq_corr_param}) to discriminate between other local minimizers of $Q_n$. We clarify here this statement which is relevant only for $n=6$ here; the global minima of $Q_7$, $Q_8$, and $Q_9$ being in fact robust as $t_0$ is varied.

For the regime analyzed here,  the ``wavy'' shape of the graph of $Q_6(\tau,T)$ may experience noticeable changes when $t_0$ varies. These changes may be manifested by the emergence of local minima that can modify substantially the location of the global minimizer and thus affect the determination of the optimal PM.

\begin{figure}[hbtp]
   \centering
\includegraphics[width=0.9\textwidth, height=0.35\textwidth]{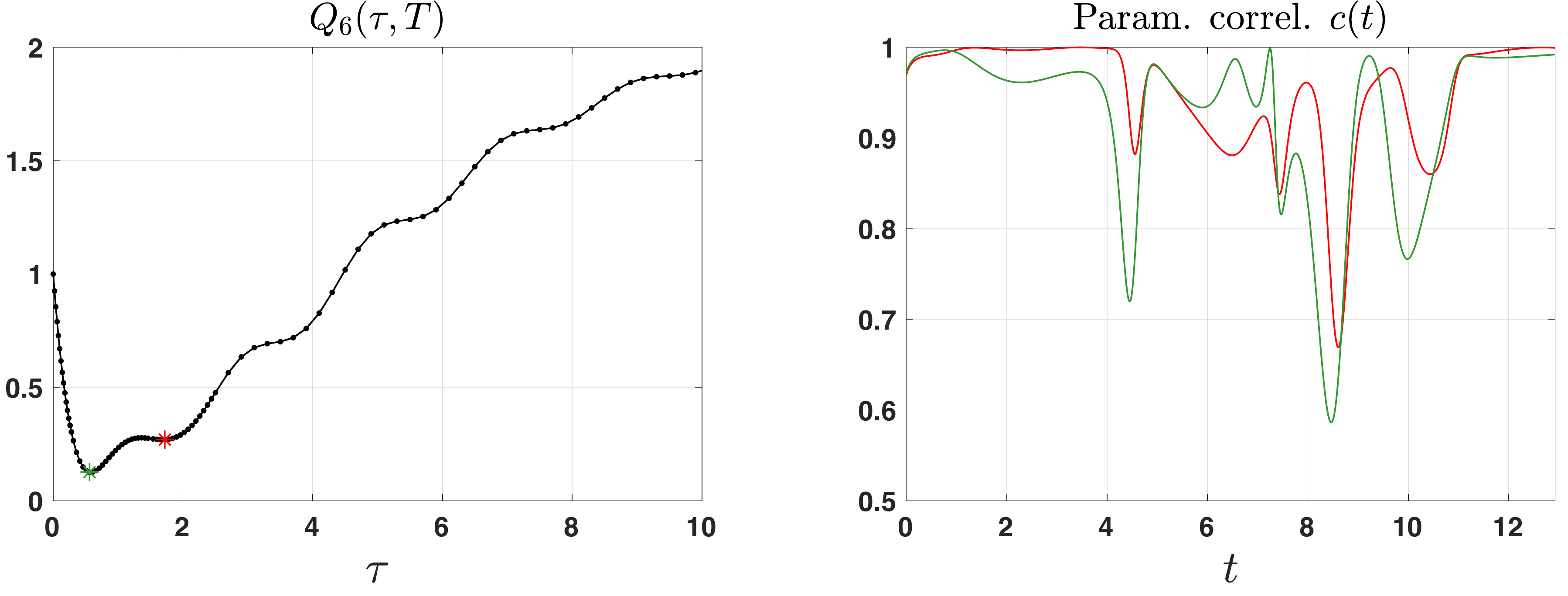}
\caption{{\footnotesize {\bf Selection of suboptimal parameterization  via parameterization correlation.}  The parameterization correlation $c(t)$ are shown in the right panel for an interval of length $T=3T_p/4$  in the period-doubling regime. {\hl Here $c(t)$ is computed from \eqref{Eq_corr_param} with $\Psi=\Phi^{(1)}_{\bftau}$ for two choices of $\bftau$. Choice 1: $\tau_n = \tau_n^*$ for all the components (green curve). Choice 2: $\tau_n = \tau_n^*$ except $\tau_6$, which is taken instead to be the local minimizer marked by the red asterisk on the left panel (red curve).}}}  
\label{Fig_Correlation_discrim}
\end{figure}

For instance the left panel of Fig.~\ref{Fig_Correlation_discrim} shows $Q_6(\tau,T)$ as obtained from another segment of the solution $y(t)$ to 
\eqref{Eq_RBC_eigenbasis}  (in the period-doubling regime), that is for another $t_0$ in \eqref{toto_Qn} than used for Fig.~\ref{fig:RB9D_Qplot_perioddoubling}. A simple visual comparison reveals that the global minimum shown for $Q_6$ in Fig.~\ref{fig:RB9D_Qplot_perioddoubling}  corresponds now to a local minimum (red asterisk), and a new global minimum  closer to $\tau=0$ has appeared (green asterisk). 

If one selects the corresponding global minimizer as $\tau_6^\ast$, the corresponding optimal closure captures only an excerpt of the dominant frequency and its harmonics (every other frequency more precisely), and the closure fails in reproducing the period-doubling.  
This issue can be easily fixed by the inspection of $c(t)$ given by \eqref{Eq_corr_param}  over [0,T]. Indeed, by using the optimal PM for which $\tau_6^\ast$ corresponds to the global minimum and the (sub)optimal PM for which $\tau_6^\ast$  corresponds to the second local minimum, we obtain two curves for $c(t)$: one associated with the optimal parameterization  (global minimum/green curve) and one associated with the suboptimal parameterization (local minimum/red curve).

The red curve is clearly closer to 1 than the green one (in average), indicating that  $\tau_6^\ast$  corresponding to the second local minimum (i.e.~the suboptimal parameterization) should be in fact retained for determining the parameterization $\Phi_n$, as indeed the corresponding PM  closure provides then similar modeling skills to those shown in Fig.~\ref{Fig_PSD_RB9D_period_doubling}.

This discrimination, made possible thanks to the parameterization correlation, $c(t)$, (prior to any simulation of \eqref{Eq_RBC_reduced}) teaches us the relevance of this non dimensional number to refine the determination of an optimal PM in practice, beyond this example and especially in presence of other local minima for a given $Q_n$ as $t_0$ is varied.   

Other tests conducted in other parameter regimes indicate that such a situation requiring the discrimination via an inspection of $c(t)$ and a selection of a suboptimal rather than optimal parameterization is rather the exception than the rule\footnote{For instance this issue is not encountered for the chaotic regime analyzed in Sec.~\ref{Sect_chaotic-RB}.}; namely the parameterization corresponding to a global minimizer of $Q_n$, provides in general the best closure results.  Nevertheless we decided to communicate on this issue subordinated to the presence of local minima as it may be encountered for other systems.

%%%%%%%%%%%%%%%%%%%%%%%%%%%%%%%%%%%
%%%%%%%%%%%%%%%%%%%%%%%%%%%%%%%%%%%%%%%%%
\er

%%%%%%%%%%%%%%%%%%%%RBC: CHAOTIC REGIME%%%%%%%%%%%%%%%%%%%%%%%%%%%%%%%%%
\subsection{Closure in a chaotic regime} \label{Sect_chaotic-RB}
We assess in this section the skills of the optimal PM closure \eqref{Eq_RBC_reduced} in a regime located right after the onset of chaos, after  the system has gone through a period doubling cascade, i.e.~for  $r=14.22$. 
We conduct also hereafter an analysis on the effect of the reduced dimension, $m$, of the reduced state space $E_\c$. 
Still this reduced state space is spanned by few dominant eigenmodes of the linear part $L$ of the perturbed system  \eqref{Eq_RBC_fluct} about the mean state $\overline{\boldsymbol{C}}$ is given by  \eqref{RB_Linear_part_perturbed}, with now the latter estimated, after removal of transient dynamics, over the training interval of length  $T=T_p$, with $T_p$ denoting the period of the solution for $r=14.1$; see previous section. 
 
Here again, the unresolved modes are parameterized by the quadratic manifold, $\Phi_n(\tau,\cdot)$, given by \eqref{Eq_Phi_tau}, optimized over the training interval $[0,T]$ by minimizing the parameterization defect  $Q_n$ given by \eqref{toto_Qn}. The distribution of energy per mode for 
this regime is shown in Table \ref{Table_Frac_energy2}. The distribution of energy is explained due to the ordering \eqref{Eq_gap} adopted here from Sec.~\ref{Sec_loc_invman}, i.e.~by ordering the modes  according to their linear rate of growth/decay; for this parameter regime again,   from the least to the most stable modes. Since $\boldsymbol{e}_4$ and $\boldsymbol{e}_5$ come in pairs (i.e.~$\textrm{Re}(\beta_4)=\textrm{Re}(\beta_5)$), we analyze hereafter the cases $m=3$, $m=5$ and $m=6$. Thus from Table \ref{Table_Frac_energy2}, the 
energy to be parameterized corresponds to $16.6\%$ of the total energy (over $[0,T]$) for the case $m=3$,  to $6.8\%$  for $m=5$, and to $2.85\%$ for $m=6$. 
%%%%%%%%%%%%%%%%%%%%%%%%%%
\begin{table}[h] 
\caption{Averaged fraction of energy over $[0,T]$: Chaotic regime}
\label{Table_Frac_energy2}
\centering
\begin{tabular}{ccccccccc}
\toprule\noalign{\smallskip}
     $\boldsymbol{e}_1$ &  $\boldsymbol{e}_2$ & $\boldsymbol{e}_3$   &  $\boldsymbol{e}_4$ &  $\boldsymbol{e}_5$ &  $\boldsymbol{e}_6$ &  $\boldsymbol{e}_7$&  $\boldsymbol{e}_8$&  $\boldsymbol{e}_9$    \\ 
\noalign{\smallskip}\hline\noalign{\smallskip}
 $37.59\%$ &  $37.59\%$ &  $8.23\%$&  $4.90\%$&  $4.90\%$&  $3.95\%$&  $0.31\%$&  $1.27\%$&  $1.27\%$\\
\noalign{\smallskip} \bottomrule 
\end{tabular}
\end{table}
%%%%%%%%%%%%%%%%%%%%%%%%%% 

 %%%%%%%%%%%%%%%%%%%%%%%%%%%%%%
\begin{table}[h] 
\caption{Optimal parameterization defects for $T=25$: Chaotic regime}
\label{table_Q}
\centering
\begin{tabular}{cccc}
\toprule\noalign{\smallskip}
      &  $m=3$ & $m=5$   &  $m=6$ \\ 
\noalign{\smallskip}\hline\noalign{\smallskip}
 $Q_4(\tau_4^\ast,T)$ & $0.09$  &    &  \\
 $Q_5(\tau_5^\ast,T)$ &  $0.09$ &    &   \\
 $Q_6(\tau_6^\ast,T)$ & $ 0.38$   &  $0.12 $ &   \\
 $Q_7(\tau_7^\ast,T)$ &  $0.22$  & $0.2$ &  $0.04$\\
 $Q_8(\tau_8^\ast,T)$ &  $0.05$ &  $0.09$ &  $0.02$ \\ 
  $Q_9(\tau_9^\ast,T)$ & $ 0.05$ &   $0.09 $  &  $0.02$\\
\noalign{\smallskip} \bottomrule 
\end{tabular}
\end{table}
%%%%%%%%%%%%%%%%%%%%%%%%%%%%%%%

Given the solution $y(t)$  of Eq.~\eqref{Eq_RBC_eigenbasis} over $[0,T]$, the minimal values  $Q_n(\tau_n^\ast,T)$ achieved by the optimal PM, $\Phi^{(1)}_{\bftau^*}$, in terms of the reduced dimension $m$ are shown in Table \ref{table_Q}. Obviously, the case $m=6$ comes with the smaller parameterization defects, while the case $m=3$ presents for the modes  $\boldsymbol{e}_6$ and  $\boldsymbol{e}_7$, values that although less than $1$ are not on the same order of magnitude than the other values of $Q_n$.  

The energy left after application of the optimal PM, represents $0.04\times 0.31 +2\times 0.02 \times 1.27=0.063\%$ of the total energy for the case $m=6$, and represents $0.765\%$ for the case $m=5$, still below $1\%$ of the total energy. To the contrary,  an amount of energy representing $5.42\%$ needs still to be parameterized after application of the optimal PM for the case $m=3$. Compared with the fraction of energy left in the corresponding unresolved modes prior parameterization, an application of the optimal PM leads to an improvement by a factor approximately equal to 45 for $m=6$, and equal to $9$ and to $3$ for respectively $m=5$ and  $m=3$. Without any surprise, the cutoff corresponding to the smallest amount of energy to be parameterized (i.e.~when $m=6$) comes with the 
best improvement in terms of parameterization when the optimal PM is used. On the other hand, the cutoff corresponding to the biggest amount of energy (i.e.~when $m=3$) comes with the poorest parameterization score in terms of energy that still needs to be parameterized after application of the optimal PM.  
Thus, one expects that an optimal PM closure should perform certainly better for $m=6$ than for $m=3$, and must show some improvements compared to the optimal PM closure for $m=5.$ 

%%%%%%%%%%%%%%%%%%PARAM ANGLE%%%%%%%%%%%%%%%%%%%%%
\begin{figure}[hbtp]
   \centering
 \includegraphics[height=.45\textwidth,width=.95\textwidth,]{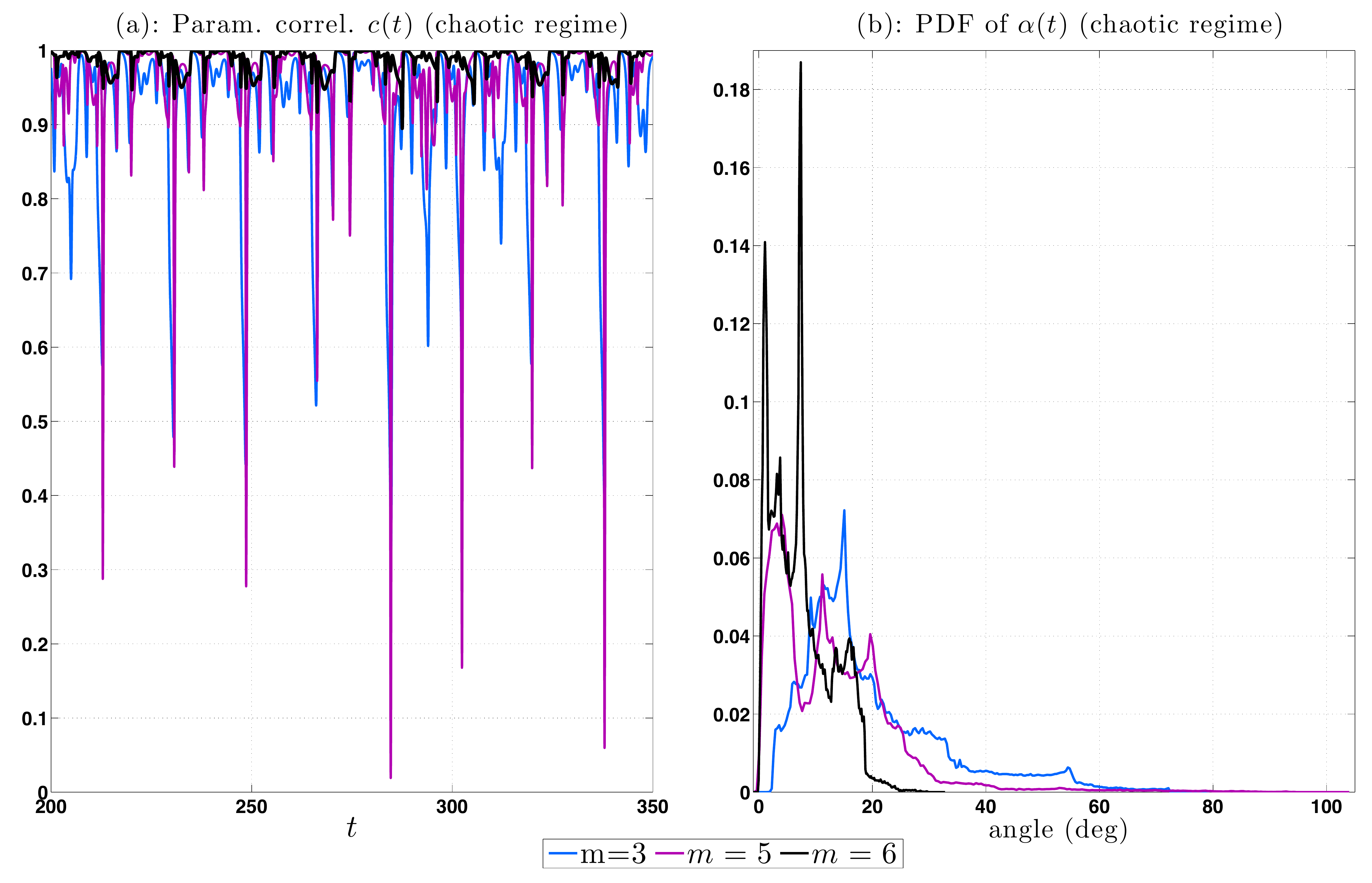}
  \caption{{\footnotesize {\bf Effect of the reduced dimension $m$: Diagnostic for $r =14.22$.} This effect is shown here on the parameterization correlation $c(t)$ (panel (a)) and the PDF of the parameterization angle $\alpha(t)$ (panel (b)) for the chaotic regime. Here $c(t)$ and  $\alpha(t)$ are respectively computed from \eqref{Eq_corr_param}  and \eqref{Eq_alpha}, with $\Psi=\Phi^{(1)}_{\bftau^*}$, the optimal PM. }}   
 \label{fig:RB9D_correlation_content}
\end{figure}
%%%%%%%%%%%%%%%%%%%%%%%%%%%%%%%%%%%%%%%% 

 %%%%%%%%%%%%%%%%%%PSDs and ACFs%%%%%%%%%%%%%%%%%%%%%
\begin{figure}[hbtp]
   \centering
 \includegraphics[height=.55\textwidth,width=.95\textwidth,]{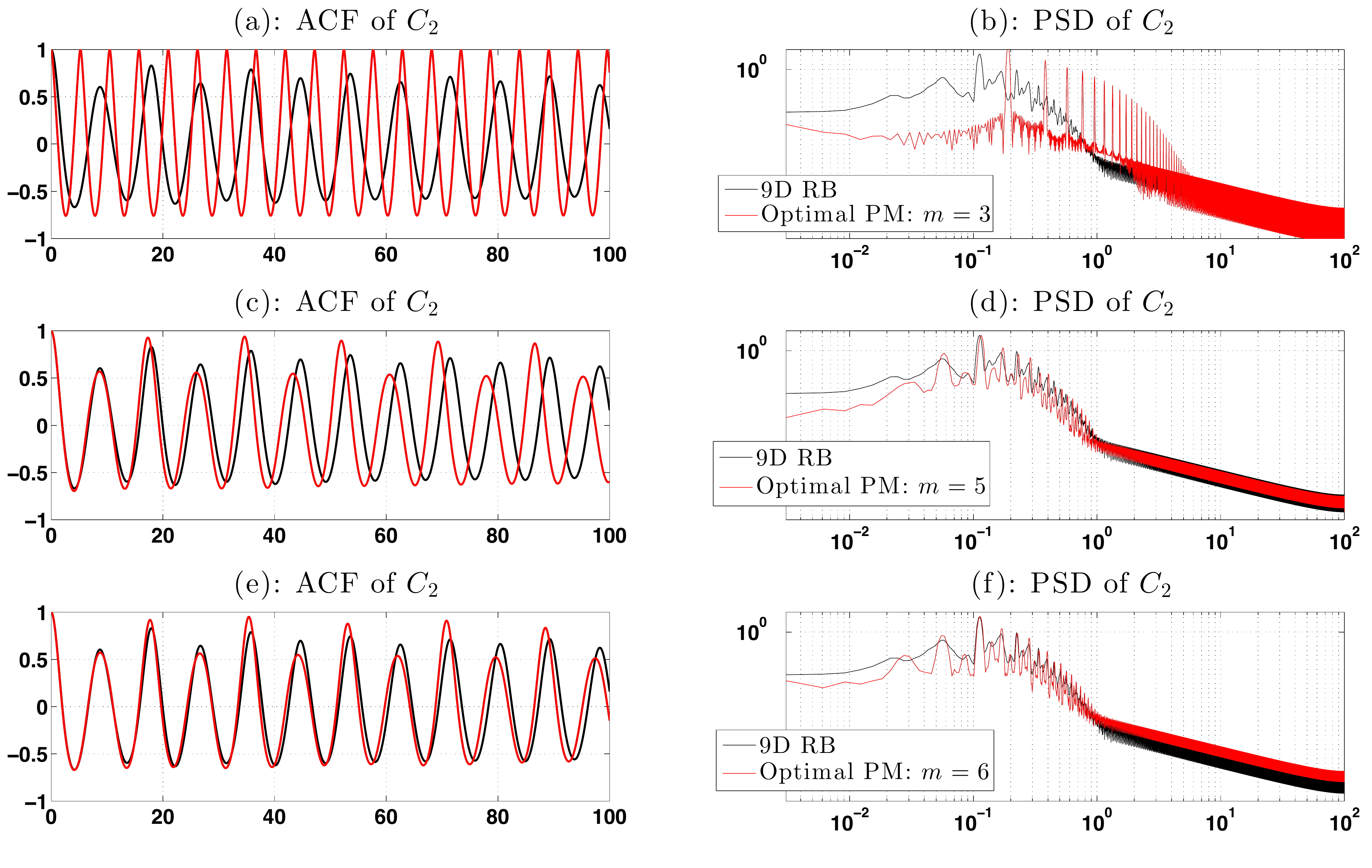}
  \caption{{\footnotesize {\bf Effect of the reduced dimension $m$: Simulation for $r =14.22$.} This effect is shown  for the chaotic regime on the ability of the optimal PM closure \eqref{Eq_RBC_reduced} to reproduce the PSD and ACF, for the second component $C_2$. A semi-log scale is used for panels (b), (d) and (f).}}   
 \label{Fig_PSDs_ACFs}
\end{figure}
%%%%%%%%%%%%%%%%%%%%%%%%%%%%%%%%%%%%%%%%
This energy budget analysis is comforted by the analysis of the parameterization correlation $c(t)$ and of the probability density function (PDF) of the parameterization angle $\alpha(t)$.  Here $c(t)$ and  $\alpha(t)$ are respectively computed from \eqref{Eq_corr_param}  and \eqref{Eq_alpha}, with $\Psi=\Phi^{(1)}_{\bftau^*}$, the optimal PM as determined for each case, $m=3$, $m=5$, and  $m=6$, from \eqref{Eq_Phi_tau_opt}, for which the optimal vector $\bftau^*$ is obtained by minimization of \eqref{toto_Qn} for the relevant $n$. As shown in panel (b) of Fig.~\ref{Fig_PSDs_ACFs}, each of these PDFs is skewed towards zero. Nevertheless the PDF that is the most concentrated (i.e.~with more mass) near zero corresponds to the case $m=6$ (black curve), then comes the PDF associated with the case $m=5$ (magenta curve), and finally the PDF for the case $m=3$ (blue curve).

These diagnostics are confirmed  when looking at the ability of the corresponding optimal PM closures \eqref{Eq_RBC_reduced}, in reproducing key statistics of the original model's dynamics such as autocorrelation functions (ACFs) and PSDs.  For the regime analyzed here ($r=14.22$), the time-variability of the chaotic dynamics is characterized by a broad band spectrum visible in each component's PSD. The black curve in either right panels of Figure \ref{Fig_PSDs_ACFs}, shows such a broad band spectrum for e.g.~the PSD of $C_2$ as estimated from integration of  Eq.~\eqref{Eq_9DRBC} after a simulation of length $T_f=1000$. Other components display similar PSDs. 

\begin{figure}[hbtp]
   \centering
\includegraphics[width=0.95\textwidth, height=0.6\textwidth]{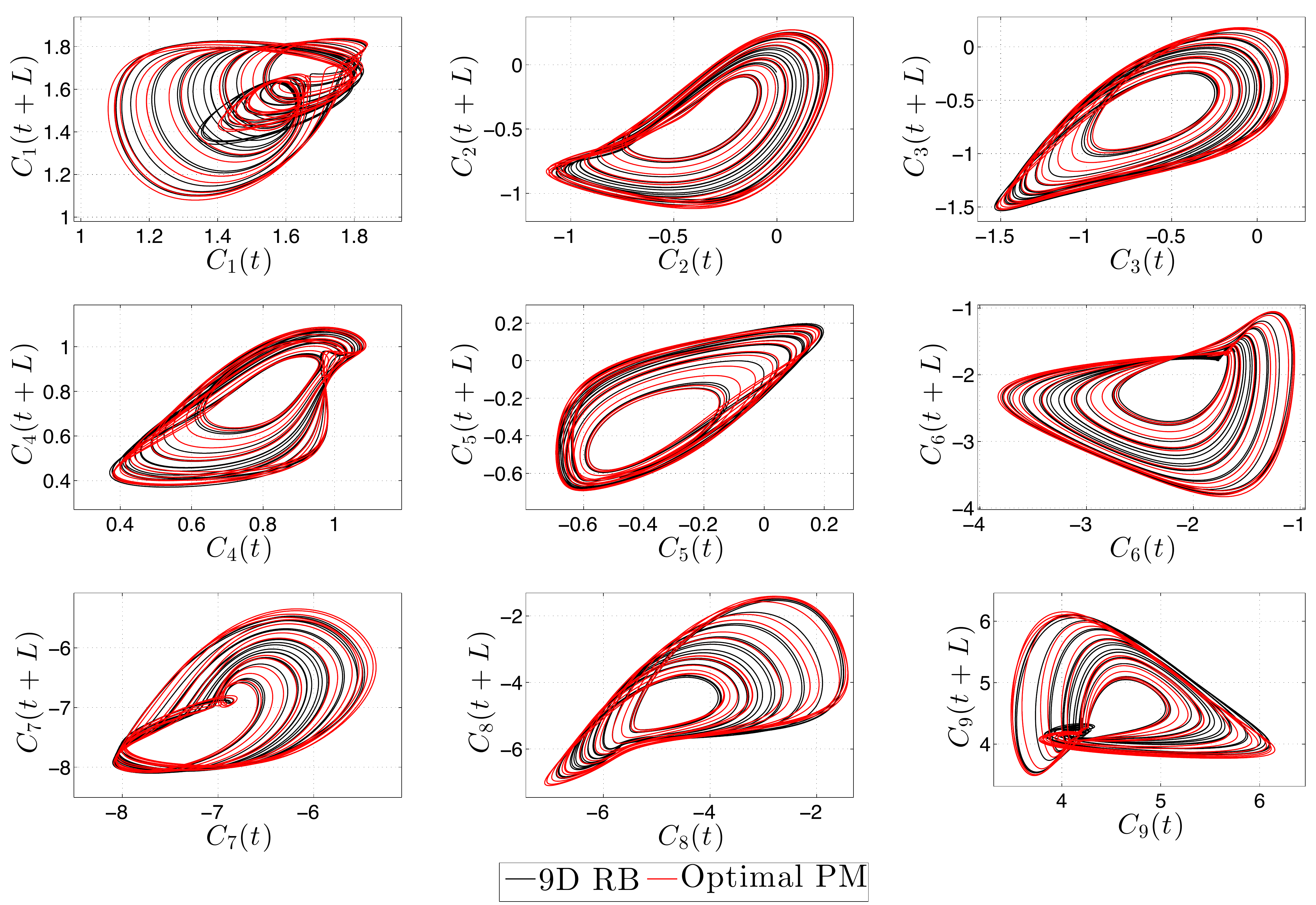}
  \caption{{\footnotesize {\bf Attractor approximation for $r =14.22$ and $m=6$}. Same as in Fig.~\ref{Fig_attractor_RB9D_perioddoubling} except $r =14.22$  (chaotic regime) and $m=6$. Here also $L=1$.}}  
   \label{Fig_RB9D_chaos_attractor}
\end{figure}

Figure \ref{Fig_PSDs_ACFs} shows clearly, as anticipated by the energy budget analysis on a short interval $[0,T]$ (with $T=17.25$) (and supported by the parameterization angle's PDF analysis), that the 5D and 6D optimal PMs provide efficient closures, with a noticeable  improvement for the ACF's reproduction of $C_2$ when the 6D optimal PM is used; see panel (e) of Fig.~\ref{Fig_PSDs_ACFs}.  
Furthermore, Fig.~\ref{Fig_RB9D_chaos_attractor} shows that the 6D optimal PM closure leads to an excellent approximation of the original model's attractor, whereas the 5D optimal PM closure although reproducing correctly most of its features fails in reproducing certain solution's large excursions in the embedded phase space (not shown). 
The 3D optimal PM fails however dramatically in the approximation of this attractor as it leads to a periodic orbit and fails thus to reproduce the time variability of the original model's chaotic dynamics; see panels (a) and (b) of Fig.~\ref{Fig_PSDs_ACFs}.

%%%%%%%%%%%%%%%%%%%ATTRACTOR%%%%%%%%%%%%%%%%%%%%%
%%%%%%%%%%%%%%%%%%%%%%%%%%%%%%%%%%%%%%%%

Based on these results, we may state that our parameterization formula of Sec.~\ref{Sect_PM_with_forcing} (i.e.~$\Phi^{(1)}_{\bftau^\ast}$ given by \eqref{Eq_Phi_tau_opt}) provides here,
seemingly, a good approximation of the optimal PM as given by the abstract Theorem \ref{Thm_variational-pb} when $m=5$ and $m=6$.  
Our optimal PM as computed for the case $m=3$, although leading to a periodic orbit, may still be a good approximation of the theoretical 
optimal parameterization \eqref{Def_h2} averaging out the unresolved variables, for the reduced state space, $E_\c=\textrm{span}\{\boldsymbol{e}_1,\boldsymbol{e}_2,\boldsymbol{e}_3\}$. It is indeed possible that  the conditional expectation as defined in Theorem \ref{Thm_variational-pb2}, gives a periodic solution for a given reduced state space. The theory of Sec.~\ref{Sect_PM_reduction} does not exclude such a scenario. 

To improve the results in the case $m=3$, stochastic parameterizations may be then superimposed to our optimal PM in order to further reduce the parameterization defect. This topic is out of the scope of the present paper but will be pursued elsewhere; see Concluding Remarks in Sec.~\ref{Sec_concluding_rmk}.

\subsection{Heat flux analysis}
We analyze here how the optimal LIA parameterization behaves in the physical domain, for the chaotic regime. We focus on the vertical heat flux, accomplished by the fluctuations around the time-averaged state that enables the system to sustain statistical equilibrium. 
Once a solution $\boldsymbol{C}(t)$ to Eq.~\eqref{Eq_9DRBC} is computed, one can evaluate the following local heat flux 
\be
  H(\boldsymbol{x},t)= w(\boldsymbol{x},t)  \theta'(\boldsymbol{x},t)-\partial_z \overline{\theta}(\boldsymbol{x}), \qquad \boldsymbol{x}=(x,y,z),
\ee
where $w$ denotes the vertical velocity, and $\theta'$ denotes the anomaly of the temperature $\theta$ with respect to  the time-mean temperature 
 $ \overline{\theta}$.  The vertical velocity $w$ and temperature $\theta$ are computed according to  Eqns.~(12) and (17) of \cite{Reiterer_al98}.

\begin{figure}[hbtp]
   \centering
\includegraphics[width=1\textwidth, height=0.35\textwidth]{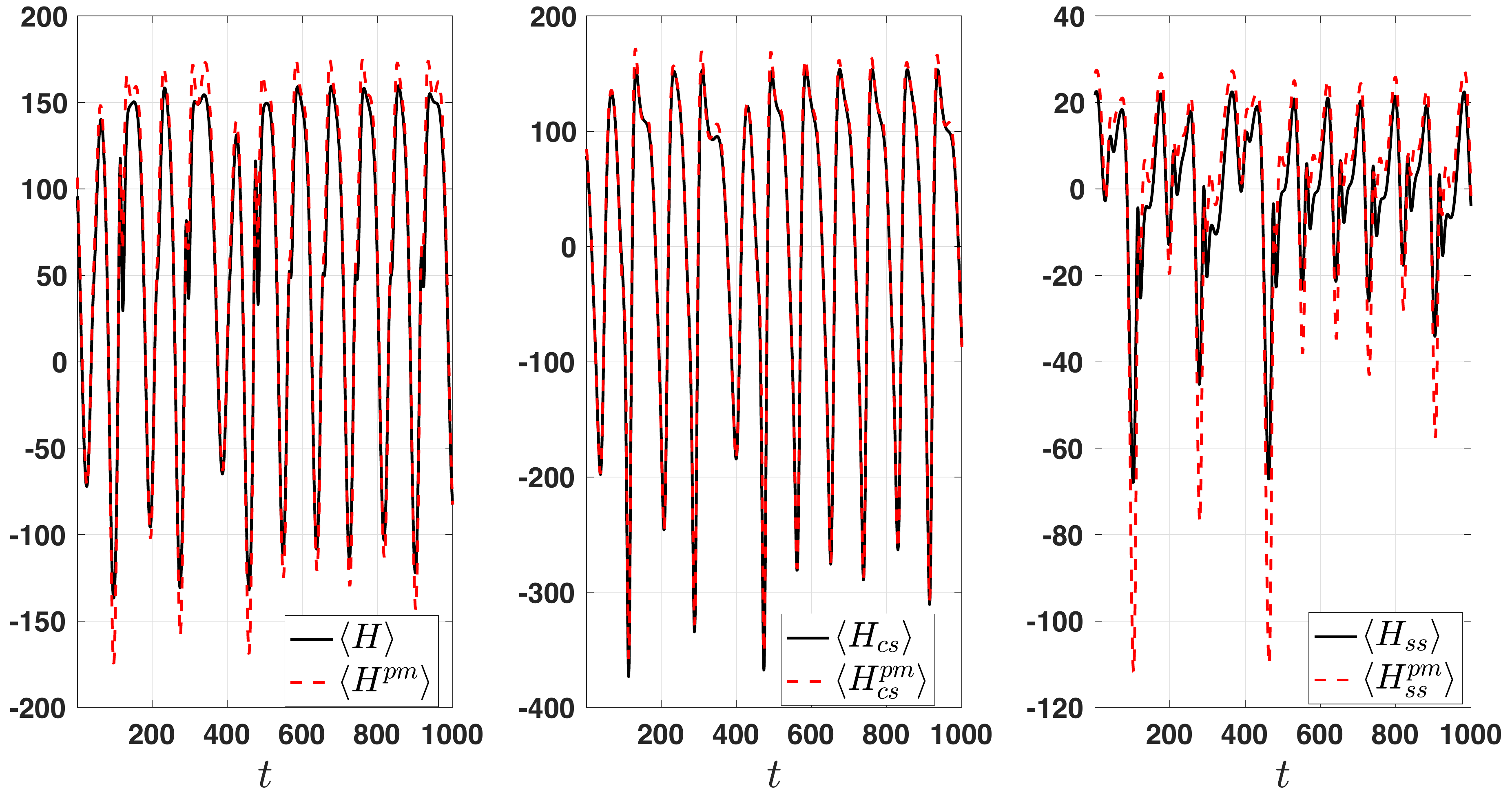}
  \caption{{\footnotesize {\bf Space-average heat fluxes for the chaotic regime}. Here the reduced state space is five-dimensional ($m=5$).}}  
   \label{Fig_RB9D_heatflux5D}
\end{figure}

\begin{figure}[hbtp]
   \centering
\includegraphics[width=1\textwidth, height=0.35\textwidth]{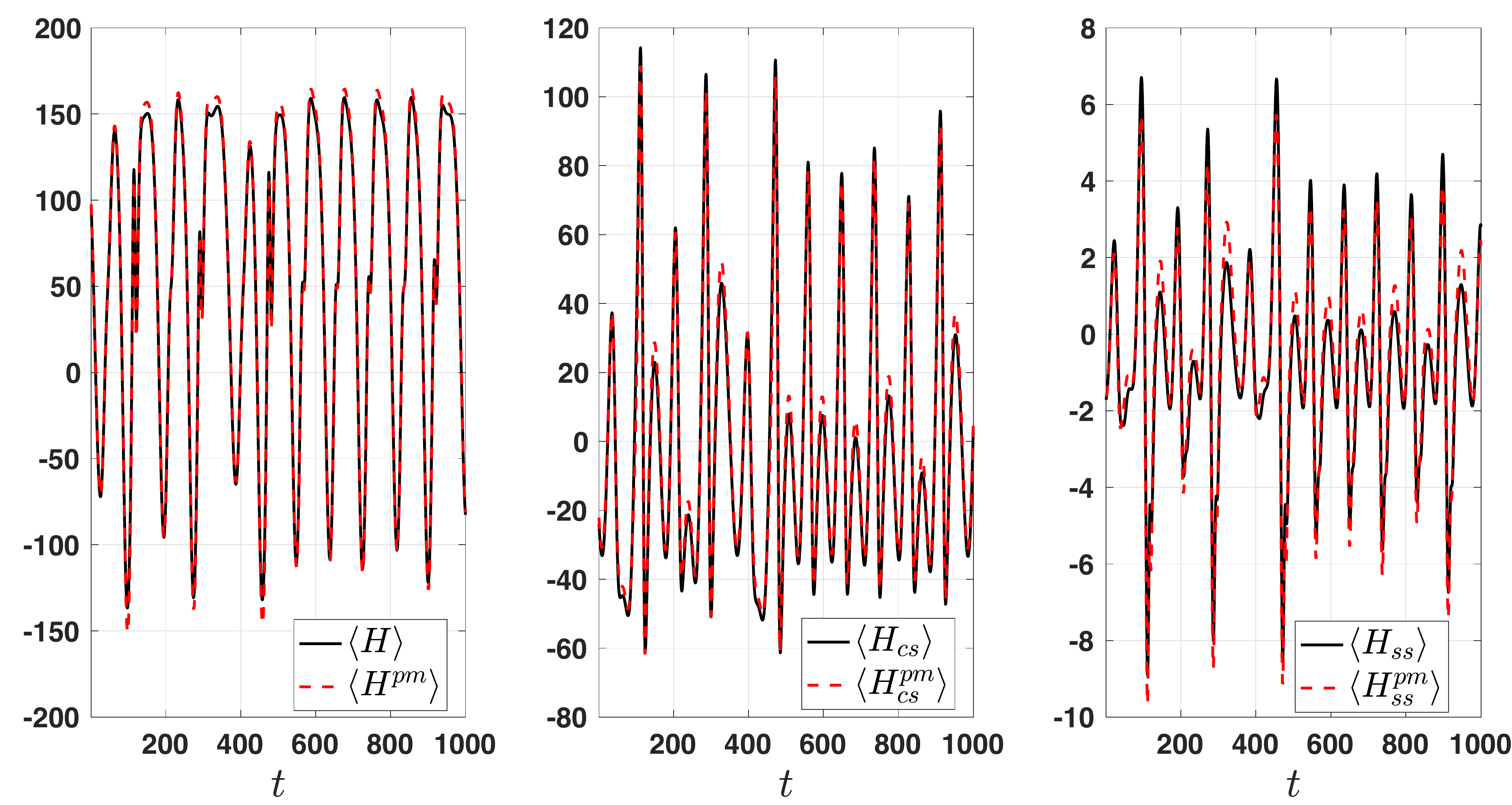}
  \caption{{\footnotesize {\bf Space-average heat fluxes for the chaotic regime}. Here the reduced state space is six-dimensional ($m=6$).}}  
   \label{Fig_RB9D_heatflux6D}
\end{figure}

Recall that our optimal PM is determined for the transformed variables, namely for Eq.~\eqref{Eq_RBC_eigenbasis}. 
In particular our splitting between low and high modes is made within the system of coordinates in the $y$-variable.   
By transforming back into the original variables we can trace the contribution of the high and low modes (defined in the transformed variables) into the original system of coordinates. By doing so, the  heat flux  $H(\boldsymbol{x},t)$ decomposes as 
\be\label{Eq_Hdecomp}
H(\boldsymbol{x},t)=H_{\c\c} (\boldsymbol{x},t)+ H_{\c \s}(\boldsymbol{x},t) +H_{\s \s}(\boldsymbol{x},t).
\ee
with 
\bea\label{Eq_Hcomp}
H_{\c\c}(\boldsymbol{x},t)&=w_{\c}(\boldsymbol{x},t)  \theta_{\c}'(\boldsymbol{x},t)-\partial_z \overline{\theta_{\c}}(\boldsymbol{x}),\\
H_{\s\s}(\boldsymbol{x},t)&=w_{\s}(\boldsymbol{x},t)  \theta_{\s}'(\boldsymbol{x},t)-\partial_z  \overline{\theta_{\s}}(\boldsymbol{x}),\\
H_{\c \s}(\boldsymbol{x},t)&=w_\c(\boldsymbol{x},t)  \theta_\s'(\boldsymbol{x},t)+w_\s(\boldsymbol{x},t)  \theta_{\c}'(\boldsymbol{x},t).
\eea
When the high-mode contribution in  \eqref{Eq_Hdecomp} and \eqref{Eq_Hcomp} is replaced by the optimal LIA parameterization derived in the previous section (chaotic regime), 
errors in the ``low-high'' and ``high-high'' interactions  to the heat flux are visible.
Table \ref{table_H_error} shows these relative errors in the $L^2$-norm in time, after space average $\langle \cdot \rangle $.    
Clearly these errors reduce as the dimension of the reduced state space (in the transformed variables) increases, but overall the reproduction of the time-variability of $\langle H\rangle $ is satisfactory, especially when $m=6$; see Figs.~\ref{Fig_RB9D_heatflux5D} and \ref{Fig_RB9D_heatflux6D}. As a comparison when only the low modes are used to approximate the heat flux like in a Galerkin truncation, the heat flux errors are substantially larger; see Table \ref{table_H_errorb}. 
Without any surprise the improvement brought by the high-mode parameterization is more pronounced when $m=5$ than when $m=6$.   
Taking volume- and time-average in \eqref{Eq_Hdecomp}, we observe that $\overline{\langle H\rangle}=54.6$. Doing the same operation in which 
the $\s$-variable is replaced by its high-mode approximation (as given by the optimal LIA) gives $\overline{\langle H^{\textrm{app}}\rangle}=61.4$ for $m=5$, and  $\overline{\langle H^{\textrm{app}}\rangle}=56.1$, for $m=6.$

 %%%%%%%%%%%%%%%%%%%%%%%%%%%%%%
\begin{table}[h] 
\caption{Heat fluxes: Relative error when  ``$\s$'' is replaced by optimal PM}
\label{table_H_error}
\centering
\begin{tabular}{ccc}
\toprule\noalign{\smallskip}
       & $m=5$   &  $m=6$ \\ 
\noalign{\smallskip}\hline\noalign{\smallskip}
 $\langle H \rangle $  &  $15 \%$  & $4.5\%$ \\
 $\langle H_{\c\s} \rangle$  &  $7.6 \%$  & $11.2\%$ \\
  $\langle H_{\s\s}\rangle $ &   $64\%$  &   $21.9\%$ \\
\noalign{\smallskip} \bottomrule 
\end{tabular}
\end{table}
%%%%%%%%%%%%%%%%%%%%%%%%%%%%%%%

 %%%%%%%%%%%%%%%%%%%%%%%%%%%%%%
\begin{table}[h] 
\caption{Relative error $\mathcal{E}_\c=|\langle H-H_{\c\c}\rangle|_{L^2}/|\langle H \rangle|_{L^2}$}
\label{table_H_errorb}
\centering
\begin{tabular}{ccc}
\toprule\noalign{\smallskip}
       & $m=5$   &  $m=6$ \\ 
\noalign{\smallskip}\hline\noalign{\smallskip}
  $\mathcal{E}_\c$ &   $132\%$  &   $35\%$ \\
\noalign{\smallskip} \bottomrule 
\end{tabular}
\end{table}
%%%%%%%%%%%%%%%%%%%%%%%%%%%%%%%

%%%%%%%%%%%%%%%%%%%Heat fluxes%%%%%%%%%%%%%%%%%%%%%
%%%%%%%%%%%%%%%%%%%%%%%%%%%%%%%%%%%%%%%%

%%%%%%%%%%%%%%%%%%%Heat fluxes%%%%%%%%%%%%%%%%%%%%%
%%%%%%%%%%%%%%%%%%%%%%%%%%%%%%%%%%%%%%%%

%%%%%%%%%%%%%%%%%%%%%%%%%%%%%%%%%%%%%%%%%%%%%%%%%%%%%
\section{Closing Kuramoto-Sivashinsky turbulence and fixing backscatter errors}\label{Sec_KS_turbulence}
%%%%%%%%%%%%%%%%%%%%%%%%%%%%%%%%%%%%%%%%%%%%%%%%%%%%%%%%%%%%%%%%%%
 In this section we show that the PM approach allows for deriving efficient closures for the Kuramoto-Sivashinsky (KS) turbulence, in strongly turbulent regimes. The closure results presented hereafter are obtained for cutoff scales placed well within the inertial range, keeping only the unstable modes in the reduced state space.  The underlying optimal PMs obtained by our variational approach are far from slaving and  allow for  remedying the excessive backscatter transfer of energy to the low modes encountered by the LIA or the QSA parameterizations in their standard forms, when they are used at this cutoff wavelength.

\subsection{Preliminaries and background}
We consider the KS equation (KSE) \cite{kuramoto1976persistent,sivashinsky1977nonlinear} posed on the domain, $\mathcal{D}=(0,L)$,  and subject to periodic boundary conditions:
\be 
\label{eq:KSE}
\partial _t u =  - \nu \partial_{x}^4 u - D \partial_x^2 u   -  \gamma u  \partial_x u, 
\ee
where $\nu,D$ and $\gamma$ are positive parameters. 
The KSE is commonly considered as a basic case study for spatio-temporal chaos.

Note that the KSE  in its formulation \eqref{eq:KSE} can be rescaled as posed on the interval $(0,2\pi)$:
\be
\label{eq:KSE3}
\partial _{\overline{t}} \overline{u} = - 4 \partial_{\overline{x}}^4 \overline{u} - \alpha \bigg(\partial_{\overline{x}}^2 \overline{u}  + \overline{u}  \partial_{\overline{x}}\overline{u}\bigg),
\ee
by using the following scaling
\be \label{Eq_conversion_I_to_III}
L =  \sqrt{\frac{\nu \alpha}{D}} \pi, \quad u = \frac{2 D^{3/2}}{\gamma \sqrt{\nu \alpha}} \overline{u},  \quad  x = \frac{\sqrt{\nu\alpha}}{2\sqrt{D}}\overline{x}, \quad  t = \frac{\nu \alpha^2}{4D^2} \overline{t}.
\ee 
Although mathematically equivalent, depending on the purpose one may prefer one formulation to the other for the closure exercises considered hereafter; see Remark \ref{Rmk_rescaling}.

We aim at closure of the KSE. Various purposes are pursued regarding what a low-dimensional closure should do and this may cause confusion when comparing methods. Among the purposes targeted in the literature concerning the closure/reduction problem of the KSE, are the following: (i) finite-time approximation error such as in AIM theory \cite{MT89,devulder1993rate} or renormalization group (RG) methods \cite{schmuck2015new}, (ii) reproduction of local and global bifurcations \cite{armbruster1989kuramoto,JKT90,brown1991minimal,jolly1991preserving}, (iii) optimal prediction of resolved variables \cite{stinis2004stochastic}, and (iv) reproduction of long-term statistics such as the energy spectrum. We follow clearly this latter path, to which we add the question of reproduction by closure of patterns and their statistical features.  For the KSE, only few works have addressed the closure in the latter sense. We refer to \cite{lu2017data} for closure aimed at reproducing long-term statistics and to \cite{stinis2004stochastic} for optimal prediction. In all these works, the regimes for which an efficient closure is sought correspond either  to specific solutions or to weakly turbulent regimes associated with a few pairs of unstable modes:  2 pairs in \cite{armbruster1989kuramoto}, up to 4 pairs of unstable modes for \cite{JKT90,brown1991minimal,jolly1991preserving}, and 3 pairs in \cite{stinis2004stochastic,lu2017data}. 

In this study, we aim at determining efficient closures for the reproduction of patterns and long-term statistics in two strongly turbulent regimes: one regime corresponding to 31 pairs (Regime A, Table \ref{Table_RegimeA}) of unstable modes and another one corresponding to 90 pairs of unstable modes (Regime B, Table \ref{Table_RegimeB}).  Our approach relies on optimal PMs that allow for approximating the conditional expectation (Theorem \ref{Thm_variational-pb2}) without assuming separation of scales and differ in that sense from averaging techniques and other RG methods.

 The reproduction of the energy spectrum of KS solutions will be one of the core metrics to assess the quality of our parameterizations.
For either formulation \eqref{eq:KSE} or \eqref{eq:KSE3}, a typical energy spectrum, $E(k)$, of a chaotic KS solution is shown as the black curve in panel (e) of Fig.~\ref{Quantized_PM}.
Four parts of this spectrum are distinguishable \cite{wittenberg1999scale}: (i)  The {\it large scale region} as $k\rightarrow 0$ which is characterized by a plateau reminiscent of a thermodynamic regime with equipartition of energy; (ii) the {\it active scale region} that contains most of the energy, with a peak corresponding to a characteristic length $l_p=L/(2\pi k_p)$ with $k_p$ that corresponds to the wavenumber of the most linearly unstable mode; (iii) a power law decay with an exponent experimentally indistinguishable from $4$ within this active region; and (iv) an exponential tail due to the strong dissipation at small scales.
 It is tempting to think of the region $E(k) \sim k^{-4}$, where production and dissipation are almost balanced ($Dk^{2}\approx \nu k^4$), as an ``inertial range.'' This latter aspect has been already discussed in the literature; see \cite{pomeau1984intrinsic}.

From a mathematical perspective, the KSE is a well-known example of PDE that possesses an inertial manifold, in the invariant space of odd functions \cite{Foias_al88KS,CFNT89}, and in the general periodic case \cite{temam1994estimates,robinson1994inertial},
but the current IM theory \cite{zelik2014inertial} predicts that the underlying slaving of the high modes to the low modes, holds when the cutoff wavenumber, $k_\c$,  is taken sufficiently far within the dissipative range, especially in ``strongly'' turbulent regimes that correspond to the presence of many unstable modes; see {\hl the Supplementary Material}. Still, as the AIM theory underlines, satisfactory closure may be expected to be derived for $k_\c$ corresponding to scales larger than what predicts the IM theory. Nevertheless, as one seeks to further decrease $k_\c$   within the inertial range, standard AIMs fail typically in providing relevant closures  and one needs to rely on no longer a fixed  cutoff but instead a dynamic one so as to avoid energy accumulation on the cutoff level \cite{debussche1995nonlinear,dubois1998incremental,dubois1998dynamic}. This situation has been already documented for the Navier-Stokes equations \cite{pascal1992nonlinear}, but is less known for the KSE.

As pointed out below, such a failure by traditional (nonlinear) parameterizations for closing the KSE when $k_\c$ is  placed low within the inertial range occurs e.g.~for Regime A considered hereafter and whose parameters\footnote{These parameters become $\alpha=4000$, $\overline{\delta t}=10^{-7}$ and $N_x=256$ when scaling \eqref{Eq_conversion_I_to_III} is applied; see Remark \ref{Rmk_rescaling}.} are listed in Table \ref{Table_RegimeA}.  For this regime, the KS flow is strongly turbulent (see Fig.~\ref{Quantized_PM}-(b)) and possesses $31$ pairs of unstable modes. We selected $k_\c$ to be the wavenumber corresponding to the smallest scale present among the unstable modes, corresponding here to $k_\c=31$ for Regime A, and making thus the reduced state space, $E_\c$, to be spanned by the unstable modes. This choice of $k_\c$ places the cutoff wavelength within the aforementioned inertial range, as one can observe in Fig.~\ref{Quantized_PM}-(d). The fraction of energy to parameterize is quite substantial for this cutoff as it represents $15.7\%$ of the total energy. For this selection of $k_\c$, the energy distribution nearby this cutoff scale is comparable to the energy $E(k)$ contained in the large scales ($k\sim 1$). Beyond $k_\c$, the energy does not drop suddenly (due to its decay following a power law) and actually takes values on a same order of magnitude compared to $E(1)$ for roughly $k_\c<k<1.5 k_\c$ while only after $k>k_1=2 k_\c$, the energy $E(k)$ drops  faster (exponentially);  see black curve Fig.~\ref{Quantized_PM}-(e).

Thus to  close the KSE at this cutoff scale, makes, a priori, the closure problem difficult because quite a few  energetic modes need to be properly parameterized. Actually, as discussed in Sec.~\ref{Sec_backscatter_fixed} below, this difficulty is manifested when using nonlinear parameterizations such as the standard QSA \eqref{Eq_QSA} that suffers from a backscattering transfer of energy particularly overwhelming for the large scales. In this case an over-parameterization of the neglected scales (i.e.~an excessive parameterization of the unresolved energy) leads to an incorrect reproduction of the backscatter transfer of energy due to nonlinear interactions between the modes, especially those near the cutoff scale. We speak of an inverse error cascade,  i.e.~errors in the modeling of the parameterized scales that contaminate gradually the larger scales and spoil the closure skills for the resolved variables.

 To illustrate such an inverse error cascade in a simple context, {\mkr we invite the reader to consult the AB-system in the Supplementary Material; see Eq.~(17) therein}. {\mkr For this system, let us} assume that an error of size $\epsilon \overline{B}$ is made on the parameterized variable $\overline{B}$ at the steady state $(\overline{A}, \overline{B})$ given by {\hl (18) in the Supplementary Material}. This error propagates {\mkr then} to the resolved variable $\overline{A}$ through nonlinear coupling
as $\overline{A}_{\mathrm{app}} = \sqrt{(\nu_2  \overline{B}_{\mathrm{app}} - \alpha \overline{B}_{\mathrm{app}}^3)/\gamma_2}$  where $\overline{B}_{\mathrm{app}} = (1\pm \epsilon) \overline{B}$.
The ($L^2$) error on the resolved variable becomes then $|\overline{A}^2 - \overline{A}^2_{\mathrm{app}}|$: of order $\epsilon$ when $\epsilon$ is small, and of order $\epsilon^3$ when $\epsilon$ is large.  This simple example shows that an error made on the parameterization may be amplified through the nonlinear interactions as it propagates to the resolved variables when the parameterization is not accurate. Such an inverse error cascade is even more pronounced as the number of nonlinear interaction terms gets large while the neglected scales contain a non-negligible amount of energy. In that respect, the parameter regimes considered here for the KSE are particularly demanding  to avoid an incorrect reproduction of the backscatter transfer of energy to the large scale.

Our purpose is to show that the parametric QSA  formulas 
\eqref{Eq_Psin}-\eqref{QSA_tau} of Sec.~\ref{Sec_FMTtau}, when optimized by solving the minimization problems \eqref{Min_formulation_h1_quadcase_b}, allow for 
fixing the backscatter transfer of energy issue encountered by the standard QSA \eqref{Eq_QSA}.  As shown hereafter, the amount of data required to determine the underlying optimal PMs (here given as optimal QSAs), is related to mixing properties such as encoded into decay of temporal correlations. Typically, the faster the decay of (temporal) correlations is, the less the amount of data (in the time direction) required, is.  The PM approach and its apparatus provides furthermore new understanding about essential variables and their interactions for closure of the KSE.

To apply the PM approach and the parameterization formulas of Sec.~\ref{Sec_FMTtau} to Eq.~\eqref{eq:KSE}  we first recall  
the spectral elements of the  operator $A=- \nu \partial_{x}^4  - D \partial_x^2 $, under periodic boundary conditions. These are given by
\be\label{Eig_KSE}
\beta_k= -\frac{16 \nu \pi^4 k^4}{L^4} + \frac{4D\pi^2 k^2}{L^2},
\ee
for the eigenvalues, and 
\be\label{Modes_KSE}
\boldsymbol{e}_k^{\ell}(x)=\begin{cases}
\sqrt{\frac{2}{L}}\cos\bigg(\frac{2\pi k x}{L}\bigg), \quad \mbox{ if } \ell =0\\
\sqrt{\frac{2}{L}}\sin\bigg(\frac{2\pi k x}{L}\bigg), \quad \mbox{ if } \ell =1,
\end{cases}
\ee
for the eigenmodes. Note that because the spatial average of our KS-solutions considered hereafter is zero (see \eqref{Eq_meanzero}), we consider  $k\geq 1$ in what follows.

Adopting the convention of Sec.~\ref{Sec_loc_invman}, and after having reordered the $\beta_k$'s in descending order,  the reduced state space is
\be
E_\c=\mbox{span}\{\boldsymbol{e}_{p(1)}^{\ell},\cdots,\boldsymbol{e}_{p(m)}^{\ell},\; \ell=0,1\},
\ee
where $p(j)$ denotes the wavenumber of the cosine/sine pair associated with the $j^{\mathrm{th}}$ largest eigenvalue.
Note that due to the distribution of the $\beta_k$'s given by \eqref{Eig_KSE}, this reordering matters only when $m<m_u$ with $m_u$ denoting the total number of pairs of unstable modes.

The projector $\Pi_\c$ onto $E_\c$ is then given by
\be\label{Pi_c_formula}
\Pi_\c u= \sum_{\ell=0}^1\sum_{j=1}^m \langle u, \boldsymbol{e}^\ell_{p(j)}\rangle \boldsymbol{e}^\ell_{p(j)}.
\ee 
Hereafter we will consider closure for $m\geq m_u$. In this case, the reduced state space is simply given by
\be\label{Eq_Hc_used}
E_\c=\mbox{span}\{\boldsymbol{e}_{1}^{\ell},\cdots,\boldsymbol{e}_{m}^{\ell},\; \ell=0,1\}.
\ee
Here the ambient space is taken to be the Hilbert space $\cH=L^2(0,L)$, and $\langle \cdot , \cdot\rangle$  denotes its natural inner product.  Hereafter we denote by 
 $\Pi_\s$ the orthogonal complement of $\Pi_\c$ in $\cH$, i.e.~$\Pi_\s=\mbox{Id}_{\cH}-\Pi_\c$. 

%%%%%%%%%%%%%%%%%%%%%%%%%%
\begin{table}[h] 
\caption{Regime A: Parameters for Eq.~\eqref{eq:KSE}}
\label{Table_RegimeA}
\centering
\begin{tabular}{cccccc}
\toprule\noalign{\smallskip}
     $\nu$ &  $D$ & $L$  &$\gamma$ &  $\delta t$ &  $N_x$    \\ 
\noalign{\smallskip}\hline\noalign{\smallskip}
 $2 \times 10^{-4}$ &  $0.2$ &  $2\pi$ & 1 &  $10^{-3}$ & 256\\
\noalign{\smallskip} \bottomrule 
\end{tabular}
\end{table}
%%%%%%%%%%%%%%%%%%%%%%%%%% 

%%%%%%%%%%%%%%%%%%%%%%%%%%
\begin{table}[h] 
\caption{Regime B: Parameters for Eq.~\eqref{eq:KSE3}}
\label{Table_RegimeB}
\centering
\begin{tabular}{ccc}
\toprule\noalign{\smallskip}
     $\alpha$  &  $\delta t$ &  $N_x$    \\ 
\noalign{\smallskip}\hline\noalign{\smallskip}
 $33000$  &  $10^{-9}$ & 2048\\
\noalign{\smallskip} \bottomrule 
\end{tabular}
\end{table}
%%%%%%%%%%%%%%%%%%%%%%%%%% 

Another regime that will be dealt with in Sec.~\ref{Sec_More_Results} below has its parameters listed in Table \ref{Table_RegimeB} for the KSE 
written under its formulation  \eqref{eq:KSE3}. This regime is even more turbulent than Regime A, as it exhibits 90 pairs of unstable modes.  
Either for Regime A or B,  the benchmark KS solution for the closure exercises conducted hereafter,  is obtained by transforming the KSE in Fourier space and by using 
 a modification of the exponential time-differencing fourth-order Runge-Kutta (ETDRK4) method proposed in \cite{kassam2005fourth} in order to solve the resulting stiff ODE system. The number of Fourier modes retained ($N_x$) and  time step used  ($\delta t$) for each regime, are listed in Tables \ref{Table_RegimeA} and \ref{Table_RegimeB}, for Regimes A and B, respectively.  We refer hereafter to a KS solution thus obtained as a Direct Numerical Solution (DNS).  The ODE closure derived hereafter are integrated with an semi-implicit Euler scheme, in which the linear terms are treated implicitly while the nonlinear ones, explicitly. These closure systems are integrated with the same time step as listed in    Tables \ref{Table_RegimeA} and \ref{Table_RegimeB}, depending on the regime. 
 
 In all our numerical experiments that follow, the KSE is integrated from the following initial datum with zero-mean
 \be\label{Eq_initKS}
 u_0(x)=\cos(x)(1+\sin(x)).
 \ee 
In such a case, since the spatial average is a conserved quantity for the KS solution $u(x,t)$, we have for all $t$,
\be\label{Eq_meanzero}
\int_0^L u(x,t) \d x =0. 
\ee

Note that compared with the original ETDRK4 proposed in \cite{cox2002exponential}, the modification in \cite{kassam2005fourth} consists of evaluating key coefficients as given by \cite[Eq.~(2.5)]{kassam2005fourth} using contour integrals rather than direct evaluation to avoid possible cancellation errors. The contours are taken to be circles of radius $\delta t$ centered around each of the eigenvalues of the discretized linear operator, and the contour integrals are approximated using trapezoid rules with $M$ equally spaced points on the circle. We have set $M=64$ for both parameter regimes considered. In our numerical calculations performed in Matlab (version R2018a), compared to the script given in \cite[Fig.~7]{kassam2005fourth}, the spatial discretization is taken to be \texttt{x = L* (0:Nx-1)'/Nx} instead of \texttt{x = L*(1:Nx)'/Nx} to suit the way the fast Fourier transform (FFT) is implemented in the Matlab built-in function \texttt{fft}.

%%%%%%%%%%%%%%%%%%%%%%%
\br\label{Rmk_rescaling}
When the  scaling \eqref{Eq_conversion_I_to_III} is performed, we find for Regime A that $\alpha=4000$ and $\overline{t}=\theta t$ with $\theta=5\times 10^{-5}$. After transient is removed, to reach the same energy level, $\|u\|_{L^2}^2$ than by integrating \eqref{eq:KSE} (with the same solver), we have found that we can decrease the time-step compared to $\delta t$ by a factor $a=10^4$, that is $\overline{\delta t}=10^{-7}$. Given an interval of length $T$ in the original time variable $t$, it corresponds to $\overline{T} = 5\times 10^{-5} T$, that is an amount of data in time that is given by $\overline{N} = \overline{T}/\overline{\delta t} = 500 T$ data points. Thus, since $N = T/\delta t=1000 T$,  we have that $\overline{N}=N/2$. Although mathematically equivalent, we can thus store twice more data (while keeping $N_x$ identical) by integrating numerically the formulation   
 \eqref{eq:KSE3} than by integrating  the formulation  \eqref{eq:KSE}, integrating the dynamics up to the same time instant (taking into account the rescaling). 
Such observations have their interest to draw statistics from long time integration. For Regime A it turns out that the simulations performed hereafter were already sufficient to draw robust statistics with the formulation   \eqref{eq:KSE}.  We use however formulation \eqref{eq:KSE3} to simulate the turbulent Regime B with a higher number of unstable modes than for Regime A.

\er

%%%%%%%%%%%%%%%%%%%%%%%%%%%%%%%%%%%%%%%%%%%%%%%%%%%%%
\subsection{Fixing the backscatter transfer of energy for KS turbulence with optimal PMs}\label{Sec_backscatter_fixed}
\begin{figure}
 \centering
  \includegraphics[height=.7\textwidth,width=1\textwidth,]{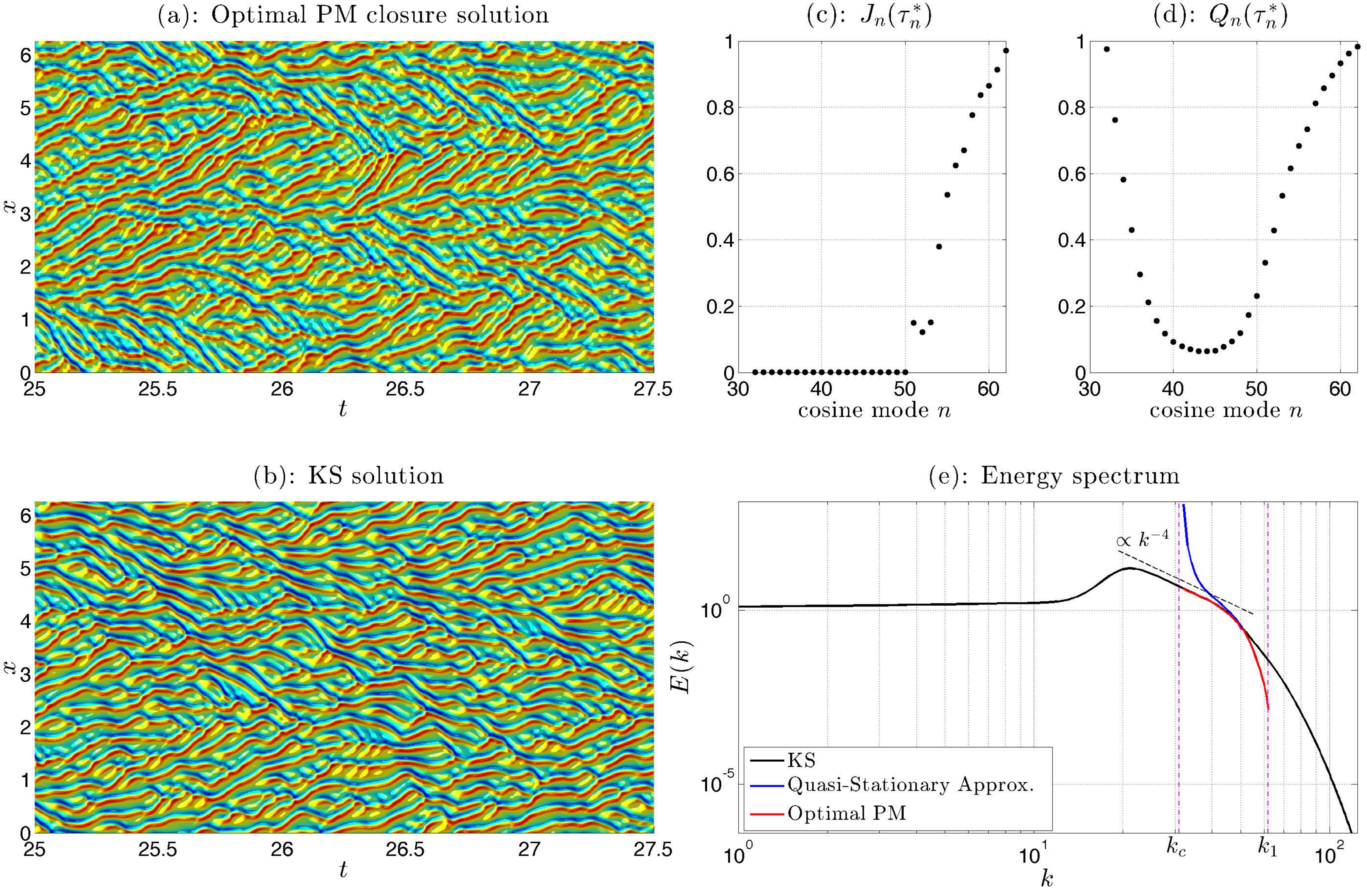}
\caption{{\footnotesize {\bf Closure and parameterization results Regime A}.  Panel (a) shows the solution obtained from the optimal PM closure \eqref{KSE_closure} with $m=31$, while panel (b) shows the KS solution as obtained from DNS of Eq.~\eqref{eq:KSE}.  
Here the optimal PM  is obtained as QSA$(\bftau^\ast)$ with $\bftau^\ast$ obtained by optimization of the cost functional $J_n$ given by \eqref{J_n_normalized} (with $t=1$ and $T=4$). The optimal values $J_n(\tau_n^\ast)$ are shown in panels (c) for the 
parameterized cosine modes. The corresponding $Q_n$-values are shown in panel (d), with $Q_n$ given by \eqref{toto_Qnb}.
The resulting optimal QSA parameterizes the wavelength band, $k_\c<k<k_1=2 k_\c$, as shown by the red curve in panel (e) on the energy spectrum $E(k)$ (log-log scale). Here $k_\c$ is the wavenumber corresponding to the smallest scale present among the unstable modes, that is $k_\c=31$.   
The blue curve shows the dramatic failure of the standard quasi-stationary approximation (QSA) \eqref{Eq_QSA} for parameterizing this  wavelength band, especially for $k$ near $k_\c$.}}
 \label{Quantized_PM}
 \vspace{-.1cm}
\end{figure}

It is known that when the cutoff wavelength is too low within the inertial range, the standard QSA  \eqref{Eq_QSA}  suffers typically from over-parameterization leading to an incorrect backscatter transfer of energy, i.e.~errors in the modeling of the parameterized (small) scales that contaminate gradually the larger scales. In the case of Regime A, when $k_\c=31$ (corresponding to $E_\c$ spanned by 31 pairs of unstable modes), the QSA leads to an over parameterization of $E(k)$ by an amount of about $5800 \%$ (in average) over the wavenumbers  $32\leq k\leq 36$; see blue curve in Fig.~\ref{Quantized_PM}-(e).  The nonlinear interactions between these modes and the unstable modes corresponding to $k\leq k_\c$ lead in this case to such an excessive backscatter transfer of energy, that a closure in which the unresolved modes are approximated by the QSA, blows up after few iterations no matter the numerical scheme used. 

 %%%%%%%%%%%%%%%%%%%%%%%%%%%%%%%%%%%%%%%%%%

As  pointed out in Sec.~\ref{Sec_FMTtau}, the parametric QSA formulas \eqref{Eq_Psin}-\eqref{QSA_tau} involve the same interaction coefficients, the $B_{ij}^n$'s given by \eqref{B_ijn} as for the standard QSA, $K(\xi)$.  However the magnitudes of the nonlinear interactions, as encapsulated in the coefficients $\delta_n(\tau)$'s given by \eqref{Eq_D_term00}, is different from the coefficients $-\beta_n^{-1}$ appearing in $K(\xi)$.  The coefficients $\delta_n(\tau)$'s enable us here to counterbalance the excess of energy in the parameterization compared to a standard QSA.  Furthermore, as explained below, these coefficients are optimized in the $\tau$-variable by solving the minimization problems \eqref{Min_formulation_h1_quadcase_b} over short training periods of length comparable to a characteristic decorrelation time of the dynamics.

In the case of the KSE, the parametric QSA  \eqref{QSA_tau}, QSA$(\bftau)$, takes  
the following form
 \be\label{Eq_optPM}
 \Psi_{\bftau}(\xi)=\sum_{\ell=0}^1 \sum_{n = m+1}^{2m} \Psi_n^{\ell} (\tau_n, \boldsymbol{\beta}, \xi) \boldsymbol{e}_n^{\ell},
 \ee
 with 
\be \label{Eq_Phi_tau2}
\Psi_n^{\ell}(\tau_n^\ell, \boldsymbol{\beta}, \xi) = \sum_{i, j = 1}^m \delta_n(\tau_n^\ell) \Big(E_{ij}^{n,\ell} \xi^{0}_{i} \xi^{0}_{j} +C_{ij}^{n,\ell} \xi^{0}_{i} \xi^1_{j}+F_{ij}^{n,\ell}  \xi^{1}_{i} \xi^{1}_{j}\Big), \qquad \xi \in E_\c.
\ee
The index $m$ in the upper bound of the sum is taken here to be equal to $k_\c=31$, which corresponds to the number of pairs of unstable modes for Regime A. The reduced state space $E_\c$ is thus $2m$-dimensional, taking into account $\ell=0,1$.

In  \eqref{Eq_Phi_tau2}, $\delta_n(\tau_n^\ell)$ is given by \eqref{Eq_D_term00} while
\be\label{E_interaction}
E_{ij}^{n,\ell}=\begin{cases}
\langle B(\boldsymbol{e}_i^{0},\boldsymbol{e}_j^{0}),\boldsymbol{e}_n^{0}\rangle, \quad \mbox{ if } \; \ell=0\\
\langle B(\boldsymbol{e}_i^{0},\boldsymbol{e}_j^{0}),\boldsymbol{e}_n^{1}\rangle, \quad \mbox{ if } \;\ell=1,
\end{cases}
\ee
\be\label{C_interaction}
C_{ij}^{n,\ell}=\begin{cases}
\langle B(\boldsymbol{e}_i^{0},\boldsymbol{e}_j^{1}),\boldsymbol{e}_n^{0}\rangle + \langle B(\boldsymbol{e}_j^{1},\boldsymbol{e}_i^{0}),\boldsymbol{e}_n^{0}\rangle \quad \mbox{ if } \; \ell=0\\
\langle B(\boldsymbol{e}_i^{0},\boldsymbol{e}_j^{1}),\boldsymbol{e}_n^{1}\rangle + \langle B(\boldsymbol{e}_j^{1},\boldsymbol{e}_i^0),\boldsymbol{e}_n^{1}\rangle \quad \mbox{ if }\; \ell=1,
\end{cases}
\ee
and
\be\label{F_interaction}
F_{ij}^{n,\ell}=\begin{cases}
\langle B(\boldsymbol{e}_i^1,\boldsymbol{e}_j^1),\boldsymbol{e}_n^0\rangle\quad \mbox{ if }\; \ell=0 \\
\langle B(\boldsymbol{e}_i^1,\boldsymbol{e}_j^1),\boldsymbol{e}_n^1\rangle \quad  \mbox{ if }\; \ell=1.
\end{cases}
\ee
%%%%%%%%%%%%%%%%%%%%%%%%%%
These  coefficients correspond to  the aforementioned interaction coefficients. They possess a simple analytic expression here given
the nonlinearity and the trigonometric  eigenfunctions. In particular,  a majority of these coefficients are actually zero for $m+1 \leq n\leq 2m$, leaving only a few of them non-zero. 

%%%%%%%%%%%%%%%%%%%%%%%%%%%%%%%%%%%%%%%%
More precisely, we have
\be\label{Sparse1}
\langle B(\boldsymbol{e}^0_i,\boldsymbol{e}^0_j),\boldsymbol{e}^0_n\rangle = 
\langle B(\boldsymbol{e}^0_i,\boldsymbol{e}^1_j),\boldsymbol{e}^1_n\rangle =
\langle B(\boldsymbol{e}^1_i,\boldsymbol{e}^0_j),\boldsymbol{e}^1_n\rangle =
\langle B(\boldsymbol{e}^1_i,\boldsymbol{e}^1_j),\boldsymbol{e}^0_n\rangle = 0, \quad \forall \; i,j,n,
\ee
\be\label{Sparse2}
\langle B(\boldsymbol{e}^0_i,\boldsymbol{e}^1_j),\boldsymbol{e}^0_n\rangle 
= \langle B(\boldsymbol{e}^1_j,\boldsymbol{e}^0_i),\boldsymbol{e}^0_n\rangle
= \begin{cases}
-\frac{ \gamma \pi n}{\sqrt{2} L^{3/2}}, &\text{ if $n = i+j$},  \vspace{0.2em}\\ 
\frac{ \gamma \pi (i-j)}{\sqrt{2} L^{3/2}}, &\text{ if $n = |i-j|$},  \vspace{0.2em}\\ 
0, & \text{otherwise},
\end{cases}
\ee
and 
\be\label{Sparse3}
\langle B(\boldsymbol{e}^\ell_i,\boldsymbol{e}^\ell_j),\boldsymbol{e}^1_n\rangle  = \begin{cases}
(-1)^\ell \frac{ \gamma \pi n}{\sqrt{2} L^{3/2}}, &\text{ if $n = i+j$, \; $\ell \in \{0,1\}$}, \vspace{0.2em}\\ 
\frac{ \gamma \pi n}{\sqrt{2} L^{3/2}}, &\text{ if $n = |i-j|$, \; $\ell \in \{0,1\}$}, \vspace{0.2em} \\ 
0, & \text{otherwise}. 
\end{cases}
\ee
 Note that formulas \eqref{Sparse1}-\eqref{Sparse3} show that the parameterization $\Psi_n^\ell$ in \eqref{Eq_Phi_tau2} is sparse, for $m+1 \leq n\leq 2m$ and identically zero for $n\geq 2m +1$.

The optimal QSA, $\Psi_{\bftau^\ast}$, is obtained by solving the minimization problems \eqref{Min_formulation_h1_quadcase_b}.
The corresponding normalized parameterization defect, 
\be\label{J_n_normalized}
J_n(t,\tau)=\frac{\bigg|\int_{t}^{t+T} [\Pi_n u(s)]^2 \d s - \int_t^{t+T}[\Psi_n(\tau,\bfbeta, u_\c(s))]^2 \d s \bigg|}{\int_t^{t+T} |\Pi_n u(s)|^2 \d s},
\ee
is shown in panel (c) of Fig.~\ref{Quantized_PM} for the $\tau=\tau_{n}^\ast$'s that correspond to the optimal values for the cosine modes, dropping here the dependence on $\ell=0$. The results for the sine modes are almost identical, and are thus not shown.
Here $t$ is chosen after the transient behavior, as measured through the  
 energy, $\|u(t)\|_{L^2}$ of the DNS for Regime A.  In our case, it corresponds to $t=1$.  The training length $T$ is chosen to be $T=4$.

Note that unlike the case dealt with in Sec.~\ref{Sect_period-doubling-RB}, the cost functional  
$J_n$ does not exhibit local minima  (in contrast with Remark \ref{Remark_training_interval}) and thus the dependence on $t$ is secondary as far as one is concerned with optimal values: $J_n(t,\tau_n^\ast)$ will be hereafter denoted by $J_n(\tau_n^\ast)$. Instead, $\tau\mapsto J_n(\tau)$ exhibits,  for $n=32$ through $n=50$, sharp gradients near the origin that lead to $\tau_n^\ast$-values close to zero for these modes.  

It is striking to observe that  $J_n(\tau_n^\ast)$ is almost identical to zero
for $n=32$ up to $n=50$ (see Fig.~\ref{Quantized_PM}-(c)), resulting by an almost perfect parameterization of the energy contained into 
the corresponding modes; compare the red curve with the black curve in Fig.~\eqref{Quantized_PM}-(e).  For instance, the corresponding optimal QSA comes with a (average) relative error of only $1.3\%$ over the wavenumbers $32\leq k\leq 36$, allowing in turn to fix the dramatic backscatter transfer of energy issue encountered by the standard QSA and even by standard Galerkin approximations with $m> k_c$; see Remark \ref{Rmk_LIA_Galerkin} below.   

\begin{figure}
 \centering
 \includegraphics[height=.4\textwidth,width=1\textwidth,]{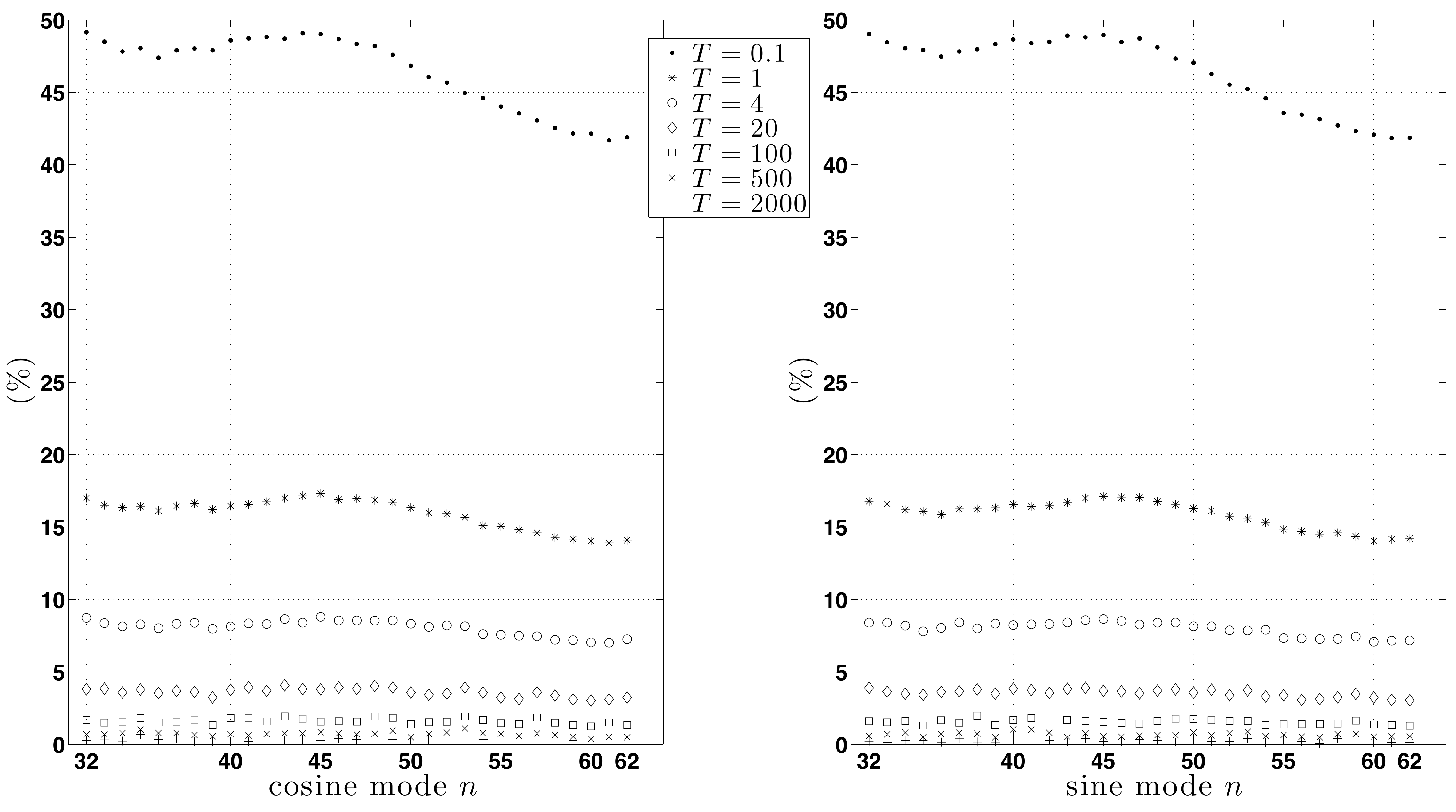}
\caption{{\footnotesize {\bf Relative error of $\frac{1}{T}\int_{t}^{t+T} [\Pi_n u(s)]^2 \d s$ compared to $E(n)$}. Here the energy contained in $E(n)$ is estimated over $4\times 10^6$ snapshots, that is for $T=4000$.}}
 \label{Fig_modewise_error}
\end{figure}

This ability of the optimal QSA to accurately reproduce the amount of energy contained in the consecutive high modes located after the cutoff scale, is even more striking when one notes that QSA$(\bftau)$ is optimized by minimizing $J_n$ on DNS data over a training length $T=4$ (corresponding to $4\times 10^3$ snapshots) whereas the energy spectrum $E(k)$ shown in  Fig.~\ref{Quantized_PM}-(e) is estimated over $T=4000$ ($4\times 10^6$ snapshots). The relative error $r$ of  $\frac{1}{T}\int_{t}^{t+T} [\Pi_n u(s)]^2 \d s$ compared to $E(n)$ is shown as $T$ evolves in Fig.~\ref{Fig_modewise_error} for the cosine and sine modes. 
For $T=4$ the average error is about 8$\%$. Even if $T=1$ (corresponding to $r\approx 16\%$) is selected to evaluate  $J_n$, the resulting optimal QSA performs similarly than that optimized with $T=4$, regarding the reproduction of the amount of energy contained in the high modes (not shown). 

%%%%%%%%%%%%%%%%%%%%%%%%%%%%%%%%%%%%%%%%%%%%%%%
%%%%%%%%%%%%%%%%%%%%%%%%%%%%%%%%%%%%%%%%%%% 

These observations show the usefulness of our variational approach: By optimizing the parameterization QSA$(\bftau)$ according to the cost functional $J_n$, one fixes the backscatter transfer of energy issue encountered by the standard QSA, while relying only on a short integration of the KSE. 
Furthermore, on a practical ground, it is worthwhile noting that one benefits greatly from the dynamically-based formulas QSA$(\bftau)$ (see \eqref{Eq_Psin}-\eqref{QSA_tau}) to operate this optimization. As a comparison, a blind regression using homogeneous polynomials of degree 2 in the $\xi$-variable, would lead in this case to $31\times 15 \times 3=1395$ coefficients\footnote{Obtained by counting the number of (distinct) monomials $\xi_i^{\ell} \xi_j^{\ell'}$, with $i,j \in \{1,\cdots, 31\}$, and $\ell, \ell' \in \{0,1\}$. } to estimate for each high mode and by taking $T=1$ or $T=4$ ($4\times 10^3$ snapshots) the resulting regression problem would be either 
underdetermined or non-robust statistically.  Instead, due to the parametric form of QSA$(\bftau)$,  only $2$ scalar parameters ($\tau_n^\ell$, $\ell=0,1$) need to be determined, for each high mode.

As a complimentary diagnosis metric, we show in Fig.~\ref{Quantized_PM}-(d), for the $\tau_n^\ast$'s obtained by minimizing \eqref{J_n_normalized}, the values of the following parameterization defect,
\be\label{toto_Qnb}
Q_n(\tau_n^\ast)= \frac{\int_{t}^{t+T} \big| \Pi_n u(s)- \Psi_n(\tau_n^\ast, \boldsymbol{\beta}, u_\c(s))\big|^2 \d s}{\int_t^{t+T} |\Pi_n u(s)|^2 \d s},
\ee
also for the cosine modes, and for $t=1$ and $T=4$. Clearly for the modes whose wavenumbers are located right above the cutoff wavelength, $k_\c$, the $Q_n$-values, although less than 1, are not as close to zero  as for the  $J_n$-values shown in Fig.~\ref{Quantized_PM}-(c). 
Remark that since the mean values of the components of our KS-solution are zero, minimizing $Q_n$ consists of minimizing the variance of the residual error, i.e.~$\overline{|u_n-f(\tau,u_\c)|^2}$, for a given parameterization $f(\tau,\cdot)$. By construction, minimizing $J_n$ consists instead of minimizing the residual error of the variance approximation, i.e.~$|\overline{|u_n|^2}-\overline{|f(\tau,u_\c)|^2}|$.

It is noteworthy that the $Q_n$-values in \eqref{toto_Qnb} differ slightly from the optimal ones that would be found by minimizing directly the $Q_n$'s in the $\tau$-variable, over the same training length. Nevertheless, the resulting differences in the corresponding minimizers matters as one would encounter an under-parameterization of about $50\%$ (in average) for the modes near the cutoff wavelength ($32\leq n\leq 36$); see Remark \ref{Rmk_LIA_Galerkin} below.

\begin{figure}
 \centering
  \includegraphics[height=.6\textwidth,width=.95\textwidth,]{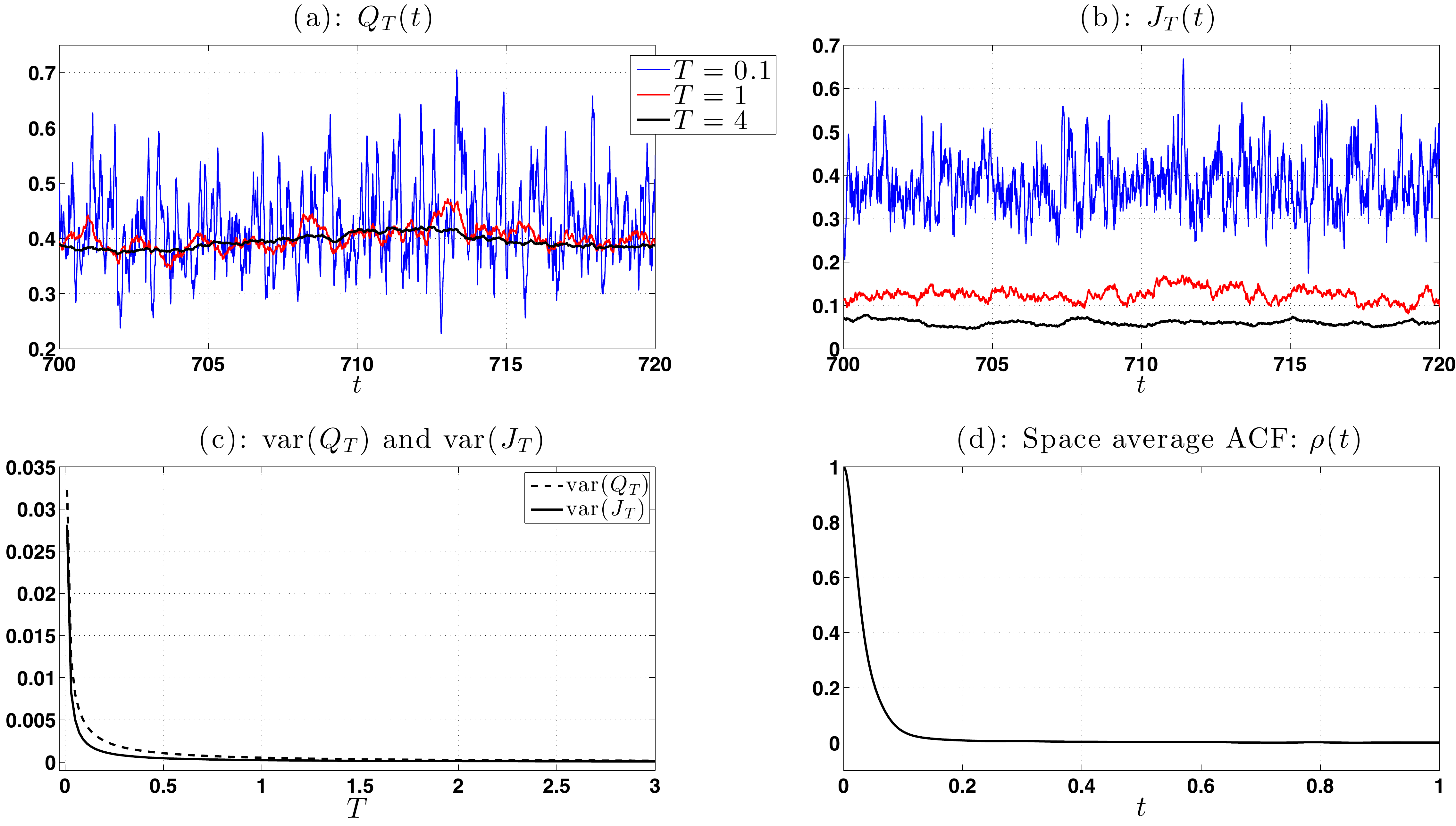}
\caption{{\footnotesize {\bf Effects of the training period, $T$, on the parameterization defects  $J_T(t)$ and $Q_T(t)$}.  Here, we observe that: (i) as $T$ is increasing, $J_T(t)$ and $Q_T(t)$ are converging towards a constant value (Panels (a) and (b)), (ii) the variance of $J_T(t)$ (resp.~$Q_T(t)$), $\mbox{var}(J_T)$ (resp.~$\mbox{var}(Q_T)$), decays to zero (Panel (c)), and (iii) the rate of decay of the latter is comparable to that of the space average ACF, $\rho(t)$, given by  \eqref{Eq_temporal_rho} (Panel (d)).}}
 \label{QT_plots}
\end{figure}

To better understand the effect of the training length $T$ (that determines the amount of data from DNS to be stored), 
we proceeded as follows.  Given a training length $T$, the optimal QSA, $\Psi_{\bftau^\ast}$, is determined by minimizing the corresponding cost functional $J_n$ given by  
\eqref{J_n_normalized} (with $t=1$), providing thus the optimal parameters, $\tau_n^\ast$'s. 
Recalling that the interaction coefficients are zero for $n\geq 2m+1$ (see \eqref{Sparse1}-\eqref{Sparse3}), we analyzed then numerically the dependence on $t$ and $T$ of the following averaged parameterization defect
\be
J_T(t, \Psi_{\bftau^\ast}) = 
\frac{ \sum_{n=m+1}^{2m} \Big| \int_{t}^{t+T} [\Pi_n u(s)]^2 \d s - \int_{t}^{t+T} [\Psi_n(\tau_n^\ast,\bfbeta, u_{\c}(s))]^2 \d s \Big|}{\sum_{n=m+1}^{2m}  \int_{t}^{t+T}  [\Pi_n u(s)]^2 \d s},
\ee
as well as of the parameterization defect $Q_T(t,\Psi_{\bftau^\ast})$ given by \eqref{Eq_PD}.  
To simplify the notations, we denote hereafter $J_T(t, \Psi_{\bftau^\ast})$ and $Q_T(t,\Psi_{\bftau^\ast})$ by $J_T(t)$ and $Q_T(t)$, respectively. Panels (a) and (b) of Fig.~\ref{QT_plots} show the dependence on $t$ of $J_T(t)$ and $Q_T(t)$, respectively.
This dependence is shown here for three values of $T$: $T=0.1$, $T=1$, and $T=4$. In each case, $Q_T(t)<1$ showing that $\Psi_{\bftau^\ast}$ is a  PM, even for the short training length $T=0.1$. Either for $Q_T(t)$ or $J_T(t)$ we observe that the amplitude of the oscillations in time is reduced as $T$ is increased. This is further confirmed by inspecting the variance of $Q_T$ and $J_T$ as $T$ is varied: both exhibit a fast convergence towards zero as $T$ grows; see panel (c) of Fig.~\ref{QT_plots}. 

The decay towards zero of these variances can be put into perspective with the
following space average temporal ACF,
\be\label{Eq_temporal_rho}
\rho(t)=\frac{1}{2\pi T} \int_{0}^{2\pi} \int_0^T u(x,s) u(x,t+s) \d s \d x.
\ee
The latter quantity informs us on how the spatio-temporal field, $u(x,t)$, decorrelates in time, after averaging over $x$. This space average ACF is shown in panel (d) of Fig.~\ref{QT_plots}. It exhibits decay of correlations on timescales comparable to those for the variances of $Q_T$ and $J_T$ supporting thus an earlier statement that the coefficients $\delta_n(\tau)$'s in \eqref{Eq_D_term00} are optimized in the $\tau$-variable by solving the minimization problems \eqref{Min_formulation_h1_quadcase_b} over short training periods of length comparable to a characteristic decorrelation time of the dynamics. For our closure results presented hereafter we selected $T=4$.

%%%%%%%%%%%%%%%%%%%%%%%%%%%%%%%%%%%%%%%%%%
%%%%%%%%%%%%%%%%%%%%%%%%%%%%%%%%%%%%%%%%%%
Thus, after minimization in the $\tau$-variable of the cost functionals, $J_n$'s, given by  \eqref{J_n_normalized}, 
(with $T=4$ and after removal of transient, $t=1$), we use the resulting optimal (and sparse) PM, QSA($\bftau^\ast$) (i.e.~$\Psi_{\boldsymbol{\tau}^\ast}$), with 
\bes
\boldsymbol{\tau}^\ast=\{\tau_{n,\ell}^{\ast}, \, : \,  m+1 \leq n\leq 2m, \; \ell=0,1\},
\ees
to construct the following optimal PM closure
\be\label{KSE_closure}
\frac{\d z_j^\ell} {\d t} = \beta_j z_j^\ell + \Big\langle B(z+\Psi_{\boldsymbol{\tau}^\ast}(z),z+\Psi_{\boldsymbol{\tau}^\ast}(z)), \boldsymbol{e}_k^{\ell}\Big\rangle, \quad 1\leq j \leq m, \; \; \ell \in \{0,1\},
\ee
where $z(x,t)=\sum_{\ell=0}^1\sum_{j=1}^m z_j^{\ell}(t) \boldsymbol{e}_j^{\ell}(x)$, for $m=31$, that, we recall,  corresponds to the number of pairs of unstable modes. 

Good closure skills are already visible with naked eyes, by simply comparing the solution patterns, $u(x,t)$, obtained by a full integration of Eq.~\eqref{eq:KSE} over $N_x$ modes (i.e.~$u$ obtained by DNS), with the patterns exhibited by the optimal PM closure solution, 
\be\label{PM_approx}
v(x,t)=z(x,t) + \Psi_{\boldsymbol{\tau}^\ast}(z(x,t)), 
\ee
obtained by resolving only $m=31$ pairs of reduced variables (i.e.~by solving system \eqref{KSE_closure}); compare panels (a) and (b) of  Fig.~\ref{Quantized_PM}. 

\begin{figure}
 \centering
  \includegraphics[height=.4\textwidth,width=1\textwidth,]{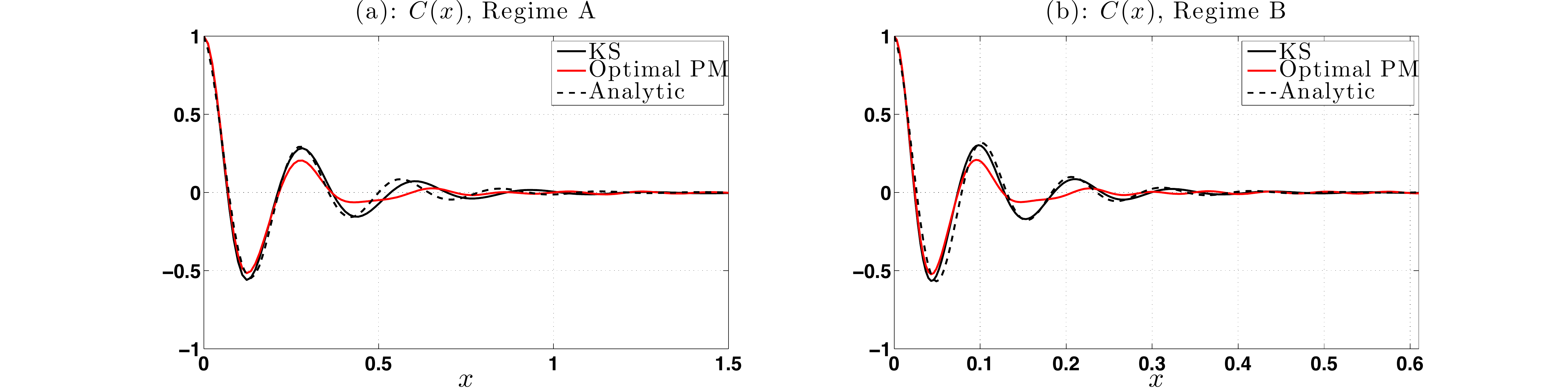}
\caption{{\footnotesize {\bf Time average spatial ACF, $C(x)$, for Regimes A and B}. In both cases, the spatial ACF, $C(x)$, is estimated  from \eqref{Eq_spatial_Cx} based on long simulations of the KSE and the optimal PM closure \eqref{KSE_closure}, with $\bftau^\ast$ minimizing the $J_n$'s given by \eqref{J_n_normalized}.  The simulation lengths correspond here, respectively, to $N=4\times 10^6$ snapshots for Regime A, and to $N=2\times 10^6$ snapshots for Regime B. These estimated ACFs are compared with the analytic formula for $C(x)$ proposed in \eqref{analytic_formula}.}}
 \label{Fig_spatial_Cx}
\end{figure}

To further assess the ability to reproduce the spatio-temporal dynamics by the optimal PM closure \eqref{KSE_closure}, we estimated  the following time average spatial ACF 
\be\label{Eq_spatial_Cx}
C(x)=\frac{1}{L T_f} \int_{0}^{T_f} \int_0^{L} u(x',t) u(x+x',t) \d x' \d t,
\ee
for $u$ as obtained from DNS and its approximation $v(x,t)$ given by \eqref{PM_approx}, both integrated up to $T_f=4000$, while we recall that the training length is $T=4$ to determine $\Psi_{\bftau^\ast}$. The results are shown in panel (a) of Fig.~\ref{Fig_spatial_Cx}. The correlation function $C(x)$ captures both the underlying oscillatory, cellular spatial structure of the KS dynamics, and the rapid spatial decorrelation reflecting the spatial disorder in the spatio-temporal chaotic regime analyzed here. These features are thus well captured by the optimal PM closure \eqref{KSE_closure}. 

Following \cite{wittenberg1999scale}, we observed that the time average spatial ACF is well modeled for the DNS by the following analytic formula, 
\be\label{analytic_formula}
C(x)\approx \cos(k_p^{-1} x ) \exp(-x/\lambda), 
\ee
with $k_p$ that corresponds to the wavelength associated with the peak in the energy spectrum $E(k)$ shown in Fig.~\ref{Quantized_PM}-(e), and $\lambda$ to a correlation length for which  
spatial coupling becomes negligible beyond a few multiples of $\lambda$. For Regime A, we found $k_p=21$ and $\lambda=0.23$.  Only for large lags in the $x$-variable, the optimal PM fails to reproduce accurately this theoretical prediction.

%%%%%%%%%%%%%%%%%%%%%%%%%%%%%%%%%%%%%%%%%%
%%%%%%%%%%%%%%%%%%%%%%%%%%%%%%%%%%%%%%%%%%

\br\label{Rmk_implicitEuler_FMT}
The QSA  \eqref{Eq_QSA} may also be obtained as the limit of the parameterization
\be\label{Eq_Ktau}
K_{\tau}(\xi)= -\tau (\mathrm{Id} + \tau A \Pi_{\s})^{-1} \Pi_{\s} B(\xi,\xi), 
\ee
obtained by using an implicit Euler method to approximate the high modes and by simplifying the nonlinear terms; see \cite{foias1988computation} and \cite[Sec.~7.1]{foias1989exponential}. In this case we have, 
\be\label{Eq_implicitEuler}
\lim_{\tau \rightarrow \infty} -\tau (\mathrm{Id} + \tau A \Pi_{\s})^{-1} \Pi_{\s} B(\xi,\xi)=- A_{\s}^{-1} \Pi_{\s} B(\xi,\xi). 
\ee
Note that in \eqref{Eq_Ktau} unlike in \cite{foias1988computation}, we consider the operator $A$ to be the full linear operator and 
not only given by the 4th-order term.  In its standard formulation,  the parameterization $K_\tau$ is not optimized and  $\tau$ is chosen to be $\lambda_{m+1}^{-1}$, where $\lambda_m=16 \nu \pi^4 m^4/L^4$ denotes the eigenvalue of $\nu \partial_x^4$.

Taking $A=\nu \partial_x^4 +D \partial_x^2$, the  analytic expression of the parameterization $K_\tau$ is the same as for QSA$(\bftau)$\eqref{Eq_Psin}, except that $\delta_n(\tau)$ therein is replaced by $\tau(1-\beta_n \tau)^{-1}$. Since $0 \le \tau(1-\beta_n \tau)^{-1} < -\beta_n^{-1}$, the range of this coefficient is the same as that of $\delta_n(\tau)$ (see discussion at the end of Sec.~\ref{Sec_FMTtau}), and the parameterization $K_\tau$ once optimized by minimizing the cost functional $J_n$  leads also to similar closure skills than those obtained by the optimal QSA.\footnote{Note that by taking $A$ to be given by $\nu \partial_{x}^4$ the resulting coefficients are bounded by $\lambda_n^{-1}$,  and since $\lambda_n^{-1}<-\beta_n^{-1}$ the optimized $K_\tau$ is not a priori of comparable parameterization defects, and in fact leads to less efficient closures.}
We see thus here that the PM approach is not limited to the QSA-class nor the LIA-class introduced respectively in Secns.~\ref{Sec_FMTtau} and \ref{Sect_PM_with_forcing}, but applies actually to any parametric family of nonlinear parameterizations.     
\er
%\ref{Rmk_implicitEuler_FMT} 

\br\label{Rmk_LIA_Galerkin}
We report briefly here on the closure skills obtained when QSA($\bftau$) is optimized by minimizing the $Q_n$'s instead of the 
$J_n$'s.  The metrics used to assess these skills are $\overline{\|u\|}_{L^2}$ (after transient removal) and its standard variation, $\texttt{std}(\|u\|_{L^2})$. The time averages are here estimated on an interval of length $T=100$ ($10^5$ snapshots).  
We observe from Table \ref{Table_energy_PM_regimeA} that the relative error of approximation for $\overline{\|u\|}_{L^2}$ is increased while that for $\texttt{std}(\|u\|_{L^2})$ is reduced, when the 62D closure  \eqref{KSE_closure} ($m=31$) is driven by the optimal QSA($\bftau^\ast$) with $\bftau^\ast$ minimizing the $Q_n$'s. Comparison with standard Galerkin approximations, show that only starting from a 118D Galerkin approximations (m=59), one starts to improve, compared to the 62D closure\footnote{Driven by the optimal QSA($\bftau^\ast$) with $\bftau^\ast$ minimizing the $J_n$'s.},  the approximation of the   
mean value of $\|u(t)\|_{L^2}$ (and comparable skills for $\texttt{std}(\|u\|_{L^2})$) although a good reproduction of the KS patterns' qualitative features, is observed for lower dimension. However this latter aspect seems to be germane to the KSE. In general, indeed,   an error in the reproduction of the right amount of energy come with failures in the reproduction of qualitative features as well, due to an incorrect reproduction of the backscatter transfer of energy. For instance, regarding the wind-driven circulation of the oceans \cite{GCS08}, the jet extension and variability \cite{dijkstra2005low} are notoriously difficult to get parameterized due to eddy backscatter \cite{berloff2005dynamically,berloff2005random}.   

%%%%%%%%%%%%%%%%%%%%%%%%%%
\begin{table}[h] 
\caption{1st and 2nd moments of $\|u\|_{L^2}$: Relative error for Regime A}
\label{Table_energy_PM_regimeA}
\centering
\begin{tabular}{cccc}
\toprule\noalign{\smallskip}
& Energy contained in $E_\s$ &  $\overline{\|u\|}_{L^2}$ &  $\texttt{std}(\|u\|_{L^2})$   \\ 
\noalign{\smallskip}\hline\noalign{\smallskip}
QSA($\bftau^\ast$)-closure \eqref{KSE_closure}, $\bftau^\ast$ minimizing the $J_n$'s & $15.7\%$ &   $3.2\%$    &  $3.8\%$ \\
\noalign{\smallskip}\hline\noalign{\smallskip}
QSA($\bftau^\ast$)-closure \eqref{KSE_closure}, $\bftau^\ast$ minimizing the $Q_n$'s & $15.7\%$ &   $6.9\%$    &  $1.3\%$ \\
\noalign{\smallskip}\hline\noalign{\smallskip}
Galerkin ($m=49$) & $0.9\%$  & $42.1\%$ & $307\%$ \\
\noalign{\smallskip}\hline\noalign{\smallskip} 
Galerkin ($m=53$) & $0.4\%$  & $16.6\%$ & $101\%$ \\
\noalign{\smallskip}\hline\noalign{\smallskip} 
Galerkin ($m=58$) & $0.2\%$  & $3.1\%$ & $5.8\%$ \\
\noalign{\smallskip}\hline\noalign{\smallskip} 
Galerkin ($m=61$) & $0.1\%$  & $0.8\%$ & $3.1\%$ \\
\noalign{\smallskip} \bottomrule 
\end{tabular}
\end{table}
%%%%%%%%%%%%%%%%%%%%%%%%%%
\er

\subsection{Closure results in presence of 90 pairs of unstable modes}\label{Sec_More_Results}
The ability of the optimal QSA to fix the backscatter transfer of energy issue, providing thus an efficient closure, is further tested by applying the PM approach to an even more turbulent regime, namely Regime B (see Table \ref{Table_RegimeB}) that exhibits 90 pairs of unstable modes. Due to the scaling \eqref{Eq_conversion_I_to_III} and the large value of $\alpha$ (see Table \ref{Table_RegimeB}) the time variable for Eq.~\eqref{eq:KSE3} evolves on a much smaller timescale than for  Eq.~\eqref{eq:KSE} and as a consequence we will often emphasize the number of snapshots that a given time instant represents rather than giving the (small) value of this time.

 %%%%%%%%%%%%%%%%%%%%%%%%%%%%%%%%%%%%%%%%%%
\begin{figure}
 \centering
  \includegraphics[height=.8\textwidth,width=1\textwidth,]{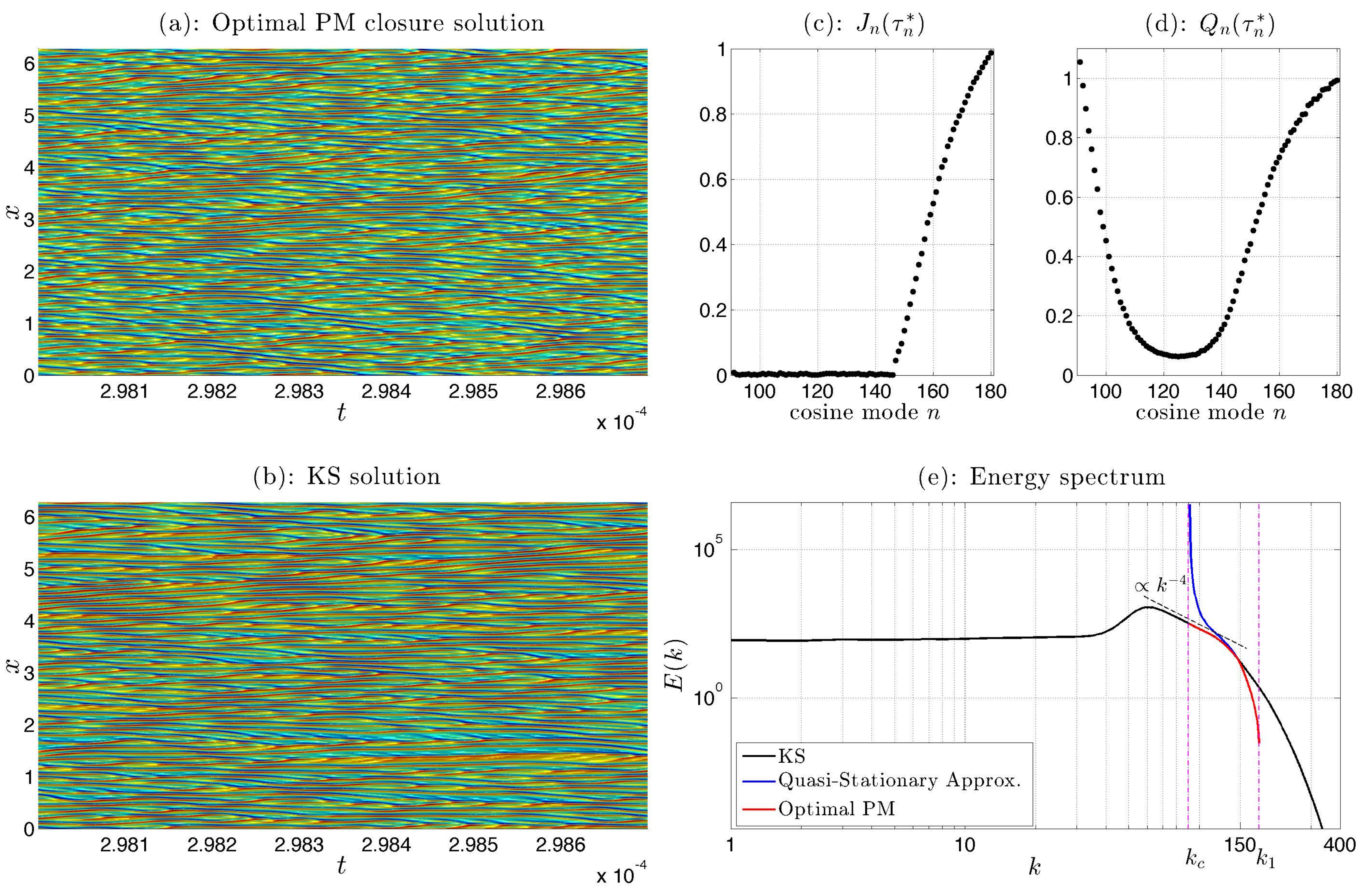}
\caption{{\footnotesize {\bf Closure and parameterization results for Regime B}.  Same as Fig.~\ref{Quantized_PM} except that $k_\c=90$, since Regime counts 90 pairs of unstable modes. The energy spectrum $E(k)$ in panel (e) is estimated over $N=2\times 10^6$ snapshots whereas the optimal QSA is determined by minimizing the cost functional, $J_n$, exploiting the first $2\times 10^4$ snapshots (after removal of transient).  
Figure \ref{Quantized_PM4} shows blowup regions of panels (a) and (b) corresponding to $2.5\leq x\leq 4$.}}
 \label{Quantized_PM3}
\end{figure}
%%%%%%%%%%%%%%%%%%%%%%%%%%%%%%%%%%%%%%%%%%

 %%%%%%%%%%%%%%%%%%%%%%%%%%%%%%%%%%%%%%%%%%
\begin{figure}
 \centering
  \includegraphics[height=.4\textwidth,width=1\textwidth,]{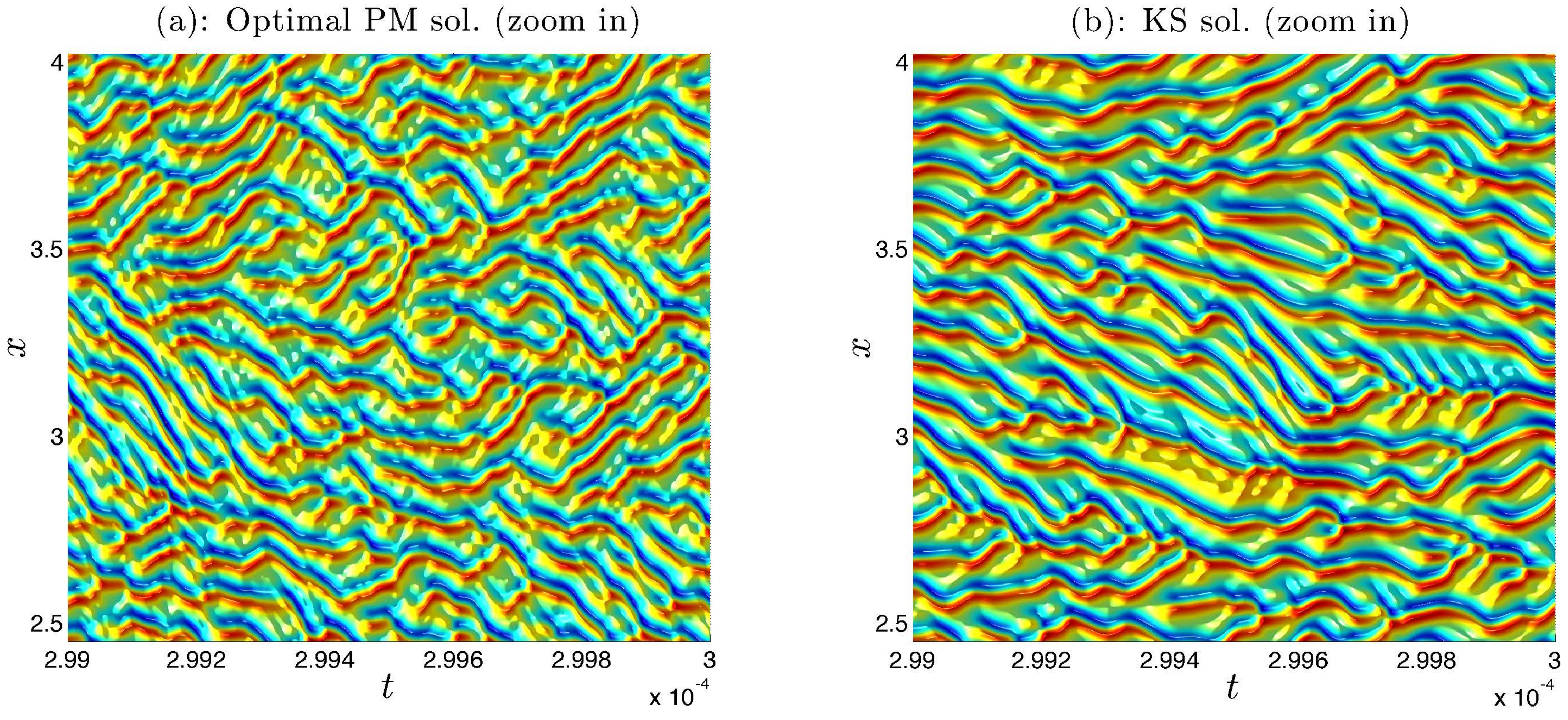}
\caption{{\footnotesize {\bf Closure results for Regime B: Patterns}. Blowup regions of panels (a) and (b) of Fig.~\ref{Quantized_PM3} corresponding to $2.5\leq x\leq 4$.}}
 \label{Quantized_PM4}
\end{figure}
%%%%%%%%%%%%%%%%%%%%%%%%%%%%%%%%%%%%%%%%%%

Here again we take the cutoff scale to be given by the smallest scale (higher wavenumber) contained among the unstable modes. Thus for Regime B, $k_\c=90$, and here also, $15.7\%$ of the total amount of energy needs to be parameterized at this cutoff scale.  For this more turbulent regime, the standard QSA fails even more dramatically than for Regime A and leads to an (ridiculous) over-parameterization of $E(k)$ by an amount of about $35 \times 10^3$ \% (in average) over the range of wavenumbers  $91\leq k\leq 121$; see blue curve in Fig.~\ref{Quantized_PM3}-(e). In contradistinction, the optimal QSA, QSA($\bftau^\ast$), obtained by minimizing  $J_n$ given in \eqref{J_n_normalized} with $T$ that corresponds to the first $2\times 10^4$ snapshots (after removal of transient)\footnote{Note that a blind regression would lead in this case to $89\times 45 \times 3=12015$ coefficients to estimate for each high mode; a number of coefficients comparable to the number of snapshots making thus the estimated coefficients by regression   
non-robust. Instead, one benefits here again greatly from the parametric (and dynamically-based) form of QSA$(\bftau)$ and only $2$ scalar parameters ($\tau_n^\ell$, $\ell=0,1$) need to be determined, for each high mode.}, leads to an average error of about $0.7\%$ over the same range of wavelengths, fixing thus here again the backscatter transfer of energy to the large scales.    As a consequence, good closure skills are obtained as shown in Fig.~\ref{Quantized_PM3} for the reproduction of KS patterns, demonstrating furthermore the robustness of our approach to even more turbulent regimes.  Note that $Q_n$ is greater than $1$ only for $n=91$ (see panel (d) of Fig.~\ref{Quantized_PM3}).  This does not affect the overall quality of the QSA$(\bftau^\ast)$-parameterization (optimized for the $J_n$'s) and we have still  
$Q_T$ given by \eqref{Eq_PD} that is strictly less than $1$, here. 

 A finer inspection of the patterns is made possible by Fig.~\ref{Quantized_PM4} which shows blowup regions of panels (a) and (b) of Fig.~\ref{Quantized_PM3}. Here, we observe that as time evolves the creation and annihilation of the humps displayed by the optimal PM closure solution is reminiscent with what can be observed for the KS solution. Statistically, the spatial correlations are also well reproduced for Regime B as shown in panel (b) of Fig.~\ref{Fig_spatial_Cx}. Only the small-scale features of the optimal PM closure solution and the spatial coherence at long-range distance require improvements, and in that respect one might pursue some ideas proposed in Sec.~\ref{Sec_concluding_rmk} below.

%%%%%%%%%%%%%%%%%%%%%%%%%
\begin{figure}
 \centering
 \includegraphics[height=.5\textwidth,width=1\textwidth,]{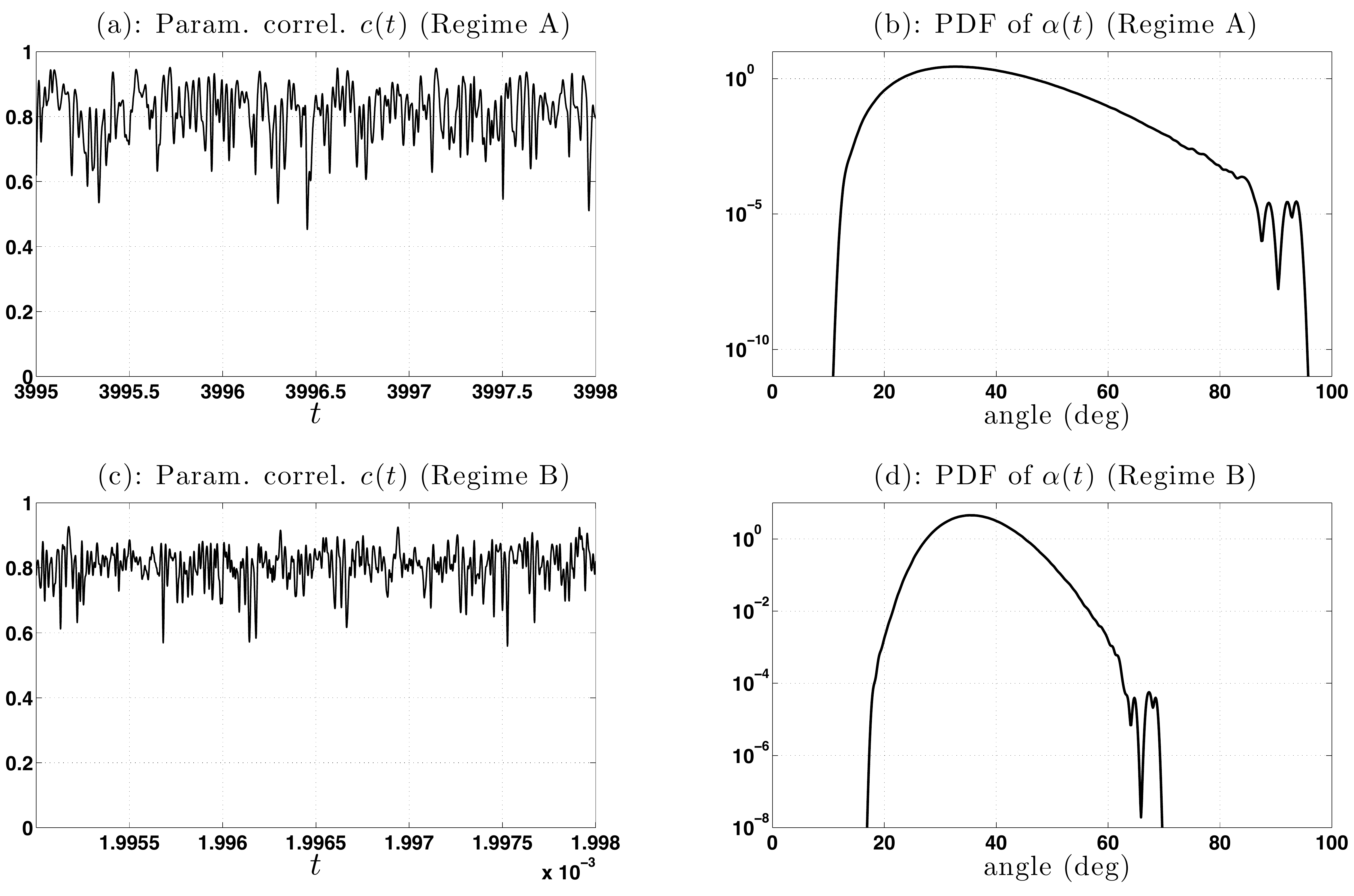}
\caption{{\footnotesize {\bf Parameterization correlation $c(t)$, and PDF of the parameterization angle $\alpha(t)$.} Here these results are obtained for the optimal QSAs, QSA($\bftau^\ast$) used in Fig.~\ref{Quantized_PM} for Regime A, and in Fig.~\ref{Quantized_PM3} for Regime B, that is with $\bftau^\ast$ minimizing the $J_n$'s with $n\geq k_\c=31$ for Regime A, and $n\geq k_\c=90$, for Regime B. A semi-log scale is used for panels (b) and (d).}} 
 \label{Fig_Q_param_correl}
\end{figure}
%%%%%%%%%%%%%%%%%%%%%%%%%

These closure and parameterization skills are put into perspective by computing for each regime, the parameterization correlation, $c(t)$, (see \eqref{Eq_corr_param}) and PDF of the corresponding parameterization angle, $\alpha(t)$ (see \eqref{Eq_alpha}).
As shown in panels (a) and (c) of Fig.~\ref{Fig_Q_param_correl}, $c(t)$ fluctuates away from 1, and $\alpha(t)$ fluctuates over a broad range of values relatively far away from zero. This situation is indicative that for both regimes, the optimal PM computed here is far from a slaving situation.

However, the distribution of $\alpha(t)$ does not seem to be consistent with the good closure results shown here and the rule of thumb pointed out in Sec.~\ref{Sec_corr_param}. The reason behind this is the large number of modes parameterized (here 90 pairs) that makes the parameterization correlation less representative of the quality of a given parameterization than for low-dimensional systems. In the same vein that we have used modewise parameterization defects (the $Q_n$'s) instead of the global parameterization defect $Q_T(t,\Psi_{\bftau^\ast})$ given by \eqref{Eq_PD}, we  inspect below a modewise version of $c(t)$ to diagnose our parameterizations.  

\begin{figure}
 \centering
  \includegraphics[height=.35\textwidth,width=1\textwidth,]{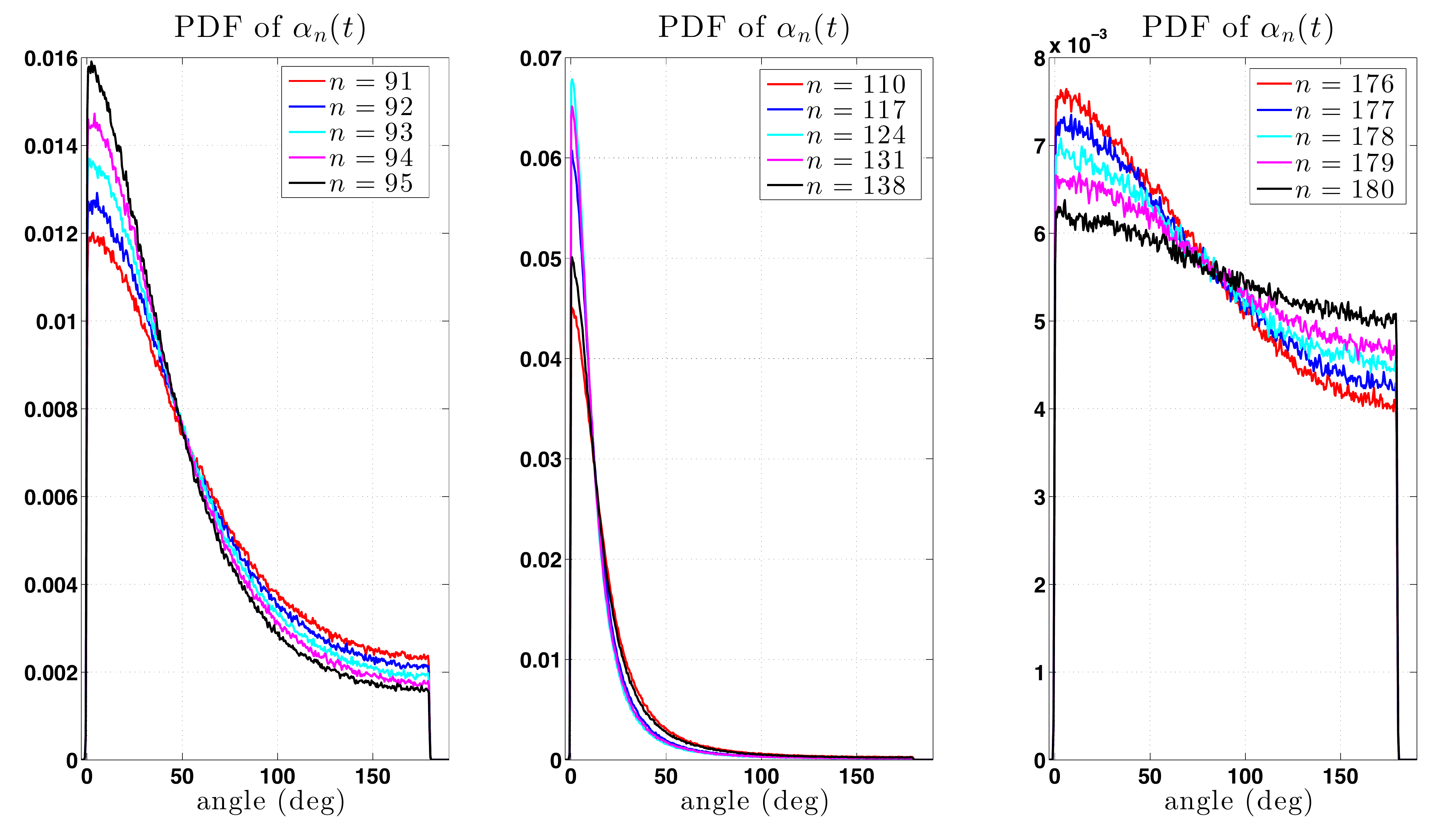}
\caption{{\footnotesize {\bf PDFs of  $\alpha_n(t)$ given by \eqref{Eq_alpha_n}.} Here the PDFs are shown in linear scale.}}
 \label{Angle_cn}
\end{figure}

In that respect, for the bidimensional real vector $\boldsymbol{f}_n(t)=(f_n^0(t),f_n^1(t))$  with  $f_n^\ell(t)=\Psi_n^{\ell}(\tau_{n,\ell}^{\ast}, y_{\c}(t))$, $\ell=0,1$, we introduce 
\be\label{Eq_corr_param_n}
c_n(t)= \frac{\langle \boldsymbol{f}_n(t), \boldsymbol{y}_n(t) \rangle}{\|\boldsymbol{f}_n(t)\| \; \|\boldsymbol{y}_{n}(t)\|}.
\ee
and the following parameterization angle,
\be\label{Eq_alpha_n}
\alpha_n(t)=\arccos (c_n(t)). 
\ee
We computed $c_n(t)$ and $\alpha_n(t)$ for $n=91$ through $n=180$. Figure \ref{Angle_cn} shows 
the results for the PDFs of $\alpha_n(t)$, as gathered into three groups: a group of parameterized modes adjacent to the cutoff scale, a group of modes (well) within the inertial range, and a group of modes corresponding to the smallest scales parameterized. Clearly the PDFs corresponding to the 2nd group of modes correspond to the best modewise parameterizations; compare middle panel of Fig.~\ref{Angle_cn} with the two other panels of the same figure. Here, we observe for this group of modes PDFs that exhibit features discussed in  Sec.~\ref{Sec_corr_param}. These PDFS are  indeed skewed towards zero with the most frequent value of $\alpha_n(t)$ also close to zero; cf.~black curve in Fig.~\ref{Fig_intro_angle}.  These features are also shared by the PDFs of the adjacent modes to the cutoff scale (left panel of Fig.~\ref{Angle_cn}) with however a fat tail towards high values of $\alpha_n(t)$.   The last group of modes corresponding to high wavenumbers  (right panel of Fig.~\ref{Angle_cn}) corresponds to the less accurate modewise parameterizations  as manifested by PDFs of $\alpha_n(t)$ that although skewed are somewhat close to a uniform distribution.

These small-scale modes are weakly energetic, they contain less than 0.6 $\%$ of the total energy for $n>150$, and  here do not spoil the parameterization noticeably.  However the fat tails of the PDFs corresponding to the adjacent parameterized modes is a determining factor responsible of pushing  the (global) parameterization correlation, $c(t)$ (given by \eqref{Eq_corr_param}), away from  1, as it can be observed by removing the contribution of these modes in the calculation of $c(t)$ (not shown). On the other hand,  these adjacent modes are important dynamically and cannot be removed for closure as they contain an amount of energy comparable to that of the modes right below the cutoff scale (i.e.~for $k<k_\c$).

%%%%%%%%%%%%%%%%%%%%%%%%%%%%%%%%%%%%%%%%%%
%%%%%%%%%%%%%%%%%%%%%%%%%%%%%%%%%%%%%%%%%%

We conclude by reporting on how the choice of the cost functional and class of parameterization impacts the closure skills.   
The metrics used to assess these skills are those used for  Table \ref{Table_energy_PM_regimeA}, namely $\overline{\|u\|}_{L^2}$ (after transient removal) and the standard variation, $\texttt{std}(\|u\|_{L^2})$. The time averages are here estimated on $2\times 10^4$ snapshots.  As Table  \ref{Table_energy_PM_regimeB} shows, minimizing the $Q_n$'s instead of the $J_n$'s leads to a deterioration in the approximation of  $\overline{\|u\|}_{L^2}$ but an improvement in the standard variation within a given class of parameterizations.

%%%%%%%%%%%%%%%%%%%%%%%%%%%%%%%%%%%%

\begin{figure}
 \centering
  \includegraphics[height=.5\textwidth,width=.85\textwidth,]{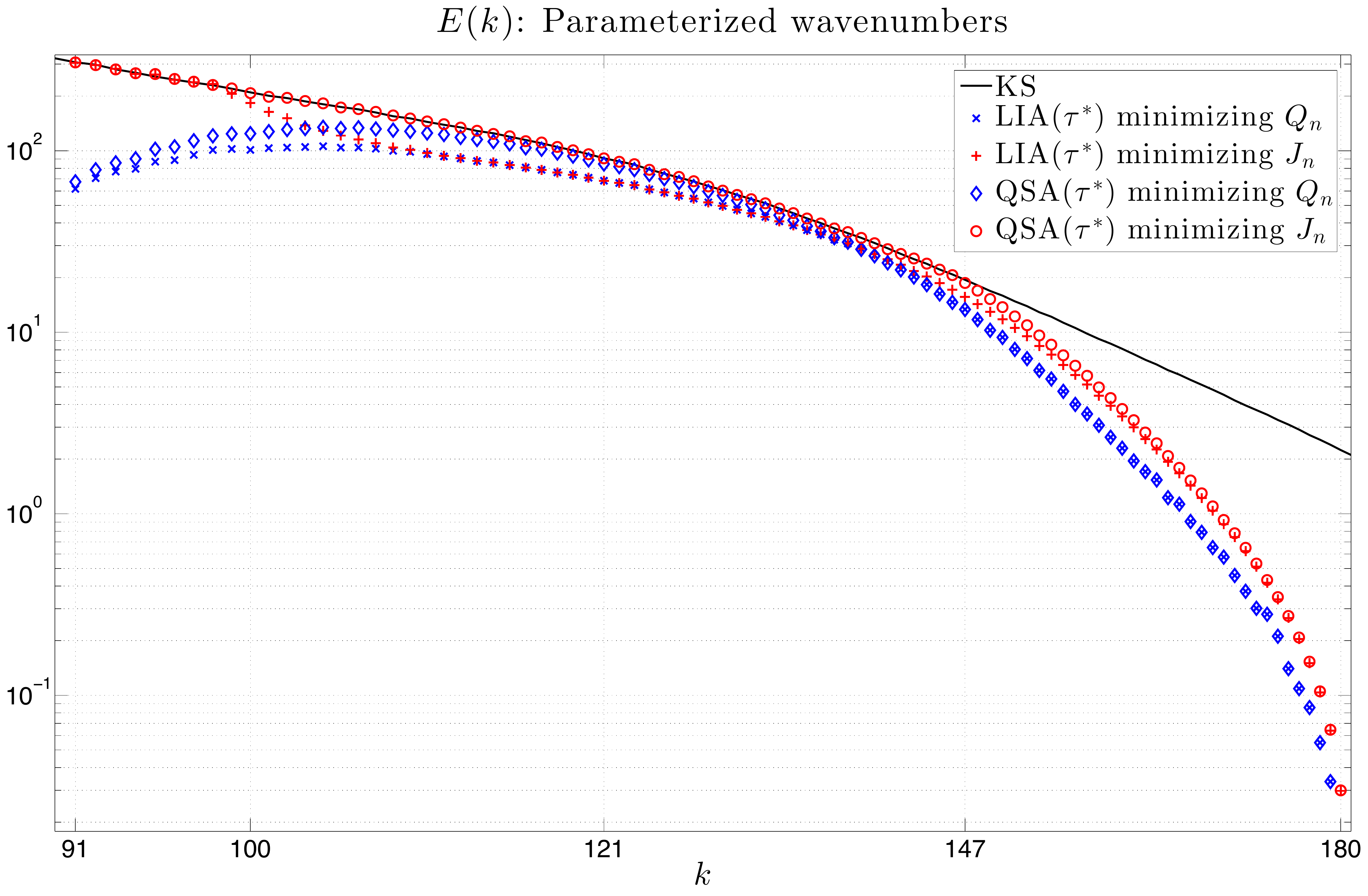}
\caption{{\footnotesize {\bf Approximations of $E(k)$ for $k_\c< k\leq k_1$ for Regime A}. Optimal LIA($\bftau^\ast$) and QSA($\bftau^\ast$) with $\bftau^\ast$ minimizing either the cost functionals $J_n$'s or $Q_n$'s. Recall that $k_\c=90$ and $k_1=2 k_\c$. A log-log scale is used here.}}
 \label{Quantized_PM3b}
\end{figure}

The portion of the energy spectrum $E(k)$ parameterized---by the optimal LIA($\bftau^\ast$) or QSA($\bftau^\ast$) with $\bftau^\ast$ minimizing either the cost functionals $J_n$'s or $Q_n$'s---is shown in Fig.~\ref{Quantized_PM3b}. As one can observe, the
QSA($\bftau^\ast$) obtained by minimizing the $J_n$'s provides the best result and an almost perfect parameterization of the energy contained in the high modes over the range of wavenumbers, $91\leq k\leq 147$, resulting thus into the good closure skills shown in Fig.~\ref{Quantized_PM3} and panel (b) of Fig.~\ref{Fig_spatial_Cx}.  We emphasize that as for Regime A, these skills are obtained from an optimal PM designed 
from a training interval over which the statistics of $|u_n|^2$ have not yet stabilized; cf.~discussion relative to Fig.~\ref{Fig_modewise_error} for Regime A.    When the $Q_n$'s are used to optimize either the LIA($\bftau$)- or the QSA($\bftau$)-parameterization, one observes an under-parameterization more pronounced near the cutoff scale $k_c=90$ and that vanishes as $k$ is increased, before re-emerging beyond wavenumbers that contain a small fraction of the total energy $E_{\textrm{tot}}$; for instance the scales beyond $k=147$, contain only $0.6\%$ of $E_{\textrm{tot}}$.  Despite this under-parameterization, the optimal  LIA($\bftau^\ast$) and  QSA($\bftau^\ast$) with $\bftau^\ast$ minimizing the $Q_n$'s, provide also closure skills comparable to those shown in Fig.~\ref{Quantized_PM3} and panel (b) of Fig.~\ref{Fig_spatial_Cx}.  The main differences are actually observed at the level of the approximation of $\overline{\|u\|}_{L^2}$ and $\texttt{std}(\|u\|_{L^2})$, as summarized in Table \ref{Table_energy_PM_regimeB}. We refer to the heuristic discussion at the end of Sec.~\ref{Sec_FMTtau} to better appreciate the nuances between the LIA- and QSA-classes of parameterizations  in regards of these numerical results.

%%%%%%%%%%%%%%%%%%%%%%%%%%
\begin{table}[h] 
\caption{1st and 2nd moments of $\|u\|_{L^2}$: Relative error for Regime B}
\label{Table_energy_PM_regimeB}
\centering
\begin{tabular}{ccc}
\toprule\noalign{\smallskip}
 &  $\overline{\|u\|}_{L^2}$ &  $\texttt{std}(\|u\|_{L^2})$   \\ 
\noalign{\smallskip}\hline\noalign{\smallskip}
QSA($\bftau^\ast$)-closure , $\bftau^\ast$ minimizing the $J_n$'s  &   $4\%$    &  $3.2\%$ \\
\noalign{\smallskip}\hline\noalign{\smallskip}
QSA($\bftau^\ast$)-closure , $\bftau^\ast$ minimizing the $Q_n$'s  &   $7.5\%$    &  $1.6\%$ \\
\noalign{\smallskip}\hline\noalign{\smallskip}
LIA($\bftau^\ast$)-closure with $\bftau^\ast$ minimizing the $J_n$'s  &   $8.9\%$    &  $1.7\%$ \\
\noalign{\smallskip}\hline\noalign{\smallskip}
LIA($\bftau^\ast$)-closure with $\bftau^\ast$ minimizing the $Q_n$'s &   $10.2\%$    &  $0.3\%$ \\
\noalign{\smallskip} \bottomrule 
\end{tabular}
\end{table}
%%%%%%%%%%%%%%%%%%%%%%%%%%

%%%%%%%%%%%%%%%%%%%%%%%%%%%%%%%%%%%%%%%%%%
%%%%%%%%%%%%%%%%%%%%%%%%%%%%%%%%%%%%%%%%%%

%%%%%%%%%%%%%%%%%%%%%%%%%%%%%%%%%%%%%%%%%%%%%%%%%%%%%
\section{Concluding remarks}\label{Sec_concluding_rmk}

Thus, the PM approach is not limited to a class of parametric parameterizations nor to a particular cost functional. As the closure exercise shows here 
in the context of KS turbulence, a good choice of the cost functional and class of parameterizations to optimize is nevertheless key to approximate certain features better than others. This is where the specificities of the problem at hand plays an important role\footnote{In that respect, we may mention the variational normal mode initialization in Meteorology,  pioneered by Daley \cite{daley1978variational}, who combined the Machenhauer \cite{machenhauer1977dynamics} non-linear normal-mode initialization within a variational procedure allowing for the adjustment of confidence weights arising in a fidelity functional $I$;  see also \cite{tribbia1982variational}. In these works, the manifold $\mathcal{M}$ is fixed a priori and it is the point on $\mathcal{M}$ nearest to the observation using the ``metric'' defined by $I$, that is sought.} 
 and where one may benefit from the flexibility of the PM approach to optimize relevant parameterizations known by the practitioner, once the underlying formulas are made parametric, i.e.~made as a function of a (collection) of (independent) scalar variable(s).

Rooted in the rigorous approximation theory of invariant manifolds (Part I), this articles  provides 
a natural framework to extend the corresponding approximation formulas as nonlinear parameterizations useful when 
slaving relations do not hold anymore, e.g., away from criticality (Part II).  The framework opens up several possible directions for future research. We outline some of these directions below.

{\it 1. Time-dependent parameterizing manifolds for non-autonomous systems.}
 As for the autonomous 
case discussed here, formulas for time-dependent PMs may be rooted in the approximation theory of time-dependent invariant manifolds \cite{PR06,PR09}. The leading order approximation, $h_2$, becomes now time-dependent and satisfies the following version of the homological equation \eqref{h1_eqn} (with $\mathcal{L}_A$ defined in \eqref{eq:ad_A}),
 \be\label{Eq_linPDEsyst}
  \Big(\partial_t+ \mathcal{L}_A\Big) h=\Pi_\s B(\xi,\xi) +\Pi_\s F(t),
 \ee
 for a system of the form 
 \be\label{Eq_ODE_gen_bis}
\frac{\d y}{\d t} = A y + B(y,y) +F(t), \qquad y\in \mathbb{C}^N.
\ee
The backward-forward method to derive parametric formulas for PMs, extends to this non-autonomous setting and provides 
a parametric family of time-dependent manifold function, $\Psi^{(1)}_\tau(t,\cdot)$, that satisfies 
for instance in the case $\Pi_\c F=0$, the following modification of Eq.~\eqref{Eq_h1_taugen}
\be \label{Eq_invariance_NDS_v2}
 \Big(\partial_t+  \mathcal{L}_A \Big)\Psi^{(1)}_\tau (t,\xi)= \Pi_{\s} B(\xi,\xi) - e^{\tau A_{\s}} \Pi_{\s} B(e^{-\tau A_{\c}} \xi, e^{-\tau A_{\c}} \xi)  +\Pi_\s F(t)- e^{\tau A_{\s}}  \Pi_{\s} F(t-\tau).
\ee

Due to the time-dependent coefficients to calculate in $\Psi^{(1)}_\tau(t,\cdot)$, the evaluation of the parameterization  defect gets more involved than in the autonomous case. Nevertheless,  the optimal value for the free parameter $\tau$ may be still obtained by minimizing this defect, leading to an optimal PM, in the $\Psi^{(1)}_\tau(t,\cdot)$-class  and thus to closures with time-dependent coefficients.  The measure-theoretic framework of Sec.~\ref{Sect_PM_reduction} may benefit here from the theory of SRB measures for non autonomous systems \cite{Young2016}. The formulas for the LIA and QSA parameterizations of Secns.~\ref{Sect_PM_with_forcing} and \ref{Sec_FMTtau} respectively, extend to this non-autonomous setting as well. The case of a stochastic forcing can be dealt with along the same lines, the backward-forward method providing in this case parametric formulas for PMs that come with non-Markovian coefficients depending on time-history of the noise (exogenous memory terms) \cite{CLW15_vol2}.

{\it 2. Combining PMs with stochastic parameterizations.}  To set the framework, we discuss stochastic improvements that can be made to the LIA class of Sec.~\ref{Sect_PM_with_forcing}, but the ideas apply to the QSA class of Sec.~\ref{Sec_FMTtau} as well.
Given a cutoff dimension $m$, the optimal PM obtained by solving the minimization problems  
\eqref{Min_formulation_h1_quadcase}, for $n\geq m+1$, is the best manifold  --- in the LIA class --- that averages out the unresolved fluctuations lying in $E_\s$.  Once the optimal PM, $\Phi_{\bftau^\ast}^{(1)}$, has been determined, we may still want to parameterize these fluctuations. 
These fluctuations are given by the residual $\eta_t$ whose components are determined after having solved \eqref{Min_formulation_h1_quadcase} for each $n\geq m+1$.
We have then 
\be\label{Eq_key}
y_\s(t)=\Phi_{\bftau^\ast}^{(1)}(y_\c(t)) +\eta_t.
\ee 
From a closure viewpoint, we are thus left with the stochastic modeling of $\eta_t$. The next step consists of seeking for a stochastic parameterization $\zeta_t$ of $\eta_t$.  Here several approaches are possible; see \cite{gottwald_crommelin_franzke_2017} for a survey.  The idea of incorporating a stochastic ingredient as a supplement to a nonlinear parameterization is not new and has been proposed in the context of two-dimensional turbulence \cite{leith1990stochastic}, atmospheric turbulence \cite{frederiksen2006dynamical} and more recently, oceanic turbulence \cite{zanna2017scale}.

Once a satisfactory stochastic parameterization $\zeta_t$ has been determined, we arrive at the following closure for the resolved variable  (in the case of bilinear system),
\be\label{Eq_closure_gen}
\frac{\d z}{\d t } = A_\c z +  \Pi_\c B\Big(z+ \Phi_{\bftau^\ast}^{(1)}(z)+\zeta_t,  z+ \Phi_{\bftau^\ast}^{(1)}(z)+\zeta_t\Big) + \Pi_\c F.
\ee
Thinking of $B$ as given by a nonlinear advective term, we see that the stochastic parameterization \eqref{Eq_key} brings new elements to the closure \eqref{Eq_closure_gen} such as stochastic advective terms compared to a closure that would be only based on the optimal PM. 
 Other recent approaches have shown the relevance of such stochastic advective terms  
to derive stochastic formulations of classical representations of fluid flows as well as for emulating suitably the coarse-grained dynamics \cite{holm2015variational,resseguier2017geophysical,resseguier2017geophysicalb,resseguier2017geophysicalc,arnaudon2018noise,cotter2019numerically}.

The selection of the best parameters (e.g.~lags for an auto-regressive process) of a given stochastic parameterization aimed at emulating the residual, $\eta_t$, can here again 
be guided by the minimization of the parameterization defect $Q_n$; the parameters of $\zeta_t$ being determined so as to minimize further $Q_n$ 
compared to when the optimal PM is used alone.  Complementarily, the parameterization correlation, $c(t)$, for which  $\Psi=\Phi_{\bftau^\ast}^{(1)}+\zeta_t$  in \eqref{Eq_corr_param}, can then be evaluated to  further revise other ingredients in the  stochastic parameterization, so that the probability distribution of the corresponding correlation angle $\alpha(t)$ gets skewed towards zero as much as possible.   In other words, one should not only parameterize properly the statistical effects of the subgrid scales but also avoid to lose their phase relationships with the retained scales \cite{mccomb2001conditional}. In that respect, the residual noise $\eta_t$ in \eqref{Eq_key} is expected to depend on the state of the resolved variable $\xi$.  The abstract formula \eqref{Def_h2} for the optimal PM suggests that subgrid-scale parameterization techniques with conditional Markov chains \cite{crommelin2008subgrid,kwasniok2012data,gottwald2016data} constitute a  consistent tool with our approach for the design of a stochastic parameterization $\zeta_t$.

3. {\it Beyond conditional expectation: Memory effects and noise.} An alternative to the inclusion of stochastic ingredients as discussed above, relies on  Theorem \ref{Thm_variational-pb2} as a starting point. The latter theorems shows that once an optimal PM is found, it provides the conditional expectation (in the case $\eta=0$). Nevertheless, as shown in Sec.~\ref{Sec_L9D_resurect}, the conditional expectation alone, let us say $\R$, is sometimes insufficient to close fully the system.  The Mori-Zwanzig formalism \cite{Mor65,zwanzig2001}
of statistical physics, instructs us then that 
a complete closure exists under the form of the following {\it generalized Langevin equation (GLE)} \cite{GKS04,Chorin_Hald-book,MSM2015,gottwald_crommelin_franzke_2017},
\begin{equation}\label{Eq_GLE}\tag{GLE}
\dot{x}=\R(x)+\int_{0}^{t}\mathbf{G}(t,s,x(s))\d s+\eta_t.
\end{equation}

Here, the integral term accounts for the nonlinear interactions between the resolved and unresolved variables that are not accounted for in $\R$; it involves the past of the macroscopic variables and conveys {\it non-Markovian (i.e.~memory) effects}.  The term $\eta_t$ accounts for effects of the unresolved variables which are uncorrelated with the resolved variables. This last term  can be thus represented by a state-independent noise that may still involve correlations in time, i.e.~of ``red noise'' type. It is well known that the analytical determination of the constitutive elements of the GLE is a difficult task in practice.  By relying on Theorem \ref{Thm_variational-pb2} and formulas of Sec.~\ref{Sect_PM_formulas}, the PM approach  can be seen as providing an  efficient way to approximate the conditional expectation $\R$ in \eqref{Eq_GLE}. However, the practical  determination of the memory and stochastic terms remains a challenge, especially for fluid flows \cite{MSM2015,gottwald_crommelin_franzke_2017}.   Various approaches have been proposed to address this aspect that include for instance short-memory approximations \cite{Chorin2002}, the $t$-model \cite{hald_stinis,stinis_euler}, formal expansions of the Koopman operator \cite{wouters2012,WL13},  NARMAX techniques  \cite{chorin2015discrete,lu2017data}, and the dynamic-$\tau$ model \cite{parish2016reduced,parish2017dynamic}. 
 See also \cite{kraichnan1959structure,kraichnan1964approximations,Herring72}, \cite{mcwilliams2012elemental}, and \cite{mana2014toward,zanna2017scale} for other reduced modeling/parameterization approaches that involve memory terms (and noise) in the context of homogeneous turbulence, shear dynamo and oceanic turbulence, respectively.

Once $\R$ is approximated from an optimal PM, the practical  determination of the memory and stochastic term could also benefit from the data-driven modeling techniques of \cite{CK17},  to model the residual, $\dot{y_\c}-\R(y_\c)$, where  $y_\c$ denotes the low-mode projection of a fully resolved solution $y$. 
As illustrated and discussed in \cite{KCB18} for a wind-driven ocean gyres model, the data-driven techniques of \cite{CK17} have been successfully applied to model the coarse-scale dynamics. 
To operate in practice,  the data-driven techniques of \cite{CK17} require observations of $y(t)$ of length comparable also to a decorrelation time of the dynamics  \cite{CK17,KCG18,KCYG18}, as for the optimization of the dynamically-based PMs of Sec.~\ref{Sect_PM_formulas}.

 4. {\it Combining modal reductions and the PM approach.} 
In many applications such as arising in turbulence, the number of ODEs  associated to a given discretization, is very large. 
This is where modes computed in the physical domain from DNS may be used to proceed to a first reduction (data compression) of the phase space.  
Among the most commonly employed modal decomposition techniques are the proper orthogonal decomposition (POD) \cite{Holmes_al12}, 
and its variants; see \cite{taira2017modal} and references therein. Of demonstrated relevance for the reduction of nonlinear PDEs are also the principal interaction patterns (PIPs) modes \cite{hasselmann1988pips,kwasniok1996reduction,kwasniok1997optimal} that find a compromise between minimizing tendency error with maximizing explained variance in the resolved modes; see \cite{kwasniok2004empirical,kwasniok2007reduced} for applications to atmospheric models, and \cite{crommelin2004strategies} for a very clear comparison between POD and PIP modes. 
In the last decade, related promising techniques such as the dynamic mode decomposition (DMD) \cite{rowley2009spectral,schmid2010dynamic,williams2015data,taira2017modal} have also emerged; see \cite{tu2013dynamic} for a discussion on the relationships between PIPs, DMD, and the linear inverse modeling \cite{penland1993prediction}. 

Also, the use of time-dependence in the basis elements ---
the so-called Dynamical Orthogonal  (DO) modes \cite{sapsis_al09,sapsis_al12} --- have been considered, as in principle it allows for the representation of the transient character of the solution using much fewer modes. A dynamical orthogonality condition leads then to a closed set of equations that allows for the evolution of the mean field, the DO modes and the corresponding (stochastic) coefficients \cite{feppon2018geometric}.  From the mean, the time-dependent patterns of the DO modes plus the distribution of the stochastic coefficients (at a certain time $t$), an approximation to the probability density function of the state vector can be obtained \cite{SD13,ueckermann2013numerical,narayanan2018probabilistic}. In terms of computational performance, there is however a trade-off between fewer modes to consider on one hand, and more equations (including interactions between the modes) to solve, on the other.

For certain problems of turbulence, even  after modal reduction, one may wish still to further reduce the dimension of the ODE approximation.
Whatever the modes used to represent the dataset at hand, one should avoid to compute parameterizations by taking the reduced state space, $E_\c$, to be spanned by only the first few modes. There are several reasons behind this caution. 
One reason is that these modes may mix the large and small spatial scales, making the distinction between $E_\c$ and $E_\s$ not obvious. 
Another reason, more technical,  is that $E_\c$ and its complement $E_\s$ are no longer invariant subspaces for the linear part of the original PDE, which introduces linear interaction terms  between the modes in $E_\c$  and $E_\s$ that have to be taken into account for the parameterization. Although one could still apply formally the backward-forward method of Sec.~\ref{Sect_PM_formulas} to derive parametric families of parameterizations, a more reasonable approach consists of proceeding directly from the Galerkin ODE systems obtained by projecting the original PDE onto these modes. 
This way, we are indeed left with the theory and techniques presented in this article, and by determining the equations for the perturbed variable about a mean state and work within the eigenbasis of the linearized operator, we can then use the dynamically-based formulas of Sec.~\ref{Sect_PM_formulas} to calculate and optimize the parameterizations.

\section*{Acknowledgments}
MDC wishes to acknowledge David Neelin for the stimulating discussions on the closure problem of convective processes in the tropical atmosphere. MDC and JCM are also thankful to  Darryl Holm for his constructive comments at the beginning of this work.   Finally,  MDC and HL are greatly indebted to Shouhong Wang for the numerous and stimulating discussions about this work over the years, and it is a pleasure to express our gratitude to Shouhong for his constant encouragement. This work has been partially supported by the National Science Foundation grants DMS-1616981 (MDC) and DMS-1616450 (HL).

\appendix
\section{Parameterization defect minimization algorithm} \label{Sect_gradient_descent}

We present in this Appendix a simple gradient-descent method to solve efficiently the  minimization problem \eqref{Min_formulation_h1_quadcase} in order to determine the  optimal $\tau$-value, $\tau^\ast$, 
for the parameterization, $\Phi_n(\tau,\boldsymbol{\beta}, \xi)$, given by \eqref{Eq_Phi_tau}. As shown below, the method allows furthermore for making apparent the dependence of the parameterization defect  on statistical moments (up to order 4) of the original system's solution.

To present the method, we first recast the parameterization defect associated with $\Phi_n$, 
\be \label{Eq_defect_Markovian}
\mathcal{Q}_n(\tau,T)= \frac{1}{T} \int_0^T \big| \Pi_n y(t)- \Phi_n(\tau, \boldsymbol{\beta}, \Pi_{\c}y(t))\big|^2 \d t,
\ee 
into a matrix format. For this purpose, we arrange the coefficients $D^{n}_{i, j}(\tau,\boldsymbol{\beta})B_{i, j}^n$ involved in the expression of $\Phi_n(\tau,\boldsymbol{\beta}, \xi)$ into an $m^2\times 1$ vector $\boldsymbol{d}(\tau)$ so that the indices $(i,j)$'s are arranged in lexicographical order; namely the $k^{\mathrm{th}}$ component of $\boldsymbol{d}(\tau)$ is given by
\be \label{Eq_Dijn_1D}
d_{k}(\tau) = D^{n}_{i, j}(\tau,\boldsymbol{\beta})B_{i, j}^n, \quad k = 1, \cdots, m^2, 
\ee
where $(i, j)$ is the unique low-mode pair of indices satisfying 
\be
(i-1)m + j = k, \quad  \text{ with } i, j  \in \{1,\cdots, m\}.
\ee
More precisely, the index pair $(i, j)$ in \eqref{Eq_Dijn_1D} is determined by: 
\be \label{Eq_lex-order-relation}
 \begin{cases}
{\displaystyle i = \frac{k-\mathrm{mod}(k,m)}{m}+1} \;\text{ and }\; j = \mathrm{mod}(k,m), & \text{ if }  \mathrm{mod}(k,m) \neq 0, \medskip \\
{\displaystyle  i = \frac{k}{m}}  \;\text{ and }\;  j = m, &\text{ otherwise}.
 \end{cases}
 \ee
 Similarly, we define an $m^2\times 1$ vector $\boldsymbol{\gamma}(\tau)$, whose components are given by
 \be
 \gamma_k(\tau) = V_{i,j}^n(\tau, \boldsymbol{\beta}) F_{j} (B_{i,j}^n + B_{j,i}^n),  \quad k = 1, \cdots, m^2. 
 \ee

Now, given the solution $y(t)$ to the underlying $N$-dimensional ODE system \eqref{Eq_quadODEs}
over $[0,T]$, we introduce
\bes
u_{k}(t) = \Pi_k y(t),\qquad k =1,\cdots, m,
\ees
where $\Pi_k$ denotes the projection onto the mode $\boldsymbol{e}_k$; see \eqref{Proj_formula}. 

We define next the column vectors $\boldsymbol{Q}_1$, $\boldsymbol{Q}_2$, $\widehat{\boldsymbol{Q}}_2$ and $\boldsymbol{Q}_3$ of size $m^2 \times 1$ as well as the matrices $\widetilde{\boldsymbol{Q}}_2$, $\widetilde{\boldsymbol{Q}}_3$ and $\boldsymbol{Q}_4$ of size $m^2 \times m^2$ as follows: 
\bea \label{Eq_Q_elements_continuous}
& (\boldsymbol{Q}_1)_{p} = \langle {\hl \overline{u}_{p_1}} \rangle_T,  \quad p = 1,\cdots, m^2, \\
& (\boldsymbol{Q}_2)_{p} = \langle {\hl \overline{u}_{p_1} \overline{u}_{p_2}} \rangle_T,  \quad p = 1,\cdots, m^2, \\
& (\widehat{\boldsymbol{Q}}_2)_{p} = \langle u_{n} {\hl \overline{u}_{p_1}} \rangle_T,  \quad p = 1,\cdots, m^2, \\
& (\boldsymbol{Q}_3)_{p} = \langle  u_{n} {\hl \overline{u}_{p_1} \overline{u}_{p_2}} \rangle_T,  \quad p = 1,\cdots, m^2,\\
& (\widetilde{\boldsymbol{Q}}_2)_{pq} = \langle  {\hl \overline{u}_{p_1}}u_{q_1} \rangle_T,  \quad p,\, q = 1,\cdots, m^2, \\
& (\widetilde{\boldsymbol{Q}}_3)_{pq} = \langle  {\hl \overline{u}_{p_1}} u_{q_1}u_{q_2} \rangle_T,  \quad p,\, q = 1,\cdots, m^2, \\
& (\boldsymbol{Q}_4)_{pq} = \langle {\hl \overline{u}_{p_1} \overline{u}_{p_2}} u_{q_1} u_{q_2} \rangle_T,  \quad p, \, q = 1,\cdots, m^2,
\eea
where {\hl $\overline{z}$ denotes the complex conjugate of $z$ in  $\mathbb{C}$}, $\langle \cdot \rangle_T$ denotes the time average over $[0,T]$, and the low-mode index pair $(p_1, p_2)$ (resp.~$(q_1, q_2)$)   relates to $p$ (resp.~$q$) according to \eqref{Eq_lex-order-relation}, namely where $p$ (resp.~$q$)  plays the role of $k$  and $(p_1, p_2)$ (resp.~$(q_1, q_2)$) that of  $(i, j)$ in \eqref{Eq_lex-order-relation}.

Besides, let us recall the constant terms given in the RHS of \eqref{Formula_GammaF} for the parameterization,  $\Phi_n(\tau,\boldsymbol{\beta}, \xi)$:
\bea \label{Formula_alphaF}
\alpha_n(\tau) =  \sum_{i, j = 1}^m  U_{i, j}^n(\tau, \boldsymbol{\beta}) B_{i, j}^n F_{i}F_{j} - \frac{1 - e^{\tau \beta_n}}{\beta_n} F_n.
\eea

Thus, we rewrite the parameterization defect $\mathcal{Q}(\tau, T)$ recalled in \eqref{Eq_defect_Markovian} as follows:  
\bea \label{Eq_defect_Markovian_recast}
\mathcal{Q}_n(\tau,T) = \boldsymbol{d}(\tau)^{\ast} \boldsymbol{Q}_4 \boldsymbol{d}(\tau) 
 & - 2 \mathrm{Re} \big( \boldsymbol{Q}_3^{\ast} \boldsymbol{d}(\tau) \big) 
  +  2 \mathrm{Re} \big( \boldsymbol{\gamma}(\tau)^{\ast}  \widetilde{\boldsymbol{Q}}_3 \boldsymbol{d}(\tau) \big) 
    +   \boldsymbol{\gamma}(\tau)^{\ast}  \widetilde{\boldsymbol{Q}}_2 \boldsymbol{\gamma}(\tau) \\
    & - 2 \mathrm{Re} \big(\widehat{\boldsymbol{Q}}_2^{\ast} \boldsymbol{\gamma}(\tau) \big) 
 + 2 \mathrm{Re} \big( \overline{\alpha}_n(\tau) \boldsymbol{Q}_2^{\ast} \boldsymbol{d}(\tau) \big) 
 + 2 \mathrm{Re} \big( \overline{\alpha}_n(\tau) \boldsymbol{Q}_1^{\ast} \boldsymbol{\gamma}(\tau) \big)   \\
 & + \langle u_n \overline{u}_n \rangle_T - 2 \mathrm{Re} \big(\overline{\alpha}_n(\tau) \langle u_n \rangle_T \big) + \alpha_n(\tau) \overline{\alpha}_n(\tau),
\eea
where $M^\ast$ denotes the conjugate transpose of a given vector or matrix $M$.

Note also 
 \bea \label{Eq_Qprime}
\frac{\d}{\d \tau}\mathcal{Q}_n(\tau,T)  = 2 \mathrm{Re} \Big( & \boldsymbol{d}(\tau)^{\ast} \boldsymbol{Q}_4 \boldsymbol{d}'(\tau) 
      - \boldsymbol{Q}^{\ast}_3 \boldsymbol{d}'(\tau) 
      + \boldsymbol{\gamma}'(\tau)^{\ast} \widetilde{\boldsymbol{Q}}_3 \boldsymbol{d}(\tau)  
            + \boldsymbol{\gamma}(\tau)^{\ast} \widetilde{\boldsymbol{Q}}_3 \boldsymbol{d}'(\tau)  \\
     & + \boldsymbol{\gamma}(\tau)^{\ast} \widetilde{\boldsymbol{Q}}_2 \boldsymbol{\gamma}'(\tau) 
     - \widehat{\boldsymbol{Q}}^{\ast}_2 \boldsymbol{\gamma}'(\tau) 
     + \overline{\alpha}'_n(\tau) \boldsymbol{Q}^{\ast}_2 \boldsymbol{d}(\tau)  
     + \overline{\alpha}_n(\tau) \boldsymbol{Q}^{\ast}_2 \boldsymbol{d}'(\tau) \\
& +\overline{\alpha}'_n(\tau) \boldsymbol{Q}^{\ast}_1 \boldsymbol{\gamma}(\tau)   
 +\overline{\alpha}_n(\tau) \boldsymbol{Q}^{\ast}_1 \boldsymbol{\gamma}'(\tau)   
       - \overline{\alpha}'_n(\tau) \langle u_n \rangle_T 
       +  \alpha'_n(\tau) \overline{\alpha}_n(\tau) \Big).
\eea

With the above expression of $\mathcal{Q}_n(\tau,T)$ and of its derivative,  the minimization of $\mathcal{Q}_n(\tau,T)$ in the $\tau$-variable can now be performed efficiently by application of a gradient-descent method as described in Algorithm~\ref{Algorithm-tau-search}.
Note that if the first moments up to the 4th order are known, then the determination of $\tau^\ast$ 
by Algorithm \ref{Algorithm-tau-search} does not require any data from direct integration of the full system. There is a vast literature about moment closure techniques and we refer to \cite{kuehn2016moment} for a recent survey on the topic.

%%%%%%%%%%%%%%%%%%%%%%%%%%%%%%%%%%%%%%%%%%%%%%%%%%%%%%%%
\begin{algorithm}[ht]
  \SetAlgoLined
  \vskip 1mm

 {\bf Setup}: Let $[0,T]$ be a training interval, and $\delta t = T/K$ with $K>0.$
 We assume that for each $k$ in $\{0,\cdots, K-1\}$, a numerical solution of Eq.~\eqref{Eq_quadODEs} is computed, which is denoted by $\boldsymbol{y}^k$. 
  \vskip 1mm
 
  {\bf Input}:   It consists of collecting the following projections of the numerical solution 
  \bes 
\underset{k = 0, \cdots, K}\bigcup (u_1^k,\cdots, u_m^k; u_n^k),
\ees
where $u^k_i= \langle \boldsymbol{y}^k, \boldsymbol{e}^*_i \rangle $ for $i = 1,\cdots, m$, and $u^k_n = \langle \boldsymbol{y}^k,  \boldsymbol{e}^*_n\rangle$, with $ \boldsymbol{e}^*_j$'s denoting the generalized eigenvectors associated with $A$ in Eq.~\eqref{Eq_quadODEs}.

  {\bf Output}: The optimal $\tau$-value,  $\tau^\ast$, that minimizes \eqref{Eq_defect_Markovian} is obtained as follow:
    \vskip 1mm
    
\nl Set parameter values for $\tau$, $\delta \tau$ and $\epsilon$,  which represent respectively the initial guess of $\tau^\ast$, the initial step size of $\tau$, and the convergence tolerance for the iteration. For instance, 
    \beas
    & \tau = 0; && \text{\texttt{\% initial guess}} \\
    & \delta \tau = 0.1; && \text{\texttt{\% initial step size of $\tau$}} \\
   % & I_{\mathrm{max}} = 10^4; && \text{\texttt{\% maximal number of iterations}} \\
    & \epsilon = 10^{-10}; && \text{\texttt{\% convergence tolerance}} \\
    \eeas
    \vskip 1mm

\nl  Compute $\boldsymbol{Q}_1$, $\boldsymbol{Q}_2$, $\widehat{\boldsymbol{Q}}_2$, $\widetilde{\boldsymbol{Q}}_2$, 
$\boldsymbol{Q}_3$, $\widetilde{\boldsymbol{Q}}_3$, and $\boldsymbol{Q}_4$ defined in \eqref{Eq_Q_elements_continuous} as well as $\langle u_n \overline{u}_n \rangle_T$ and $\langle u_n \rangle_T$ appearing in \eqref{Eq_defect_Markovian_recast} by using a standard numerical quadrature.

  \vskip 1mm
 
\nl Evaluate $\mathcal{Q}' = \frac{\d}{\d \tau}\mathcal{Q}(\tau, T)$ by using \eqref{Eq_Qprime}\;
 
 \While{$|\mathcal{Q}'| > \epsilon$}{
 \vskip 1mm
 Set $\tau_\delta = \tau - \mathrm{sgn}(\mathcal{Q}') \delta \tau$\;
 \vskip 1mm
  
 Compute $\mathcal{Q}'_\delta = \frac{\d}{\d \tau}\mathcal{Q}(\tau_\delta, T)$ by using \eqref{Eq_Qprime}\;
  \vskip 1mm
  
  \uIf{$|\mathcal{Q}_\delta'| > \epsilon \;\; \mathrm{ and } \;\; \mathrm{sgn}(\widetilde{\mathcal{Q}}') \neq \mathrm{sgn}(\mathcal{Q}')$}{
              \vskip 1mm
             $\delta \tau = \delta \tau/2$;
              \vskip 1mm
       }
       \Else{
              \vskip 1mm
              $\tau = \widetilde{\tau}$;
               \vskip 1mm
               
               $\mathcal{Q}' = \widetilde{\mathcal{Q}}'$;
               \vskip 1mm               
       }
 }
\vskip 1mm               
%$\tau^\sharp = \tau$.
% 
  \caption{Find the optimal $\tau$ for the minimization problem \eqref{Min_formulation_h1_quadcase} using a gradient-descent}
  \label{Algorithm-tau-search}
\end{algorithm}
%%%%%%%%%%%%%%%%%%%%%%%%%%%%%%%%%%%%%%%%%%%%%%%%%%%%%%%%
%\end{appendices}

%\bibliographystyle{amsalpha}
\newcommand{\etalchar}[1]{$^{#1}$}
\providecommand{\bysame}{\leavevmode\hbox to3em{\hrulefill}\thinspace}
\providecommand{\MR}{\relax\ifhmode\unskip\space\fi MR }
\providecommand{\MRhref}[2]{%
  \href{http://www.ams.org/mathscinet-getitem?mr=#1}{#2}
}
\providecommand{\href}[2]{#2}

\end{document}